\documentclass{article}
\usepackage[utf8]{inputenc}
\usepackage[margin=1in]{geometry}
\usepackage{relsize}
\usepackage{svg}
\usepackage{biblatex}
\usepackage[colorlinks=true, urlcolor=blue, linkcolor=red]{hyperref}
\usepackage{amsmath, amsfonts, amsthm, amssymb}
\usepackage{enumerate}
\usepackage{graphicx}
\usepackage{rotating}
\usepackage{floatpag}
\usepackage{fancyhdr}
\fancypagestyle{floatpage}{%
  \fancyhf{}
  \fancyfoot[C]{\makebox[\textwidth][r]{%
    \smash{\raisebox{\dimexpr\footskip+.5\textheight}{\rotatebox{90}{\thepage}}}}}%
}
\usepackage{comment}
\usepackage{color}
\usepackage{tikz-cd}
\newcommand{\R}{\mathbb{R}}
\newcommand{\Z}{\mathbb{Z}}
\newcommand{\kk}{\mathbf{k}}
\newcommand{\C}{\mathbb{C}}
\DeclareMathOperator{\ind}{ind}
\newcommand{\bS}{\mathbb{S}}
\renewcommand{\SS}{\mathbb{S}}
\newcommand{\CC}{\mathcal{C}}

\newcommand{\Grad}{\nabla}
\newcommand{\M}{\mathcal{M}}
\newcommand{\cM}{\overline{\mathcal{M}}}
\newcommand{\tensor}{\otimes}
\DeclareMathOperator{\colim}{colim}
\DeclareMathOperator{\Fred}{Fred}
\DeclareMathOperator{\vdim}{vdim}
\DeclareMathOperator{\coker}{coker}
\DeclareMathOperator{\dom}{dom}
\DeclareMathOperator{\im}{im}
\theoremstyle{definition}
\newtheorem{definition}{Definition}
\newtheorem{theorem}{Theorem}
\newtheorem{proposition}{Proposition}
\newtheorem{remark}{Remark}
\newtheorem{lemma}{Lemma}
\title{Cyclotomic Structures in Symplectic Topology}
\author{Semon Rezchikov}
\date{\today}

\begin{document}

\maketitle
\abstract{We extend the Cohen-Jones-Segal construction of stable homotopy types associated to flow categories of Morse-Smale functions $f$ to the setting where $f$ is equivariant under a finite group action and is Morse but no longer Morse-Smale. This setting occurs universally, as equivariant Morse functions can rarely be perturbed to nearby equivariant Morse-Smale functions. The method is very general, and allows one to do equivariant Floer theory while avoiding all the complications typically caused by issues of equivariant transversality. The construction assigns a (genuine) equivariant orthogonal spectrum to an equivariant framed virtually smooth flow category. Using this method, we construct, for a compact symplectic manifold $M$, which is symplectically atoroidal with contact boundary, and is equipped with an equivariant trivialisation of its polarization class, a cyclotomic structure on the spectral lift of the symplectic cohomology $SH^*(M)$. This generalizes a variant of the map which sends loops to their $p$-fold covers on free loop spaces to the setting of general Liouville domains, and suggests a systematic connection between Floer homology and $p$-adic Hodge theory. }

\section{Introduction}
It is a recurrent theme in mathematics that more symmetric representatives of a class of objects are forced to be more singular. This occurs already in the simplest geometric figures: it is easy to draw three lines in the plane such that the figure has a $\Z/6$ symmetry, but then no perturbation of the lines preserving the full $\Z/6$ symmetry will resolve the resulting triple point into double points.  In this paper, we will be concerned with \emph{equivariant Morse theory}, and with the problem of how to deal with the singular behavior of Morse functions which are invariant under the action of a group $G$.  It turns out that taking a homotopical perspective gives a universal solution to this problem.

Recall that a function $f: M \to \R$ on a manifold $M$ is \emph{Morse} if all of its critical points are nondegenerate, and that it is \emph{Morse-Smale} if all the stable manifolds $W^s(x)$ and unstable manifolds $W^u(y)$ intersect transversely for each pair of critical points $(x,y)$ of $f$. It is a basic fact that a generic function $f$ is both Morse and Morse-Smale. Moreover (if $M$ is e.g. closed) one can compute the homology of $M$ from counts of isolated gradient flow lines of $f$ via the \emph{Morse-Witten complex} $CM_*(f)$, and more generally one can compute more complex homotopy-theoretic invariants of $M$ using elaborate constructions based on gradient flows. One terminus of this perspective is the Cohen-Jones-Segal construction \cite{cohen1995floer, cohen2007floer}, which recovers the stable homotopy type $\Sigma^\infty_+ M$ from the data of all gradient flow lines of $f$, packaged together into the \emph{flow category} $\mathcal{C}(f)$. 

Now, if a finite group $G$ acts on $M$ and we require $f$ to be $G$-invariant, then, while Wasserman proved that a generic $G$-invariant $f$ is Morse \cite{wasserman1969equivariant}, then exist open neighborhoods of many equivariant Morse functions $f$ which \emph{contain no Morse-Smale representatives}. As such, it is impossible to straightforwardly generalize Morse-theoretic constructions to the equivariant setting: one must either have access to a non-generic equivariant function $f$ which is Morse-Smale, or one must perform some Borel-type construction to make the $G$-action free, which limits the range of equivariant homotopy invariants that one can study.
In this paper, we will see how to generalize these constructions to the setting where we require that $f$ is equivariant, by \emph{dropping the requirement that $f$ is Morse-Smale}. We will then apply these constructions to the setting of Floer theory, for which the finite-dimensional set-up of Morse theory is a model case.

Specifically, we focus on the setup of Hamiltonian Floer homology, where we have symplectically aspherical symplectic manifold $M$, and a Hamiltonian function $H: M \times S^1 \to \R$, where $S^1$ is a time coordinate, and we wish to investigate the \emph{dynamics} of the time-dependent vector field $X_H$ associated to $H$ via Hamilton's equations by studying the time $k$-periodic points of the flow of $X_H$, i.e. the periodic points of the time-1 flow $\phi^1_H$ of $X_H$. To gain insight into these periodic points,  in symplectic topology we consider the formal gradient flow of the symplectic action functional $\mathcal{F}: \mathcal{X} \to \R$, where $\mathcal{X} = LM$ is the space of contractible loops in $M$, and $\mathcal{F} = \mathcal{A}_H$ is the symplectic action functional, whose critical points correspond to time-$1$ periodic points of $X_H$. By speeding up time and rescaling $H$, we produce a modified sequence of action functionals $\mathcal{A}_{H^{\#k}}$ whose critical points are now in bijection with the time-$k$ periodic points of $X_H$, which include the time $1$-periodic points as a subset. Moreover, the action of $\Z/k\Z$ on $L_0X$ via \emph{loop rotation} makes $\mathcal{A}_{H^{\#k}}$ into a $\Z/k\Z$-equivariant action functional, and the fixed points of the $\Z/m\Z \subset \Z/k\Z$ subgroup for each $m$ dividing $k$ are the time $k/m$ periodic points of the flow of $X_H$.

The Morse homologies of $\mathcal{A}_{H^{\# k}}$ are the \emph{Hamiltonian Floer cohomology} groups $HF^*(H^{\# k})$, and the study of the underlying chain complexes and certain operations on filtrations these groups \emph{as $k$ varies} has led to many results on Hamiltonian dynamics (e.g. \cite{ginzburg2010conley, ginzburg2018hamiltonian, Shelukhin2022, prssrevisited}). Often, it is crucial to utilize the $\Z/k\Z$ symmetry of $\mathcal{A}_{H^{\#k}}$ in some way to establish these dynamical results.

In the symplectically atoroidal case (i.e. the integral of the symplectic form over any map from a torus to $X$ is zero) and under a topological assumption on $X$ the Hamiltonian Floer cohomologies of $H$ admit \emph{lifts} to the stable homotopy category 
\[CF^\bullet(H^{\# k}, \SS) \in Spectra \] in the sense that $CF^\bullet(H^{\# k}, \SS) \wedge H\Z$ is quasi-isomorphic to the Hamiltonian Floer chain complex \cite{cohen2007floer, large2021spectral, abouzaid2024foundation}. In this paper, we incorporate the $C_k$ symmetry into this construction. In the following, $M$ will always be compact and symplectically atoroidal, and $H$ will be admissible for $M$ (see Section \ref{sec:hamiltonian-floer-review}).

\begin{theorem}
\label{thm:equivariant-floer-spectra-exist}
Suppose that the equivariant polarization class $[\rho] \in KO_{C_k}^1(LX)$ (see Section \ref{sec:index-theory}) vanishes, and we have chosen a trivialization $\tilde{\rho}$ of the polarization map in the sense of Definition \ref{def:trivialization-of-polarization-map}. Then there 
    exist orthogonal $C_k$-spectra 
\begin{equation}
    \label{eq:spectral-lift-equivariant-hf-1}
    CF_\bullet(H^{\# k}, \SS) \in C_k-Spectra
\end{equation}
depending on $\tilde{\rho}$, which are well-defined up to canonical isomorphism in the homotopy category,  and which refine Floer homology in the sense that
\[ H_*(CF_\bullet(H^{\#k}, \SS), \Z) = HF_*(H^{\#k}, \Z). \]
\end{theorem}

One should think of the objects $CF_\bullet(H^{\#k}, \SS)$ in the above theorem as containing all possible relations between all conceivable $C_k$-equivariant Floer homology invariants.
Crucially, \emph{this construction cannot be done without the new methods of this paper} because no small $C_k$-equivariant perturbation of $\mathcal{A}_{H^{\# k}}$ will produce a functional with enough transversality to apply previously-existing methods in Floer homotopy theory. This obstruction exists for dynamical reasons, due to jumps in the Conley-Zehnder indices of \emph{elliptic} periodic points under iteration. We will explain this issue later in the introduction.

\begin{remark} In fact, the method of proof of Theorem \ref{thm:equivariant-floer-spectra-exist} shows that there is a weak equivalence of non-equivariant $H\Z$-module spectra
\[ H\Z \wedge CF_\bullet(H^{\#k}, \SS) \simeq CF_\bullet(\widetilde{H^{\#k}}, \Z) \]
where  $\widetilde{H^{\#k}}$ is a small regular admissible perturbation of $H^{\#k}$, and we view the Floer chain complex of $\widetilde{H^{\#k}}$ as an $H\Z$-module spectrum via the standard comparison map. 
\end{remark}

Now, part of the data of a genuine $C_k$-equivariant spectrum $Y$ is the data of the \emph{geometric fixed point spectra} $Y^{\Phi C_m}$ for every $m$ dividing $k$, which, in the case that $Y = \Sigma^\infty_+ Y_0$ is the suspension spectrum associated to a $G$-CW complex $Y_0$, capture the spectra associated to the actual fixed points on the $G$-action on $Y_0$, i.e. 
\[ (\Sigma^\infty_+ Y_0)^{\Phi C_m} = \Sigma^\infty_+ Y_0^{C_m}.\]

Using this construction, we prove the following relation between the equivariant stable homotopy type of a Hamiltonian and that of its iterate. Suppose that $H=H_0^{\#m}$ is an $m$-th iterate of another Hamiltonian, and otherwise we are in the setting of Theorem \ref{thm:equivariant-floer-spectra-exist}. 
\begin{theorem}
\label{thm:fixed-point-thm}
  There is a weak equivalence
    \[ (CF_\bullet(H^{\# k}, \SS))^{\Phi C_k} \simeq CF_\bullet(H,\SS) \in C_{m}-Spectra.\]
\end{theorem}
This is a Floer-theoretic manifestation of the obvious geometric homeomorphism 
\[ \phi_k: LX \simeq LX^{C_k} \subset LX, \gamma \mapsto \gamma^{\#k}\]
sending a loop to its $k$-fold iterate,
together with the identity of Floer action functionals 
\[ \mathcal{A}_{H} \circ \phi_k = k\mathcal{A}_{H_0}.\]

\begin{remark}
    By a standard argument this result immediately re-proves the localization theorem for Hamiltonian Floer homology that holds in this setting (which is proven in \cite{Shelukhin2021} in general; the setting of this paper is closer to that of \cite{seidel2010localization} specialized to the setting of the diagonal Lagrangian and the graph of the the Hamiltonian diffeomorphism associated to $H^{\#k}$).   Understanding in what sense a genuine equivariant lift of \cite{Shelukhin2021} holds over spectra in the exact setting (or over Kuranishi bordism \cite{pardon2022enough}, in the general setting) is an interesting topic for future work.
\end{remark}

Moreover, taking $H$ to be a Hamiltonian with appropriate dynamics at infinity for $X$, there is a standard definition of an invariant which measures something about the time $k$-periodic points for all $k$, called the \emph{symplectic cohomology}:
\[ SH^*(M) = \colim_{k} HF^*(H^{\# k}).\]
This invariant also admits a spectral lift $SH^*(M; \SS)$ when the polarization class vanishes \cite{cohen2007floer}. Like symplectic cohomology, this invariant is independent of of $H$. 

One can obtain Floer `cohomotopy' types $CF^\bullet(H; \SS)$ by considering the derived function spectra $F(CF_\bullet(H), \mathbb{S})$, where $\mathbb{S}$ is the sphere spectrum. By working in the category of genuine equivariant spectra, Theorem \ref{thm:equivariant-floer-spectra-exist} lets us make sense of $CF^\bullet(H^{\#k}; \SS)$ as an object of $C_k$-spectra. 

Under the conditions of Theorem \ref{thm:fixed-point-thm}, we then prove:
\begin{theorem}
\label{thm:sh-is-cyclotomic}
    Assume that we are in the setting of Theorem \ref{thm:fixed-point-thm}. Then there is an object $SH^\bullet(M, \SS)$ in the $\infty$-category of \emph{genuine $p$-cyclotomic spectra} \cite{hesselholt2004rham, nikolaus-scholze} such that the corresponding object $SH^\bullet(M, \SS)_0$ in the category of $\infty$-spectra has homology groups 
    \[ \pi_*(SH^\bullet(M, \SS))= SH^{-*}(M).\]
    In particular, there are lifts of $SH^\bullet(M, \SS)_k$ of $SH^\bullet(M, \SS)_0$ to the $\infty$-categories of $C_{p^k}$ spectra for all $k$ such forgetting the $C_{p^k}$ action on $SH^\bullet(M, \SS)_k$ to $C_{p^{k-1}}$ gives an object which is equivalent (as a $C_{p^{k-1}}$-spectrum) to $SH^\bullet(M, \SS)_k$, and similarly 
    \begin{equation}
    \label{eq:cyclotomic-equivalences-for-sh}
        (SH^\bullet(M, \SS)_k)^{\Phi C_k} \simeq (SH^\bullet(M, \SS)_{k-1}. 
    \end{equation}
\end{theorem}

This is a Floer-theoretic analog of the statement that the fixed points of the $C_k$ action on the free loop space are homeomorphic to the free loop space itself via the map which takes loops to their $k$-fold covers; moreover, this is an $(S^1 \simeq S^1/C_k)$-equivariant homeomorphism. Indeed, symplectic cohomology is formally the cohomology of a Morse function on the free loop space of $M$, and thus carries many structures like those carried by the homology of free loop spaces \cite{abbondandolo2006floer}. 

\begin{remark}
Theorems $1$, $2$, and $3$ continue to hold for symplectically aspherical manifolds without contact boundary, provided that one only considers Floer homology or symplectic cohomology of contractible fixed points.
\end{remark}

\begin{remark}
    The methods of this paper can be generalized to construct a `genuine $\mathbb{Q}/\mathbb{Z}$-cyclotomic spectrum $SH^\bullet(M, \SS)$; to build the full $S^1$ action, one requires the extension of the methods of this paper to the setting of \emph{Morse-Bott} functions. The latter correspond to time-independent Hamiltonians. 
\end{remark}

\begin{remark}
    One may use the construction of C\^{o}t\'{e}-Kartal \cite{cote-kartal-1} to associate a genuine equivariant $S^1$-spectrum to symplectic cohomology, and then forget the action to $C_k$ to produce genuine equivariant $C_k$-spectra with the same homology groups as symplectic cohomology. However, these latter $C_k$-spectra will \emph{not} agree with the spectra produced in this paper. Instead, they are a \emph{Borel analog} of the spectra produced in this paper: the $S^1$-actions (and thus $C_k$-actions) on the spectra in \cite{cote-kartal-1} are \emph{free}, and thus their geometric fixed points under \emph{any} nontrivial subgroup of $S^1$ are simply the \emph{zero} spectrum, in contrast to \eqref{eq:cyclotomic-equivalences-for-sh}. We expect that that by taking the smash product of the spectra of this paper with $ES^1$ recovers the spectra equivalent to those in  \cite{cote-kartal-1}. We will refer to the spectra contructed by C\^{o}t\'{e}-Kartal `Borel-equivariant spectral symplectic cohomology' below.
\end{remark}

\begin{remark}
A genuine $p$-cyclotomic spectrum is in particular a genuine $C_{p^\infty}$-spectrum, where $C_{p^\infty} \subset \mathbb{Q}/\mathbb{Z}$ is the subgroup of $p$-power torsion elements. As such, one can take geometric (or categorical, etc) $C_{p^k}$ fixed points, for any $k$, of such an object. However, even though $C_{p^\infty}$ is dense in $S^1$, one cannot extract the $S^1$-fixed points from the underlying genuine $C_{p^\infty}$ spectrum of a genuine $S^1$-spectrum \cite{degrijse2023proper}. 
\end{remark}

This construction turns out to give a connection between equivariant string topology and certain algebraic structures that arise in the context of $p$-adic Hodge theory. Fully explaining this connection is outside the scope of this paper, but we will return to it later in this introductory section.

\paragraph{Morse theory and equivariant transversality.} 

Let us go back to the more concrete setting of equivariant Morse theory. 

The basic obstruction to equivariant transversality in the setting of Morse theory comes from equivariant index theory. Imagine for example that $G = \Z/2$ acts isometrically on a closed manifold $M$ with Riemannian metric $g$. Then the fixed locus $M^{G}$ a manifold. Suppose furthermore that the $G$-invariant Morse function $f$ has a pair of invariant critical points $(x,y)$ with $\{x,y\} \subset  M^G$. We will write $\ind_{M}(x)$ for the index of $x$ as a critical point of $f$, and $\ind_{M^G}(x)$ the index of $x$ thought of as a critical point of $f|_{M^G}: M^G \to \R$. Suppose further that $f|_{M^G}$ is Morse-Smale, as will hold $C^\infty$-generically among $G$-invariant functions $f$; and moreover that $\ind_{M^G}(x) - \ind_{M^G}(y) = 1$, and that there is a single isolated gradient flow line $\gamma$ of $f|_{M^G}$ from $x$ to $y$. What may occur in this setting is that we have $\ind_M(x) - \ind_M(y) < 1$: the flow line $\gamma$ is an \emph{obstructed} flow line of the the gradient flow of the ambient Morse-function $f$. Thus, if we could make $f$ Morse-Smale by a $C^\infty$-small perturbation of $f$ among $G$-equivariant functions, we would expect that there would no longer be any flow lines of the perturbed function $\tilde{f}$ from $x$ to $y$. However, this is impossible: small equivariant perturbations of $f$ give small perturbations of $f|_{M^G}$, which is already Morse-Smale, and an isolated flow-line of a Morse-Smale function persists under all sufficiently small perturbations. Thus, such a function $f$ has no small equivariant perturbations which will make it Morse-Smale. 

\begin{figure}[h!]
    \centering
    \includegraphics[width=0.8\textwidth]{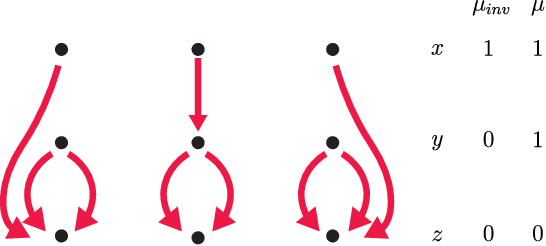}
    \caption{\emph{Flow category of an $C_2$-equivariant Morse function $f$ which is Morse but not Morse Smale.}. The flow lines of such a function function, depicted by the middle diagram on the left of the figure, as an invariant flow line from $x$ to $y$ and a pair of flow lines, forming a free $C_2$ orbit, from $y$ to $z$. All three critical points $x,y,z$ are invariant, and the indices thought of as critical points of $f^{C_2}$ are denoted $\mu_{inv}$, while their indices in the ambient space are denoted $\mu$. In particular, the flow line is regular from $x$ to $y$ is regular as a flow line in the invariant locus but not regular as a flow line in the ambient space, and no small perturbation of $f$ is Morse-Smale in the ambient space because $\mu(x) - \mu(y) < 1$. Thus, upon a perturbation (corresponding to moving from the middle to the left or to the right side diagrams of the left half of the figure) this flow line disappears; however, it is forced to glue to exactly one of the regular flow lines from $y$ to $z$ to produce a single flow line from $x$ to $z$, regardless of perturbation. (See Figure \ref{fig:virtually-smooth-category-obstructed-example} for further analysis). Note that looking at `half' of this figure we get the usual transformation of a Morse function through a handle-slide; thus the function depicted in the middle diagram is `in the middle of two handle-slides'.  }
    \label{fig:obstructed-c2-example}
\end{figure}

What goes wrong in this example is that, if we look at each invariant critical point $x \in M^G$ of $f$, and we view the the negative eigenspaces $T^-_xM \subset T_xM$ of the Hessians of $f$ as $G$-representations, then the spaces $T^-_xM \cap (TM^G)^\perp$ of eigenvectors normal to the invariant locus may jump in dimension as the invariant critical point $x$ varies. The dimensions of these spaces are simply $\ind_M(x) - \ind_{M^G}(x)$, and agree with the sum of the dimensions of the nontrivial isotypic components of $T^-_xM$ viewed as a $G$-representation. The fact that the numbers  $\ind_M(x) - \ind_{M^G}(x)$ can jump is the first topological obstruction to achieving equivariant transversality for small perturbations of $f$. While this obstruction is in general nontrivial, in the finite-dimensional context one can construct functions $\tilde{f}$ for which this obstruction vanishes, and which are indeed equivariant and Morse-Smale. When $M$ is a closed manifold and $G=\Z/2$, the method to construct such $f$ is as follows:  first,  chooses a Morse-Smale function $f|_{M^G}$ on $M^G$, then \emph{extends} $f$ $f|_{M^G}$ to a neighborhood $U$ of $M^G$ by adding a large positive $G$-invariant quadratic form along the normal bundle to $M^G$ in $M$, and then extends $f|_U$ to the rest of $M$ \cite{seidel2010localization}. For such a function, we will have that $\ind_{M^G}(x) = \ind_{M}(x)$ for all invariant critical points $x$. This construction can be generalized to the setting of general finite groups $G$, and shows that \emph{in the finite-dimensional context}, given a $G$-invariant Morse function $f$, one can perform a \emph{large} perturbation to $f$ to produce a $G$-equivariant Morse-Smale function. 

Considering the infinite dimensional context, the difference between the index of a fixed point of a Hamiltonian and the index of its iterate (which is an invariant fixed point of the iterated Hamiltonian) can jump in complicated ways, leading to the difficulties with equivariant transversality described in this section.

\paragraph{Floer theoretic approaches to equivariant transversality.} However, in the infinite-dimensional setting of Floer homology, the large perturbations of Morse functions outlined in the previous setting are not available. Recall that in Floer theory, one aims to define Morse complex $CM_*(\mathcal{F})$ of some infinite-dimensional Morse function $\mathcal{F}: \mathcal{X}\to \R$, typically arising from a gauge theory or a $\sigma$-model. The basic phenomenon distinguishing Morse theory from Floer theory is that the gradient flows of $\mathcal{F}$ are badly behaved, and do not have long-time existence from most initial initial conditions. As such, when one succeeds at defining $CM_*(\mathcal{F})$, this homology does \emph{not} compute the homology of the ambient space $\mathcal{X}$, but instead a different and \emph{more useful} invariant. In this setting, due to analytic challenges, one only has access to \emph{small} perturbations of $\mathcal{F}$ which can be constructed via applications of the Sard-Smale theorem. The obstruction that $\ind_M(x) - \ind_{M^G}(x)$ is independentof $x$ then becomes an obstruction of equivariant Atiyah-Singer family index theory, and specifically can be expressed as the nontriviality of a class in $\rho^{\perp G} \in KO^1_G(\mathcal{X}^G)$. When this class is nontrivial, there is no hope of achieving transversality for small perturbations of $\mathcal{F}$. 

However, because this problem arises recurrently in Floer homology, various approaches have been introduced to deal with this problem. One approach is to simply understand multiple-cover contributions `by hand'; this is the perspective taken in Wendl's work in superrigidity \cite{wendl2023transversality} and subsequent works \cite{doan2021gopakumar, doan2021castelnuovo}. A related approach is to modify $\mathcal{X}$ by performing a \emph{blow-up} along $\mathcal{X}^G$, replacing $\mathcal{X}$ with a related space $\tilde{\mathcal{X}}$ on which the $G$-action is free and such that the obstructions to small perturbations coming from equivariance no longer arise in the modified problem. This is the approach taken by Kronheimer-Mrowka in the definition of Seiberg-Witten Floer homology \cite{kronheimer2007monopoles}, and have subsequently been used by researchers in gauge theory \cite{lin2014morse, miller2019equivariant} and symplectic topology \cite{large2019equivariant} to define Floer-theoretic invariants in the presence of symmetry. The complexity of this approach scales with the complexity of the group $G$, but it tends to have a concrete and analytic flavor. 

A radically different approach can be found in the work of Bauer-Furuta \cite{bauer2004stable} and Manolescu \cite{manolescu2003seiberg}, which define certain Floer-theoretic invariants as morphisms or objects in the \emph{genuine equivariant stable homotopy category}.  These approaches  bypass the tension between equivariance and  transversality by constructing the underlying theoretic objects directly, without perturbing the defining data. The idea behind these methods is a generalization of the following: given any map of spheres $a: S^n \to S^n$, one can define the degree of $a$ by perturbing to a nearby $\tilde{a}$ and counting the number of preimages of a generic point with sign, \emph{or} one can define the degree directly via homotopy theory as the corresponding equivalence class of maps $\Sigma^\infty_+ [a]$ in the stable homotopy category, \emph{without} perturbing $a$. One produces the map $a: S^n \to S^n$ from Floer-theoretic equations via finite-dimensional reduction technique; in the equivariant setting, the map $a$ is automatically an \emph{equivariant} mp of spheres, and so in the homotopy-theoretic approach one automatically constructs a \emph{more refined invariant} than the constructions arising from the Kronheimer-Mrowka method by considering the map as a morpism in the \emph{equivariant stable homotopy category}. Because this homotopical approach does not involve perturbations of the resulting action functional, it sometimes involves less analytic and geometric detail than the Kronheimer-Mrowka blow-up construction.  However, the homotopical approach in gauge theory relies crucially on the technique of \emph{global finite-dimensional approximation}, which in the end requires some kind of global linear structure on $\mathcal{X}$. In the case of Seiberg-Witten theory, the linear structure is available because $\mathcal{X}$ is the set of connections on a manifold. Analogously, the technique of \emph{generating functions} can be used construct homotopy theoretic lifts of Floer-theoretic invariants in symplectic topology \cite{abouzaid2016immersion}, when some variant of a global linear structure is available,such as when when $\mathcal{X}$ is the space of free loops in the cotangent bundle to a closed manifold $Q$. However, in many cases, such as when $\mathcal{X}$ is the space of free loops in a Liouville domain, no global linear structure is available to make global finite-dimensional approximation techniques possible. 

\paragraph{Virtual Equivariant Morse Theory.} In this paper, we extend the homotopical techniques of \cite{bauer2004stable, manolescu2003seiberg} to general setting where global finite-dimensional approximation techniques are not available. Specifically, given a function $\mathcal{F}$ of Floer theoretic type, one is (in many cases) able to construct a \emph{flow category} $\mathcal{C}_\mathcal{F}$ with objects the critical points of $\mathcal{F}$ and morphism spaces $\mathcal{F}(x,y)$ given by a compactification of the space of gradient flows of $\mathcal{F}$ from $x$ to $y$. When $G$ acts on $\mathcal{X}$ and $\mathcal{F}$ is $G$-equivariant, the flow category $\mathcal{C}_\mathcal{F}$ becomes $G$-equivariant: there is a $G$-action on objects and a $G$-action on morphisms which commutes with the source and target maps. When one can achieve transversality for the gradient flow of $\mathcal{F}$, the flow category $\mathcal{C}_\mathcal{F}$ becomes \emph{smooth}. Namely, when $\mathcal{C}_\mathcal{F}$ is smooth, the morphism spaces of $\mathcal{C}_\mathcal{F}$ are smooth manifolds with corners, and one can apply the Cohen-Jones-Segal (CJS) construction \cite{cohen1995floer} to associate a stable homotopy type to $\mathcal{F}$ in a way which recovers the stable homotopy type of $M$ when $\mathcal{X}=M$ finite-dimensional. However, due to the tension between equivariance and transversality, one generally cannot produce $\mathcal{F}$ such that $\mathcal{C}_\mathcal{F}$ is a smooth $G$-equivariant flow category.

\begin{figure}[h!]
    \centering
    \includegraphics[width=\textwidth]{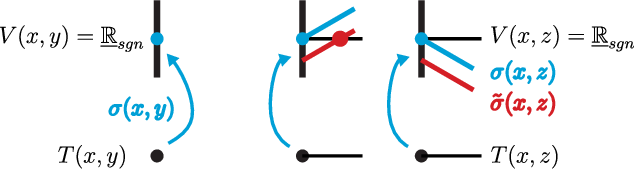}
    \caption{\emph{Analysis of the $C_2$-equivariant function of Figure \ref{fig:obstructed-c2-example} via virtually smooth equivariant flow categories.} (All other virtual smoothings must have similar features). Writing $\mu_\nu= \mu - \mu_{inv}$, since the normal index difference $\mu_\nu(x) - \mu_\nu(y) = 1$ we have that the obstruction space $V(x,y)$ associated to the single flow line from $x$ to $y$ is a copy of the sign representation of $C_2$; the obstruction section is a map from a point to zero in this vector bundle, denoted by the blue arrow. The two flow lines from $y$ to $z$ are regular, so they do not carry an obstruction bundle; thus the Kuranishi chart around the flow lines from $x$ to $z$ is forced to look as in the right hand side of this figure; the two pieces of the graph of $\sigma(x,z)$ have opposite slopes because of the requirement that $\sigma(x,z)$ is $C_2$ equivariant, and the $C_2$-action permutes the two components of $T(x,z)$ corresponding to the two thickenings around the two broken flow lines $x \to y \to z$. Choosing a compatible family of non-equivariant perturbations of the Kuranishi section (red) must move $\sigma(x,y)$ above or below zero (in this case below), and thus forces the perturbation of $\tilde{\sigma}(x,z)$ of $\sigma(x,z)$ to be such that one has a zero on one component of the thickening but not the other. This corresponds to the fact that upon a perturbation of the Morse function to a (non-invariant) Morse-Smale function, the flow line from $x$ to $y$ will glue on to \emph{one} of the flow lines from $y \to z$, as in Figure \ref{fig:obstructed-c2-example}.
   Note that in this example, the Morse complex of any small perturbation of this function will satisfy $dx = z \;(mod\, 2)$, even there is no invariant flow line from $x$ to $z$. In particular, the Morse complex of $f^{C_2}$ cannot be computed from a neighborhood of invariant part of the virtually smooth flow category, but interacts with the rest of the flow category in a global fashion. This is a subtlety in the behavior of equivariant Morse functions which are not Morse-Smale with respect to the theory of equivariant localization in for $C_2$-actions.}
    \label{fig:virtually-smooth-category-obstructed-example}
\end{figure}

However, one when $\mathcal{F}$ is Morse but not Morse-Smale, in many settings one can equip $\mathcal{C}_\mathcal{F}$ with the structure of a \emph{virtually smooth flow category}: one can write each morphism space $\mathcal{C}_\mathcal{F}(x,y)$ as the zero-set of a section $\sigma(x,y)$ of a vector bundle $V(x,y)$ on a manifold with corners $T(x,y)$, in a way that is compatible with composition in the category $\mathcal{C}_\mathcal{F}$. In particular, this has been achieved when $\mathcal{F}$ is symplectic action functional via \cite{rezchikov2022integral, bai2022arnold}, following \cite{AMS}. The first observation of this paper is to note that the CJS construction can be \emph{modified} to a new construction, which we call the \emph{virtual CJS construction}, which takes as input a \emph{virtually smooth flow category} rather than a smooth flow category, and reduces to the CJS construction when the virtually smooth flow category is smooth. The second observation of this paper is that, while it is difficult to perturb $\mathcal{F}$ among equivariant functionals to some equivariant $\tilde{\mathcal{F}}$ such that $\mathcal{C}_{\tilde{\mathcal{F}}}$ is smooth, there is automatically a $G$-action already on the virtually smooth flow category associated to $\mathcal{F}$: there are $G$-actions on $T(x,y)$ and $V(x,y)$ such that $V(x,y)$ is a $G$-equivariant vector bundle, $\sigma(x,y)$ is a $G$-equivariant section, and the $G$-actions are compatible with composition. The third observation is that virtual CJS construction easily adapts to produce genuine equivariant stable homotopy types from equivariant virtually smooth flow categories. This construction encapsulates very complicated combinatorics about the way in which obstructed invariant trajectories can combine with non-invariant trajectories to contribute nontrivially to the differential of the Morse complex of a non-Morse-Smale $G$-equivariant Morse-function, as we illustrate via a simple example in Figure \ref{fig:obstructed-c2-example} (which describes the virtual smoothing of the flow category described in Figure \ref{fig:virtually-smooth-category-obstructed-example}).

\begin{remark}
    A slogan for the ideas of this paper is: 
    \emph{A Morse function $f$ that is not Morse-Smale does have a `Morse complex', but does still define a filtered homotopy type; the condition that $f$ is Morse-Smale simply means that this filtration is a CW filtration}. Thus the spectral sequence of a CW filtration (the `Morse complex') collapses at the second page; in general there many pages of this spectral sequence (corresponding to `contributions of broken Morse trajectories to the differential'). 
\end{remark}

The general general constructions regarding virtuallly smooth equivariant flow categories can be found in Sections \ref{sec:flow-categories}, \ref{sec:flow-category-constructions}, and Section \ref{sec:floer-homotopy}. However, to \emph{produce} these categories, one must do some geometry, consisting essentially of the construction of the appropriate equivariant virtual smoothings (Section \ref{sec:global-charts}) and their equivariant framings (Section \ref{sec:producing-framings}). The management of the compatibility conditions associated to a framing of an equivariant virtually smooth flow category takes up a significant portion of the work of this paper. In particular, to make sure that the conditions of Theorem \ref{thm:sh-is-cyclotomic}, we prove the following 
\begin{theorem}
\label{thm:condition-for-vanishing-of-equivariant-polarization-class}
    The equivariant polarization class $[\rho] \in KO^1_{S^1}(LM)$ is zero whenever the tangent bundle of $M$ is \emph{polarized}, i.e. it is a complexification of a real bundle. 
\end{theorem}
When $M$ is a Weinstein domain, this is the condition under which it is known to be arborealizable \cite{alvarez2020positive}; it would be interesting to see if this class always vanishes on arborealizable Weinstein domains. The vanishing of the image of this class in $KO^1(LM)$ is expected to be the requirement for the definition of the spectral Fukaya category of this Weinstein domain \cite{cohen2007floer, large2021spectral}.   Part of the contribution of this paper is to clarify the technical aspects of the existing constructions of framings of framed flow categories in Floer homology (in part to verify that equivariance does not pose a technical challenge); such technical aspects which are remarked on throughout the paper.

\paragraph{Topological Hochschild Homology.}
As explained earlier in the introduction, we apply the equivariant virtual CJS construction to the setting of Hamiltonian Floer homology to produce the equivariant spectra of Theorem \ref{thm:sh-is-cyclotomic}. These equivariant structures on spectral symplectic cohomology are meant to give geometric models of certain structures associated to the Fukaya category of $M$ over the sphere spectrum. 

Indeed, when $M$ is aspherical and the polarization class $[\rho] \in KO^1(LM)$ vanishes, there is conjecturally a spectrally-enriched category $Fuk_\SS(M)$ which recovers the full subcategory of the Fukaya category $Fuk_\kk(M)$ consisting of those Lagrangian submanifolds withich are ``oriented with respect to $\SS$'' defined over a field $\kk$ by smashing with $H\kk$. On the level of stable homotopy groups, this category has been studied by Large \cite{large2021spectral}. Similarly, there should be spectral analogs $WF_\SS(M)$ of the wrapped Fukaya category $WF_\kk(M)$ \cite{abouzaid2011cotangent},  and many of the usual ideas about generation criteria \cite{abouzaid2010geometric, ganatra2020covariantly}  are expected to have analogs in the spectral setting. A constructible sheaf model of this latter category has already been constructed \cite{jin2017brane}.

In particular, one expects that there is an \emph{open-closed map}
\begin{equation}
\label{eq:spectral-open-closed-map}
    \mathcal{OC}: THH(WF_\SS(M)) \to \Sigma^n SH^\bullet(M, \SS)
\end{equation}
from the \emph{Topological Hochschild Homology} \cite{hesselholt1997k} to the spectral symplectic cohomology, paralleling the corresponding construction on the level of homology \cite{abouzaid2010geometric}. This map should be a map of cyclotomic spectra, where one uses the standard cyclotomic structure on $THH$ for the domain and the cyclotomic structure of Theorem \ref{thm:sh-is-cyclotomic} for the codomain. In good cases such as the setting of Weinstein domains, the spectral analog of the generation criterion \cite{abouzaid2010geometric} should be satisfied, and so $SH^\bullet(M, \SS)$ with its cyclotomic structure should give a geometric model of $THH(WF_\SS(M))$. Making sense of this idea is the fundamental motivation of this paper.
\begin{remark}
    The putative genuine $S^1$-action on symplectic Floer cohomotopy is a subtle invariant. In particular, the author believes that the map \eqref{eq:spectral-open-closed-map} will not be realizable as a map of map of genuine $S^1$ spectra. Indeed, the genuine $S^1$-action on a putative genuine $S^1$-variant of $SH^\bullet(M, \SS)$ should satisfy 
    \[ SH^\bullet(M, \SS)^{\Phi S^1} = F(\Sigma^\infty_+ M, \SS).\]
    As such this spectrum knows the underlying stable homotopy type of the given symplectic manifold. In contrast, the $S^1$-fixed points of the geometric realization of a cyclic set are always discrete. The object $THH(\CC)$ for a stable $\infty$-category $\CC$ should instead be thought of as defining an object in the $\infty$-category of $\mathbb{Q}/\Z$-spectra equipped with a compatible homotopy-action of $S^1$ (we thank Andrew Blumberg for this comment). We expect that the open-closed map \eqref{eq:spectral-open-closed-map} can be refined to intertwine all of these latter structures. 
\end{remark}

\begin{remark}
    The expected incompatibility of the open closed map with the expected genuine $S^1$-spectrum structure on symplectic cohomology explains in part the following well-known phenomenon. Subcritical surgeries on $M$ do not change its Fukaya category up to weak equivalence, while they do change the stable homotopy homotopy type of $M$. If one believes that subcritical surgeries do not change the  spectral wrapped Fukaya category as well, then one is forced to conclude  (conditional on the equation above one) that   \eqref{eq:spectral-open-closed-map} cannot be made $S^1$-equivariant. In contrast, the methods of \cite{cote2024recovering} rely on the action filtration on Borel-equivariant spectral symplectic cohomology, which is not preserved under any operations on $M$.    
\end{remark}

Given this correspondence, one expects that the cyclotomic structure on spectral symplectic cohomology should have rich algebraic properties. For example, on the level of cohomology, under the condition that $c_1(M) = 0$, $SH^*(M)$ is a BV algebra, with the $S^1$ action associated to the $S^1$ action on Hochschild homology via the homological analog of the open-closed map \eqref{eq:spectral-open-closed-map}, and the Gerstenhaber algebra structure associated to the Gerstenhaber algebra structure on $HH^*(Fuk_{\kk}(M))$ via a \emph{closed-open map} from $SH^*(M)$. The compatibility of the two structures is a reflection of the \emph{Calabi-Yau property} of the Fukaya category \cite{ganatra2019cyclic}. Analogs of the Calabi-Yau property should exist for $WF_\SS(M)$ \cite{cohen2019twisted} under strong topological conditions on $M$ (e.g. when $M= T^*Q$ for a closed framed manifold $Q$, see below) and the cyclotomic structure should satisfy intriguing compatibility relations with an $E_2$-ring structure on $SH^\bullet(M, \SS)$ \cite{abouzaid2022framed}.

One can study these compatibility relations from multiple perspectives. 

\paragraph{String topology.}
Let us specialize the open-closed map \eqref{eq:spectral-open-closed-map} to to the case where $M = T^*Q$ with $Q$ a closed manifold. This is the setting of \emph{string topology} \cite{cohen2005notes}, which studies the inter-relationship between homological and cohomological operations on the free loop space $LQ$ based on geometric operations involving loops in $Q$; algebraic structures on the Hochschild (co)homology of the $dg$-rings $C_*(\Omega Q)$ and $C^*(Q)$; and symplectic operations on symplectic cohomology $SH^*(T^*Q)$. The general pattern that occurs is that there is a precise correspondence between structures in all three settings, and symplectic topology allows one to generalize the setting of cotangent bundles to the general setting of Weinstein domains (which can be thought of as `cotangent bundles of singular spaces'). 

By analogy to the situation over a field, in the case that $Q$ is a stably framed manifold, we should have that
\[ WF_\SS(T^*Q) = Perf_\SS(\Sigma^\infty_+ \Omega Q),\]
where $\Sigma^\infty_+ \Omega Q$ is an $\mathbb{E}_1$-algebra and $Perf_\SS$ denotes the $\infty$-category of perfect modules over this algebra. (Outside the case that $Q$ is stably framed, one may need to appropriately twist the suspension spectrum of the based loop space, by analogy to the conditions on the spin-ness of $Q$ that arise in \cite{abouzaid-sh}.) Thus, in this setting, \cite{hesselholt1997k}
\begin{equation}
    \label{eq:thh-of-free-loop-space}
    THH(WF_\SS(T^*Q)) = \Sigma^\infty_+ LQ.
\end{equation}
The cyclotomic structure on $THH$ in this setting really is induced from the map 
\begin{equation}
\label{eq:p-fold-cover-map}
    \phi_p:\Sigma^\infty_+(LQ) \to \Sigma^\infty_+(LQ)
\end{equation} 
sending loops to their $p$-fold covers. The equivalence \eqref{eq:thh-of-free-loop-space} should factor through the open-closed map \eqref{eq:thh-of-free-loop-space} by composing the latter with a spectral version of the Abbondondalo-Schwartz map 
\cite{abbondandolo2006floer}
\begin{equation}
    AS: SH^\bullet(T^*Q, \SS) \to \Sigma^{-n}\Sigma^\infty_+(LQ). 
\end{equation}
This latter map, which is known to exist over $\Z$, is known to have a rich interaction with the string-topology operations on $\Sigma^{-n}\Sigma^\infty_+(LQ)$; in particular, which should make the latter into a framed $E_2$ algebra \cite{cohen2002homotopy}. 

From the categorical perspective (e.g. \eqref{eq:thh-of-free-loop-space} or \eqref{eq:spectral-open-closed-map}), while $THH(C)$ does not have an algebra structure when $C$ is a spectral category, $THC(C)$, the \emph{topological Hochschild cohomology of $C$}, does. The latter is isomorphic to $THH(C; P)$ for $P$ some `twisting bimodule' up to a shift when $C$ is a twisted, smooth Calabi-Yau category; in the case at hand since $Q$ is stably framed this twisting bimodule is the trivial bimodule \cite{cohen2019twisted}. Thus the $E_2$ structure on $THC$ and the $S^1$ action on $THH$ combine to a framed $E_2$-algebra structure, and the open closed map together with the dual `closed-open map' \cite{abouzaid2010geometric}
\[ \mathcal{CO}: SH^\bullet(M, \SS) \to THC(WF_\SS(M))\]
should intertwine a framed $E_2$-algebra structure on symplectic cohomology with each of these structures (this is known after taking smash products with $H\Z$, \cite{abouzaid2022framed}).

A description, even in this relatively simple setting, of the compatibility relations between the cyclotomic structure (which is induced by the very geometric $p$-fold cover map \eqref{eq:p-fold-cover-map}) and the framed $E_2$-algebra structure on $\Sigma^{-n}\Sigma^\infty_+ LQ$, is not known to the author, and we leave this to future work. We expect that that there are general relations between the cyclotomic structure on symplectic cohomology described in this paper and the framed $E_2$ algebra structure that hold on all Weinstein domains for which these structures can be defined; articulating this when $M = T^*Q$ is a basic problem in algebraic topology. It is known
that on the level of rational string topology of simply connected manifolds, the $p$-fold cover map \eqref{eq:p-fold-cover-map} induces an \emph{addtional grading} on the string topology algebra $H_{-n}(LQ, \mathbb{Q})$ \cite{felix2008rational}. Moreover, this map has a rich history of applications in dynamics \cite{gromoll1969periodic} and in Viterbo's original appliations of symplectic cohomology \cite{viterbo1997exact}. This paper can be motivated by the geometrically natural problem of generalizing the existence of the $p$-fold cover map \eqref{eq:p-fold-cover-map} as much as possible to the setting of all of general symplectic manifolds. Due to the fact that symplectic cohomology enjoys convenient functoriality properties, understanding the relation between string topology operations and the cyclotomic structure can likely be applied to  conventional problems in symplectic topology.

\paragraph{Symplectic topology and $p$-adic Hodge theory.}
From a different perspective, we can think of the cyclotomic structure on spectral symplectic cohomology as allowing us to access \emph{arithmetic aspects of mirror symmetry}. Indeed, in mirror symmetry, the symplectic topology of $M$ is often related to the algebraic geometry of some mirror variety $M^\vee$. In many cases, $M^\vee$ is not an arbitrary variety over $\C$, but is in fact naturally defined over a finitely-generated, arithmetically interesting ring $R$. Moreover, certain arithmetic aspects of the mirror can be seen from the symplectic topology of $M$, in special cases \cite{kontsevich-vologodsky, arithmetic-greene-plesser, Lekili2016}. From the perspective of symplectic topology, we can partially explain this phenomenon by noting that for geometric reasons, the Fukaya category of many manifolds is naturally defined over $\Z$, or over an explicit localization of $\Z$, rather than over $\C$ \cite{seidel2008fukaya, sheridan2021homological}. On the level of enumerative geometry, mirror symmetry connects the quantum or symplectic cohomology of $M$ with the Hodge theory of $M^\vee$. However, when $M^\vee$ is e.g. defined over $\Z$, its Hodge theory also has $p$-adic aspects \cite{2005.07919}, which can be in part understood by reference to Topological Hochschild Homology, as explained by Hesselholt \cite{hesselholt1997k, hesselholt2004rham}, Kaledin \cite{kaledin2010motivic}, and Bhatt-Morrow-Scholze \cite{bhatt2019topological}. As such, we expect that the cyclotomic structure of Theorem \ref{thm:sh-is-cyclotomic}  (or rather, its expected generalizations to general symplectic manifolds) should allow us to assign analogs of ``$p$-adic Hodge theoretic'' invariants to general symplectic manifolds, which should match with the corresponding invariants in cases where mirror symmetry applies. 


A proper discussion of conjectures and computations about the cyclotomic structure on symplectic cohomology and its relation to $p$-adic Hodge theory (or algebraic $K$-theory) is far outside the scope of this paper. We hope to return to these ideas in future works. In this paper, we establish the foundational geometric techniques that make discussion of this structure possible. We will end the introductory discussion with a precise foundational conjecture, which we do not know how to resolve, and has a significant impact on any future theory of ``arithmetic aspects of symplectic topology'':

\textbf{Question:} Let $M$ be a Weinstein domain and suppose that we are in the setting of Theorem \ref{thm:sh-is-cyclotomic}. Does the composed map 
\begin{equation}
    SH^\bullet (M, \SS) \to (SH^\bullet(M, \SS))^{\Phi C_p} \to  (SH^\bullet(M, \SS))^{tC_p},
\end{equation}
where the second map in the composition is the canonical map from geometric fixed points to tate fixed points
admit a \emph{Frobenius lift}, i.e. a lift to the homotopy fixed points $(SH^\bullet(M, \SS)_{\hat{p}})^{hC_p}$? The same question may be asked for the corresponding map obtained by $p$-completing $SH^\bullet(M, \SS)$.

Answering this question in the affirmative would give a direct analog of the $p$-cover map \eqref{eq:p-fold-cover-map} to general symplectic manifolds by composing with the standard map $(SH^\bullet(M, \SS))^{hC_p} \to (SH^\bullet(M, \SS))$ obtained by evaluating a map from $EC_p$ at the base-point. 

\paragraph{Outline of paper and its techniques.} We begin this paper with the pedagogical Section \ref{sec:pedagogical-intro}. We first explain in full detail the special case of the virtual CJS construction, and its equivariant analog, in the case where the flow category has only two objects, both of which are invariant. We then recall enough about the CJS construction to explain why the virtual CJS construction agrees with the CJS construction in this case. 

Then, in Section \ref{sec:equivariant-homotopy-review}, we review the very basics of equivariant homotopy theory and the equivariant spanier-whitehead category; this section can be skipped by those familiar with the material. In Section \ref{sec:flow-categories}, we explain what is a virtually smooth flow category, how to frame such a category, and how to associated a stable homotopy type to such a category given certain auxiliary data. A significant amount of additional informal discussion can be found at the end of Section \ref{sec:virtual-cjs}, after the definition of the Floer homotopy type itself. 
Subsequently, in Section \ref{sec:equivariant-flow-categories}, we explain the equivariant analog of the previous constructions. 

The technical Section \ref{sec:flow-category-constructions} is focused on the management of certain linear algebra data which organizes compatibility conditions between framings in an equivariant virtually smooth flow category. Because we wish to work with many possible groups in a compatible way, we cannot use the ordering on the coordinates on $\R^n$, as \cite{cohen2007floer, lipshitz2014khovanov} do to organize the framings. To deal with this we introduce the convenient notions of free and semi-free parameterizations (see Remark \ref{rk:explanation-of-free-parameterizations} for the significance of these apparently trivial notions), and show that stabilization of virtually smooth framed flow categories by such parameterizations lets us construct the additional data that allows us to define the Floer homotopy type associated to such a category. We must also deal with another kind of stabilization having to do with changing the stratifications associated to our virtually smooth flow categories, another phenomenon that only arises when working with irregular Floer data (in the regular case, the stratification is controlled by the Maslov index). 

We review Hamiltonian Floer Homology and fix conventions in Section \ref{sec:hamiltonian-floer-review}. We construct virtually smooth equivariant flow categories associated to Hamiltonian Floer Homology in Section \ref{sec:global-charts}, as well as for the flow categories associated to continuation maps and homotopies of compositions of continuation maps. As part of this construction, we verify all necessary compatibility conditions involving from passing from a Hamiltonian $H$ to its iterate $H^{\#k}$.

Section \ref{sec:index-theory} reviews index theory and its equivariant analog, and framings, compatible with all relevant symmetries, of the virtually smooth flow categories of Section \ref{sec:global-charts}. It also proves Theorem \ref{thm:condition-for-vanishing-of-equivariant-polarization-class}. This section was written in part to avoid deficiencies in the literature around arguments for how to construct framings of flow categories, and in order to verify that framings of our virtually smooth flow categories can be constructed in a fashion compatible with symmetry. 

Finally, Section \ref{sec:homotopy-theory-background} reviews the basics of equivariant orthogonal spectra and some elements of model categories, and Section \ref{sec:floer-homotopy} uses the geometric results established in the earlier sections to prove the theorems stated in this introductory section.  

\begin{remark}
    This paper does not prove that for an admissible Hamiltonian (which has very small positive linear slope at infinity if $M$ is Liouville) that the (possibly $C_k$-equivariant) spectrum $CF_\bullet(H, \SS)$ agrees with $\Sigma^\infty_+ M$ (equipped with the trivial $C_k$ action). Nor does it prove that for a $G$-Morse-Smale function $f$ on a closed $G$-manifold $M$, the Floer homotopy type associated to the flow category of $f$ agrees with the spectrum $\Sigma^\infty_+ M$ (which now has a nontrivial $G$ action. Neither of these are difficult to establish with the methods of this paper. 

    A more significant deficiency is that we do not work out in full detail a proof that the cyclotomic lift of Theorem \ref{thm:sh-is-cyclotomic} is well-defined up to equivalence in the category of genuine $p$-cyclotomic spectra; however, see the end of Section \ref{sec:nikolaus-scholze-comparison} for a sketch of an argument for this invariance.
\end{remark}

\paragraph{Acknowledgements.} This project has been gestating in the author's notes since 2019. I would like to thank Mohammed Abouzaid and Andrew Blumberg for encouragement and discussion throughout this project. I would also like to thank Paul Seidel for his encouraging interest in my earlier attempts to construct the equivariant $p$-cover map \eqref{eq:p-fold-cover-map}; in particular, he first pointed out to me the difficulties with more naive approaches to the definition of this map. I would like to thank Egor Shelukhin for extensive discussions about a Kronheimer-Mrowka variant of Hamiltonian Floer homology, which occurs nowhere in this paper but informed my understanding of the subject. Finally, I would like to thank Kristen Hendricks for pointing out that certain related statements in Lagrangian Floer cohomology have proven difficult to establish without the methods of this paper.

\section{Virtual Equivariant Morse Theory: Elementary Examples}
\label{sec:pedagogical-intro}
In this informal, pedagogical section, we describe the full construction of the genuine equivariant stable homotopy type associated to an virtually smooth equivariant flow category in the simplest nontrivial case: when the flow category has two objects, each of which are fixed by the $G$-action. The general construction presented in this paper is simply a generalization of the construction in this section to the setting of manifolds-with-corners and to flow categories with many objects.

\paragraph{The Cohen-Jones-Segal construction.} 
Let $\mathcal{C}$ be a smooth framed flow category with two objects $x, y$ of degrees $|x| = n > m=|y|$. Thus, $\CC(x,y)$ is a smooth closed manifold of dimension $n-m-1$ with a stable framing, namely, a trivialization of $T\CC(x,y) \oplus \R^k$ for some $k$. Embeddings $\iota$ of $\CC(x,y)$ into $\R^{n-m-1+k'}$ are unique up to isotopy for large $k'$. Morever, composing such an embedding with the embedding $\R^{n-m-1+k'} \subset \R^{n-m-1+k'+(n-m-1+k)}$ gives an embedding $\bar{\iota}$ with a trivial normal bundle: we have 
\begin{equation}
    \label{eq:trivializing-the-normal-bundle}
    \eta_{\bar{\iota}} = \eta_{\iota} \oplus \underline{\R^{n-m-1+k}} \simeq \eta_{\iota} \oplus T\CC(x,y) \oplus \underline{\R^k} \simeq \underline{\R^{n-m-1+k' + k}}.
\end{equation}
Here and henceforth, given a vector space $V$, we write $\underline{V}$ for the vector bundle $V \times B$ (over some implied base $B$ suppresed from the notation) with fiber $V$. The second isomorphism in the above equation \eqref{eq:trivializing-the-normal-bundle} arises from the stable framing of $\CC(x,y)$, and the last isomorphism comes from the original embedding $\iota$ into $\R^{n-m-1+k'}$. 

Writing $K = k' + k$, Pontrjagin-Thom collapse map associates to the stably framed manifold $\CC(x,y)$ well defined stable homotopy classes of maps $[S^{n-m-1+K}, S^K]$ which are preserved by the stabilization functor. In particular, there is an associated map of stable spheres $f: \SS^{n-1} \to \SS^m$. (See the subsequent Section \ref{sec:equivariant-homotopy-review} for a very short introduction to stable homotopy theory.) The Cohen-Jones-Segal construction associates a stable homotopy type $|\CC|$, which in this case is given by $[\bS^{n-1} \to \bS^m] = Cone(f)$. This is a spectrum with a $2$-step filtration, with the lowest filtration level given by $\bS^m$, and the quotient of the next filtration by this one given by $\bS^n$; in other words, it is a stable homotopy analog of a CW complex with cells corresponding to the two objects of the flow category.

\begin{remark}
    Technically, in the formulation of Cohen-Jones-Segal, $\CC(x,y)$ would be a $\langle n -m-1\rangle$-manifold with all boundary faces empty, and we would look for embeddings into $\R^{n-k-1}_+ \times \R^K$ with trivial normal bundle. However, the image of $\CC(x,y)$ under these emnbeddings would miss all boundary faces of $\R^{n-k-1}_+ \times \R^K$, and so the Pontrjagin-Thom collapse map would collapse all boundary faces of $\R^{n-k-1}_+ \times \R^K$ to the base-point, giving the same map $f$ as in our construction. We will be using the language of $\langle k \rangle$-manifolds extensively in our construction in Section \ref{sec:flow-categories}.
\end{remark}

\paragraph{The Virtual CJS construction.}
The first observation is that the CJS construction naturally generalizes to the case of \emph{virtually smooth} framed flow categories. A virtually smooth flow category is, informally, speaking, a flow category equipped with smooth Kuranishi charts with no isotropy groups (Definition \ref{def:flow-category}). A virtually smooth flow category can also be \emph{framed} (Definition \ref{def:framing-of-flow-category}). We will spell out the meaning of these notions completely in the setting of a a virtually smooth framed flow category $(\CC, \CC')$ with two objects $x,y$, and for simplicity we will assume that the \emph{integral action} (Definition \ref{def:integral-action}), an auxiliary piece of combinatorial data, satisfies $\bar{A}(x,y)=1$. This virtually smooth flow category is fully specified by the data of the morphism object $\CC'(x,y)= (T(x,y), V(x,y), \sigma(x,y))$, which is a derived manifold. In other words, $\CC'(x,y)$ consists of the following data:
\begin{itemize}
    \item $T(x,y)$, a smooth manifold, called the \emph{thickening}; 
    \item $V(x,y)$,  a smooth vector bundle over $T(x,y)$ called the \emph{obstruction bundle}, and 
    \item $\sigma(x,y): T(x,y) \to V(x,y)$ a continuous section called the \emph{obstruction section}.
    \item Moreover, the data of $\CC'$ is said to be a \emph{virtual smoothing} of $\CC$ if there is a homeomorphism $\sigma(x,y)^{-1}(0) \simeq \CC(x,y)$, where we assume that $\CC(x,y)$ is a compact Hausdorff space, and if $\dim T(x,y) - \dim V(x,y) = n-m-1$. We say that $\CC'(x,y)$ is a \emph{smooth Kuranishi chart} on $\CC(x,y)$, or that $\CC'(x,y)$ is a \emph{virtual smoothing} of $\CC(x,y)$.
\end{itemize}
A framing in this case is a trivialization of the virtual vector bundle $TT - V$ over $T$, i.e. a pair of isomorphisms $TT \oplus \underline{\R^k} \simeq \underline{\R^{n-m-1+k_1}}, V(x,y) \oplus \underline{\R^k} \simeq \underline{\R^{k_1}}$ for some $k$.

There is a certain equivalence relation on smooth Kuranishi charts called \emph{stabilization}, where one replaces $\CC'(x,y)$ with $\CC'(x,y)_W$ for some vector bundle $\pi_W: W \to T(x,y)$, where $\CC'(x,y)_W = (W, \pi_W^*V(x,y) \oplus \pi_W^*W, \pi_W^*\sigma(x,y) \oplus id_W)$. Two smooth Kuranishi charts that differ by a stabilization, or by the replacement of $T(x,y)$ by an open subset $U \subset T(x,y)$ containing $\sigma(x,y)^{-1}(0)$, should be thought of as equivalent. In particular, after a stabilization by $\underline{\R^k}$, we can arrange such that $\CC'(x,y)$ is framed in the sense that there are isomorphisms $T(x,y) \simeq \underline{\R^{n-m-1+k_1}}, V(x,y) \simeq \underline{\R^{k_1}}$.

Proceeding as in the CJS construction, we can find an embedding of $T(x,y)$ into some $\R^{K'}$ such that the normal bundle to the embedding is trivialized using the framing of of $\CC(x,y)$. Choosing a tubular neighborhood of this normal bundle and using this trivialization, after stabilizing $\CC'(x,y)$ by the normal bundle of this embedding, we can identify $T(x,y)$ as an open subset of $U \subset \R^{K'}$. The framing of $\CC'(x,y)$ and the section $\sigma(x,y)$ thus defines a map $\sigma: U \to \R^{K'-(n-m-1)}$. Using the compactness of $\CC(x,y)$, we can restrict  $\sigma$ to a precompact open neighborhood of $\sigma^{-1}(0)$ and then extend $\sigma$ to a map $\tilde{\sigma}: (\R^{K'})_* \to (\R^{K'-(n-m-1)})_*$ with $\sigma^{-1}(0) = \CC(x,y)$, and with $X_*$ denoting the one-point compactification of a space $X$.

It is clear from this construction that the corresponding stable homotopy class of maps $f: \bS^{n-1} \to \bS^{m}$ is independent of all choices made, and only depends on the original virtual smoothing $\CC'$ of $\CC$ and its framing. The virtual CJS construction defines $|\CC'| = Cone(f)$. 

\begin{figure}[h!]
    \centering
    \includegraphics[width=0.8\textwidth]{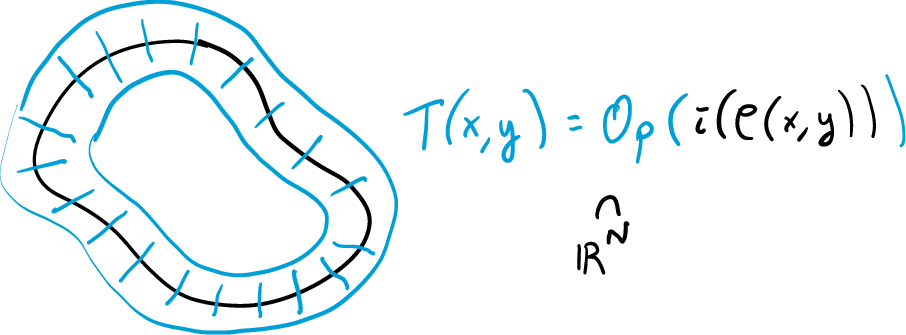}
    \caption{The Pontrjagin-Thom construction associates a map between stable spheres to a cobordism class of framed manifolds. One can view the Pontrjagin-Thom construction as first stabilizing the manifold $\CC(x,y)$, thought of as a Kuranishi space with trivial obstruction bundle, to a Kuranishi space such that the tangent bundle of the thickening is trivial, and then subsequently embedding the thickening into $\R^N$ and subsequently extending the obstruction section (corresponding to the trivialization of the normal bundle of $\CC(x,y)$ in $\R^N$) to a one-point-compactification of the ambient space. The basic observation of this paper is that this construction manifestly continues to make sense when one starts with a general Kuranishi space, and admits an equivariant generalization. This allows us to avoid all equivariant transversality issues arising in Floer homology.}
    \label{fig:virtual-pt}
\end{figure}

We now note that the virtual CJS construction  is a strict extention of the CJS construction. Indeed, a smooth manifold has an associated tautological Kuranishi chart on itself with $V$ a trivial bundle. Embedding this manifold into $\R^{n-m-1+k}$ and noting that the normal bundle is trivial exactly corresponds stabilizing this tautological chart by a trivial vector bundle and subsequently embedding the thickening of the stabilized chart. The map to the sphere in the Pontrjagin-Thom collapse map is, up to homotopy, the same as the map constructed by our prescription for this embedding associated to the stabilized Kuranishi chart. 

The idea that Floer homotopy theory should allow us to work with non-transverse situations is not new. Indeed, the same phenomenon arises in the Bauer-Furuta construction in Seiberg-Witten theory, where, instead of \emph{perturbing} the Seiberg-Witten equation to get a well defined count of zeros, they instead consider the unperturbed equation and note that it gives a well-defined  map of stable spheres. 
\paragraph{The Virtual Equivariant CJS Construction.}
We now explain the construction in the case of framed smooth virtual equivariant flow categories with two $G$-fixed objects. Let $G$ be a compact Lie group. We say that $\CC'$ is $G$-equivariant if 
\begin{itemize}
    \item $G$ acts on the objects preserving their degrees, which means it must fix $x$ and $y$, and 
    \item $G$ acts on $T(x,y)$ and $V(x,y)$, making $V(x,y)$ into an equivariant vector bundle and $\sigma(x,y)$ into an equivariant section, and
    \item finally the homeomorphism $\sigma(x,y)^{-1}(0) \to \CC(x,y)$ is $G$-equivariant.
\end{itemize}
Given a $G$-vector bundle $W$ over $T(x,y)$, we can stabilize $\CC'(x,y)$ by $W$ as before. An \emph{equivariant framing} is, in this case, a trivialization of the virtual equivariant vector bundle $TT(x,y) - V(x,y)$ up to such stabilization; in other words, it is a pair of isomorphisms of equivariant vector bundles 
\[ TT(x,y) \oplus \underline{W} \simeq \underline{W^0}, V(x,y) \oplus \underline{W} \simeq \underline{W^1}\]
for some triple of orthogonal $G$-representations $(W, W^0, W^1)$. Here, as before, $\underline{W}$ denotes the $G$-vector bundle (over an implied $G$-equivariant base $B$) with fiber an orthogonal $G$-representation $W$.  Clearly, one can stabilize a framed equivariant Kuranishi chart to a framed equivariant Kuranishi chart for which the representation $W$ in the above definition of a framing is empty. 

This natural definition of an equivariant framing of a virtually smooth equivariant flow category suggests that the objects of an equivariant flow category $\CC'$, just as the spheres of the equivariant stable homotopy category, should be labeled by a \emph{virtual $G$-representation} rather than by an integer. Indeed, if $\CC$ was the flow category of a $G$-equivariant Morse-Smale function $f:M \to \R$, then every critical point $x \in Crit(f)$ has a negative eigenspace $W_x^-$ of the Hessian of $f$ at $x$, and this negative eigenspace is \emph{preserved} by the $G$-action on $TM_x$. If there was no $G$-action, then the dimension of $W_x^-$ -- its isomorphism class among virtual vector spaces -- would give the index of $x$ in the flow category. Thus, we instead set $|x| = [W_x^-] \in RO(G)$ as a \emph{virtual $G$-representation} the equivariant setting. While these virtual $G$-representations happen to be $G$-representations in the setting of Morse theory, in the context of Floer theory they are produced using equivariant index theory of Fredholm operators (Section \ref{sec:index-theory}) and will be genuinely virtual. 

This modification of the grading theory requires one more modification of the condition on an equivariant framing.The modification we make is that the framings are compatible with the grading: in the case with two objects, we simply require that the equivariant framing satisfies $[W_0] - [W_1] = |x|-|y|-[\R]$ of virtual vector bundles, were $[\R] \in RO(G)$ denotes the trivial $1$-dimensional $G$-representation.

We now use the fact that the equivariant immersion  theory works as before on open manifolds. The theorem needed is a theorem of Wasserman, and is used in Lashof's work on equivariant stable smoothing theory \cite{lashof2006stable}.

\begin{theorem}[Wasserman \cite{wasserman1969equivariant}]
There is a finite-dimensional orthogonal $G$-representation $W'$ such that $T(x,y)$ admits an equivariant embeding into $W'$. 
\end{theorem}

As such, an argument identical to one in the CJS construction shows that after large stabilization by a $G$-representation, we can find a codimension zero equivariant embedding of $T(x,y)$ into a finite-dimensional $G$-representation $W'$, and after this stabilization, $V(x,y)$ is equivariantly isomorphic to $\underline{W''}$ for an orthogonal $G$-representation $W''$. As such, we can view $\sigma(x,y)$ as an equivariant map from an open subset $U$ of $W'$ to $W''$, with $\sigma(x,y)^{-1}(0) = \CC(x,y)$.  Thus, after restricting to a precompact $G$-invariant subset of $U$ containing $\sigma(x,y)^{-1}(0)$, the equivariant Tietze extension theorem shows that this map extends to an equivariant map of \emph{representation spheres} $S^{W'} \to S^{W''}$, which are simply the one point-compactifications of the corresponding $G$-representations. Up to stabilization by $G$-represntations on both sides, this map defines a well-defined map $\tilde{f}_{eq}: \SS^{W'} \to \SS^{W''}$ for a pair of \emph{equivariant stable spheres}. 

Tracing though this process, the compatibility of the equivariant framing of $\CC'(x,y)$ with the grading implies that $[W'] - [W''] = |x| - |y|-[\R] \in RO(G)$. 
With this modification, it is immediate that we suspend $\tilde{f} _eq$ by the stable equivariant sphere $\SS^{|y|-W''}$, we will get a a map $f_eq: \bS^{|x|-1} \to \bS^{|y|}$. 

The genuine equivariant Floer homotopy type of $\CC'$ is the cone on this map $f^{virt}_eq$ in the equivariant stable homotopy category (which we partially describe in the next section). It is elementary to show that this homotopy type is only dependent on the original framed virtually smooth equivariant flow category $(\CC, \CC')$.

\pagebreak

\section{Review of Equivariant (Stable) Homotopy Theory}
\label{sec:equivariant-homotopy-review}
In this section, we give a lightning review of the basic ideas of equivariant homotopy, stable homotopy theory, and their interaction. An excellent review of the basics of the subject can be found in Adams \cite{adams2006prerequisites}, and the first few pages of \cite{may1996equivariant} contain more equivariant homotopy than we need. A clear discussion of a convenient category of orthogonal spectra can be found in \cite{schwede2022orthogonal}. This section is only for motivation and review; throughout, we will be working with orthogonal spectra, which are reviewed in \ref{sec:orthogonal-spectra-detailed-review}. Further references may be found in this latter section.

Let us fix a \emph{finite} group $G$; the notions recalled in this section can be generalized to the case when $G$ is a compact Lie group, but we will not need this generality in this paper. Let $Rep(G)$ denote the set of finite-dimensional real $G$-representations up to isomorphism, and $RO(G)$ the group-completion of the monoid $Rep(G)$ under the monoidal operation of direct sum. 

A \emph{$G$-space} is a space $X$ with an action of $G$ by homeomorphisms. The category $Top_G$ has objects $G$-spaces and morphisms the $G$-equivariant continuous maps. A \emph{$G$-homotopy} between morphisms $f_0, f_1: X \to Y$ in $Top_G$ is a morphism $f: X \times [0,1] \to Y$ in $Top_G$, where the $G$ action on $X \times [0,1]$ is the product of the $G$-action on $X$ and the trivial action on $[0,1]$, such that $f|_{X \times \{i\}} = f_i$ for $i=0,1$. Thus, a $G$-homotopy equivalence $f: X \to Y$ is a morphism such that there is a morphism $g: Y \to X$ with $f \circ g$ and $g \circ f$ both $G$-homotopic to identity maps. A \emph{weak $G$-equivalence} $f: X \to Y$ in $Top_G$ is a morphism such that for all subgroups $H \subset G$, the induced map on $H$-fixed points $f^H: X^H \to Y^H$ is a weak equivalence. 

\begin{definition} A $G$-CW complex structure on a $G$-space $X$ is a CW complex structure on $X$ such that the $G$-actions preserves the skeleta and such that the $n$-skeleton $X_n$ is obtained from the $(n-1)$-skeleton $X_{n-1}$ by attaching $I \times D^n$ along a $G$-equiuvariant map $I \times S^n \to X_{n-1}$, where the $G$-action on $I \times D^n$ is the product of a $G$-action on the set $I$ and the trivial action on $D^n$.
\end{definition}

\begin{definition}
The $G$-homotopy category $hTop_G$ is the category obtained from $Top_G$ by inverting the weak $G$-equivalences. 
\end{definition}

We now state two fundamental results:
\begin{theorem}[Corollary I.3.3 \cite{may1996equivariant}]
A morphism $f: X \to Y$ in $Top_G$ between $G$-CW complexes is a $G$-homotopy equivalence if and only if it is a weak $G$-equivalence. 
\end{theorem}

\begin{theorem}[Chapter I.3 \cite{may1996equivariant}]
    The category $hTop_G$ is equivalent to category with objects the $G$-CW complexes in $Top_G$ and morphisms given by $G$-homotopy-equivalence classes of morphisms in $Top_G$. 
\end{theorem}

Given a pair of spaces $X, Y \in Top_G$, we denote by $[X,Y]_G$ the set of $G$-homotopy-equivalence classes of morphisms in $Top_G$ from $X$ to $Y$. When $G=1$ we drop $G$ from the notation.

The category of based $G$-spaces $Top_{G, *}$ is the category of pairs $(X, b_x)$, with $X \in Top_G$ and $b_x \in X$, with morphisms $Top_{G, *}((X, b_x), (Y, b_y)) = \{f: Top_G(X,Y) : f(b_x) = b_y)\}$. One defines based $G$-homotopies and based $G$-CW complexes in the obvious way. The smash product of based $G$-spaces is the smash product of the underlying spaces with the action on equivalence classes of elements $[(x,y)] \in X \wedge Y = X \times Y / \sim$ given by  $g [(x,y)] = [(gx, gy)]$. We will write $[X,Y]_G$ when $X$ and $Y$ are pointed $G$-spaces to refer to homotopy-equivalence classes of pointed $G$-maps. 

\paragraph{Equivariant Stabilization.} Just as the spheres $S^n$ are distinguished objects in the category of topological spaces from the perspective of homotopy theory, equivariant homotopy theory has the \emph{representation spheres}, which are spheres with $G$-actions induced from $G$-representations.
\begin{definition}
    Let $V$ be an orthogonal $G$-representation. Then the representation sphere $S^V \in Top_G$ is the one-point compactification of $V$ viewed as a $G$-space. 
\end{definition}
By the equivariant triangulation theorem \cite{illman1983equivariant}, $S^V$ admits the structure of a $G$-CW complex. 

Recall now the definition of the \emph{Spanier-Whitehead category}:
\begin{definition}
\label{def:sw-category}
    The \emph{Spanier-Whitehead category} $SW$ is the category with objects $\Sigma^{-k} \Sigma^\infty X$, for $X$ a pointed CW complex and $k \geq 0$ a nonnegative integer, with morphisms defined as
    \[ SW(\Sigma^{-k} \Sigma^\infty X,\Sigma^{-\ell} \Sigma^\infty Y) = \colim_n [\Sigma^{n-k} X, \Sigma^{n-\ell} Y]\]
    where the colimit is over all natural $n \geq k + \ell$ and the maps between the terms in the diagram are induced by the suspension functor. 
\end{definition}
In the Spanier-Whitehead category, the suspension functor $\Sigma: SW \to SW$, which sends $\Sigma^k \Sigma^\infty X$ to $\Sigma^k \Sigma^\infty \Sigma X$, has an inverse functor $\Sigma^{-1}: SW \to SW$, which sends $\Sigma^k \Sigma^\infty X \to \Sigma^{k-1} \Sigma^\infty X$. Moreover, $SW$ is monoidal, with $\Sigma^k \sigma^\infty X \wedge \Sigma^\ell \Sigma^\infty Y = \Sigma^{k+\ell} \Sigma^\infty (X \wedge Y)$, with $\mathbb{S}^0 = \Sigma^0 \Sigma^\infty S^0$ the unit of the monoidal structure. This makes the Spanier-Whitehead category into a symmetric monoidal category with duals, with $\mathbb{S}^{-n} = \Sigma^{-n} \Sigma^\infty S^0$ a monoidal inverse of $\mathbb{S}^n = \Sigma^0 \Sigma^\infty S^n$ for every $n > 0$. 

An analogous notion, the $G$-Spanier-Whitehead category, exists in the realm of equivariant homotopy theory. Before giving one, we recall a useful notion:
\begin{definition}
    A complete $G$-universe $\mathcal{U}$ is a countably-infinite-dimensional real orthogonal $G$-representation containing countably many copies of the regular representation $\mathbb{R}[G]$ of $G$. 
\end{definition}
Any pair of complete $G$-universes are isomorphic as orthogonal $G$-representations, and a canonical choice of $G$-universe given by $\mathbb{R}[G] \tensor_\R \R^{\oplus \infty}$. Note that a related orthogonal $G$-representation, $\mathbb{R}[G] \tensor_\R \mathcal{H}$ where $\mathcal{H}$ is a separable Hilbert space, is not an example of a complete $G$-universe, but instead arises as the completion of a complete $G$-universe under a natural topology. We now proceed with the definition of the equivariant Spanier-Whitehead category:

\begin{definition}
\label{def:g-sw-category}
    Fix a complete $G$-universe $\mathcal{U}$.
    The $G$-Spanier-Whitehead category $SW_G$ is the category with objects $\Sigma^{-V} \Sigma^\infty X$ for $X$ a pointed $G$-CW complex and $V$ a finite-dimensional real orthogonal $G$-subrepresentation of $\mathcal{U}$. The morphisms are defined by
    \[ SW_G(\Sigma^{-V} \Sigma^\infty X, \Sigma^{-W} \Sigma^\infty Y) = \colim_{V \subset \mathcal{U}} [\Sigma^{V'\ominus V}X, \Sigma^{V'\ominus W} Y]-G\]
    where 
    \begin{itemize}
        \item the colimit is taken over the category with objects 
    \[\{ V' \subset \mathcal{U} \text{ a subrepresentation }\}\] such that $V'$ contains $V + W$ (note that $V \cap W$ may be empty)
    and morphisms are given by inclusions of subrepresentations of $\mathcal{U}$; and 
    \item The notation $V'\ominus W$ denotes $W^\perp \subset V'$, and 
    \item The map in the diagram associated to an inclusion $(V', \iota) \to (V'', \iota')$ is induced by suspension by $S^{V'' \ominus V'}$ where $V'' \ominus V'$ is the orthogonal complement to $V'$ in $V''$. 
    \end{itemize}
\end{definition}
One can show the following properties of the $G$-Spanier Whitehead category \cite{schwede2022orthogonal}:
\begin{itemize}
    \item The category $SW_G$ is independent of the choice of complete $G$-universe $\mathcal{U}$ up to equivalence of categories. 
    \item If $V \simeq W$ as orthogonal $G$-representations, there is a canonical isomorphism  $\Sigma^{-V} \Sigma^\infty X \simeq \Sigma^{-W} \Sigma^\infty X$ in $SW_G$;
    \item For any finite-dimensional orthogonal $G$-representation $V$, the functor $\Sigma^V: SW_G \to SW_G$ which acts on objects as $\Sigma^V(\Sigma^{-W}\Sigma^\infty X) = \Sigma^{-W}\Sigma^V \Sigma^\infty X$ has an inverse functor $\Sigma^{-V}$ which acts on objects as $\Sigma^{-V}(\Sigma^{-W} \Sigma^\infty X) = \Sigma^{-V \oplus W}\Sigma^\infty X$; and finally
    \item The operation $\Sigma^{-V} \Sigma^\infty X \wedge \Sigma^{-W} \Sigma^\infty Y = \Sigma^{-V \oplus W} \Sigma^\infty(X \wedge Y)$ makes $SW_G$ into a symmetric monoidal category with duals, with a monoidal inverse to $\mathbb{S}^V = \Sigma^0 \Sigma^\infty S^V$ given by $\mathbb{S}^{-V} = \Sigma^{-V} \Sigma^\infty S^0$. 
\end{itemize}

\paragraph{(Equivariant) Stable Homotopy categories}
The category $SW$ embeds as a full subcategory into the \emph{stable homotopy category}, which we denote by $hSpectra$. Helpful references for the construction and properties of this category can be found in \cite{schwede2022orthogonal}; it is, loosely speaking, a homotopy-ocompletion of the Spanier-Whitehead category $SW_G$. One convenient construction of $hSpectra$ is the symmetric monoidal category of orthogonal spectra $Sp^O$ \cite{mandell2002equivariant}, which is equipped with a set of weak equivalences such that the homotopy category of $Sp^O$ defined using these weak equivalences is $hSpectra$. 

Similarly, in equivariant stable homotopy theory, the Spanier-Whitehead category embeds fully faithfully into $hSpectra_G$, which can be defined as a localization of $Sp^O_G$, the category of \emph{orthogonal $G$-spectra}. For the discussions of these notions, as in Definition \ref{def:g-sw-category}, we first fix a complete $G$-universe, which we will suppress from the notation.

First, in our notation for $SW$ and $SW_G$, we will write $\Sigma^\infty X$ for $\Sigma^0 \Sigma^\infty X$. There a functor 
\[ \Sigma^\infty: hTop_G \to hSpectra_G\]
which, on finite $G$-CW complexes $X$, extends the functor $X \to \Sigma^\infty X$ (in the notation of Definition \ref{def:sw-category}) to $SW_G$ by composing with the inclusion of $SW_G$ into $hSpectra_G$. Objects in the image of $\Sigma^\infty$ are called \emph{suspension spectra}.

Second, there are forgetful functors 
\[ hSpectra_G \to hSpectra_H\]
for any subgroup $H \subset G$ which ``only remember the $G$-action''; in particular, for any $X \in hTop_G$, they send $\Sigma^{-V} \Sigma^\infty X$ to $\Sigma^{-V} \Sigma^\infty X$, where now $X$ is seen as an $H$-CW complex and $V$ is seen as an $H$-representation.

Third, there is a functor called the \emph{geometric fixed point functor} 
\begin{equation}
    (\,\cdot\,)^{\Phi G}: hSpectra_G \to hSpectra
\end{equation}
which behaves on suspension spectra as 
\begin{equation}
\label{eq:defining-property-geometric-fixed-points}
    (\Sigma^\infty X)^{\Phi G} = \Sigma^\infty (X^G),
\end{equation}
is monoidal in the sense that $(X \wedge Y)^{\Phi G} = X^{\Phi G} \wedge Y^{\Phi G}$, and commutes with homotopy colimits (see \cite{schwede2022orthogonal}). This functor is one of three standard fixed point functors on the stable homotopy category of which the other two are the \emph{categorical fixed points} $X^G$ and the \emph{homotopy fixed points} $X^{hG}$. While property \eqref{eq:defining-property-geometric-fixed-points} may seem exceedingly natural, applying \eqref{eq:defining-property-geometric-fixed-points} gives, for any orthogonal $G$-representation $V$, that
\begin{equation}
    (\mathbb{S}^V)^{\Phi G} = (\Sigma^\infty S^V)^{\Phi G} = \Sigma^\infty S^{V^G} = \mathbb{S}^{V^G}
\end{equation}
which combined with the monoidal property implies that 
\begin{equation}
\label{eq:geometric-fixed-points-negative-representation-sphere}
    (\mathbb{S}^{-V})^{\Phi G} = \mathbb{S}^{-V^G}.
\end{equation}
In particular, if we think of $\dim \mathbb{S}^V = \dim V$, then the geometric fixed point functor does not have to decrease the dimension of a representation sphere. 

\begin{remark}
    All the equivariant spectra constructed by geometric methods in this paper will lie in the some category of equivariant Spanier-Whitehead spectra, although we will leave these categories in Section \ref{sec:cyclotomic-structure-cofibrancy}. Thus, we will be dealing with reasonably concrete geometric objects throughout.

    As a Floer theorist, one should think of $\mathbb{S}^V$, with $V$ a virtual $G$-representation as the cell of the equivariant Floer homotopy spectrum associated to a critical point $x$ of the Floer action functional with $V$ the given by the difference between the positive eigenspaces and the negative eigenspaces of the Hessian of the Floer action functional at $x$, thought of as  $G$-representations. One can make rigorous sense of this notion using spectral flow, as we do in Section \ref{sec:producing-framings}.
\end{remark}

\section{Flow Categories: Definitions}
\label{sec:flow-categories}
In this section, we go into more detail about the theory of flow categories. 

\paragraph{Notational conventions}
This paper contains a relatively elaborate construction and so is somewhat burdened by notation. We will sometimes use multiple notations to denote the same object in order to improve readability. Below we state several general notational conventions used in the paper.

Given a category $\CC$, we write $x \in \CC$ for an object of $x$. We write $x \geq y$ if the morphism object $\CC(x,y)$ is not the zero object, and $x > y$ if $x \geq y$ and $x \neq y$. Whenever this notation is used in the paper, the objects of the category being referred to assemble themselves into a poset (which we also denote by $\CC)$ under this relation.

Throughout this paper, given an object with two induces $\alpha_{xy}$ we will sometimes write $\alpha_{x,y}$ for clarity; these denote the same notion. 

All finite dimensional inner product spaces will be real vector spaces equipped with a nondegenerate positive definite inner product, unless otherwise specified.

We will write $Op(S)$ for `an open neighborhood of $S$'. 

Finally, we will sometimes consider a datum $\{a_w\}_{w \in W}$ where the $a_w$ are objects in some category and $W$ is some indexing set. When the indexing set is clear from the context we will suppress it from the notation and write $\{a_w\}$ to refer to $\{a_w\}_{w \in W}$. 

\subsection{Stratified spaces}
\label{sec:stratified-spaces}
First, we introduce a combinatorial language for discussing the kinds of stratifications on spaces and manifolds that we will need to manipulate in order to generalize the virtual CJS construction to manifolds with corners. Definitions essentially equivalent to the notions discussed in this section arise in \cite{cohen2007floer} as well as in \cite{laures2000cobordism} and \cite{albin2011resolution}.

Recall that a \emph{stratification of a space $X$ by a poset $\mathcal{P}$} is an assignment of a closed subset $X(p)$ for each $p \in \mathcal{P}$ such that $X(p) \subset X(q)$ if $p \leq q$. 

Let $[n] = \{1, \ldots, n\}$ for $n \geq 1$, and let $[0] = \emptyset$.

\begin{definition}
\label{def:k-space}
A \emph{$\langle k \rangle$-space} is a space $X$ stratified by the power set $2^{[k]}$ ordered by inclusion, such that 
\begin{itemize}
\item Writing $X_j = X([k] \setminus j)$, we have 
\begin{equation}
\label{eq:boundary-defining-hypersurfaces-equation}
X(T) = \cap_{j \notin T} X_j,
\end{equation}
and 
\item The top open stratum $X \setminus \cup_j X_j$ is nonempty.
\end{itemize}
\end{definition}

With the definition above, the empty space is not a $\langle k \rangle$-space for any $k$.
\begin{definition}
We define the unique $\langle \bullet \rangle$-space to be the empty space $X = \emptyset.$
\end{definition}

The positive orthants $\R^k_+ = \{x \in \R^k | x_i \geq 0, i = 1, \ldots, k\}$ are a collection of $\langle k \rangle$-spaces with 
\[\R^k_+(S) = \{x \in \R^k_+ | x_i = 0 \text{ if } i \notin S\}.\]
They are canonically smooth manifolds with corners, and form the model spaces for topological and smooth $\langle k \rangle$-manifolds. Before we define these, we first introduce some convenient notation related to $\langle k \rangle$-spaces.  

Let $X$ be a $\langle k \rangle$-space. 
For $x \in X$ we write $S(x) = \min \{S \in 2^{[k]} \, | x \in X(S)\}$.
Given $S \in 2^{[k]}$ write $[k] \setminus S = \{r'_1, \ldots, r'_\ell\}$. Then there is a unique isomorphism of posets 
    \begin{equation}
    	r_S: \{T \subset [k]\; | \: T \supset S\} \to 2^{[k \setminus |S|]} \text{ such that } r_S(S \cup \{s'_j\}) = \{j\}.
    \end{equation}

\begin{definition}
       A \emph{$n$-dimensional $\langle k \rangle$-manifold chart} (centered at $x \in X$) is, for some $S \subset [k]$, a homeomorphism
    \begin{equation}
    	\psi: \R^{[k] \setminus S}_+  \times \R^{n - |k \setminus S|} \supset U \to V \subset X
    	\end{equation}
    from an open subset $U \ni (0, 0)$ to an open subset $V$ such that
    \begin{itemize} 	
    	\item For $y \in V$, we have that $S(y) \geq S(\psi(0,0))$, and moreover
    	\item For $[k] \supset T \supset S$, have that $\psi^{-1}(X(T)) \subset \R^{[k] \setminus S}_+(r_S(T)) \times \R^{n - |k \setminus S|}$; and finally
    	\item we say that the chart $\psi$ is centered at $x$ if $\psi(0, 0) = x$. 
    \end{itemize}
\end{definition}

\begin{definition}
\label{def:topological-k-manifold}
        A  \emph{topological $\langle k \rangle$-manifold} of dimension $n$ is a second-countable Hausdorff space $X$ stratified by $2^{[k]}$ such that for every $x \in X$, there exists an $n$-dimensional $\langle k \rangle$-manifold chart centered at $x$, and such that the interior of $X([k]) \setminus \cup_{j \in [k]} X_j$ is nonempty. A \emph{smooth} $\langle k \rangle$-manifold is a topological $\langle k \rangle$-manifold equipped with a maximal subset (an \emph{atlas})of $n$-dimensional $\langle k \rangle$-charts such that transition maps between charts are smooth maps between open subsets of $\R^{\ell}_+  \times\R^{n-\ell}$.
\end{definition}

\begin{remark}
\label{rk:restratify-1}
    In particular, a $\langle 0 \rangle$-manifold is just a manifold, and a $\langle 1 \rangle$-manifold is a manifold with boundary. Just as in those settings, the property of being a topological $\langle k \rangle$-manifold is a \emph{condition} on a $\langle k \rangle$-space, while a smooth $\langle k \rangle$-manifold is a topological $\langle k \rangle$-manifold with \emph{additional structure}.
\end{remark}

\begin{remark}
    We will not state the dimension $n$ of a $\langle k \rangle$ manifold when it is implied. We declare a countable disjoint union of topological $\langle k \rangle$-manifolds of varying dimension to be a topological $\langle k \rangle$-manifold. 
\end{remark}

\begin{remark}
\label{rk:stratification-language-1}
Definitions \ref{def:k-space} and \ref{def:topological-k-manifold} all have counterparts ($\langle S\rangle$-spaces, $\langle S\rangle$-manifolds) where one uses the power set $2^{S}$ to stratify a manifold, for some arbitrary finite set $S$. However, when the set $S$ is not specified, we will use the term $\langle k \rangle$-manifold to refer to any one of these objects; thus we may say that $\R_+^1, \R_+^2, \ldots$  are all $\langle k \rangle$-manifolds. When $k$ is being used to refer to a specific quantity, this should be clear from the context.  
\end{remark}
\begin{remark}
\label{rk:stratification-language-2}
In keeping with typical notation in differential topology,  even though an $\langle S \rangle$-space is a really a pair $(S, X)$, we will often simply refer to the $\langle S \rangle$-space $(S,X)$ by $X$ while suppressing $S$ from the notation. When using this notation, we will use $S_X$ to refer to the set controlling the stratification of $X$; thus, $X$ is an $\langle S_X\rangle$-space. In keeping with this notation, we will write $S_\emptyset=\bullet$.
\end{remark}

\begin{remark}
    The definition of a topological $\langle k \rangle$-manifold $X$ forces the underlying stratified space to satisfy condition \eqref{eq:boundary-defining-hypersurfaces-equation}. The $X_j$ are referred to as the \emph{faces} of $X$ in \cite{cohen1995floer}, and as boundary hypesurfaces in \cite{albin2011resolution}. Any subspace of a topological $\langle k \rangle$-manifold with the induced stratification will still satisfy conditoin \eqref{eq:boundary-defining-hypersurfaces-equation}; thus, while a previous usage of the term $\langle k \rangle$-space \cite{rezchikov2022integral} does not have condition \eqref{eq:boundary-defining-hypersurfaces-equation}, the usage of the term is consistent with \cite{rezchikov2022integral}.
\end{remark}

\begin{remark}A smooth manifold with corners may be a $\langle k \rangle$-manifold in many different ways. For example, while a closed manifold with boundary $X$ is canonically a $\langle 1 \rangle$-manifold, it can also be thought of as a $\langle 2 \rangle$-manifold in $2$ different, by setting either $X(\{1\}) = \emptyset$ or $X(\{2\}) = \emptyset.$ We will use this freedom of stratification later in the paper.
\end{remark}

\subsection{Derived $\langle k \rangle$-manifolds.}

As discussed in the introduction, when studying equivariant Morse- and Floer-theory, we are forced to deal with moduli spaces of flows which are not cut out transversely, and are thus no longer manifolds. To deal with this issue, we have to work with \emph{virtually smooth flow categories}. Just as smooth flow categories are flow categories enriched in $\langle k \rangle$-manifolds, virtually smooth flow categories are flow categories enriched in \emph{derived $\langle k \rangle$-manifolds}, which we now describe.

\begin{definition}
    A \emph{derived $\langle k \rangle$-manifold} $\mathcal{T}$  is a triple $(T, V, \sigma)$ where 
    \begin{itemize}
    \item $T$ is a smooth $\langle k \rangle$-manifold,
    \item $V$ is a vector bundle over $T$, and 
    $\sigma: T \to V$ is a continuous section.
    \end{itemize}
    The \emph{underlying $\langle k \rangle$-space} of $\mathcal{T}$ is $\sigma^{-1}(0) \subset T$ with the induced stratification. The space $T$ is the \emph{thickening} of $\sigma^{-1}(0)$, with $V$ the \emph{obstruction bundle} and $\sigma$ the \emph{Kuranishi section} or the \emph{obstruction section}. A \emph{virtual smoothing} of a $\langle k \rangle$-space $X$ is a stratum-preserving homeomorphism from the underlying $\langle k \rangle$-space of a derived $\langle k \rangle$-manifold to $X$. Moreover, we require that every connected component of $T(S)$ contains a point in $\sigma^{-1}(0) \cap T(S)$.

    A \emph{topological derived $\langle k \rangle$-manifold} has the same definition but without the smooth structure on $T$. Similarly, derived $\langle k \rangle$ is the same notion but $T$ is only required to be a $\langle k \rangle$-space. 
\end{definition}

\begin{remark}
We will continue to use the linguistic conventions of Remarks \ref{rk:stratification-language-1} and \ref{rk:stratification-language-2} with derived $\langle k \rangle$-manifolds. 
\end{remark}

When a $\langle k \rangle$-space $X$ has a virtual smoothing $\mathcal{T}$, it naturally comes with a series of other virtual smoothings, the \emph{stabilizations} of $\mathcal{T}$, which should be thought of as equivalent to $\mathcal{T}$.
\begin{definition}
\label{def:stabilize-derived-manifold}
    A \emph{stabilization} of a derived $\langle k \rangle$ manifold $\mathcal{T} = (T, V, \sigma)$ is any derived $\langle k \rangle$-manifold $\mathcal{T}_W = (W, \pi^*_V \oplus \pi^*W, \pi^*\sigma \oplus id)$ associated to a vector bundle $\pi: W \to T$, where $W$ is given the $\langle k \rangle$-manifold structure $W(S') = \pi^{-1}(T(S'))$ for $S' \subset S_T$, and $id$ denotes the tautological section $W \to \pi^*W$.
\end{definition}

\begin{remark}
    Replacing a virtual smoothing $\mathcal{T} = (T, V, \sigma)$ with $\mathcal{T'} = (T', V_{T'}, \sigma|_{T'})$ for $T' \subset T$ an open subset of $\sigma^{-1}(0)$ should also be thought of as an equivalence, and we may perform such replacements during some constructions without explicitly pointing them out. 
\end{remark}

\subsection{Virtually smooth flow categories}
\label{sec:flow-categories-the-def}

\begin{definition}
\label{def:kspc'}
Let $\langle k \rangle-Spc'$ be the category with objects all $\langle S\rangle$-spaces (including the empty space, which is a $\langle \bullet \rangle$-space), with morphisms from $(S_X, X)$ to $(S_Y, Y)$ given by pairs $(f, \bar{f})$, where 
\begin{itemize}
\item if $S_X = \bullet$ or $S_Y = \bullet$ then there are no morphisms, unless $S_X = S_Y = \bullet$ in which case there is the identity morphism. Otherwise,
\item $f: S_X \to S_Y$ is an injective map of sets, and 
\item $\bar{f}: X \to Y$ is a map of spaces which is a homeomorphism onto its image $Y(f(S_X))$. Moreover, $\bar{f}$ induces homeomorphisms $X(S') \to Y(f(S'))$ for $S' \subset S_X$. 
\end{itemize}
\end{definition}

The category $\langle k \rangle-Spc'$ is symmetric monoidal with product $\times$ on nonempty $\langle k \rangle$-spaces given by 
\[ ((S_X, X), (S_Y, Y)) \mapsto ((S_X \cup S_Y), X \times Y))\]
(with the stratification given by the product stratification), extended to include the empty $\langle k \rangle$-space via $\emptyset \times X = X \times \emptyset = \emptyset$ in the natural way. The unit object of this monoidal structure is the $\langle 0\rangle$-space with one point.

\begin{definition}
\label{def:flow-category}
    A  flow category $\CC$ is a category enriched in $\langle k \rangle-Spc'$ such that the endomorpisms of any object are the unit object in $\langle k \rangle-Spc'$, and such that the compositions \emph{cover the boundaries of the morphism spaces}, namely that any morphism in $\partial \CC(x, y)$, for any pair of distinct objects $(x, y)$ of $\CC$, is the composition of some pair of morphisms in $\CC$. Here, we say that the elements of the $\langle k \rangle$-space $\CC(x,y)$ are the morphisms of $\CC$, and the compositions of morphisms are given by the maps $\bar{f}$ in Definition \ref{def:kspc'}. 
\end{definition}

\begin{definition}
    A \emph{grading} on a flow category $\CC$ is a map $\mu: Ob(\CC) \to \Z$. Or pairs of objects $(x,y)$ of $\CC$, we define
    \begin{equation}
        \mu(x,y) = \mu(x) - \mu(y) - 1.
    \end{equation}
\end{definition}

\begin{definition}
\label{def:integral-action}
    An \emph{integral action} on a flow category $\CC$ is a map $\bar{A}: Ob(\CC) \to \{ n \in \Z : n \geq -1\}$ 
    such that for every pair of objects $(x,y)$ of $\CC$ such that $\CC(x, y) \neq \bullet$ and $x \neq y$, we have
    \[|S_{\CC(x,y)}| = \bar{A}(x) - \bar{A}(y) - 1;\]
    together with isomorphisms 
    \[ \beta_{x,y}: S_{\CC(x,y)} \to [\bar{A}(x) - \bar{A}(y)-1]\]
    such that for the composition map $(f, \bar{f})$ where 
    \[f: S_{\CC(x,y)} \cup S_{\CC(y,z)} \to S_{\CC(x,z)}\]
    one has 
    \[ \beta_{x,z}  \circ f \circ \beta_{x,y}^{-1}(a) = a\]
    and 
    \begin{equation}
    \label{eq:standard-inclusion}
    \beta_{x,z}  \circ f \circ \beta_{y,z}^{-1}(a) = a+\bar{A}(x) -\bar{A}(y).
    \end{equation}
    More invariantly, we can think of the integral action as a choice of lift of the sets $S_{\CC(x,y)}$ to totally ordered sets such that the maps $f: S_{\CC(x,y)} \cup S_{\CC(y,z)} \to S_{\CC(x,z)}$ associated to composition in $\CC$ are order-preserving, where we give the domain of $f$ the ordering where every element of the first term in the disjoint union is less than every element of the second term.
    
    In particular, there are only morphisms from $x$ to $y$ if $\bar{A}(x) > \bar{A}(y)$, since these conditions force $S_{\CC(x, y)}$ to be equal to $\bullet$ if $\bar{A}(x) \leq \bar{A}(y)$.

    When $\bar{A}$ is an integral action, we write 
    \begin{equation}
        \bar{A}(x,y) = \bar{A}(x) - \bar{A}(y) - 1.
    \end{equation}
\end{definition}

\begin{remark}
    The integral action is a helpful device for organizing the combinatorics of flow categories which arises naturally in existing constructions. Later in the paper, we will be forced to modify flow categories by \emph{changing their integral actions}, modifying the stratifications of the morphism spaces $\CC(x, y)$ from $\langle k \rangle$-manifold to $\langle k + \ell\rangle$-manifolds, as in Remark \ref{rk:restratify-1}.
\end{remark}

\begin{remark}
    The reason for the condition on the range of $\bar{A}$ because of the need to fix such a convention when working with the Floer homotopy type. Note that if $y \in \CC$ achieves the minimal value $\bar{A}(y) = -1$ for the integral action then $\bar{A}(x,y) = \bar{A}(x)$; this is the reason for using this particular convention.

    If $\bar{A}$ is an integral action then so is $\bar{A} + 1$; this operation will turn out not to change the Floer homotopy type up to a canonical weak equivalence. 
\end{remark}

\begin{remark}
\label{rk:integral-actions-and-ordering}
    Without loss of generality, we can think of the nonempty $\langle k \rangle$-spaces in a flow category with an integral action as being as satisfying $S_{\CC(x,y)} = [\bar{A}(x,y)]$, with all maps $\beta_{x,y}$ of Definition \ref{def:integral-action} being identity maps. The existence of an integral action it being possible to ehnance the $S_{\CC(x,y)}$ to \emph{totally ordered} sets, such that the maps of Definition \ref{def:kspc'} become maps of \emph{ordered} sets, where we define the monoidal structure on ordered sets to be 
    \[(S, T) \mapsto S \cup T, \text{ where } s \leq t \text{ for all } s \in S, t \in T.\]
\end{remark}

\begin{definition}
\label{def:integral-action-set}
    Given an integral action $\bar{A}$ on $\CC$, we can identify the sets $S_{\CC(x,y)} \simeq [\bar{A}(x)]$ with subsets of the set
    \begin{equation}
        \label{eq:global-integral-action-set}
        S_{\bar{A}} = \{0, 1, \ldots, \max_{x \in \CC} \bar{A}(x)-1\}
    \end{equation}
    via 
    \[n \mapsto \beta_{xy}(n)+\bar{A}(y).\] 
    The set $S_{\bar{A}}$ will be helpful later when organizing certain framing data and reasononing about operations which change the stratification on $\CC$.

    Note that the set 
    \begin{equation}
        \label{eq:excess-set}
        S_{\bar{A}}^- = \{0, 1, \ldots, \min_{x \in \CC} \bar{A}(x)\}
    \end{equation} 
    is nonempty whenever $\bar{A}$ never takes on the value $-1$; this set preciely the complement of all the images of all of the $S_{\CC(x,y)}$.
\end{definition}

\begin{remark}
    When $\CC$ is equipped with an integral action, for the general discussion regarding operations on flow categories and Floer homotopy types, there is no loss in identifying elements of $S_{\CC}(x,y)$ with their images under $\beta_{xy}$; we will do so freely in the discussion. 
\end{remark}

\begin{definition}
    \label{def:poset}
    The objects of $\CC$ form the objects of a poset, with the partial order defined by the property that $x > y$ if $\CC(x,y) \neq \emptyset$, $x \neq y$. We will say that a  pair $(x,y)$, $x > y$, is \emph{primitive} if there are no objects $z$ such that $x > y > z$; in other words, $(x,y)$ is primitive if $G(x,y) = \emptyset$, where 
    \begin{equation}
            G(x,y) = \left\{ (x_0, x_1) \in \CC \times \CC\, | \, \CC(x, x_0) \neq \emptyset, \CC(x_0, x_1) \neq \emptyset, \CC(x_1, y) \neq \emptyset \right\} \setminus \{(x,y)\}.
    \end{equation}
A \emph{primitive path} from $x$ to $y$ is a sequence of elements $x_0 = x, \ldots, x_r = y$ such that $x_i > x_{i+1}$ for $i =1, \ldots, r$, and each pair $(x_i, x_{i+1})$ is primitive.
\end{definition}

\begin{definition}
The category $Man'$ is the category with objects $\langle k \rangle$-manifolds, and morphims given by morphisms $(f, \bar{f})$ of the corresponding objects in $\langle k \rangle-Spc'$ such that $\bar{f}$ is a smooth map of manifolds with corners.
\end{definition}

\begin{definition}
    A \emph{smooth flow category} is a pair $(\CC, \CC')$, where $\CC$ is a flow category with a grading $\mu$ and an integral action $\bar{A}$, together with a category $\CC'$ enriched in $Man'$ which gives $\CC$ upon applying the forgetful functor from $Man'$ to $\langle k \rangle-Spc'$. Moreover, we must have that 
    \[ \dim \CC'(x,y) = \mu(x,y) \]
    for all distinct pairs of objects $(x,y)$ of $\CC$ for which $\CC'(x,y)$ is nonempty.
    
    As in differential topology, when the data of $\CC'$ is clear from the content, we will denote a smooth flow category $(\CC, \CC')$ simply by $\CC$, and think of $\CC(\tilde{x}, \tilde{y})$ as smooth $\langle k \rangle$-manifolds.
\end{definition}

\begin{remark}
    This definition of a smooth flow category is essentially equivalent to the definition given in Cohen \cite{cohen1995floer} when the integral action is \emph{equal to} the grading. The decoupling of the stratification from the grading in the virtual setting is a necessary complication forced by the geometry.
\end{remark}

\begin{remark}
We will now define a category of derived $\langle k \rangle$-manifolds and morphisms given by isomorphisms onto boundary strata of the codomain, in order to make sense of a derived analog of a smooth flow category. It is straightforward to work with topological or smooth analogs of such a definition. However, most naive constructions of compatible families of charts for moduli spaces in Floer homology do now have the property that that given $x > y > z$ objects of a flow category with $T(x_0,x_1)$ denoting the thickening of the moduli space of floer tracjectories $\cM(x,y)$, that 
\[ \partial T(x,z) = \cup_{x > y > z} \text{ in } T(x,y) \times T(y,z).\]
Instead, one has that boundary strata of $T(x,z)$ are unions of \emph{vector bundles} over the products $T(x,y) \times T(y,z)$, and more generally that the boundaries of the derived $\langle k \rangle$-manifold charts are \emph{stabilizations} of the product charts associated to the boundary strata. The reason for this is because unless one works carefully, one needs \emph{more} potential perturbations of the Floer equation to make the linearized operators all Floer trajectories in $\cM(x,z)$ surjective than one gets from the choices previously made for the boundary strata.

However, in Section \ref{sec:global-charts}, where we construct a virtual smoothing of the flow category of Hamiltonian Floer trajectories, we in particular explain how to carefully choose perturbation data so as to avoid these ``stabilizing vector bundles'' that show up in previous constructions \cite{bai2022arnold, rezchikov2022integral}. It is possible to develop a version of a theory of flow categories where such stabilizing vector bundles are \emph{allowed}; however, constructing and manipulating framings, embeddings and extensions of embeddings for such categories (see Definition \ref{def:framing-of-flow-category} and below) is significantly more difficult, and we leave the study of such more general virtually smooth flow categories to future work. 
\end{remark}

\begin{definition}
\label{def:derman'_nostab}
    The category $DerMan'_{sm, no-stab}$ is the category with objects given by derived $\langle k \rangle$-manifolds, and morphisms from the derived $\langle S_\mathcal{T}\rangle$-manifold $\mathcal{T}$ to the derived $\langle S_{\mathcal{T}'}\rangle$-manifold $\mathcal{T}' = (T', V', \sigma')$ to be tuples $\mathfrak{f} = (f, \widetilde{f}) = (f, (\bar{f}, \bar{f}'))$, where 
    \begin{itemize}
        \item $f: S_\mathcal{T} \to S_{\mathcal{T}'}$ is an injective map;
        \item $\bar{f}: T \to T'$ is a diffeomorphism onto its image $T'(f(S_\mathcal{T}))$, such that $\bar{f}$ induces diffeomorphisms $T(S') \to T'(f(S'))$ for each $S' \subset  S_\mathcal{T}$; 
        \item $\bar{f}': V \to V'|_{T'(f(S_\mathcal{T}))}$ is a bundle isomorphism covering $\bar{f}$, such that 
        \[ \bar{f}' \circ \sigma = \sigma' \circ \bar{f}.\]
    \end{itemize}

    We define $DerMan'_{top, no-stab}$ to be the same category where $\mathcal{T}$ is a topological derived $\langle k \rangle$-manifold and $\bar{f}$ is only a homeomorphism. Similarly, $DerSp'$ is the same category where $\mathcal{T}$ is only a derived $\langle k \rangle$-space and $\bar{f}$ a homeomorphism. 
\end{definition}

\begin{definition}
\label{def:virtually-smooth-flow-category}
    A \emph{virtually smooth flow category} (without stabilizations) is a pair $(\CC, \CC')$, where $\CC$ is a flow category with grading $\mu$ and integral action $\bar{A}$, and $\CC'$ is a category enriched in $DerMan'_{sm, no-stab}$ equipped with an isomorphism to $\CC$ upon applying the forgetful functor $DerMan'_{sm, no-stab}\to \langle k \rangle-Spc'$, satisfying two further conditions:
    \begin{itemize}
    \item Compatibility with the grading: 
    \[ \vdim \CC'(x, y) = \mu(x,y)\]
    for all distinct pairs of objects $(x,y)$ of $\CC$; 
    \item Writing $\CC'(x,y) = (T_{\CC(x,y)}, V_{\CC(x,y)}, \sigma_{\CC(x,y)})$, any element of $\partial T_{\CC(x,y)}$ is in the image of the maps $\bar{f}$ (Definition \ref{def:derman'_nostab}) associated to some composition of morphism objects in $\CC'$.
    \end{itemize}

    A \emph{virtually topologically-smooth flow category} is the same definition with $DerMan'_{sm, no-stab}$ replaced with $DerMan'_{top, no-stab}$. A \emph{virtual topological flow category} is the same definition but with $DerMan'_{top, no-stab}$ replaced with $DerSp'$. 
\end{definition}

\begin{remark}
Virtually (topologically-)smooth flow categories together with extra data be used to construct stable homotopy types; in this paper we will explain the construction for virtually smoooth flow categories (equipped with extra data), leaving the topologically smooth case to later work. In contrast, virtual topological flow categories do not carry homotopical meaning, but are simply a convenient method for packaging data produced in the process of producing virtually smooth flow categories (as in the proof of Proposition \ref{prop:cyclotomic-compatibility-all-perturbation-data-choices}). 
\end{remark}

\subsection{Framings of flow categories}
\label{sec:flow-categories-framings}
As discussed in Section \ref{sec:pedagogical-intro}, any derived $\langle k \rangle$-manifold $\mathcal{T} = (T,V, \sigma)$ carries a class in $KO^0(T)$ given by the virtual vector bundle $[TT]-[V]$. Moreover, using index theory \ref{sec:producing-framings} one may be able to arrange for this class  to equal to a class of the form $[\R^{\vdim \mathcal{T}}]$, or in other words (as in Section \ref{sec:pedagogical-intro}) for there to be some vector bundle $W$ such that the stabilization $\mathcal{T}_W = (T_W, V_W, \sigma_W)$ of $\mathcal{T}$ by $W$ has both $TT_W \simeq \underline{\R}^N$ and $V_W \simeq \underline{\R}^{N - \vdim \mathcal{T}}$ for some $N$. In this section, we spell out precisely precisely is a compatible family of such pairs of isomorphisms across all morphisms in a flow category. First however we define some notation that is helpful for discussing framings for $\langle k \rangle$-manifolds.
\paragraph{Tangent bundles of $\langle k \rangle$-manifolds.}

We now recall some basic structures associated to the tangent bundles of $\langle k \rangle$-manifolds. The model $\langle k \rangle$-manifolds $\R^k_+ \times \R^\ell$, have trivial tangent bundles $T(\R^k_+ \times \R^\ell)$ which, when restricted to $(\R^k_+ \times \R^\ell)(S) = \R^k_+(S) \times \R^\ell$, contain the trivial subbundles $T(\R^k_+ \times \R^\ell)(S)$. The perpendiculars to these sub-bundles are canonically trivial: for each $j \in [k] \setminus S$, there is a vector in $T(\R^k_+ \times \R^\ell)(S)^\perp$ that is the \emph{inwards-pointing normal vector} to $(\R^k_+ \times \R^\ell)_j$ (in the notation of Definition \ref{def:k-space}). 

Similar structures exist on the tangent bundles of a general $\langle k \rangle$-manifold $X$, namely, subbundles $TX(S)$ of $TX|_{X(S)}$. By choosing a convenient of Riemannian metric on $X$ (Definition \ref{def:convenient-metric}), one can show that on any $\langle k \rangle$-manifold $X$, there exist decompositions $TX|_{X(S)} = TX(S) \oplus \underline{\R}^{[k] \setminus S}$ which satisfy natural compatibility conditions between the strata; choosing a convenient Riemannian metric on $X$ will determine these decompositions uniquely. 

Given a $\langle k \rangle$-manifold $X$, write $\underline{T\R}^k_+$ to be the trivial $\R^k$ bundle over $X$, together with the distinguished subundles $\underline{\R}^{[k] \setminus S} \subset \underline{T\R}^k_+|_{X(S)}$ for all $S \subset [k]$. The natural notion of a trivialization of the tangent bundle $TX$ of $X$ is a vector bundle isomorphism $\eta: TX \to \underline{T\R}^k_+ \times \underline{\R}^{\dim X - k}$ which induces isomorphisms $TX(S) \to \underline{T\R}^k_+(S) \times \underline{\R}^{\dim X - k}$ for all $S \subset [k]$. Of course, this notion immediately generalizes to $\langle S\rangle$-manifolds for general finite sets $S$.

\paragraph{How to frame a flow category.}
We now define a framing of a virtually smooth flow category:

\begin{definition}
\label{def:framing-of-flow-category}
A \emph{framing} of a virtually smooth flow category without stabilizations $(\mathcal{C}, \mathcal{C}')$ is a sequence, for all pairs $x, y \in \CC$ with $x>y$,   of vector bundle isomorphisms 
\[ \eta_{xy}: TT(x,y) \to \underline{T\R}_+^{\bar{A}(x,y)} \times \underline{V^0(x,y)} \]
\[ \eta'_{xy}: V(x,y) \to \underline{V^1(x,y)}\]

Here the vector spaces $V^0(x,y)$, $V^1(x,y)$, and $W^0(x,y,z)$ are choices of finite-dimensional vector spaces that are part of the data of the framing. The remaining data of a framing are isomorphisms of inner product spaces
 \begin{equation}
    \label{eq:little-linear-iso}
    \begin{gathered}
        \tau^0_{xyz}: V^0(x,y) \oplus V^0(y,z)  \to V^0(x,z), \\
        \tau^1_{xyz}: V^1(x,y) \oplus V^1(y,z) \to V^1(x,z).
        \end{gathered}
    \end{equation}

Moreover, the maps $\{(\eta_{xy}), (\eta'_{xy}), (\eta_{xyz}), (\tau^0_{xyz}), (\tau^1_{xyz})\}$ are required to satisfy the following compatibility conditions. 
\begin{itemize}
    \item The isomorphism $\eta_{xy}$ are trivializations of $TT(x,y)$ as the tangent bundle of a $\langle k \rangle$-manifold; thus, using the isomorphisms $S_{\mathcal{C}(x,y)} \simeq [\bar{A}(x) - \bar{A}(y) -1]$ of Definition \ref{def:integral-action}, we have that $\eta_{xy}$ induces isomorphisms 
    \[TT(x,y)(S) \to \underline{T\R}_+^{\bar{A}(x) - \bar{A}(y) -1}(S')\times \underline{V^0(x,y)} \]
    for all $S' \subset S$;
    \item For all triples of objects $(x,y,z)$ of $\CC$,
    we require that the following diagram commutes
    \begin{equation}
    \label{eq:framing-diagram-1}
    \begin{tikzcd}
       TT(x, y) \times TT(y, z)\ar[r, "d\bar{f}_{xyz}"] \ar[d, "\eta_{xy} \boxplus \eta_{yz}"]     
       & TT(x,z) \ar[d, "\eta_{xz}"]           
           \\
        \underline{\R^{\bar{A}(x,y)}} \oplus \underline{V^0(x,y)} \times \underline{\R^{\bar{A}(y,z)}} \oplus \underline{V^0(y,z)} \times \underline{W^0(x,y,z)} \ar[r, "r"] 
        &
        \underline{\R^{\bar{A}(x,z)}} \times \underline{V^0(x,z)}
        \end{tikzcd}
    \end{equation}
    where the bottom map $r$ is the bundle map covering $\bar{f}_{xyz}$ induced from the composition of the reordering and inclusion map
    \[\R^{\bar{A}(x,y)} \oplus V^0(x,y) \oplus \R^{\bar{A}(y,z)} \to \R^{\bar{A}(x,y)} \oplus \R^{\bar{A}(y,z)}\oplus V^0(x,y)  \subset \R^{\bar{A}(x,z)}\times V^0(x,y),\]
    with  $\tau_{xyz}^0$ \eqref{eq:little-linear-iso}. 
   
    \item Finally, continuing with the notation of the previous item, the diagrams
    \begin{equation}
    \label{eq:framing-diagram-2}
    \begin{tikzcd}
         V(x,y) \boxplus V(y,z)  \ar[r, "\bar{f}'_{xyz}"] \ar[d, "\eta'_{xy} \boxplus \eta'_{yz})"] &  V(x,z) \ar[d, "\eta'_{xz}"] \\
        \underline{V^1(x,y)} \oplus \underline{V^1(y,z)} \oplus \ar[r] & \underline{V^1(x,z)}
    \end{tikzcd}
    \end{equation}
    where the bottom horizontal arrow is the map of trivial bundles covering $\bar{f}_{xyz}$ induced from $\tau^1_{xyz}$ \eqref{eq:little-linear-iso}.
\end{itemize}

\end{definition}

\paragraph{Parameterizations of framings. }
The data of the maps $\{(\tau^0_{xyz}), (\tau^1_{xyz})\}$ associated to a framing of a flow category is essentially `homotopically auxiliary'. We will not make this statement completely precise in this paper. We will refer to the data $(\{V^0(x,y)\}, \{V^1(x,y)\}, \{(\tau^0_{xyz}), (\tau^1_{xyz})\})$ as the \emph{parameterization} of the framing; specifically, it is an $F'_2$-parameterization (Lemma \ref{lemma:framed-flow-categories-give-F_2'-params}). The data of the  framing forces the existence of many compatibility conditions between the maps $\{(\tau^0_{xyz}), (\tau^1_{xyz})\}$, which are spelled out in this lemma.  We discuss these data and the relevant terminology further in Section \ref{sec:embedding-framings}. 

\begin{definition}
    A framed flow category is \emph{semi-freely parameterized} when the associated $F'_2$-parameterization is an $F^s_2$-parameterization. This notions and related notions are defined in Section \ref{sec:embedding-framings}.
\end{definition}

\paragraph{Canonical parameterizations of framings.}

Now, there are canonical maps \begin{equation}
\label{eq:reordering-map}
\begin{gathered}
    r_{a,a',b, b'}: \R^a_+ \times \R^{a'} \times \R_+^b \times \R^{b'} \to \R^{a+b+1} \times \R^{a' + b'} \\
    r_{a,a', b, b'}(x,y,z,t) = ((x,0, z), (y, t)). 
\end{gathered}
\end{equation}

\begin{definition}
    We say that a framing is \emph{canonically parameterized} if all vector spaces $V^0(x,y), V^1(x,y)$ are simply euclidean vector spaces $\R^n$ of the appropriate dimension and the maps in the bottom of the diagrams \eqref{eq:framing-diagram-1}, \eqref{eq:framing-diagram-2} are induced from canonical isomorphisms \eqref{eq:reordering-map}.
\end{definition}

\begin{remark}
\label{rk:why-canonical-parameterizations-are-not-enough}
In previous work on flow categories \cite{cohen2007floer, lipshitz2014khovanov}, one considers framed flow categories without stabilizations for which the framings are canonically parameterized. An $F_2^s$ parameterization is a generalization of a canonical parameterization which allows for a direct sum of a canonical parameterization with a a 'free parameterization' (an $F_2$-parameterization). 

We are forced to introduce more general parameterizations of framings in our definition of a framing of a virtually smooth flow category in order to deal with framings of $G$-equivariant flow-categories. Specifically, there is an analog of a canonically-paramerized $G$-equivariant flow category, but unfortunately, when one views this flow category as an $H$-equivariant flow category for $H \subset G$, its framing will no longer be canonically parameterized. For further discussion on the kinds of parameterizations of framings that we use, see Section \ref{sec:embedding-framings}.
\end{remark}

\paragraph{Embeddings of framings.}
\begin{definition}
\label{def:embedding-of-framing}
    An \emph{embedding} of a framing of a virtually smooth flow category is a sequence of finite-dimensional inner product spaces $V^0(x), V^1(x)$, such that 
  together with isometric  embeddings 
    \[\bar{\alpha}^{0}_{xy}: V^0(x,y) \oplus V^0(y) \to V^0(x), \bar{\alpha}^1_{xy}: V^1(x,y) \oplus V^1(y) \to V^1(x)\]
    such that the following diagram commutes    for $i=0,1$:
    \begin{equation}
    \label{eq:embedding-of-framing-compatibility}
    \begin{tikzcd}
        V^i(x,y) \oplus V^i(y,z) \oplus V^0(z) \ar[d] \ar[r] & V^i(x,y) \oplus V^i(y) \ar[d]\\
        V^i(x,z) \oplus V^i(z) \ar[r] & V^i(x).
    \end{tikzcd}
    \end{equation}
A framing with a chosen embedding is said to be embedded. We write $\alpha^i_{xy}: V^i(x,y) \to V^i(x)$ $, i=0,1$, for the restriction of $\bar{\alpha}^i_{xy}$, and $\hat{\alpha}^i_{xy}$ for the restriction of $\bar{\alpha}^i_{xy}$ to $V^i(y)$. 

We will refer to each of the data $(\{V^i(x)\}, \{\hat{\alpha}^i_{xy}\})$ for $i=0,1$ as \emph{$E'$-parameterizations}, and the union of these data for $i=0$ and $i=1$ as an $E'_2$-parameterization. An embedding of a framing the same as the data of an $E'_2$-parameterization together with an isomorphism of the associated $F'_2$-parameterization with the one underlying the framing. We discuss these notions further in \ref{sec:embedding-framings}.

\end{definition}

\begin{lemma}
    Semi-freely parameterized framings of virtually smooth flow categories can be embedded. 
\end{lemma}
\begin{proof}
    This is a tautology by Definitions \ref{def:free-parameterizations}, \ref{def:canonical-parameterization}, and \ref{def:semi-free-parameterization}.
\end{proof}

\begin{definition}
    An embedded framing of a graded virtually smooth flow category $\CC$ is \emph{adapted}, or \emph{compatible with the grading}, if 
            \begin{equation}
    \label{eq:index-condition-for-embedding-of-framing}
        \dim V^0(x) - \dim V^1(x) = \mu(x) - (\bar{A}(x) + 1).
    \end{equation}
    for all $x \in \CC$. 
\end{definition}

\begin{remark}
    Implicitly, we can use the formula above to \emph{define} $\mu$ given a framed virtually-smooth flow category with an embedded framing; with this choice of $\mu$ the category $\CC$ would be adapted. From this perspective the grading is actually auxiliary data, and the embedding of the framing is a geometric enrichment of the concept of a grading.  
\end{remark}

\subsection{Auxiliary data for flow categories.}
\label{sec:flow-categories-auxiliary-data}
\begin{definition}
    A flow category $\CC$ is \emph{hom-proper} if the morphism spaces $\CC(x,y)$ are compact for all pairs of objects $x,y \in \CC$.
\end{definition}
\begin{definition}
    A flow category $\CC$ is integral action $\bar{A}$ is proper if it is hom-proper and the sets 
    \[ \{ x \in \CC | \bar{A} \leq C\}\]
    are finite for all real numbers $C$.
\end{definition}

\begin{definition}
A flow category $\CC$ is \emph{finite} if it is hom-proper and has a finite number of objects.
\end{definition}

\paragraph{Embeddings of flow categories.}
\begin{definition}
    An \emph{embedding} $\mathbf{E}$ of a finite virtually smooth flow category $(\CC, \CC')$ with an  embedded framing $(\eta_{x,y}, \eta'_{x,y})$ is a sequence of proper codimension-0 embeddings (the \emph{embedding maps})
    \begin{equation}
        \iota_{xy}: T(x,y) \to \R_+^{\bar{A}(x,y)} \times V^0(x,y).
    \end{equation}
    such that the following conditions are satsified. First, the following diagrams commute:
    \begin{equation}
    \label{eq:embeddings-commute}
    \begin{tikzcd}
    T(x,y) \times T(y, z) \ar[d, "\bar{f}_{xyz}"] \ar[r, "\iota_{xy} \times \iota_{yz}"] &  \R_+^{\bar{A}(x,y)} \times V^0(x,y) \times \R_+^{\bar{A}(y,z)} \times V^0(y,z) \ar[d, "r"]\\
    T(x,z) \ar[r, "\iota_{xz}"] &   \R_+^{\bar{A}(x,z)} \times V^0(x,z)
    \end{tikzcd}
    \end{equation}
    where $r$ is the composition of a reordering map and the map $\tau^0_{xyz}$ of \eqref{eq:little-linear-iso}, as in \eqref{eq:framing-diagram-1}. Second, one has that for each pair of objects $x>y$ of $\CC$, one has that 
        $d\iota_{xy}|_{TT(x,t)_j^\perp}$ maps to the inwards pointing normal vector of $(\R^{\bar{A}(x,y)}_+)_j$. 

    Finally, we require that for each $x \in \CC$ and each $a < \bar{A}(x)$, writing 
    \[\{y_1, \ldots, y_r\} = \{ y \in \CC : x > y, \bar{A}(y) = a\} \]
    we have that the images of 
    \[ T(x,y_i) \times V^0(y_i) \text{ in } \R_+^{\bar{A}(x,y)} \times V^0(x)\text{ under } (id \times \bar{\alpha}_{xy_i}) \circ (\iota_{xy_i} \times id) \]
    are disjoint as $i$ ranges over $1, \ldots, r$. 
\end{definition}

\paragraph{One-point compactifications and wedge products.}
If $X$ is a non-compact space without basepoint, its one-point compactification $X_*$ is compact space with a basepoint at infinity. One immediately verifies that for a pair of spaces $(X,Y)$ without basepoint, one has a canonical identification
\begin{equation}
\label{eq:canonical-smash-homeo}
    (X \times Y)_* = X_* \wedge Y_*. 
\end{equation}
(where if $X$ is compact we use $X_*$ to denote the addition of a disjoint basepoint $*$). 

The one-point compactifications of Euclidean spaces are all spheres with a base-point, i.e. 
\begin{equation}
    \R^n_* = S^n,
\end{equation}
and one has the canonical identification 
\begin{equation}
    (\R^1_+ \times X)_* = c(X_*) = \{ [0,1] \times X\}/((0, x) \sim (t, *) \text{ for any pair} (t,x) \in [0,1] \times X);
\end{equation}
here $c(X)$ is the \emph{cone} on the pointed space $X_*$. Thus, we think of $(\R^k_+ \times \R^\ell)_*$ as the \emph{iterated cone} on the sphere $S^\ell$; these spaces arise in the Cohen-Jones-Segal construction \cite{cohen2007floer}.

\paragraph{Extensions of embeddings of flow categories}

\begin{definition}
    \label{def:extension}
    Let us assume now that we have an embedding $(\iota_{xy})$ of a finite virtually smooth flow category $\CC$ an embedded framing.  An \emph{extension}  $\bar{\mathbf{E}}$ of the embedding $\mathbf{E}$ is the data of \emph{extension maps}
\begin{equation}
\label{eq:extension-maps}
    e_{xy}: (\R_+^{\bar{A}(x,y)} \times V^0(x,y))_* \to (V^1(x,y))_*  
\end{equation}
 These maps are required to satisfy five conditions: 
\begin{enumerate}[(A)]
    \item There exists an open subset $U_{xy} \subset \R_+^{\bar{A}(x,y)} \times V^0(x,y)$ such that 
    \[ U_{xy} \supset \overline{V_{xy}} = \overline{Im \iota_{xy}} \]
    which moreover deformation retracts to $U_{xy}$, and such $e_{xy}|_{U_{xy}^c}$ is the constant map to the basepoint; 
    \item For each $x$ and each $A_r < \bar{A}(x)$, we have that the aforementioned sets
    \[\{\overline{U_{xy}}_{x > y, \bar{A}(y) = A_r}\]
    are pairwise disjoint;
    \item On $V_{xy}$ we have that $e_{xy} = \check{\eta}'_{xy}\sigma_{xy} \circ \iota_{xy}^{-1}$; here $\check{\eta'}_{xy}$ is the composition of the framing $\eta'_{xy}$ of the obstruction bundle with the projection to the fiber $V^1(x,y)$.
    \item The maps $e_{xy}$ are compatible with composition. This means, first, that for every $x > y > z$ in $\CC$ we have that the image of $U_{xy} \times U_{yz}$ under the map 
    \[ \R_+^{\bar{A}(x,y)} \times V^0(x,y) \times \R_+^{\bar{A}(y,z)} \times V^0(y,z) \simeq \R_+^{\bar{A}(x,y)} \times \{0\} \times  \R_+^{\bar{A}(y,z)} \times V^0(x,y) \times V^0(y,z) \xrightarrow{i \times \tau^0_{xyz}} \R^{\bar{A}(x,z)} \times V^0(x,z) \]
    lies inside $U_{xz}$. Second, we also have that, viewing $U_{xy} \times U_{yz}$ under the inclusion above, we have that $e_{xz}|_{U_{xy} \times U_{yz}} = \tau^1_{xyz}(e_{xy} \times e_{yz})_*$, where $\tau^1_{xyz}(\cdot)_*$ denotes the homeomorphism
    \[ \tau^1_{xyz}(\cdot)_*: V^1(x,y)_* \vee V^1(y,z)_* \to V^1(x,z)_*\]
    induced by monoidal property of the $1$-point-compactification functor.
     \item The preimage of $0$ under $e_{xy}$ is exactly $\iota_{xy}(\sigma^{-1}_{xy}(0))$. Finally,
    \item  Recall that our conditions on virtually smooth flow categories require that if $\sigma^{-1}_{xy}(0) \cap T(x,y)(S) = \emptyset$ then $T(x,y)(S) = \emptyset$. We require that if $T(x,y)(S) = \emptyset$ then $U_{xy} \cap \R^{\bar{A}(x,y)}_+(S) \times V^0(x,y) = \emptyset$; thus, $e_{xy}$ sends strata not containing any points sent to zero  by $\sigma_{xy}$ \emph{entirely} to the basepoint.  
\end{enumerate}
\end{definition}

\subsection{The Virtual CJS construction.}
\label{sec:virtual-cjs}

Suppose we are given an adapted  freely framed virtually smooth flow category $\CC'$ with an extension of of an embedding, as in Definition \ref{def:extension}. 

Set $V^1_{max}$ to be the shift space associated to the parameterization $\{V(x)\}$ of the framing of the obstruction bundles of $\CC'$ (Definition \ref{def:shift-spaces}); thus $V^1(x)$ is canonically identified as a subspace of $V^1_{max}$, and there is a canonical orthogonal complement $(V^1(z))^\perp$. In particular, we have canonical isomorphisms 
    \begin{equation}
        \label{eq:canonical-iso-on-perps}
        V^1(x, z) \oplus V^1(x)^\perp \simeq V^1(z)^\perp
    \end{equation} 
    for any $x > z$ in $\CC$.

\begin{definition}
\label{def:floer-homotopy-type}    
    The stable homotopy type associated to an extension of an embedding of a virtually smooth flow category $\CC$ is given by 
    \begin{equation}
    \label{eq:floer-homotopy-type-new}
    \begin{gathered}
    |\CC| = \Sigma^{-(V^1_{max})}\Sigma^\infty \{\CC\}, \text{ where } \\
    \{\CC\} = \left(\bigvee_{z \in \CC}(\R_+^{\bar{A}(x)+1} \times V^0(x) \oplus V^1(x)^\perp)_* / \sim\right).
    \end{gathered}
    \end{equation}
    The relations are defined below. Each term $(\R_+^{\bar{A}(x)+1} \times V^0(x) \times V^1(x)^\perp)_*$ is homeomorphic to a disk with base-point on the boundary; the relations be given by inductively attaching these disks to one another along a diagram organized by the poset $\CC$. More explicitly, writing $\{\CC\}_{\mathcal{P}} \subset \{\CC\}$ for the subset which is the image of the part of the wedge product under the quotient by the relations of all $z \in \mathcal{P}$, then if $\mathcal{P}$ is a downwards closed subposet of $\CC$ and $\mathcal{P}_x = \mathcal{P} \cup x$ where there is no $y \in \CC \setminus \mathcal{P}$ such that $x > y > z$ with $z \in \mathcal{P}$, then 
    \begin{equation}
        \label{eq:fundamental-attaching-map}
        \{\CC\}_{\mathcal{P}_x} \text{ is the cone on } \tilde{e}_x: \partial (\R_+^{\bar{A}(x)+1} \times V^0(x) \times V^1(x)^\perp)_* \to\{\CC\}_{\mathcal{P}}. 
    \end{equation} 

    To define the relations we will define these maps. First off, the relations for $\mathcal{P} = M(\CC)$ collapse the boundaries of the corresponding terms in the wedge product; in other words, if $\bar{A}(z) > -1$, we identify $\partial (\R_+^{\bar{A}(x)+1} \times V^0(x) \oplus V^1(x)^\perp)_* \subset (\R_+^{\bar{A}(x)+1} \times V^0(x) \oplus V^1(x)^\perp)_*$ with the basepoint. (In other words, we set the maps $\tilde{e}_z$ to be the maps to the basepoint for all $z \in M(\CC)$, whenever the domain of $\tilde{e}_z$ is nonempty, i.e. when $\bar{A}(z) > -1$. 
    Thus 
    \[ \{\CC\}_{M(\CC)} = \bigvee_{z \in M(\CC)}(\R_+^{\bar{A}(x)+1} \times V^0(x) \oplus V^1(x)^\perp)_*/\partial.\]
    
    Now, note that  $\partial(\R_+^{\bar{A}(x)+1} \times V^0(x) \times V^1(x)^\perp)_*$ can be divided into $\bar{A}(x) +1 $ pieces  $\partial_i(\R_+^{\bar{A}(x)+1} \times V^0(x) \times V^1(x)^\perp )_*$ for $i=1, \ldots, \bar{A}(x)+1$, with the $i$-th piece given by the subset of the boundary where the $i$-th coordinate of the $\R_+^{\bar{A}(x)+1}$ factor is set to zero. We have a canonical homeomorphism  
    \[ \partial_i(\R_+^{\bar{A}(x)+1} \times V^0(x))_* \simeq (\R_+^{i-1} \times V^0(x) \times V^1(x)^\perp) \wedge (\R_+^{\bar{A}(x)+1-i})_*. \] 

    We will define maps 
    \[ \tilde{e}_{x, i}:  \partial_i(\R_+^{\bar{A}(x)+1} \times V^0(x) \times V^1(x)^\perp)_* \to \{\CC\}_{\mathcal{P}_{x>}}\]
    where we write $\mathcal{P}_{x>}$ to be the set of elements of $\CC$ which are strictly less than $x$. The map $\tilde{e}_x$ is defined such that  
    \[ \tilde{e}_{x,i} = \tilde{e}_x |_{\partial_i(\R_+^{\bar{A}(x)+1} \times V^0(x))_*};\]
    this inductive definition of the attaching maps allows us to define the full Floer homotopy type by induction along $\CC$. 

    Let us now define $\tilde{e}_{x,i}$.,  
 fixing $i \in 1, \ldots, \bar{A}(x)+1$. Essentially, $\tilde{e}_{x,i}$ is the unique map we can define using the data of the $\tilde{e}_{x, y}$ as $y$ ranges over $\{ y \in \mathcal{P}_x : \bar{A}(y) = i\}$, by ``superimposing'' these maps on one another. For each such $y$,  we set
\[ U'_{xy} = \bar{\alpha}_{xy}^0(U_{xy}\times V^0(y)), U''_{xy} = U'_{xy} \times (V^1(x)^\perp \times \R_+^{\bar{A}(x)+1-i})_*. \]
In the complement of all $U''_{xy}$,  we set $\tilde{e}_{x,i}$ to map to the basepoint; and in each $U''_{xy}$ we set $\tilde{e}_{x,i}$ to be the composition 
\begin{equation}
    \label{eq:define-attaching-map-on-one-stratum} 
    \begin{gathered}U'_{xy} \times (V^1(x)^\perp \times \R_+^{\bar{A}(x)+1-i})_* \xrightarrow{(e_{xy} \wedge id) \circ(\bar{\alpha}^0_{xy})^{-1} \wedge id} (V^1(x,y) \times V^0(y) \times V^1(x)^\perp \times  (\R_+^{\bar{A}(x)+1-i})_* \simeq \\
    \simeq (\R_+^{\bar{A}(x)+1-i} \times V^0(y) \times V^1(y)^\perp   )_*\to \{\CC\}_{\mathcal{P}_{x>}}.
    \end{gathered}
\end{equation} 
Here the isomorphism from the right hand side of the first line to the left hand side of the second line is just rearranging factors and using the isomorphism \eqref{eq:canonical-iso-on-perps}, as well as the fact that $\bar{A}(y) = \bar{A}-i$ so the dimensions of the $\R_+$ factors agree. The left hand side of the bottom line is manifestly one of the terms of the wedge product in \eqref{eq:floer-homotopy-type-new}, and so the bottom map is the inclusion into the quotient by the previously defined relations. Notice the map $\tilde{e}_{x,i}$ is the composition of of a manifestly well-defined map 
\[ \tilde{e}'_{x,i}: \partial_i(\R^{\bar{A}(x)+1} \times V^0(x) \times V^1(x)^\perp)_* \to \bigvee_{y \in \mathcal{P}_{x>}: \bar{A}(y) = \bar{A}(x)-i}(\R^{\bar{A}(y)+1} \times V^0(y) \times V^1(y)^\perp)_* \]
sending the boundary the domain to the boundary of the codomain, with the inclusion of the codomain of $\tilde{e}'_{x,i}$ into $\{\CC\}_{\mathcal{P}_{x>}}$. 
\end{definition}

\begin{lemma}
    The attatching maps $\tilde{e}_{x}$ defined above are well-defined. 
\end{lemma} 
\begin{proof} (Sketch.)
    The non-overlap condition in the definition of an extension of an embedding means that the maps $\tilde{e}_{x,i}$ are each well defined. So the only question is whether these definitions agree on the boundaries; let's consider whether they agree on the intersection $\partial_i(\R_+^{\bar{A}(x)+1} \times V^0(x) \times V^1(x)^\perp )_*$ and $\partial_j(\R_+^{\bar{A}(x)+1} \times V^0(x) \times V^1(x)^\perp )_*$, for $i \neq j$. For notation let us call these strata $\partial_iC$ and $\partial_j C$. Without loss of generality suppose $i <j$. Let $U_{x,y_1}, \ldots, U_{x, y_r}$ be the open sets arising when defining $\tilde{e}_{x,j}$ and similarly let $U_{x, w_t}$ be the open sets arising when defining $\tilde{e}_{x,i}$. the compatibility condition for the extension of the embedding tells us that when the associated open sets $U''_{x,y_r} \subset   \partial_j C$, when restricted to $\partial_i C \cap \partial_j C$, become images of sets of the form 
    \[ U_{x, z_{\ell r}} \times U_{z_{\ell r}, y_r}  \times V^0(y_r) \times (\R^{\bar{A}(x)+i-j}_+)_* \]
    where $\bar{A}(z_{\ell r}) = \bar{A}(U_{w_t}) = c$ for some $c$ independent of $\ell, r, t$. 
    Dropping the point at infinity, we see that the disjointness condition for extensions of embeddings embeddings but now applied to the $U_{x, z_{\ell r}} \times V^0(z_{\ell r}) \subset \R^{\bar{A}(x, c)}_+ \times V^0(x)$ as well as the $U_{x, w_t} \times V^0(w_t) \subset \R^{\bar{A}(x,c)} \times V_0(x)$ implies that two of these sets overlap if and only if $z_{\ell r} = w_t$ for some $\ell, r, t$. But in that case, we see that for any point $\tilde{u}$ in $U''_{x, w_t} \cap \partial_j C$, $e_{i, x}(\tilde{u})$ will be identified with a point which is the image of the boundary sphere of the wedge factor associated to $w_t = z_{\ell r}$ in $\{\CC\}_{\mathcal{P}_x}$. The condition that the extension of the emebdding is compatible with composition then implies that $\tilde{e}_{j, x}(\tilde{u})$ will be sent to the same point in $\{\CC\}_{\mathcal{P}_x}$. 
\end{proof}

Thus, an extension of an embedding of a framed virtually smooth flow category $\CC'$ defines a canonical associated space $\{\CC\}$ (implicitly depending on the framing, embedding, and extension) and a corresponding spectrum $|\CC|$. We will see in the subsequent Section \ref{sec:equivariant-flow-categories} that it is elementary to extend all these definitions to the equivariant setting. The space $\{\CC\}$ is built out of disks one for each object of $\CC$, attached inductively along the poset $\CC$ with attaching maps determined by the chosen embeddings of the thickenings. In the $G$-equivariant setting, these disks become compactifications of $G$-representations and the attaching maps (and the poset $\CC$ itself) are all equivariant as well. The desuspension in the definition of $|\CC|$ is chosen such the spectrum $|\CC|$ is filtered by subspectra with associated graded given by stable spheres, with the stable sphere corresponding to $x \in \CC$ being of the form $\SS^{\mu(x)}$ where $\mu$ is the grading on $\CC$ (Section \ref{sec:comparison-with-floer-homology}). Thus, $\{\CC\}$ has the homotopy type of a CW complex, and $|\CC|$ defines an object in the Spanier-Whitehead category (Section \ref{sec:equivariant-homotopy-review}). It is a stable homotopy analog of a finite $CW$ complex, and thus a stable homotopy refinement of the finite-dimensional chain complex produced by Floer homology. In the equivariant setting, the stable spheres $\SS^{\mu(x)}$ are replaced by equivariant stable spheres $\SS^{[\mu](x)}$, where $[\mu](x)$ is a virtual $G$-representation. In the Floer theoretic situation, this virtual $G$-representation should be thought of as the stable equivariant spheres corresponding to the virtual representation given by the difference between the positive and the negative eigenvalues of the Floer action functional, (=or, more rigorously, the equivariant spectral flow of a path from the self-adjoint operator given by Hessian at a given critical point to a base-point self-adjoint operator. The latter notion is made completely formal in Section \ref{sec:index-theory} and and specialized to the discussion of Hamiltonian Floer Homology in Section \ref{sec:producing-framings}.

\begin{figure}[h!]
    \centering
    \includegraphics[width=\textwidth]{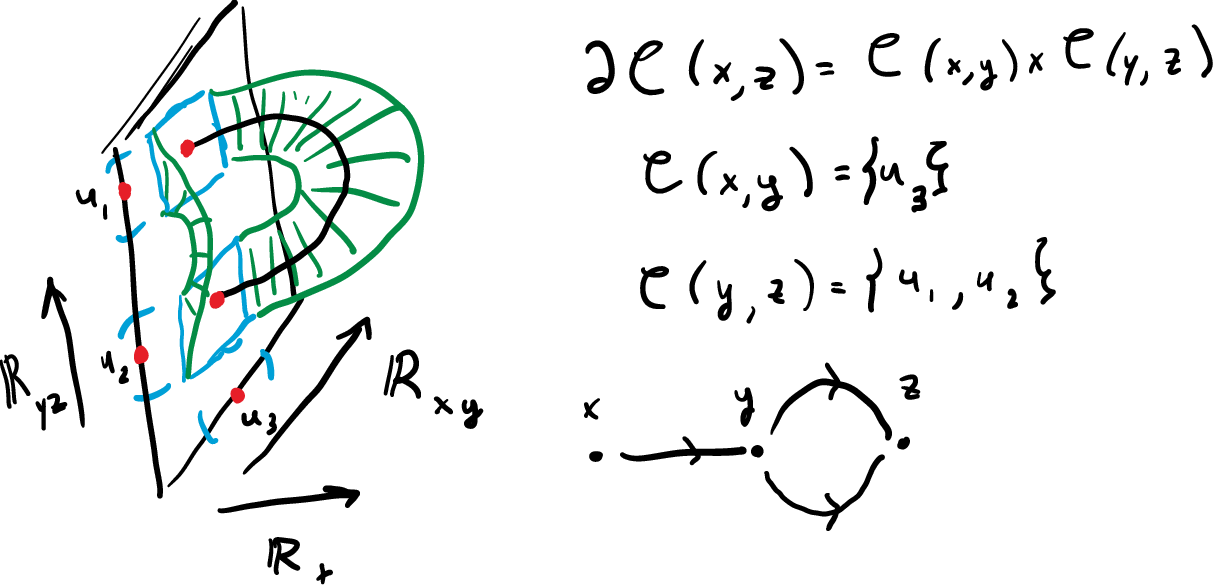}
    \caption{\emph{Flow category with composable morphisms.} This figure depicts (on the left) the embedding of a flow category $\CC$ with post depicting as on the bottom right, with all arrows depicting regular isolated flow lines of a $C_2$-equivariant Morse function. the embedding of the thickening of $\CC(x,y)$ is sitting in $\R_{xy}$, and the embedding of the thickening of $\CC(x,z)$ is sitting in $\R_{yz}$, as depicted on the left. These embeddings determine the embedding of the boundary of the thickening of $\CC(x,z)$ (which is the product of the previously mentioned thickenings -- the blue squares in the figure), which must be extended to the interior of $\R_{yz} \times \R_{xy}$ (green). Note how the framing will be forced to twist along this embedding of the thickening of $\CC(x,z)$, leading to the conclusion that $d^2=0$ in the resulting Morse complex. In the final construction of the Floer homotopy type (Definition \ref{def:floer-homotopy-type}), the cell associated to $x$ is essentially the `cone' on the figure on the right. Note that the embedding for the thickening of the virtual smoothing of Figure \ref{fig:virtually-smooth-category-obstructed-example} looks very similar, except for two modifications. First,  in order to embed $\CC(y,z)$ we must stabilize the previously zero-dimensional Kuranishi charts to $1$-dimensional Kuranishi charts by a copy of the sign representation (so that equivariant an embedding exists!); this makes making the obstruction bundle of $\CC(x,z)$ have fiber $\R_{sgn}^2$, and $\R_{yz}$ has a nontrivial $C_2$ action. The other difference is that $\R_{xy}$ in this latter case is zero dimensional (so one slices the picture vertically at the point $u_3$). The extension of the obstruction section to $\R_{xy} \times \R_+$ (corresponding to the extension of the embedding) has no zeros on the interior, and upon perturbation will have an isolated finite set of zeros on the interior.The final Floer homotopy space is inductively built by attaching disk but the corresponding filtration does not me the space into a CW complex, as the attaching map from the boundary of the $3$-dimensional cell corresponding to $x$  to the $3$-dimensional cell corresponding to $y$ (here the extra dimension comes from $V(y)^\perp$ in the formula for $\{\CC\}$)is not the zero map (although it is zero in the homotopy category).  }
    \label{fig:flow-category-with-multiple-layers}
\end{figure}

Thus, after the discussion in Section \ref{sec:equivariant-flow-categories} of the equivariant situation, we will know how to produce a genuine equivariant stable homotopy type from an equivariant analog of the above data. The following Figures \ref{fig:flow-category-with-multiple-outputs} and \ref{fig:flow-category-with-multiple-layers} explain what the equivariant Floer homotopy type looks like in simple examples; the general case is essentially an amalgamation of the cases described in those figures. The subsequent Section \ref{sec:flow-category-constructions} explains how to actually \emph{produce} framings and embeddings of framings of a framed virtually smooth (equivariant) flow category; and the actual virtually smooth flow categories associated to Hamiltonian Floer homology are produced in Section \ref{sec:global-charts} and framed in Section \ref{sec:index-theory}. Using this data one can produce the equivariant stable homotopy types associated to Hamiltonian Floer Homology; this is done in Section \ref{sec:floer-homotopy}. Similarly, geometric constructions described in these aforementioned sections allow us to relate these flow categories using the continuation map equation and its homotopies; the pattern is quite general, and fairly general tools regarding the manipulation of homotopy types using flow categories are described in Section \ref{sec:comparison-with-floer-homology}. 

In the following remarks, we comment on several aspects of the above construction.
\begin{remark}
\label{rk:embeddings-dont-really-matter}
    The Floer homotopy type $|\CC|$ may seem to depend not just on the framed virtually smooth flow category, but on its embedding and the extension of its embedding. However, Section \ref{sec:flow-category-constructions} below develops certain tools which allow one to relate the Floer homotopy type of $\CC$ to that of a \emph{stabilization} of $\CC$, or a \emph{restratification} of $\CC$. Neither of these turn out to change the (genuine equivariant) homotopy type up to canonical weak equivalence; moreover, while we do not write this out in this paper, the tools of Section \ref{sec:flow-category-constructions} and an elementary parameterized homotopy theory argument can establish that \emph{the floer homotopy type $\CC$ does not depend on the (extension of) the embedding of the framed virtual smoothing $\CC'$, but only on the (semi-freely) framed virtual smoothing $\CC'$ itself,} up to canonical weak equivariant homotopy equivalence. Thus, we should think of the framed virtually smooth category $\CC'$ as representing a (possibly equivariant) stable homotopy type.
\end{remark}

\begin{remark}
\label{rk:embeddings-dont-really-matter}
    Analogously to Definition \ref{def:stabilize-derived-manifold} and the discussion about the Pontrjagin-Thom construction in Section \ref{sec:pedagogical-intro}, the subsequent Section \ref{sec:flow-category-constructions} shows that the floer homotopy type $|\CC'|$ actually only depends on the framed virtually smooth flow category $\CC'$ up to stabilization. Similarly, the results of that section show that the data of the integral action $\bar{A}$ is essentially auxiliary (and this will be used to show invariance of the equivariant floer homotopy type $|\CC'(H, J)|$ in Section \ref{sec:invariance-of-equivariant-floer-homotopy}).  
\end{remark}

\begin{remark}
\label{rk:how-do-framings-of-thickenings-play-a-role}
    One sees by observation that while the data of the framings of the obstruction bundle and all aspects of the parameterization of the framing of $\CC'$ plays a role in the definition of the Floer homotopy type, the actual framings of the tangent bundles of the thickenings of the the morphism objects of $\CC'$ play no role in the definition of the Floer homotopy type. However, the process of \emph{producing} embeddings of $\CC'$ (see Proposition \ref{prop:embeddings-exist}) involves heavy use of the Pontrjagin-Thom construction, which involves the framings of the thickenings. One can verify that if one changes the framing of the thickenings one can also thus affect the eventual Floer homotopy type produced, in a sense that can be seen in the following example. Given a virtually smooth flow category $(\CC, \CC')$ with two objects $x$ and $y$, the Floer homotopy type $|\CC'|$ produced by applying Proposition \ref{prop:embeddings-exist} and the subsequent disscussion to produce the required extension of an embedding, will manifestly agree with the vitual Pontrjagin Thom construction on $\CC'(x,y)$, performed as described in in the introdution. Twisting the framing of the thickening of $\CC'(x,y)$ will thus, in general, change the Floer homotopy type. 
\end{remark}

\begin{remark}
    Because the Floer homotopy type is really an invariant of the framed virtual smoothing $\CC'(x,y)$ up to stabilization, one should think of the \emph{space} of framings as colimit of the spaces of framings of all the virtually smooth flow categories $\CC'(x,y)_{\mathcal{E}}$ produced by stabilizing $\CC'(x,y)$ by an arbitrarily large semi-free parameterization $\mathcal{E}$. This space is highly nontrivial; in the case of a flow category with $2$ objects and with $\CC'(x,y)$ simply a Kuranishi chart on a point with no obstruction bundle, this space can be identified with $O= \colim_n O(n)$, (and in the case of a $G$-action with its equivariant analog described in Appendix \ref{sec:classifying-spaces-for-k-theory}). The fact that this space is nontrivial even in the case of such a trivial flow category has to do with \emph{parameterized Floer theory}: one can imagine a \emph{parameterized family of flow categories} with fiber isomorphic to the given one and the framings (after stabilization of the entire family) twisting as one changes the parameter in the family, giving rise to a parameterized homotopy type that is not the same as the produce of the usual homotopy type with the parameter space. We do not develop the parameterized theory of virtually smooth flow categories in this paper, and we leave this to further works.
\end{remark}

\begin{remark}
    The filtration on the homotopy space $\{\CC\}$ by the poset underlying $\CC$ is in general \emph{not} a CW filtration, and the space is generally \emph{not} a CW complex, but only a space built inductively by attaching cells.  Such a phenomenon arises whenever one is studying the setting where there are flow lines at are \emph{obstructed}, or \emph{not regular}, as in the example of Figure \ref{fig:obstructed-c2-example}. The caption of Figure \ref{fig:flow-category-with-multiple-layers} explains how the corresponding Floer homotopy space looks like for that example. In contrast, in all previous constructions of a Floer homotopy type, the latter is built by initially building a space \emph{which is a CW complex} with the CW filtration refining to a filtration by $\CC$. The new phenomenon in our case arises because we stratify by the \emph{auxiliary integral action} $\bar{A}$ rather than by the index $\mu$, \emph{unlike all previous constructions}. This in turn is forced on us by the need to study obstructed Floer trajectories, which force the introduction of additional strata which are disallowed in the regular case. 
\end{remark}

\begin{remark}
    An analogous situation to the previous remark arises very clearly in the setting of Morse theory. A  Morse function (say, with distinct critical values)  on a closed manifold $Y$ defines a filtered space $Y^{\leq \bullet}$ by the filtration by the sublevel sets of its critical values, even when it is not Morse-Smale. however, when the Morse function is Morse-Smale, the corresponding filtered space \emph{is a CW complex} \cite{qin2021application}. As one moves between Morse functions through Morse-smale functions, the filtrations change continuously, moving through a non-CW filtration in between two different CW complex structures on $Y$. The phenomenon in the previous remark is an abstracted analog of this example; an essential insight of this paper is that \emph{a non-regular floer action functional defines a well defined filtered stable homootopy type, even as it does not define a well-defined chain complex!}
\end{remark}

\begin{figure}[h!]
    \centering
    \includegraphics[width=\textwidth]{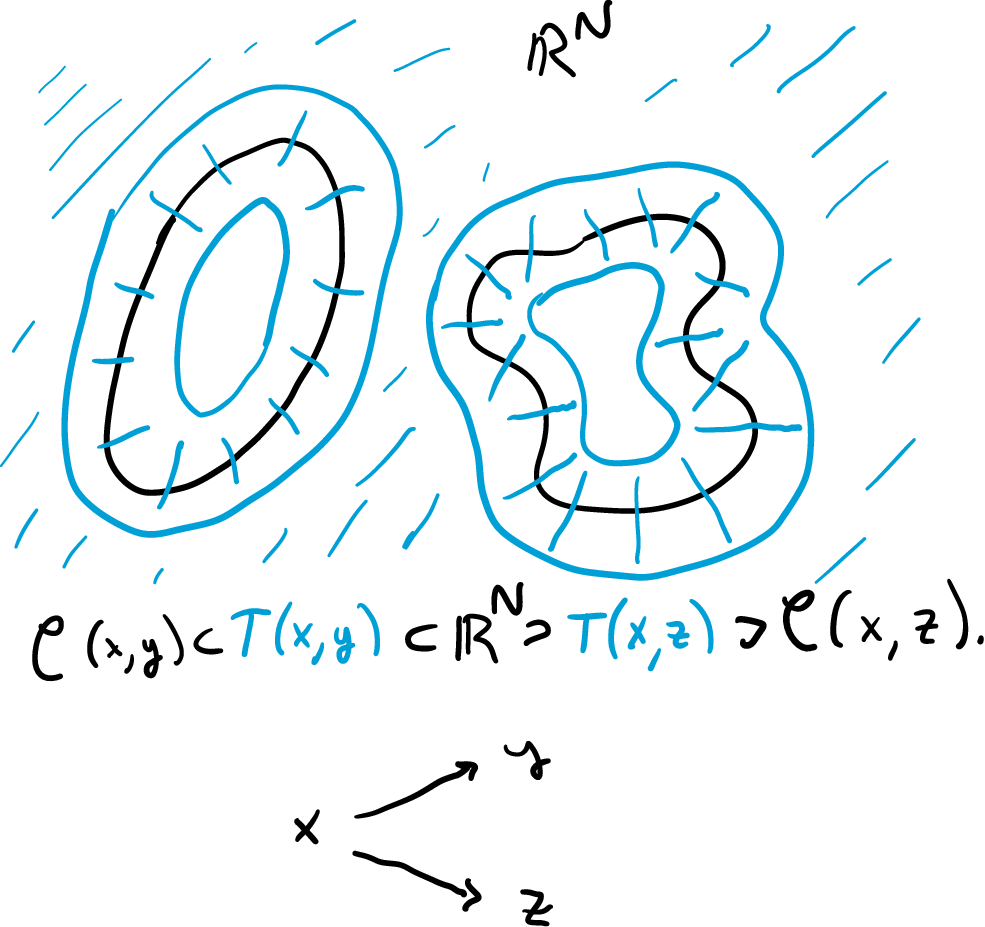}
    \caption{\emph{Flow category with multiple outputs.} Given a flow category with three objects $x,y,z$, with nontrivial morphisms $x \to y$ and $x \to z$, one produces the resulting homotopy type by finding disjoint (equivariant) embeddings of (stabilizations of) $T(x,y)$ and $T(x,z)$ into a vector space, and then extending the defining sections (which can be thought of as taking values in a wedge of spheres, each sphere coming from stabilizing $V(x,y)$ or $V(x,z)$ appropriately and compactifying) to the rest of the one-point compactification by requiring that on the complement of a slightly larger thickening of the images of $T(x,y)$ and $T(x,z)$ one maps to the basepoint of the corresponding wedge of spheres. If the maps in the stable homotopy category associated to $x \to y$ and $x \to z$ are $f_{xy}: \SS^{V(x)} \to \SS^{V(y)}$ and $f_{xz}: \SS^{V(x)} \to \SS^{V(z)}$, respectively, then the resulting map represents the canonical map $(f_{xy} \vee 0) + (0 \vee f_{yz}): \SS^{V(x)} \to \SS^{V(y)} \vee \SS^{V(z)}.$. See Lemma \ref{lemma:model-for-summing-maps}. }
    \label{fig:flow-category-with-multiple-outputs}
\end{figure}

\subsection{Equivariant Flow Categories}\label{sec:equivariant-flow-categories}
In this section, we generalize the earlier constructions of Section \ref{sec:flow-categories} to the setting to the equivariant setting. 

\begin{definition}
    Let $G$ be a finite group. A $G$-equivariant flow category $\CC$ is a flow category $\CC$ together with an action of $G$ on $Ob(\CC)$, and for each $g\in G$ and each pair $(x,y)$ of objects of $\CC$ such that $x > y$, an isomorphism in $\langle k \rangle$-Spc'
    \[ i^g_{x,y} = (\hat{i}^g_{x,y}: S_{\CC(x,y)} \to S_{\CC(gx, gy)}, \bar{i}^g_{x,y}: \CC(x,y) \to \CC(gx, gy))\] 
    \[ g_{xy}: \CC(x,y) \to \CC(gx, gy)\] 
    satisfying 
    \begin{itemize}
    \item $g'_{gx,gy} g_{x,y} = (g'g)_{x,y}$ \text{ for all } $x> y$ in $\CC$ and for each $g', g \in G$, such that
    \item $(id)_{xy}$ is the identity map, and 
    \item The action is compatible with composition: for every $x>y>z$ objects of $\CC$, writing $(f_{x,y,z}, \bar{f}_{x,y,z})$ for the composition $\CC(x,y) \times \CC(y,z) \to \CC(x,z)$, we have that $g_{xz}\bar{f}_{x,y,z} = \bar{f}_{gx,xy,gz} \circ (g_{xy} \times g_{yz}).$
    \end{itemize}
\end{definition}

\textbf{Notation.} Given a group $G$ acting some set $X$, we will denote the stabilizer of $x \in X$ either by $G_x$ or by $Stab_x$. 

\begin{definition}
    An integral action $\bar{A}$ on a $G$-equivariant flow category is the same as the correspondig data for the underlying flow category, with the additional requirements that 
    \[ \bar{A}(gx) = \bar{A}(x) \text{ and }\beta_{x,y}\circ \hat{i}^g_{xy} = \beta_{gx,gy}\text{ for any } g \in G, x,y \in Ob(\CC).\]
\end{definition}

\begin{definition}
\label{def:equivariant-grading}
    A \emph{grading} on a $G$-equivariant flow category $\CC$ is an assignment of a pair of vector spaces $(V^0(x), V^1(x))$ to each object $x \in \CC$, together with linear isomorphisms of inner product spaces $g^0_x: V^0(x) \to V^0(gx)$, $g^1_x: V^1(x) \to V^1(gx)$, for each $g \in G$ and $x \in \CC$, such that 
    \[ (id)^0_x = id, (id)^1_x = id, (g')^0_{gx} g^0_x = (g'g)^0_x, (g')^1_{gx} g^1_x = (g'g)^1_x.\]

    An \emph{isomorphism} of gradings 
    \[\{(V^0(x), V^1(x)), (g^{0}_x, g^{1}_x)\} \simeq \{(\hat{V}^0, \hat{V}^1); (\widehat{g^{0}_x}, \widehat{g^{1}_x})\}\] 
    is an assignment of a vector spaces $W(x), \widehat{W}(x)$ for each $x \in Ob(\CC)$ together with isomorphisms $g^{w}_x: W(x) \to W(gx)$,  $\widehat{g^{w}_x}: \widehat{W(x)} \to \widehat{W(gx)}$ for each $x \in \CC, g \in G$ satisfying 
    \[ (id)^{w}_x = \widehat{(id)^{w}_x}= id, (g')^{w}_{gx} g^w_{x} = (g'g)^{w}_{x}, \widehat{(g')^{w}_{gx}} \widehat{g^w_{x}} = \widehat{(g'g)^{w, i}_{x}},\]
    together with isomorphsms 
    \[ h^0_x: V^0(x) \oplus W(x) \to \widehat{V^0(x)} \oplus \widehat{W(x)}\]
    \[ h^1_x: V^1(x) \oplus W(x) \to \widehat{V^1(x)} \oplus \widehat{W(x)}\]
    such that 
    \[ (\widehat{g^{0}_x} \oplus \widehat{g^{w}_{x}})\circ h^0_x = h_{gx}(g^{0}_x \oplus g^{w}_{x}), \]
    \[(\widehat{g^{1}_x} \oplus \widehat{g^{w}_{x}})\circ h^1_x = h_{gx}(g^{1}_x \oplus g^{w}_{x}).\]

    Isomorphism classes of gradings are specified by an assignment of an element $[V^0(x)] - [V^1(x)] \in RO(H_x)$ of a virtual $H$-representation for every orbit $[x] \in Ob(\CC)/G$ with stabilizer (up to conjugacy) $H_x \in G$. We will denote the function 
    classifying isomorphism classes of gradings by 
    \[[\mu]: Ob(\CC) \to \bigsqcup_{H \subset G} RO(H).\]
    We have that $[\mu](x) \in RO(Stab_x)$, and $[\mu](gx) = g_* ([\mu](x)) \in RO(g Stab_x g^{-1})$.
\end{definition}
\begin{definition}
\label{def:equivariant-pre-grading}
    A \emph{pre-grading} has the same definition as a grading. However, we will use pre-gradings and gradings differently, and the comparison between the two will be through the (essentially auxiliary) integral action. Specifically, given a pre-grading 
    $\{(V^0(x), V^1(x)), (g^{0}_x, g^{1}_x)\} \simeq \{(\hat{V}^0, \hat{V}^1)$ on a finite $G$-equivariant flow category with integral action $\bar{A},$ there is an associated grading $\{(\widehat{V^0(x)},\widehat{V^1(x)}), (\widehat{g^{0}_x}, \widehat{g^{1}_x})\}$ well defined up to canonical isomorpism. Choosing some $N$ such that $\bar{A}(x) + N \geq 0$ for all $x \in \CC$, the associated grading is defined via 
    \[ \widehat{V^0(x)} = \R^{\bar{A}(x) + N} \oplus V^0(x), \widehat{V^0(x)} = \R^{\bar{A}(x) + N} \oplus V^0(x),\widehat{g^{0}_x} =id \oplus g^{0}_x, \widehat{g^{1}_x} =id \oplus g^{1}_x. \]
\end{definition}

\begin{definition}
    A smooth $G$-equivariant flow category is a smooth flow category such that the maps $i^g_{x,y}$ are maps in $Man'$.
\end{definition}

\begin{definition}
    A virtually smooth $G$-equivariant flow category is a virtually smooth flow category $(\CC, \CC')$ (without stabilizations), where $\CC$ is a $G$-equivariant flow category, such that we also have isomorphisms $i^g_{x,y}$ in  $DerMan'_{sm, no-stab}$ from $\CC'(x,y)$ to $\CC'(gx, gy)$, which we write as
    \[i^{g}_{x,y} = (\hat{i}^g_{x,y}: S_{\CC(x,y)} \to S_{\CC(gx,gy)}, (\bar{i}^g_{x,y}: T(x,y) \to T(gx,gy), (\bar{i}')^g_{x,y}: V(x,y) \to V(gx,gy))\]
    which satisfy
    \begin{itemize}
    \item $i^{g'}_{gx,gy}i^g_{x,y} = i^{g'g}_{x,y}$ for all $g, g' \in G$;
    \item $i^{id}_{x,y}$ is the identity map, and 
    \item The action $i^g_{x,y}$ is compatible with composition: writing $(f_{xyz},\bar{f}_{xyz}, \bar{f}'_{xyz})$ for the composition $\CC(x,y) \times \CC(y,z) \to \CC(x,z)$, we have that
    \[ \bar{f}_{gx gy gz} \circ (\bar{i}^g_{xy} \times \bar{i}^g_{yz}) = i^g_{x z}\bar{f}_{x y z},  \bar{f}'_{gx gy gz} \circ ((\bar{i})'^g_{xy} \times (\bar{i})'^g_{yz}) = (\bar{i}')^g_{x z}\bar{f}_{x y z}.\]
    \end{itemize}
\end{definition}

\begin{definition}
\label{def:equivariant-framing}
    A \emph{framing} of a virtually smooth $G$-equivariant flow category $(\CC, \CC')$ is the the data $\{(\eta_{xy}, \eta'_{xy}, \eta_{xyz}), (\tau^0_{xyz}, \tau^1_{xyz})\}$ of a framing of the underlying virtually smooth flow category, together with vector space isomorphisms 
    \[ g^0_{xy}: V^0(x,y) \to V^0(gx,gy)\]
    \[g^1_{xy}: V^1(x,y) \to V^1(gx,gy)\]
    satisfying the following conditions. First, 
    \begin{itemize}
        \item Write $\underline{g^0_{xy}}:\underline{V^0(x,y)} \to \underline{V^0(gx,gy)}$ for the map of vector bundles covering $\bar{i}^g_{xy}$ which acts by $\check{\lambda}^g_{xy}$ $g^0_{xy}$ on the fiber;
        \item Write $\underline{id}: \underline{T\R_+^{\bar{A}(x) - \bar{A}(y)-1}} \to \underline{T\R_+^{\bar{A}(gx) - \bar{A}(gy)-1}}$ for the map of vector bundles covering $\bar{i}^g_{xy}$ which acts by identity on the fiber;
        \item Write $\underline{g^1_{xy}}$ for the map of vector bundles $\underline{V^1(x,y)} \to \underline{V^0(gx,gy)}$ covering $\bar{i}^g_{xy}$ which acts by $g^1_{xy}$ on the fiber; and finally
    \end{itemize} 
    Then the following conditions must be satisfied:
    \begin{equation}
    \eta_{gx gy} \circ d\bar{i}^g_{xy} = (\underline{id} \boxplus \underline{g^0_{xy}}) \circ \eta_{x y};
    \end{equation}
    \begin{equation}
    \mu^g_{xy} \circ \eta'_{xy} = \eta'_{gx gy} (\bar{i}')^g_{xy};
    \end{equation}
    and finally that the maps $\tau^i_{xyz}$ commute with the $G$-action in the natural sense; these conditions are explicitly spelled out in Definition \ref{def:equivariant-F_2-framing}, making the data of the underlying $F_2$-parameterization and the compatible $G$-action on the $V^i(x,y)$ into the (equivariant) \emph{underlying $F_2$-parameterization} of the framed equivariant flow category (thought of as a $G$-poset) 

    The compatibility conditions
    \begin{equation}
        (g')^0_{gx gy} g^0_{x y} = (g'g)^0_{xy}; (g')^1_{gx gy} g^1_{xy} = (g'g)^1_{xy}
    \end{equation}
    are actually forced by the previous conditions. 
\end{definition}

\begin{definition}
    A framing of an equivariant flow category is compatible with a pre-grading if the pre-grading can be enhanced to some $G-E_2$-parameterization with associated $G-F_2$-parameterization isomorphic to the underlying $G-F_2$-parameterization of the category. A framing of an equivariant flow category is compatible with a grading $\mathcal{G}$ if it is compatible with some pre-grading $\mathcal{G}_{pre}$ with associated grading  $\mathcal{G}$.  
\end{definition}
\begin{definition}
    An embedding of a framing of a $G$-equivariant virtually smooth flow category is a $G$-$E_2$-parameterization with associated $G-F_2$-parameterization isomorphic to the one underlying the framing of the category. In other words, an embedding of a framing is the \emph{same thing} as a compatible (pre)-grading. We will largely use the language of embeddings of framings in the rest of the paper, but this relation is essential in the interpretation of the construction.
\end{definition}

\begin{definition}
    An embedding of a framing of a $G$-equivariant virtually smooth flow category is adapted if the underlying framing is adapted.
\end{definition}

\begin{definition}
    A \emph{framed embedding} of a finite, adapted, virtually smooth $G$-equivariant flow category without stabilizations is an embedding of the underlying flow category such that 
    \begin{equation}
    \iota_{gx gy}\circ \bar{i}^g_{xy} = (id \times \hat{\lambda}^g_{xy}) \circ \iota_{xy}.
    \end{equation}
\end{definition}

\begin{definition}
    An \emph{extension} of an embedding of a $G$-equivariant flow category $\CC$ is an embedding of the underlying flow category $\CC$ such that 
    \[e_{gx} \circ (id \times g^0_{x})_* = \left( \bigvee_{z \in M(\CC)} (g^0_z \oplus g^1_{x,z})_*\right) \circ e_x. \]
\end{definition}

We now wish to define the genuine equivariant stable homotopy type associated to an extension of a framed embedding of a $G$-equivariant flow category. 

\begin{definition}
\label{def:equivariant-stable-homotopy-type}
    The equivariant stable homotopy type asociated to an extension of an embedding of a virtually smooth $G$-equivariant flow category is given by the same formula as the definition of the stable homotopy type of the extension of the embedding the underlying virtually smooth flow category, using the fact that $G$ acts on 
   \[\bigvee_{z \in \CC}  (\R_+^{\bar{A}(x)+1} \times V^0(x) \oplus V^1(x)^\perp)_*\]
    via 
    \[ \bigvee_{z \notin M(\CC)}  (id \times g^0_x \oplus (g^1_x)^\perp)_*\]
    where $V^1_{max}$ is naturally an orthogonal $G$-representation and 
    \[ g^1_x: V^1(x)^\perp \to V^1(gx)^\perp\]
    is induced by restricting the $G$ action to $V^1(x)^\perp \subset V^1_{max}$.
\end{definition}

It is clear from the construction that 
\begin{proposition}
The spectrum $|\CC|$ underlying the equivariant stable homotopy type of $\CC$ is the stable homotopy type of $\CC$ thought of as a non-equivariant flow category.
\end{proposition} 

\begin{definition}
    Let $H \subset G$ be a \emph{normal} subgroup and let $\CC$ be a $G$-equivariant flow category. The fixed point flow category $\CC^H$ the $G/H$-equivariant flow category with objects $Ob(\CC)^H$ and morphisms $\CC^H(x,y) = \CC(x,y)^H$.

    If $(\CC, \CC')$ is a virtually smooth $G$-equivariant flow category then $(\CC^H, (\CC')^H)$ is a virtually smooth $G/H$-equivariant flow category, where 
    \[ (\CC')^H(x,y) = (T(x,y)^H, V(x,y)^H, \sigma(x,y)^H)\]
    where we define $\sigma(x,y)^H(c) = \frac{1}{|H|}\sum_{h \in H} h\sigma(x,y)(c)$ for $c \in T(x,y)^H$; this is the canonical $H$-equivariant projection $V(x,y)_c \to V(x,y)_c^H$. 
\end{definition}
\begin{definition}
    If $\mu = (V^0(x), V^1(x), g^0_x, g^1_x)$ is a grading or a pre-grading on a $G$-equivariant flow category $\CC$, then for any $x \in \CC^H$, there is an associated grading or pre-grading on $\CC^H$ defined via 
    \[(V^0(x)^H, V^1(x)^H, [g]^0_x := (g^0_x)|_{V^0(x)^H}, [g]^1_x:= (g^1_x)|_{V^1(x)^H}).\]
\end{definition}

\begin{lemma}
\label{lemma:framing-of-invariant-category}
   Given a framing of a virtually smooth $G$-equivariant flow category $\CC$, there is an associated framing of $(\CC^H, (\CC')^H)$ defined by 
   \[ V^0_{(\CC')^H}(x, y) = V^0_{(\CC')}(x, y)^H,  V^1_{(\CC')^H}(x, y) = V^1_{(\CC')}(x, y)^H,  W^0_{(\CC')^H}(x, y, z) = W^0_{(\CC')}(x, y, z)^H\]
   together with 
   \[ ([g]^0_{xy})_{(\CC')^H} = (g_{xy})^H,  ([g]^1_{xy})_{(\CC')^H} = (g^1_{xy})^H\]
   where we use the fact that for $x > y > z \in \CC^H$, 
   $(g^0_{xy},g^1_{xy})_{g \in H}$  make $(V^0(x, y), V^1(x,y), W^0(x,y,z))$ into $H$-representations as $(x,y,z) = (gx,gy,gz)$.
\end{lemma}

\begin{theorem}
\label{thm:geometric-fixed-points-of-flow-category}
The geometric fixed points of the equivariant stable homotopy type associated to a flow category agree with the stable homotopy type of fixed point flow category:
    \begin{equation}
    |\CC|^{\Phi H} = |\CC^H|. 
    \end{equation}
\end{theorem}

This follows from a point-set computation together with
\begin{proposition}
\label{prop:homotopy-type-is-G-CW-complex}
    The space $Y_{\CC}$ of Definition \ref{def:equivariant-stable-homotopy-type} has the homotopy type of a $G$-CW complex. 
\end{proposition}
\begin{proof}
    Note that $(\R_+^k)_* = c(S^{k-1}) = (\R_+)_* \wedge S^{k-1}$. Using the partial ordering on $Ob(\CC)$ given by $y \leq x$ if $\CC(x,y) \neq \emptyset$, we thus see that $Y_\CC = Y^r$, where $Y^0 = *$ and $Y^r$ is built from $Y^{r-1}$ by forming the cone on a map 
    \[f_r: \wedge_{x \in [x]} S_x^{V'} \to Y^{r-1}\]
    where $[x]$ is an orbit of the $G$-action on $Ob(\CC)$, $H$ is the stabilzer of some arbitrary base-point $x_0 \in [x]$, $V'$ is some $G/H$ representation, $x$ is simply a label on the sphere, and $g \in G$ acts on $\wedge_{x \in [x]} S_x^{V'}$ by 
    \[ \wedge_{x \in [x]} (g_*: S^{V'}_x \to S^{V'}_{gx}), \]
    i.e. by permuting spheres as $G$ permutes objects and acting on the corresponding representation spheres as in the representation $V'$.   Mapping cones on $G$-homotopy equivalent maps are $G$-homotopy equivalent; by induction we can assume that $Y_r$ has the $G$-homotopy type of a $G$-CW complex, in which case the equivariant whitehead theorem implies that $f_r$ is homotopy equivalent to a $G$-CW map and thus $Y_r$ is $G$-homotopy equivalent to a $G$-CW complex as well.
\end{proof}

\section{Flow Categories: Constructions}
\label{sec:flow-category-constructions}

\subsection{Equivariant embedings of $\langle k \rangle$-manifolds}

\begin{definition}
\label{def:convenient-metric}
A \emph{convenient metric} on a $G-\langle k \rangle$-manifold $X$ is a $G$-invariant metric for which all the boundary hypersurfaces $X_j$ are totally geodesic and intersect totally orthogonally in the sense that for any $S \in [k]$ and any $x \in X(S)$, any sequence of nonzero vectors $v_j \in T_xX(S)^\perp \subset T_xX(S \cup \{j\})$ ranging over $j \notin S$ is pairwise orthogonal in $T_xX$.
\end{definition}

\begin{lemma}
Any $G-\langle k \rangle$-manifold $X$ admits a convenient metric $g$. Given convenient metrics $g_i$ on $X_i$ such that $g_i|_{X_i \cap X_j} = g_j|_{X_i \cap X_j}$ (defining a metric $g_0$ on $\partial X$) the convenient metric $g$ can be chosen to be an extension of all $g_i$. 

\end{lemma}

\begin{proof}
We employ the \emph{doubling construction} on $X$. This construction is a generalization  to $\langle k \rangle$-manifolds of doubling construction on manifolds with boundary $Y$, which puts a smooth structure on $Y_{\cup \partial Y} Y$ making the latter into a closed $\Z/2$-manifold, such that the inclusion of $Y \to Y \cup_{\partial Y} Y$ is a smooth embedding. The analog for $\langle k \rangle$-manifolds $X$ puts a smooth structure on the union 
\[ \mathcal{D}(X) = \bigsqcup_{i \in 2^k} X/\sim  \]
of $2^k$ copies of $X$ glued appropriately along their strata $X(S)$, such that $\mathcal{D}(X)$ is a $G\times(\Z/2)^k$-manifold and the inclusion of one copy of $X$ is the fundamental domain of the $1 \times (\Z/2)^k$ action. The details of this construction are worked out in \cite{rezchikov2022integral}, although it is implicit in \cite{albin2011resolution}. 

We then apply an equivariant embedding theorem 
\cite{wasserman1969equivariant} that (for any compact Lie group 
$G$) every $G$-manifold $Y$ can be isometrically embedded into 
an orthogonal $G$-representation $V$, and if such an embedding 
was specified on a closed subset then it can be extended to the 
whole manifold after enlarging $V$ to a larger $G$-
representation. Moreover, the $G$-representation $V$ can be 
chosen to be \emph{subordinate} to the $G$ action on $Y$, which 
means that $V$ is a direct sum of irreducible $G$-
representations $V_1, \ldots, V_r$ such that for each $m \in M$ 
with stabilizer $G_m \subset G$, the manifold $G/G_m$ embeds 
equivariantly into $(\oplus_{i=1}^r V_i)^\ell$ for some $\ell$, 
as does the $G_m$-representation $T_xM/Lie(G_m)$ 
\cite[Proposition 1.7]{wasserman1969equivariant}. The 
assumptions of the theorem give a $G \times (\Z/2)^k$-
equivariant metric $g_0$ on the $(\Z/2)^k$-orbit $\mathcal{D}
(\partial X) \subset \mathcal{D}(X)$ of $\partial X \subset X$. Repeatedly applying the previous theorem 
allows us to define an embedding of $\mathcal{D}(X) $ into some $(G \times \Z/2)^k$-representation $V$ which is an isometry on $\mathcal{D}(\partial X)\subset \mathcal{D}(X)$. Pulling the metric on $V$ back onto $\mathcal{D}(X)$ defines an equivariant extension of $g_0$, and restricting this metric to $X$ defines the metric required in the theorem. Indeed, the boundary hypersurfaces $X_j$ will be totally geodesic in $X$ because the images of $X_j$ in $V$ are the the intersections of the image of $X$ with the totally geodesic submanifolds given by the subrepresentations $V_j \subset V$ on which the $j$-th $\Z/2$ factors act trivially. The condition on the orthogonality holds because it holds for any $Stab(x) \subset G \times (\Z/2)^k$-invariant inner product on $T_x\mathcal{D}(X) = T_xX$. 
\end{proof}

One would like to use the equivariant embedding theorem quoted above and the doubling construction to produce embeddings of $G-\langle k \rangle$ manifolds into $\R^k_+ \times \hat{V}$ for some $G$-representation $V$ . Unfortunately, the proof of the equivariant embedding theorem does not make it immediately clear that one can choose the representation on which $(\Z/2)^k$ acts nontrivially to be precisely the representation $\R^k$ on which the $i$-th $\Z/2$ in $(\Z/2)^k$ acts by flipping the $i$-th coordinate. As such, we give a direct proof of this fact.

\begin{lemma}
\label{lemma:embed-a-single-equivariant-manifold}
    Let $\hat{V}$ be a $G$-representation, and let $Y$ be a $G-\langle k \rangle$ manifold with a convenient metric. Suppose that one has equivariant immersions 
    \[\iota_j: Y_j  \to (\R^k_+)_j \times \hat{V} \] such that for any $S' \subset [k]$, the restrictions of $\iota_j$ for all $j \in [k] \setminus S'$ to $Y(S)$ agree and factor through $(\R^k_+)(S) \times \hat{V}$.  Suppose one also has a $G$-equivariant framing 
    \[ \eta: TX \simeq \underline{T\R^k_+ \times V}\]
    for some $G$-representation $V$. 
    
    Then there exists a $G$-representation $\hat{V}'$ together with an immersion 
    \[ \iota: Y \to \R^k_+ \times \hat{V} \times \hat{V'}\]
    which maps the top stratum of $Y(S')$ to the top stratum of $\R^k_+(S') \times \hat{V} \oplus \hat{V'}$ for all $S' \subset [k]$, is an embedding on the top open stratum of $Y$, and satisfies $\iota|_{Y_j} = \iota_j \times 0$. 

    Moreover, if $X$ has a convenient metric, then if the composition of $\iota_j$ with the projection to the $r$-th factor of $\R^k_+$ is the distance on $Y_j$ from $Y_j \cap Y_r$ in a neighborhood of the latter for every $j,r$, and and the projection of $d\iota_j$ onto $T(\R^k_+)_j$ is full rank everywhere, then we can choose the composition of $\iota$ with the projection to the $r$-th factor of $\R^k_+$ to be the distance from $Y_r$ for every $r$, with the differential full rank everywhere.
    
    This latter condition implies that
    \[d\iota(TY)^\perp \cap T(\R^k_+)_j \times T\hat{V} = (d\iota_j)(TY_j)^\perp.\]
\end{lemma}
\begin{proof}
    In a small neighborhood of $Y_j$, we choose the map to the $j$-th coordinate of $\R^k_+$ be the distance to $Y_j$, and the other coordinates are given by composition to the projection to $Y_j$ with the corresponding component of $\iota_j$. This defines a submersion to $\R^k_+$ on an open neighborhood of $\partial Y_j$, which we then extend to a submersion mapping the rest of $Y$ to the interior of $\R^k_+$ using the relative $h$-principle for holonomic approximation \cite{eliashberg2002introduction} (this uses the condition regarding the existence of the framing $\eta$). We define the map to $\hat{V}$ via an arbitrary smooth extension of the compositions of $\iota_j$ with the projection $\hat{V}$. To define the map to $\hat{V}'$, we follow the standard proof of the equivariant embedding theorem. 
\end{proof}

\subsection{Parameterizations of framings}
\label{sec:embedding-framings}

When defining the framing of a virtually smoooth flow category  (Definition \ref{def:framing-of-flow-category}), we are forced to introduce maps $\tau_{xyz}$ \eqref{eq:little-linear-iso} which relate the vector spaces serving as the targets of the trivializations of the tangent bundles $TT(x,y)$ various thickenings $T(x,y)$, in order to make sense of the compatibility conditions between the framings of these tangent bundles. Moreover, we will see that when working with virtually smooth flow categories, we will be forced to \emph{stabilize} the Kuranishi charts comprising the flow category in a \emph{compatible} way. This linear algebra data organizing the compatibility conditions can be challenging to work with, especially as one changes groups and performs various stabilizations. In particular, when comparing two \emph{different} virtually smooth flow categories, one is forced to \emph{compare} this linear algebra data, which in general is a combinatorically elaborate process.  In this section, we introduce some terminology regarding different ways of organizing this linear algebra data, and prove some helpful results. In particular, we introduce the notion of a \emph{free parameterization} (Definition \ref{def:free-parameterizations}) which is \emph{fundamental} to the arguments used in this paper. Free parameterizations are completely controlled by their associated numerical data (Lemma \ref{lemma:free-embeddings-are-controlled-numerically}), and freeness is preserved under group restriction, making free parameterizations easy to work with. A mild extension of this notion is a \emph{semi-free} parameterization, which is a direct sum of a free parameterization and a canonical parameterization, and has the same convenient properties as free parameterizations. \emph{All virtually smooth flow categories in this paper will be have framings parameterized by semi-free parameterizations.} When writing this paper, we found that the systematic use of (semi-)free parameterizations dramatically simplified various constructions with flow categories, avoiding many technical difficulties that arose with other approaches. 

\paragraph{Representations of posets.}\mbox{}\\
Let $\Gamma$ be finite poset. We will write $x \in \Gamma$ to denote that $x$ is an object of $\Gamma$.
\begin{definition}
\label{def:representation-of-poset}
 A representation $\mathcal{E}$ of $\Gamma$ is a assignment of inner product spaces $\{V(x)\}_{x \in \Gamma}$, together with isometric embeddings $\hat{\alpha}_{xy}: V(y) \to V(x)$ for each $x > y$ such that if $x > y > z$ then the embedding of $V(z) \to V(x)$ agrees with the embedding $V(x) \to  V(y) \to V(z)$.   In other words, a representation $\mathcal{E}$ of $\Gamma$ is a functor from $\Gamma$ to the category of inner product spaces and isometric embeddings of such. A morphism of representations is a natural transformation of functors.
 
 We will often refer to a representation of $\Gamma$ via the notation $\{V(x)\}_{x \in \Gamma}$, and supressing the mention of the maps $\hat{\alpha}_{xy}$.
\end{definition}

\begin{definition}
A \emph{free representation} of $\Gamma$ is a representation of $\Gamma$ satisfying the following condition. For $x \in \Gamma$, write $V(x)^* \subset V(x)$ to the perpendicular of the sum of the images of all $V(y)$ for $y<x$. A representation of $\Gamma$ is free when 
\begin{equation}
    \label{eq:decomposition-of-a-free-poset-rep}
    V(x) = \oplus_{x \geq y} V(y)^*
\end{equation}
where we have identified $V(y)^*$ with its image in $V(x)$ under the embedding $V(y) \to V(x)$. 
\end{definition}

\begin{definition} Given a representation $\{V(x)\}$ of $\Gamma$, a stabilization of $\{V(x)\}$ is a representation isomorphic to some representation $\{V'(x)\}$, where 
\begin{equation}
    \label{eq:stabilization-of-poset-rep} 
    V'(x) = V(x) \oplus \bigoplus_{y \leq x} \tilde{V}(y)
\end{equation}
for some arbitrary choice of vector spaces $\{\tilde{V}(y)\}_{y \in \Gamma}$, with the inclusions $V'(y) \to V'(x)$ for $x>y$ given by the direct sum of the inclusion map $V(x) \to V(y)$ with the natural inclusions of the $\tilde{V}(y)$ into the direct sum. When the vector spaces $\tilde{V}(y)$ are nonzero exactly for $y \in \Gamma' \subset \Gamma$ we say that this is a stabilization \emph{at} $\Gamma'$. 
\end{definition}

Let $G$ be a finite group.
\begin{definition}
    A $G$-poset is a poset $\Gamma$ with an action of $G$ on the objects of $\Gamma$, such that if $x \geq y$ then $gx \geq gy$.
\end{definition}

\begin{definition}
\label{def:equivariant-F_2-framing}
    A representation of a $G$-poset is a representation $\{V(x)\}_{x \in \Gamma}$ of the underlying poset, together with isometries 
    $g_x: V(x) \to V(gx)$ for every $g \in G$ and $x \in \Gamma$ such that 
    \begin{equation}
    \label{eq:equivariance-conditions-for-poset-rep-1}
    h_{gx} g_x = (hg)_{x} \text{ for } g, h \in G, x \in \Gamma, id_x = id,  \text{ and } 
    \end{equation}
    \begin{equation}
        \label{eq:equivariance-conditions-for-poset-rep-2}
     g_x \hat{\alpha}_{xy} = \hat{\alpha}_{gx, gy} g_y \text{ for } g \in G, x,y \in \Gamma, \text{ and } x > y,
    \end{equation}
    where $\hat{\alpha}_{xy}: V(y) \to V(x)$ is the map associated to the representation of the underlying poset.
\end{definition}

\begin{definition}
    A representation of a $G$-poset is free when it is free as a representation of the underlying poset. A stabilization of a representation $\{V(x)\}$ of a $G$-poset is defined in the same way as a stabilization of the representation of the underlying poset, but for which the action of $g_x$ on $V'(x)$ resticts to the action of $g_x$ on $V(x)$ and sends $\tilde{V}(x)$ isometrically to $\tilde{V}(gx)$ in the decompositions of \eqref{eq:stabilization-of-poset-rep}.
\end{definition}

\begin{definition}
    A direct sum $\mathcal{V} \oplus \mathcal{W}$ of representations $\mathcal{V} = \{V(x)\}, \mathcal{W} = \{W(x)\}$ of a $G$-poset $\Gamma$ assigns to $x \in \Gamma$ the vector space $V(x) \oplus W(x)$ and the $x \geq y$ the direct sum of the embeddings associated by $\mathcal{V}$ and by $\mathcal{W}$.  
\end{definition}

The following lemma is immediate:
\begin{lemma}
    Stabilization of a representation of a $G$-poset $\Gamma$ is the same as direct sum with a corresponding free representation of $\Gamma$.     Free representations of a $G$-poset are exactly the stabilizations of the zero representation. 
$\blacksquare$
\end{lemma}

Finally, there is a convenient way of constructing stabilizations of a $G$-poset:

\begin{lemma}
\label{lemma:how-to-stabilize-equivariantly}
    Given a $G$-poset $\Gamma$, a choice of $x \in \Gamma$, and a $G_x$-representation $\tilde{V}(x)$, there is a free representation of $\Gamma$ which is, as a representation of the underlying poset, the stabilization of the zero representation exactly at the $G$-orbit of $x$. This free representation of the $G$-poset is unique up to isomorphism. 
\end{lemma}
\begin{proof}
    All elements of the $G$-orbit of $x$ are necessarily incomparable. The construction of the desired free representation then follows from the construction of induced representations: we choose representatives $g_1, \ldots, g_r$ of $G/G_x$, and we ask for the the representation to assign $g_i \tilde{V}(x)$ to $g_ix$. If $g_i x = x'$ and $gx' = g_j x$, then set the maps $g_{x'}: g_i\tilde{V}(x) \to g_j \tilde{v}(x)$ for $x' \in G_x$ to be such that $\oplus_{x' \in Gx} g_{x'}$ is the action of $g$ on the induced representation $\ind_{G_x}^G \tilde{V}(x)$. The uniqueness up to isomorphism of induced representations translates into the correspondig statement for this free representation. 
\end{proof}

\paragraph{Parameterizations of framings.}

\begin{definition}
    An $E'$-parameterization associated to an (equivariant) flow category $\CC$ is a representation of the underlying poset of $\CC$. This $E'$-parameterization is an $E$-parameterization when the corresponding representation is free.  An $E'_2$-parameterization $\mathcal{E}$ for $\CC$ is a pair of $E'$-parameterizations $(\mathcal{E}_1, \mathcal{E}_2)$ for $\CC$, and similarly an $E_2$-parameterization is a pair of $E$-parameterizations.
\end{definition}

Below, we characterize the abstract linear algebra data parameterizing codomains of the framings of the tangent spaces of thickenings and obstruction bundles of a framed flow category. 

\begin{definition}
\label{def:parameterization-of-a-half-framing}
    Let $\Gamma$ be a $G$-poset. An $F'$-parameterization of $\Gamma$ is an assignment of inner product spaces $V(x,y)$ for each pair of elements $x>y$ of $\Gamma$, as well as linear isometries
    \begin{equation}
        \tau_{xyz}: V(x,y) \oplus V(y,z) \to V(x,z)
    \end{equation}
    such that for every quadruple of objects $x>y>z>t$ of $\Gamma$, and for every triple $(a,b,c) \in V(x,y) \oplus V(y,z) \oplus V(z,t)$ we have   
    \[\tau_{xyt}(a, \tau_{yzt}(b,c,0), 0) = \tau_{xzt}(\tau_{xyz}(a,b, 0), c, 0). \]
    Moreover, there for all $g \in G$, there are isometries
    \begin{equation}
        g_{xy}: V(x,y) \to V(gx,gy) \text{ for all }x,y \in \Gamma, x>y, 
    \end{equation}
    satisfying 
    \[ g'_{gx,gy}g_{x,y} = (g'g)_{x,y} \text{ for } g, g' \in G,\]
    \[ (id)_{xy} = id, id_{xyz} = id,\]
    \[ \tau_{gx,gy,gz}(g_{xy} \oplus g_{yz}) = g_{xz}\tau_{xyz} \text{ for all } g \in G,\]
    where in the above we either refer to conditions ranging over all pairs $x>y$ of $\Gamma$ or over all triples $x>y>z$ of $\Gamma$.

    An isomorphism of $F'$-parameterizations with underlying vector spaces $(\{V(x,y)\}, \{W(x,y)\})$ is a family of isomorphisms $V(x,y) \to W(x,y)$ commuting with the maps $\tau_{xyz}$ and $g_{xy}$ in the obvious way. 

    An $F'_2$-parameterization $\mathcal{F}$ of $\Gamma$ is a pair $(\mathcal{F}_1, \mathcal{F}_2)$ of $F'_2$-parameterizations of $\Gamma$. 

    A direct sum $\mathcal{F}_1 \oplus \mathcal{F}_2$ of $F'_2$-parameterizations $(\mathcal{F}_1, \mathcal{F}_2)$ with underlying vector spaces $(\{V(x,y)\}, \{W(x,y)\})$ is the $F'$-parameterization assigns to $x>y$ the vector space $V(x,y) \oplus W(x,y)$ and maps maps $\tau_{xyz}$ and $g_{xy}$ the direct sum of the corresponding maps for $\mathcal{F}_1$ and $\mathcal{F}_2$. A direct sum of $F'_2$-parameterizations $\mathcal{F}_1 \oplus \mathcal{F}_2$ with $\mathcal{F}_a = (\mathcal{F}_{a,0}, \mathcal{F}_{a,1})$ for $a=1,2$, is $(\mathcal{F}_{1, 0} \oplus \mathcal{F}_{2, 0}, \mathcal{F}_{1, 1} \oplus \mathcal{F}_{2,1})$. 
\end{definition}

\begin{lemma}
\label{lemma:framed-flow-categories-give-F_2'-params}
    Let $\CC'$ be a framed $G$-equivariant flow category. Then, in the notation of Definition \ref{def:equivariant-framing}, the tuples $(\{V^0(x)\}, \{W^0(x,y,z)\}, \{\tau^0_{xyz}\})$ and $(\{V^1(x)\}, \{W^0(x,y,z)\}, \{\tau^1_{xyz}\})$ each define an $F_2'$-parameterization. $\blacksquare$
\end{lemma}

\begin{definition}
    An embedding an $F'$-parameterization $\mathcal{F}$ of $\Gamma$ is an $E'$-parameterization  $\mathcal{E} = \{V(x)\}$ of $\Gamma$ together with a set of isometries 
    \[ \bar{\alpha}_{xy}: V(x,y) \oplus V(y) \to V(x), \alpha_{xy} = \bar{\alpha}_{xy}|_{V(x,y)}, \hat{\alpha}_{xy} = \bar{\alpha}_{xy}|_{V(y)}\]
    for every pair of elements $x>y$ of $\Gamma$  (where we are using the notation of Definitions \ref{def:parameterization-of-a-half-framing} and \ref{def:representation-of-poset}), which extend the embeddings $\hat{\alpha}_{xy}: V(y) \to V(x)$ of the $E'$-parameterization, and such that that these embeddings satisfy the diagram \eqref{eq:embedding-of-framing-compatibility} with $V^i$ replaced by $V$ and $W^0$ replaced by $W$. We say that this is an embedding of $\mathcal{F}$ into $\mathcal{E}$. 

    An embedding of an $F'_2$-parameterization $(\mathcal{F}_0, \mathcal{F}_1)$ is is a pair of embeddings of $\mathcal{F}_i$ into $\mathcal{E}_i$ for $i=0,1$, respectively. We say this embedding has underlying $E'_2$-parameterization $(\mathcal{E}_0, \mathcal{E}_1)$. 
    
\end{definition}

\begin{lemma}
\label{lemma:poset-rep-determines-half-framing}
    An $(E', E'_2)$-parameterization $\mathcal{E}$ determines an associated $(F', F'_2)$-parameterization $\mathcal{F}$ together with an embedding of $\mathcal{F}$ into $\mathcal{E}$.  
\end{lemma}
\begin{proof}
    For each $x>y$ objects of $\Gamma$, define $V(x,y)$ to be the perpendicular of the image of $V(y)$ in $V(x)$. For each $x>y>z$ objects of $\Gamma$, set $W(x,y,z)$ to be zero. For each such triple, define 
    \[ \tau_{xyz}: V(x,y) \oplus V(y,z) \to V(x,z) \]
    to be map that acts by the injective isometry $V(y) \to V(x)$ on the second factor and then adds the two elements. For $g \in G$ and $x>y$ objects of $\Gamma$, define $g_{xy}: V(x,y) \to V(gx,gy)$ to be the restrictions of the maps $g_x: V(x) \to V(y)$ (this maps sends $V(x,y)$ onto $V(gx,gy)$ because $g_y$ sends $V(y)$ onto $V(gy)$ and because $g_x \hat{\alpha}_{xy} = \hat{\alpha}_{gx,gy}g_y$). It is straightforward to see that these maps satisfy the required conditions.
\end{proof}

\begin{definition}
\label{def:free-parameterizations}
    We say that an $(F', F'_2)$-parameterization is an $(F, F_2)$-parameterization, or alternatively that it is \emph{free}, if it is isomorphic to the $(F', F'_2)$-parameterization associated to \emph{some} $(E, E_2)$-parameterization (note the lack of primes!) 
\end{definition}

\begin{definition}[Restriction and Fixed Points.]
\label{def:restriction-and-fixed-points-of-parameterization}
Let $\Gamma$ be a $G$-poset and let $H \subset G$ be a subgroup. Then $\Gamma$ is an $H$-poset under the inclusion $H \subset G$. 

Let $\Gamma^H$ is the poset of elements fixed by $H$. If $H$ is normal, this is naturally a $G/H$-poset. 

The \emph{restriction} of an $(E',E'_2,F', F'_2)$-parameterization of $\Gamma$ from $G$ to $H$ is the corresponding representation of $\Gamma$ thought of as an $H$-poset under $H \subset G$ obtained by forgetting the maps $g_{xy}$ for $g \notin H$. 

The \emph{$H$-fixed points} of an $(E',E'_2,F', F'_2)$-parameterization  $\mathcal{G}$ of $\Gamma$  are the $(E',E'_2,F', F'_2)$-parameterization $\mathcal{G}^H$ of $\Gamma^H$ obtained by taking $H$-fixed points of all assigned vector spaces and maps. If $H$ was normal, this is naturally an $(E',E'_2,F', F'_2)$-parameterization of $\Gamma^H$ as a $G/H$-poset. 
\end{definition}

\begin{lemma}
    Freeness of all $(E',E'_2,F', F'_2)$-parameterizations is preserved by restriction and taking fixed points. $\blacksquare$
\end{lemma}

\begin{lemma}
Given an embedding of $(F', F'_2)$-parameterization $\mathcal{F}$ into a $(E', E'_2$) parameterization $\mathcal{E}$, the $(F', F'_2)$ parameterization produced by Lemma \ref{lemma:poset-rep-determines-half-framing} from $\mathcal{E}$ is isomorphic to $\mathcal{F}$.   
\end{lemma}
\begin{proof}
    This follows because the condition that the embedding is adapted means that for every $x>y$ objects of $\Gamma$, the map $\bar{\alpha}_{xy}$ restricted to $V(x,y)$ must send $V(x,y)$ isometrically onto the perpendicular to the image of $V(y)$ under the isometric embedding associated to the poset representation, the latter being simply the restriction of $\bar{\alpha}_{xy}$ to $V(y)$.
\end{proof}

By \eqref{eq:decomposition-of-a-free-poset-rep}, the vector spaces underlying an $F$-parameterization can be written like 
\[ V(x,y) = \oplus_{x \geq z, y \ngeq z} V^*(y).\]
where $V^*$ are produced from the corresponding $\mathcal{E}$-representation as in \eqref{eq:decomposition-of-a-free-poset-rep}
with the maps $\tau_{xyz}$ for $x>y>z$ objects of $\Gamma$ given by the natural inclusion of direct sums.

\begin{definition}
    Given a representation $\{V(x)\}$ of a $G$-poset $\Gamma$, one can associate a map 
    \[ \mu(x): \Gamma \to \cup_{H \subset G} RO(H)\]
    such that $\mu(x) \in RO(Stab_x)$ gives the class of $V(x)$. This map is called the \emph{numerical invariant} of the representation, and will have the property that $\mu(gx)$ is the image of $\mu(x)$ under the map $RO(Stab_x) \to RO(Stab_{gx})$ induced by conjugation the isomorphism $Stab_x \to Stab_{gx}$ given by conjugation by $g$.
\end{definition}
\begin{lemma}
\label{lemma:free-embeddings-are-controlled-numerically}
    A pair of $E_2$-parameterizations $\{V(x)\}, \{V'(x)\}$ with the same numerical invariants $\mu(x), \mu'(x)$ are isomorphic. 
\end{lemma}
\begin{proof}
First note that for each each $x$, $\oplus_{y \in Gx} V(y)$ is a $G$-representation via the action of $\oplus_{y \in Gx} g_y$. Moreover, this $G$-representation is isomorphic to the induction of $V(x)$ as a $Stab_x$-representation from $Stab_x$ to $G$. This representation has the property, in particular, that for any other subgroup $H \subset G$, $\oplus_{y \in Hx} V(y)$ is an $H$-representation.

Fix $x \in \Gamma$. Assume (by inductive hypothesis) that for every $G$-orbit $o$ of $\Gamma$ which contains some element that is strictly less than $x$ in $\Gamma$, we have fixed a distinguished element $y \in o$ (which may or may not be less than $x$); denote the set of distinguished elements as $Y \subset \Gamma$. Assume further that for each distinguished element $y \in Y$ that there are $Stab_y$-equivariant isomorphisms $\gamma_y: V^*(y) \simeq (V')^*(y)$ for all $y\in \Gamma$ such that $y<x$.

Using this $\gamma_y$, we can, for any other element $y' \in Gy$, define a map 
\[ \gamma_{y'}: V^*(y') \to (V')^*(y'), g_y\gamma_{y'} = \gamma_{y'}g^{-1}_{y'} \text{ where } y'=gy \text{ for some } g \in G.\]
These maps satisfy the following properties:
\begin{itemize}
    \item They are well defined, i.e. the above formula defines a map that is independent of the choice of $g \in G$.
    \item The direct sum map 
    \[ \bigoplus_{y' \in Gy} \gamma_{y'}: \bigoplus_{y \in G'} V^*(y') \to \bigoplus_{y' \in Gy} (V')^*(y')\]
    is $G$-equivariant, and
    \item For any $H$ subgroup $H$ of $G$, letting $o'$ be an $H$-orbit in $Gy$ (which may not contain $y$),
    the map 
    \[\bigoplus_{y' \in o'}  \gamma_{y'}: \bigoplus_{y' \in o'} V^*(y') \to \bigoplus_{y' \in o'} (V')^*(y')\]
    is an $H$-equivariant isomorphism.
\end{itemize}

For any $y < x$, writing 
\begin{equation}
\label{eq:equivariant-decompostion}
V(y) = \bigoplus_{z \in \Gamma, y > z} V^*(z) = \bigoplus_{o \in \{z \in \Gamma, y > z\}/Stab_y} \bigoplus_{z' \in o} V^*(z)
\end{equation}
and similarly for $V'(y)$, we see that the corresponding direct sum of  maps $\bigoplus_{z \in \Gamma, y > z} \gamma_{z}$ defines a $Stab_y$-equivariant isomorphism $V(y) \to V'(y)$.

We now turn to $x$. If we have already chosen a distinguished element $x' \in Gx \cap Y$ and found an isomorphism $\gamma_{x'}$, then the above discussion defines $\gamma_x$ and thus a $Stab_x$-equivariant isomorphism $V(x) \to V'(x)$. If we have not, we will add $x$ to $Y$ and show that a map $\gamma_x$ as above exists. Indeed, the decomposition \eqref{eq:equivariant-decompostion}, now applied to  $V(x)$, and the properties listed above of the maps $\alpha_y$ for $y<x$  shows that $[V^*(x)] = [(V')^*(x)]$ in $RO(Stab_x)$. But then a choice of $\alpha_x$ exists; we pick any one and continue the induction.
\end{proof}

\begin{remark}
    \label{rk:explanation-of-free-parameterizations}
    Lemma \ref{lemma:free-embeddings-are-controlled-numerically} is the fundamental reason we always use (semi-)free parameterizations throughout this paper; it is the minimal formalization of the statement that the parameterizations of the framings, at least when they are taken to be (semi-)free, are entirely homotopically auxiliary data (which is suggested in the discussion of Section \ref{sec:flow-categories-framings} about parameterizatoins). If one does not impose a condition like semi-freeness such that some analog of Lemma  \ref{lemma:free-embeddings-are-controlled-numerically} holds, it becomes exteremely difficult to work with the resulting homotopy types, and the \emph{space of parameterizations} begins to play a very complicated role in all constructions. The author is entirely unsure, for example, of how to show that two parameterizations are deformation equivalent after sufficient stabilization; while Lemma \ref{lemma:free-embeddings-are-controlled-numerically} makes this completely trivial.  Beyond that, statements like the existence of shift spaces (see Definition \ref{def:shift-spaces} below) completely fail to be true for parameterizations in the general case (one can easily find an counterexample where $\Gamma$ has four objects), making it very difficult to add new objects to a flow category. The convenient language of (semi-)free parameterizations and the resulting (elementary) theory completely avoids all these technicalities hiding in the background. 
\end{remark}

\begin{remark}
    The importance of (semi)-free parameterization is the fact that they admit shift spaces (discussed below) together with the fact that they are preserved under group restriction \emph{and} under taking invariants; simple generalizations of canonical parameterizations to the equivariant setting (obtained e.g. by choosing orderings on irreducible representations of $G$) generally do not satisfy all of these properties at once.
\end{remark}

\begin{lemma}
\label{lemma:numerical-condition-for-existence-of-e-param}
    Given vector spaces $\{V(x)\}_{x \in \Gamma}$ and isometries $g_x: V(x) \to V(gx)$ satisfying \eqref{eq:equivariance-conditions-for-poset-rep-1}, define the classes $V^*(x) \in RO(G_x)$ inductively as follows. If $x$ is a minimum of $\Gamma$, set $V^*(x) = V(x)$. Otherwise,  
    If $V(y')$ is defined as $y'$ ranges over $y' \in Gy$, write
     \[\bigoplus_{y'\in G_x y} V^*(y') \in RO(G_x) \] for the corresponding virtual $G_x$-representation isomorphic to $Ind_{G_{xy}}^{G_x} Res^{G_y}_{G_{xy}}$.
     
    If $V^*(y)$ is defined for all $y < x$, define
    \begin{equation}
        \label{eq:numerical-condition-for-existence-of-e-param}
        V^*(x) = V^0(x) - \bigoplus_{y < x} V^*(y) \in RO(G_x).
    \end{equation}

    Suppose that $V^*(x)$,
    which is a-priori virtual, is represented by  an actual $G_x$-representation, for every $x$. Then the vector spaces $\{V(x)\}$ and operators $g_x$ are the underlying vector spaces and operators of an equivariant $E$-parameterization of $\Gamma$.
\end{lemma}
\begin{proof}
    We inductively define the embeddings $\hat{\alpha}_{xy}$. We induct upwards on the poset $\Gamma/G$. For the minima of $\Gamma/G$ there is nothing to do. Suppose we have defined an $E$-parameterization of the pre-image $\Gamma'$ of a downwards closed sub-poset $\Gamma'/G$ of $\Gamma/G$ under $\Gamma \to \Gamma/G$, which has underlying vector spaces and operators the ones specified.  Pick an element $[x] \in \Gamma/G$ such that all elements less than this element are in $\Gamma'/G$, and a representative $x \in [x]$. Then by the condition we can define a $G_x$-equivariant embedding of $\bigoplus_{y' < x}/G V^*(y')$ into $V(x)$. Fix the representative of the class \eqref{eq:numerical-condition-for-existence-of-e-param} to be the perpendicular of this image. Restricting these embeddings to the sums of the $V^*(y')$ contained in $V(y)$ for $y<x$ defines the embedding $\hat{\alpha}_{xy}$. One then sets $\hat{\alpha}_{gx gy} = g_x \hat{\alpha}_{xy} g^{-1}_{y}$.
\end{proof}

\begin{definition}[Shift spaces.]
\label{def:shift-spaces}
To an $E'$-parameterization $\mathcal{E}$ of a $G$-poset $\Gamma$ can associated a canonical $G$-representation called the \emph{shift space} of $\mathcal{E}$, characterized by the following property. Write $\Gamma' = \Gamma \cup m$ for the poset with $m$ a new object that is strictly larger than all previous objects of $\Gamma$, and otherwise the same partial order on the objects of $\Gamma \subset \Gamma'$; with the group action fixing $m$. The \emph{shift space} of $\mathcal{E}$, when it exists, is the orthogonal $G$-representation assigned by an extension $\mathcal{E}_m$ of $\mathcal{E}$ up to $\Gamma_m$ to thhe object $m$ which is initial among such extensions $\mathcal{E}_m$ in the category of $E'$-parameterizations of $\Gamma'$ with morphisms $\{V(x)\} \to \{W(x)\}$ given by object-wise injective natural transformations $f_x: V(x) \to W(x)$ of functors commuting with the $G$ action i.e. such that $f_{gx}g_x=g_xf_x$ for all $g \in G$ and $x \in \Gamma$. This $G$-representation is manifestly unique up to unique isomorphism in the category of orthogonal $G$-representations, but it may or may not exist.
\end{definition}

\begin{remark}
    \label{rk:shift-spaces-for-categories}
    In the construction of the Floer homotopy type (Definition \ref{def:floer-homotopy-type}) of a virtually smooth flow category $\CC'$, we use the shift space of the parameterization of the framing of the obstruction bundles and we desuspend the Floer homotopy space by this shift space to get ther Floer homotopy type. We refer to this shift space as the shift space of $\CC'$. 
\end{remark}

\begin{remark}
Shift spaces will exist for semi-free $E'$-representations, and are computed in Lemma \ref{lemma:shift-space-exists-for-free} and Lemma \ref{lemma:shift-space-exists} below; whenever we refer to \emph{the} shift space of a $\mathcal{E}'$ representation we use the formulae in those lemmas as a preferred choice of shift space.
\end{remark}

\begin{lemma}[Shift spaces for free parameterizations.]
\label{lemma:shift-space-exists-for-free}
    Suppose $\mathcal{E} = \{V(x)\}$ is a free (equivariant) $\Gamma$-representation. Then the shift space of $\Gamma$ exists and is given by 
    \begin{equation}
        V_{\mathcal{E}} =  \oplus_{x \in \Gamma} V^*(x). \blacksquare
    \end{equation}
\end{lemma}

\paragraph{Canonical and semi-free parameterizations.}

\begin{definition}
\label{def:canonical-parameterization}
    Let $\Gamma$ be a $G-$poset equipped with an \emph{integral action}, i.e. a map $\bar{A}: \Gamma \to \{n \in \Z : n \geq -1\}$ such that if $x > y$ in $\Gamma$ then $\bar{A}(x) > \bar{A}(y)$. To this data, we can associate the set $S_{\bar{A}}$ as in Definition \ref{def:integral-action-set}. 

    A \emph{canonical parameterization} of $\Gamma$ (thought of as being equipped with $\bar{A}$) is an assignment of a natural number $f(n)$ to each element of $S_{\bar{A}}$. This defines an $E'$-parameterization of $\Gamma$ via
    \[ x \mapsto \R^{f(\bar{A}(x))} \oplus \R^{f(\bar{A}(x)-1)} \oplus \ldots \oplus \R^{f(0)}. \]
    Here, the $G$-action is everywhere \emph{trivial}, and the maps associated to $x > y$ in $\Gamma$ are just natural inclusion into the right-most factors. We call this $E'$-parameterization or its associated $F'$-parameterization \emph{canonical}; we call an ($E'_2$ or $F'_2$) parameterization for which each element of the pair is canonical, canonical. 

    When a parameterization has a specified isomorphism with a canonical parameterization, we will also refer to it as canonical. 
\end{definition}

The proof of the following useful lemma is trivial:
\begin{lemma}
    The direct sum $\mathcal{E}_1 \oplus \mathcal{E}_2$ of a pair of canonical $E'$-parameterization $\mathcal{E}_1 = \{V_1(x)\}$ and  $\mathcal{E}_2 = \{V_2(x)\}$ defined by functions $f_1$ and $f_2$, respectively, is isomorphic to the canonical $E'$-parameterization $\mathcal{E}_{12}$ with defining function $f = f_1 + f_2$, via the reordering maps
    \[ \left(\R^{f_1(\bar{A}(x))} \oplus \cdots \oplus \R^{f_1(0)}\right)\oplus  \left(\R^{f_2(\bar{A}(x))} \oplus \cdots \oplus \R^{f_2(0)}\right) \to  \left(\R^{f_1(\bar{A}(x))} \oplus \R^{f_2(\bar{A}(x))} \oplus \cdots \oplus \R^{f_1(0)} \oplus\R^{f_2(0)}\right)\]
    where grouping pairs of factors on the right hand side of course gives $\mathcal{E}_{12}$. Whenever we take a direct sum of a pair of canonical parameterizations, we will view it as a canonical parameterization via the isomorphism induced by the one specified above. $\blacksquare$
\end{lemma}

\begin{definition}
\label{def:semi-free-parameterization}
    The direct sum $\mathcal{E}_f \oplus \mathcal{E}_c$ of a canonical $E'$-parameterization $\mathcal{E}_c$ with a free $E$-parameterization $\mathcal{E}_f$ will be referred to as an $E^s$ parameterization; the associated $F'$-parameterization will be referred to as an $F^s$-parameterization, and an $F'_2$ parameterization for which each element of the pair is an $F^s$-parameterization will be called an $F^s_2$-parameterization.  We will use `semi-free parameterization' to refer to any of an $E^s-, F^s-$, or $F^s_2$-parameterization.
\end{definition}

\begin{lemma}
\label{lemma:shift-space-exists}
    Any semi-free representation $\mathcal{E}_f \oplus \mathcal{E}_c$ (where $\mathcal{E}_c$ is a canonical representation defined by the function $f$) has a shift space given by 
    \[ V_{\mathcal{E}} = V_{\mathcal{E}_f} \oplus \R^{\sum_{i \in S_{\bar{A}}} f(i)} \]
    with $V_{\mathcal{E}_f}$ as in Lemma \ref{lemma:shift-space-exists-for-free} and the $G$-action is trivial on the second factor above. $\blacksquare$
\end{lemma}

\begin{remark}
    Note that the factors $T\R_+^{\bar{A}(x,y)}$ in the definition of a framing of a virtually smooth category can be viewed as the terms of a canonical $F'$-parameterization, and thus if the flow category is freely parameterized then the terms $T\R_+^{\bar{A}(x,y)} \times V^0(x,y)$ giving the targets of the framings of the thickenings $T(x,y)$ can be viewed as the terms of an $F^s$-parameterization.  Thus, even if the parameterization associated to a framing is free, in some sense the notion of a framing almost forces the appearance of semi-free parameterizations in the general discussion. However, as mentioned in Remark \ref{rk:why-canonical-parameterizations-are-not-enough}, we cannot stick to the use of canonical parameterizations, as previous works do, because if the importance of being able to easily change the symmetry group $G$. 
\end{remark}

\
\begin{definition}[Stabilization of canonical parameterizations.] We can \emph{stabilize} a canonical $E'$-parameterization $\mathcal{E}_1$ with defining function $f_1$ another canonical $E'$-parameterization $\mathcal{E}_2$ with defining function $f_2$ by specifying maps $\R^{f_1(n)} \to \R^{f_2(n)}$ for each $n \in S_{\bar{A}}$ which include the domain as the set where a chosen subset of the $f_2(n)$ coordinates of size $f_2(n)-f_1(n)$ are required to be zero. It is easy to see that $\mathcal{E}_2$ is the direct sum of $\mathcal{E}_1$ with the canonical $E'$-parameterization with defining function $f_2 - f_1$; the inclusions specified above fix the isomorphism with the direct sum. This notion defines the corresponding notion of stabilization of canonical $F'$ , $E'_2$,and $F'_2$ parameterizations. 

We will refer to a stabilization of an any semi-free parameterization $M$ as one that is obtained by first stabilizing the free part of $M$ and subsequently stabilizing the canonically parameterized part  (of course, these stabilization operations can be performed in either order, up to canonical isomorphism). 
\end{definition}

\paragraph{Linear parameterization of stabilizations.}

In this paper, we are largely able to avoid discussion of general kinds of stabilizations of flow categories. However, we do have to perform one particularly simple kind of stabilization, which we describe below.

\begin{definition}
\label{def:stabilization-by-f-parameterization}
    A stabilization of a framed virtually smooth flow category $(\CC, \CC')$ by an $F'$-parameterization $\mathcal{F}$ of $\CC$ (which is required to be $G$-equivariant if $\CC$ is) with underlying vector spaces $\{V(x,y)\}$ is the framed virtually smooth flow category $\CC'_\mathcal{F}$, which has the same objects as $\CC'$, satisfies $\CC'_\mathcal{F}(x,y) = \CC'(x,y)_{V(x,y)}$ (Definition \ref{def:stabilize-derived-manifold}) and with the composition maps 
    \[(f_\mathcal{F}, \bar{f}_\mathcal{F}, \bar{f'}_\mathcal{F}): \CC'_\mathcal{F}(x,y) \times \CC'_\mathcal{F}(y,z) \to \CC'_\mathcal{F}(x,z)\]
    induced by the maps $\tau_{xyz}$ covering the compositions
    \[(f, \bar{f}, \bar{f'}): \CC'(x,y) \times \CC'(y,z) \to \CC'(x,z).\]
    If the underlying $F'_2$-parameterization of the framing of $\CC'$ was $\mathcal{F}_{\CC'} = (\mathcal{F}_0, \mathcal{F}_1)$, then the underlying $F'_2$-parameterization of $\CC'_\mathcal{F}$ is 
    \[ \mathcal{F}_{\CC'} \boxplus \mathcal{F} :=  (\mathcal{F}_0 \oplus \mathcal{F}, \mathcal{F}_1 \oplus \mathcal{F}). \]

    If $\mathcal{F}_{\CC'}$ was an $F'_2$-parameterization arising from the $E'_2$-parameterization $\mathcal{E}_{\CC'}(\mathcal{E}_0, \mathcal{E}_1)$, and $\mathcal{F}$ is an $F'$-parameterization arising from the $E'$-parameterization $\mathcal{E}$, then  $\mathcal{F}_{\CC'} \boxplus \mathcal{F}$ is an $F_2$-parameterization arising from the $E_2$-parameterization
    \[ \mathcal{E}_\CC \boxplus \mathcal{E} = (\mathcal{E}_0 \oplus \mathcal{E}, \mathcal{E}_1 \oplus \mathcal{E}). \]
    In this case we will denote the stabilized virtual smooothing $\CC'_\mathcal{F}(x,y)$ by $\CC'_{\mathcal{E}}$ as well. 
\end{definition}

\subsection{Stabilization of embeddings and Floer homotopy types}

We have now explained how to define spaces $\{\CC\}$ associated to extensions of embeddings of framed flow categories. Certain desuspensions of these spaces are then the associated floer homotopy types $|\CC|$.  In this section, we describe the sense in which the spectra $|\CC|$ do not change (up to \emph{canonical} weak equivalance) upon \emph{stabilizations} of the flow categories and their associated embeddings. This turns out to be very helpful when trying to ``extend'' embeddings from flow subcategories to larger flow categories,  as is necessary when when defining continuation maps or any other operations on Floer homotopy types. 

Because of the way that we set up Definitions \ref{def:extension} and \ref{def:floer-homotopy-type}, we can use a very strict notion of a stabilization of an embedding of $\CC'$ and consequently prove a strict invariance statement for the floer homotopy type. Imagine $\CC'_1$ is a virtually smooth flow category with an extension of an embedding, and $\CC'_2$ is produced from $\CC'_1$ by stabilizing the category and the associated embedding. In our setup we will establish a \emph{canonical homeomorphism} of $\{\CC'_2\}$ with a suspension of $\{\CC'_1\}$. The fact that this is a homeomorphism, rather than a homotopy equivalence defined up to contractible choice, allows us to simplify some of the coherence data we need to keep track of later when constructing the cyclotomic structure on symplectic cohomology. This is possible because when defining the Floer homotopy type we use \emph{canonical} maps on the the $V^0(y)$ factors in the definition of the $\tilde{e}_{x, i}$ in \eqref{eq:define-attaching-map-on-one-stratum}. If instead we used maps which were \emph{homotopy equivalent} to these canonical maps, some some of the other constructions (e.g. when constructing or extending embeddings, as in Proposition \ref{prop:embeddings-exist} and \ref{prop:relative-variant-of-embedding-construction} below) would have been simpler, at the cost of a loss of the canonical nature of the homeomorphisms we establish below.

In any case, we proceed with the constructions. 

\begin{lemma} [Stabilizing embeddings by $E^s$-parameterizations.]
\label{lemma:stabilizing-embeddings-by-parameterizations}
    Let $\CC'$ be a flow category with an extension $\bar{\mathbf{E}} = \{e_{xy}\}$ of an embedding $\mathbf{E}=\{\iota_{xy}\}$ produced by applying Proposition \ref{prop:embeddings-exist} to $
    \CC'$. Let $\mathcal{E}$ be an $E$-parameterization of $\CC$ with associated $F$-parameterization $\mathcal{F}$. Then there is a canonical embedding $\mathbf{E}_\mathcal{E} = \{\iota_{xy}^\mathcal{E}\}$ of $\CC'_{\mathcal{F}}$ which is the \emph{stabilization} of the embedding $\mathbf{E}$ by $\mathcal{E}$, as well a corresponding extension $\bar{\mathbf{E}}_{\mathcal{E}} = \{ e_{xy}^{\mathcal{E}}\}$ of $\mathbf{E}^{\mathcal{E}}$ which we also refer to as the stabilization of $\bar{\mathbf{E}}$ by $\mathcal{E}$. These are characterized by the property that, writing $V_\mathcal{F}(x,y)$ for the vector space assigned by $\mathcal{F}$ to $x>y$ in $\CC$, we have that 
    \[ \iota_{xy}^{\mathcal{E}} = \iota_{xy} \times id_{V_\mathcal{F}(x,y)}, e^{\mathcal{E}}_{xy} = e_{xy} \wedge id_{(V_{\mathcal{F}}(x,y))_*}. \]
\end{lemma}
\begin{proof}
The only nontrivial condition to check is the disjointness condition for embeddings and for extensions of embeddings. Write $V_\mathcal{E}(x)$ for the vector space associated by $\mathcal{E}$ to $x \in \CC$. To see that the disjointness condition still holds holds, we note that the embedding maps for $\mathcal{F}$ give canonical isomorphisms 
\[ V_{\mathcal{F}}(x,y) \times V_\mathcal{E}(y) \to V_{\mathcal{E}}(x);\]
as such, the embedding maps maps in the disjointness condition for $\bar{\mathbf{E}}_{\mathcal{E}}$ look like 
\begin{equation}
\label{eq:nice-decomposition} 
(T_{\CC'}(x,y_i) \times V_\mathcal{F}(x,y_i)) \times (V^0_{\CC'}(y_i) \times V_\mathcal{E}(y_i)) \subset \R_+^{\bar{A}(x,y_i)} \times V^0_{\CC'}(x) \times V_{\mathcal{E}}(x).
\end{equation}
thus the disjointness condition still holds because the new embeddings are just the products of the old embeddings with $V_{\mathcal{E}}(x)$, so they are still disjoint as $y_i$ ranges over the relevant set of objects of $\CC$. 
\end{proof}

\begin{lemma}
\label{lemma:effect-of-stabilizing-by-parameterization-on-floer-homotopy-type}
    Given a framed virtually smooth flow category $\CC'$ with an extension $\bar{\mathbf{E}}$ of an embedding $\mathbf{E}$, and an $E^s$-parameterization $\mathcal{E}$ of $\CC$, writing $\{\CC\}$ and $|\CC|$ for the Floer homotopy space/type associated to this data and writing 
    $\{\CC_{\mathcal{F}}\}, |\CC_{\mathcal{F}}|$ for the Floer homotopy space/type associated to the extension $\bar{\mathbf{E}}_{\mathcal{E}}$ of of the embedding $\mathbf{E}_{\mathcal{E}}$ of $\CC_{\mathcal{F}}$, we have an equivariant homeomorphism
    \begin{equation}
        \label{eq:suspension-homoeo-floer-homotopy-space}
        \Sigma^{V^{\mathcal{E}}_{max}} \{\CC\} \simeq \{\CC_{\mathcal{F}}\}
    \end{equation} 
    where $V^{\mathcal{E}}_{max}$ the shift space assigned to $\mathcal{E}$ (see Definition \ref{def:shift-spaces}). This homeomorphism is characterized by the following property: write $(\{V^0(x)_{\CC'}\}, \{V^1(x)_{\CC'}\})$ for the $E_2^s$-parameterization associated to the framing of $\CC'$, and $V^1_{\CC', max}$ for the shift space assigned to $\{V^1(x)_{\CC'}\}.$ Then we manifestly have 
    \[V^1_{\CC'_\mathcal{F}, max} = V^1_{\CC', max} \oplus V^{\mathcal{E}}_{max}, V^0_{\CC'_{\mathcal{F}}}(x) = V^0_{\CC'}(x) \oplus V^{\mathcal{E}}(x), (V^1_{\CC'_\mathcal{F}(x)})^\perp = V^1_{\CC'}(x)^\perp \oplus V^{\mathcal{E}}(x)^\perp \]
    where $V^{\mathcal{E}}(x)^\perp$ is the perpendicular to the image of $V^{\mathcal{E}}(x)$ in $V^{\mathcal{E}}_{max}$. From this we see that there is a homeomorphism
    \[   \Sigma^{V^\mathcal{E}_{max}}(\R_+^{\bar{A}(x)+1} \times V^0_{\CC'}(x) \times V^1_{\CC'}(x)^\perp)_* \simeq (\R_+^{\bar{A}(x)+1} \times V^0_{\CC'_{\mathcal{F}}}(x) \times V^1_{\CC'_{\mathcal{F}}}(x)^\perp)_*\]
    The homeomorphism \eqref{eq:suspension-homoeo-floer-homotopy-space} is defined by requiring that it restricts to the composition of the above homeomorphism with the inclusion of its codomain into $\{\CC_{\mathcal{F}}\}$ on the image of its domain in $\Sigma^{V^\mathcal{E}_{max}}\{\CC\}$. 
\end{lemma}
\begin{proof}
    One simply needs to check that the definition of the maps $\tilde{e}_x$ for $\CC'_{\mathcal{F}}$ is the smash product of the same maps for $\CC'$ with the sphere $S^{V^{\mathcal{E}}_{max}}$. One can check by checking this for each $\tilde{e}_{x, i}$, and subsequently on  on each open $U''_{xy}$ where the definition of the map may be nontrivial. Now there are canonical isomorphisms 
    \[ V^{\mathcal{E}}(x) \oplus V^{\mathcal{E}}(x)^\perp \simeq V_{\mathcal{E}}(y) \oplus V^{\mathcal{E}}(x,y) \oplus V^{\mathcal{E}}(x)^\perp \simeq V^{\mathcal{E}}(y) \oplus V^{\mathcal{E}}(y)^\perp\]
    which commute with the isomorphisms of the left and right hand sides of this equation with $V^{\mathcal{E}}_{max}$ described in statement of the lemma; the maps on each $U''_{xy}$ for $\CC'_{\mathcal{F}}$ simply built out of the maps on $U''_{xy}$ by summing with the middle term of the above equation, analogously to \eqref{eq:nice-decomposition}, and the map is just the smash product with (the 1-point compactification) of the second map in the above equation. This verifies that the attaching maps indeed agree with the smash product of the previous attaching maps with the identity map on $V^{\mathcal{E}}_{max}$. 
\end{proof}

\subsection{Linear algebra of the Pontrjagin-Thom construction}
\label{sec:linear-algebra-of-pontrjagin-thom}

In this section, we review certain elementary aspects of the Pontrjagin-Thom construction, as preparation for the construction of embeddings of flow categories in Section \ref{sec:embeddings-of-flow-categories} below. 

Traditionally, in the Pontrjagin-Thom construction, given a manifold $M^n$ and a framing $TM \simeq \underline{\R^n}$, we embed $M$ into $\R^{m}$ and then subsequently into $\R^{m+N}$ for some large $N >> 0$ by composing with the \emph{canonical} embedding $\R^m \subset \R^{m+N}$, then then argue that the normal bundle of $M \subset \R^{m+N}$ is canonically trivialized. In fact, the normal bundle is already canonically trivialized with $N = n$. For every $x \in M$, we have a map $A: T_xM \to \R^{m}$, and an isomorphism $B: \R^n \to T_xM$. Writing $\bar{A}: T_xM \to \R^{m} \oplus \R^n$ for $\bar{A} = A \oplus 0$, we have that $\coker \bar{A} = \coker A \oplus \R^n$. But there is a canonical isomorphism $\R^m \simeq  \dom AB \oplus \coker AB$; since $\coker AB = \coker A$ and $\dom AB = \R^n$, as $A$ and $B$ vary continuously, so does the above isomorphism. The collections of these isomorphisms over $x \in M$ thus define a trivialization of the normal bundle of $M \subset \R^{m + n}$. . 

An important aspect of the Pontrjagin-Thom construction is that, in a certain precise sense that is not always completely articulated, this trivialization of the normal bundle of $M$ is compatible with stabilization. Let us contemplate replacing the embedding $\iota_{m+n}: M^n \to \R^{m+n}$ with the embedding $\iota_{m+N}: M^n \to \R^{m+N}$, $N > n$, given by composing with the canonical inclusion $\R^{m+n} \to \R^{m+N}$ into the first $m+n$ coordinates. The normal bundle $\nu_{\iota_{m+N}}$  is canonically isomorphic to $\nu_{\iota_{m+n}} \oplus \underline{\R}^{N-n}$, and we would like to say that there is a trivialization of $\nu_{\iota_{m+N}}$ which is under this isomorphism simply the trivialization of $\nu_{\iota_{m+n}}$ summed with the identity map on $\underline{\R}^{N-n}$. \emph{However, this happens only if we keep track of the image of $\R^n \subset \R^{m+n} \subset \R^{m+N}$, and use \textbf{this} copy of $\R^n$ to define the trivialization of  $\nu_{\iota_{m+N}}$}. 

We could think of the relevant copy of $\R^n$ as the `first $n$ coordinates of $\R^N = 0 \oplus \R^N \subset \R^m +N$. This uses the \emph{ordering} on the natural numbers, and this convention is one of the usual ones in textbook discussions of the Pontrjagin-Thom construction. However, when one wishes to consider \emph{equivariant} embeddings of $G$-manifolds, one will have to replace the vector spaces involved with various $G$-representations, and there is no longer a canonical ordering that can be used to track this data. 

\begin{remark}The issue above with the use of the \emph{ordering} on the natural numbers, which does not immeidately generalize to the equivariant context, is closely related to the problems related to the ordering of coordinates when defining smash products of spectra \cite{adams-blue-book}. Carefully tracking certain decompositions of vector spaces while performing Pontrjagin-Thom constructions, as we do below, is heuristically analogous to the solution to this ordering problem pursued in the construction of orthogonal spectra. 
\end{remark}

We now describe the abstract linear algebra setup of the Pontrjagin-Thom construction. Given 
\begin{itemize}
    \item inner product spaces $V, X, Y$;
    \item an inclusion $\bar{A}: V \to X$, and 
    \item an orthogonal decomposition $X \simeq Y \oplus V$ such that $Y \supset \im A$, 
\end{itemize}   
there is a canonical isomorphism $\coker A \simeq Y$ given by writing $A: V \to Y$ for the factorization of $\bar{A}$, and identifying $\coker \bar{A} \simeq (\coker A)\oplus Y \simeq Y$. Indeed, we can canonically identify $\coker A\simeq (\im A)^\perp$, and we can write the above isomorphism as $(\im \bar{A})^\perp = (\im A)^\perp \oplus Y \simeq Y$, where the last map sends $(v, w) \to v + Aw$. 

Thus, we can describe the above construction by stating that \emph{the normal bundle to $\bar{A}$ is $Y$-framed}. This latter isomorphism clearly values continuously with $A$ so long as $\im A$ continues to stay in $Y$. Moreover, if all inner product spaces and maps are upgraded to orthogonal $G$-representations and $G$-equivariant maps, the isomorphism $\coker A \simeq Y$ becomes $G$-equivariant, as well.

In the next three paragraphs, we will discuss how this map behaves under several standard operations. 

\textbf{Renaming the framing.}

Clearly, if we choose an isomorphism $Y' \simeq Y$ then the normal bundle to $\bar{A}$ will by $Y'$-framed by composing with the inverse of this isomorphism.

\textbf{Behavior under stabilization of $X$.} 

If we have an (equivariant, isometric) inclusion $X \subset X'$, then we can replace $\bar{A}$ with the composed embedding $\bar{A}': V \to X'$; we have that $(\im \bar{A})^\perp = (\im \bar{A}) \oplus X^\perp$, and replacing $Y$ with $Y' = Y \oplus X^\perp$, we have that the decomposition $X' \simeq Y' \oplus V$ still holds and thus there is an isomorphism $(\im \bar{A}')^\perp \simeq Y'$.  
(Here we use $\oplus$ to denote direct sum of subspaces rather than abstract direct sum; of course, one is canonically isomorphic to the other.) In fact, this isomorphism is just the direct sum of the original isomorphism $(\im \bar{A})^\perp \simeq Y$ with the identity map $X^\perp$.

\textbf{Products of embeddings.}
Finally, suppose that we are given a pair of embeddings $\bar{i}_1: M_1^{n_1} \subset X_1$ and $\bar{i}_2: M_2^{n_2} \subset X_2$, as well as a pair of decompositions $X_i = Y_i \oplus V_i, i=1,$ as in the previous paragraph, as well as trivializatoins $TM_i \simeq \underline{V_i}$. An isometry $P_{12}: X_1 \oplus X_2 \to X$ defines an embedding $\bar{i}_{12}: M_1 \times M_2 \to X$ by product and composition, and a decomposition $\bar{P}_{12}: \eta\bar{i}_{12} \simeq \eta d\bar{i}_1 \oplus \eta d\bar{i}_2$ by applying $P_{12}^{-1}$ to the orthogonal of the image of $d\bar{i}_{12}$. Define $Y = P_{12}(Y_1 \oplus Y_2)$. An isomorphism $P'_{12}: V \to V_1 \oplus V_2$ then defines a trivialization $T(M_1 \times M_2) \simeq \underline{V}$ by precomposition. The resulting isomorphism $\eta \bar{i}_{12} \simeq Y$ then satisfies the condition that when compositing with $P'_{12}$ and with $\bar{P}^{-1}_{12}$ it gives the sum of the trivializations of $\eta\bar{i}_1$ and $\eta\bar{i}_2$. 

\begin{remark}
There is another kind of natural kind of stabilization operation that can be performed on in this linear algebra setup, where one replaces stabilizing the domain $V$ of $X$ by a stabilized domain vector space $V \oplus W$ for some vector space $W$. In the setting of Pontrjagin-Thom constructions, this operation would arise if we wish to understand the interaction of the Pontrjagin-Thom construction on manifold $M$, thought of a smooth derived manifold with a zero dimensional obstruction bundle, with the procedure of stabilizing the manifold (and the embedding) by a framed vector bundle. However, we are able to systematically avoid the discussion of this latter operation in this paper, except in the case of codimension zero embeddings, where the situation is rather simple.
\end{remark}

\subsection{Embeddings of flow categories.}
\label{sec:embeddings-of-flow-categories}
\begin{proposition}
    \label{prop:embeddings-exist}
    Let $(\CC, \CC')$ be a finite virtually smooth framed $E'_2$-parameterized  $G$-equivariant flow category equipped with a very convenient metric. Then $\CC'$ can be stabilized by an $E$-parameterization $\mathcal{E}_{stab}$ to $\CC''$ such that a restriction $\CC'''$ of $\CC''$ admits an embedding. 
\end{proposition}

This proposition is an analog of the equivariant Pontrjagin-Thom construction for compatible families of manifolds with corners. Moreover, the proof of this proposition uses a certain \emph{equivariant double induction} strategy which is used repeatedly in the rest of the paper. 

\begin{proof}[Proof of Proposition \ref{prop:embeddings-exist}]
This follows via an \emph{equivariant double induction}. We explain the strategy in abstract for future reference. 

\paragraph{Abstract description of equivariant double induction.}

The outer loop of the induction is on downwards-closed $G$-subposets $\mathcal{P} \subset \CC$. The base case of this induction is taken to be the empty poset. 

For each step of the outer loop, one takes a \emph{minimal element} $x$ of the poset $\CC \setminus \mathcal{P}$, and writes $\mathcal{P}_x = \mathcal{P} \cup \{x\}$, $\mathcal{P}' = \mathcal{P} \cup Gx$. Write also $\mathcal{P}_{\bar{x}} \subset \mathcal{P}_x$ for the subset of elements $y$ such that $x \geq y$. After completing the inner loop, one performs an operation (generally involving \emph{transport of structure} by a natural $G/Stab(x)$-action) such that the inductive hypothesis is satisfied for $\mathcal{P}'$, and then continues the induction. 

The inner loop of the induction is on upwards-closed $G_x := Stab(x)$-subposets of $\mathcal{Q} \subset \mathcal{P}_x$. The base case is taken to be $\mathcal{Q} = \{x\}$. At each step, one chooses a maximal element $z \in \mathcal{P}_x \setminus \mathcal{Q}$, and writes $\mathcal{Q}_{z} = \mathcal{Q} \cup \{z\}$ and $\mathcal{Q}' = \mathcal{Q} \cup G_xz$. If $z \notin \mathcal{P}_{\bar{x}}$ one does nothing; otherwise, one performs an operation on the full subcategory of the current state of $\CC''$ and associated data over the objects of $\mathcal{Q}_z$; then, usually using transport of structure under a $G_x/G_{xz}$ action (where $G_{xz} = Stab(x,z) \subset G$) one achieves the inductive hypothesis for $\mathcal{Q}' \supset \mathcal{Q}$. Because something is only done for $z \in \mathcal{P}_{\bar{x}}$, we may think of this as an induction over $z \in \mathcal{P}_{\bar{x}}$, or more appropriately over $\mathcal{P}_{\bar{x}}/G_x$. 

Both inductions terminate because of finiteness of the number of objects of $\CC$. 

\paragraph{Application of inductive strategy to embeddings}.

Write $(\{V^0(x)\}, \{V^1(x)\})$ for the vector spaces the $E'$-parameterization of the framing of $\CC'$. With notation as above, during the induction we will define $E$-parameterizations $\mathcal{E}_{\mathcal{P}}$ such that if $\mathcal{P}_2 \subset \mathcal{P}$ is a downwards-closed $G$-sub-poset of $\mathcal{P}$ then $\mathcal{E}_{\mathcal{P}}|_{\mathcal{P}_2} = \mathcal{E}_{\mathcal{P}_2}$. We will write $\bar{V}_{\mathcal{P}}(x)$ for the $G_x$-representation that $\mathcal{E}_{\mathcal{P}}$ assigns to $x$ for any $\mathcal{P} \ni x$ as in the induction. Similarly, for any $y > x$ in such a $\mathcal{P}$, we will write $\delta_{xy}: \bar{V}(y) \to \bar{V}(x)$ for the corresponding embedding defined by $\mathcal{E}_{\mathcal{P}}$. 

The inductive assumption to be verified for $\mathcal{P}$ is that we have defined, for every pair of elements $x > y$ in $\mathcal{P}$, smooth embeddings
\[\iota_{xy}: T(x,y) \to \R_+^{\bar{A}(x,y)} \times im \delta_{xy}^\perp \] 
Moreover, these maps are required to satisfy the following  properties A, B, C, D defined below. 
Write $\bar{f}: T(x,y) \times T(y,z) \to T(x,z)$ for the map associated to the composition $\CC'(x,y) \times \CC'(y,z) \to \CC'(x,z)$ for objects $x > y > z$ of $\CC'(x,y)$. 
Property A is compatibility of the embeddings $\iota_{xy}$ with composition: we must have that for $x > y > z$ objects of $\mathcal{P}$, we have that 
\begin{equation}
    \label{eq:compatibility-condition-for-embeddings-of-manifolds}
    \iota_{xz} \circ \bar{f} (a, b)= \iota_{xy}(a) + \delta_{xy} \iota_{yz}(b).
\end{equation}
(where we make sense of the sum after rearranging the $\R_+$ factors to be adjacent). Property B is the equivariance condition: 
\begin{equation}
    (id \times g_x) \iota_{xz} = \iota_{xz} g_{xz}.
\end{equation}
Property $C$ is that the maps to $\R^{\bar{A}(x,y)}_+$ induced by composing the $\iota_{xy}$ with the projection have full rank differential, and near $T(x,y)_j$ the corresponding component of the map is the distance from $T(x,y)_j$. 

Finally, one requires that the analog of the disjointness condition on embeddings is satisfied for $\mathcal{P}$: for every $x \in \mathcal{P}$ and any $a < \bar{A}(x)$, we have that the closures of the sets $Im \iota_{xy} \times V^0(y)$ in $V^0(x)$  are pairwise disjoint, as $y$ runs over $\{y \in \mathcal{Q} : \bar{A}(y) = a\}$.

Assume that we have satisfied the inductive assumption for $\mathcal{P}$. As in the framework for the double induction, choose a minimum $x$ of $\CC \setminus \mathcal{P}$ and write $\mathcal{P}_x$, $\mathcal{P}_{\bar{x}}$, and  $\mathcal{P}'$ as above. The inner loop of the double induction runs over upwards-closed $G_x$-subposets $\mathcal{Q} \subset \mathcal{P}_{x}$. For each $\mathcal{Q}$, we define an $E$-parameterization $\mathcal{E}_{\mathcal{P}, \mathcal{Q}}$ of $\CC|_{\mathcal{P}_{\bar{x}}}$ by stabilizing $\mathcal{E}_\mathcal{P}$ by a $G_x$-representation $V^*(x)$ at $x$ (Lemma \ref{lemma:how-to-stabilize-equivariantly}) which depends on $\mathcal{Q}$ (although we suppress this dependence from the notation).  We will initially set $V^*(x)$ to be the zero vector space; at each step of the inner loop we will enlarge $V^*(x)$ as a $G_x$-representation. 

The inductive assumptions on the inner loop are that that for all $x > y$ in $\mathcal{Q} \cup \mathcal{P}$ we have chosen maps $\iota_{xy}$ which satisfy property $B$ and $C$ above as well as property $A$ for all $x > y > z$ in $\mathcal{Q} \cup \mathcal{P}$, together with property D: that for $a < \bar{A}(x)$, the images of $\iota_{xy} \times V^0(y)$ in $V^0(x)$ are pairwise disjoint, as $y$ runs over $\{y \in \mathcal{Q} : \bar{A}(y) = a\}$.  

To proceed with the induction on the inner loop, we choose a $y \in \CC$ such that 
\begin{itemize}
    \item the element $y$ is a maximum of $\mathcal{P}_{\bar{x}} \setminus \mathcal{Q}$, and 
    \item $\bar{A}(y) = \max_{y \in \mathcal{P}_{\bar{x}} \setminus \mathcal{Q}} \bar{A}(y). $
\end{itemize}
We  write $\mathcal{Q}_y = \mathcal{Q} \cup \{y\}$, $\mathcal{Q}' = \mathcal{Q} \cup Gy$.

Condition A for $\mathcal{Q}$ lets us define a map
\[  \iota_{xy}|_{\partial T(x,y)}: \partial T(x,y) \to \R_+^{\bar{A}_{x,y}} \times (\im \delta_{xy})^\perp\]
which automatically satisfies the equivariance condition. We would like the image of $ \iota_{xy}|_{\partial T(x,y)} \times V^0(y)$ in $\R_+^{\bar{A}(x,y)} \times V^0(x)$ to be disjoint from the images of all $\iota_{xy'} \times V^0(y')$ as $y'$ runs over $\{ y' \in \mathcal{Q} : \bar{A}(y) = \bar{A}(y')\}$. If there are no such other $y'$ then this condition is vacuous. If there are such $y'$ then it turns out to be satisfied due to the inductive assumptions: if the images of these two sets for $y$ and $y'$ overlap, then there must be a $z \in\mathcal{Q}$ such that $x > z > y$, $z > y'$ (by considering the projections of these images to $\R_+^{\bar{A}(x,a')} \times V^0(x)$ over various $a'$ with $\bar{A}(x) > a' > a$, and using condition $D$ for $\mathcal{Q}$). But then by condition $D$ for $\mathcal{P}$ and condition $A$ for $\mathcal{P}$ and $\mathcal{Q}$ these images will not overlap because the corresponding images of $\iota_{zy} \times V^0(y)$ and $\iota_{zy'} \times V^0(y')$ in $\R_+^{\bar{A}(z, y)} \times V^0(z)$ will not overlap.

We then enlarge $V^*(x)$ sufficiently such that there is a $G_{xy}$-equivariant extension $\iota_{xy}$ of  $\iota_{xy}|_{\partial T(x,y)}$ satisfying condition $C$, which satisfies the corresponding analog $D$ of condition $D$, namely that 
\begin{itemize}
    \item the $G_x$-orbit of the image of $\iota_{xy} \times V^0(y)$ is disjoint from the images of $pr_x^*\iota_{xy'}$ as $y'$, runs over the minima of $\mathcal{Q}$; and 
    \item Choosing representatives $\{g_1, \ldots, g_r\}$ of $G_x/G_{xy}$, the sets $g_i(\text{Im }\iota_{xy} \times V^0(y)) = g_i(\text{Im} \iota_{xy}) \times V^0(g_i y)$ are disjoint as $i=1, \ldots, r$. 
\end{itemize}
Such an extension exists by Lemma \ref{lemma:embed-a-single-equivariant-manifold}; we extend $\iota_{xy}$ such that the projection of the image of $\iota_{xy}|_{T(x,y) \setminus \partial T(x,y)}$ to the $G_{xy}$-representation $V^*(x)$ we are enlarging $V^*(x)$ by does not intersect $0$, and then enlarge $V^*(x)$ by $Ind_{G_{xy}}^{G_x} V^*(x)$ after choosing this embedding. 

For all $y' \in G_x y'$, we then define $\iota_{xy'}$ by requiring that the equivariance condition $B$ holds; that these can be chosen to satisfy condition $D$ is possible by the second part of condition $D'$ above. This completes the inductive step for the inner loop. After completing the inner loop (so $\mathcal{Q} = \mathcal{P}_{\bar{x}}$) we set $\mathcal{E}_{\mathcal{P}'}$ to be the stabilization of $\mathcal{E}_{\mathcal{P}}$ by the free representation which stabilizes by the final value of $\tilde{V}(x)$ at $x$ (as in Lemma \ref{lemma:how-to-stabilize-equivariantly}).  Similarly, to complete the inductive step for the outer loop, for all $y \in \mathcal{P}_{\bar{x}}$ and all $g \in G$,  one defines $\iota_{gx, gy}$ such that the equivariance condition holds.  

Note that we have now constructed  an $\mathcal{F}$-parameterization $\mathcal{F}^{new}$ with underlying vector spaces $W(x,z) = \im \delta_{xz}^\perp$ embedded in the the $E$-parameterization $\mathcal{E}_\CC$ with underlying vector spaces $\{\bar{V}(x)\}$ (Lemma \ref{lemma:poset-rep-determines-half-framing}). We set $\CC'' = \CC'_{\mathcal{E}_\CC}$ (Definition \ref{def:stabilization-by-f-parameterization}).

We now think of the maps $\iota_{xy}$ as maps 
\[ \iota_{xy}: T_{\CC'}(x,y) \to \R_+^{\bar{A}(x,y)} \times  V^0_{\CC''}(x,y).\]
This is possible because $V^0_{\CC''}(x,y)$ is canonically identified with $(im (\hat{\alpha}^0_{xy})_{new})^\perp$, which contains 
$(im \delta_{xy})^\perp$. Because the metric on the underlying virtual smoothing is very convenient, we have that at every point $p \in im \iota_{xy}$, projection of $(d\iota_{xy})_p(\eta^{-1}_{xy}( 0 \oplus V^0(x,y)))$ to $V^0_{\CC''}(x,y)_p$ is full rank (since the projection of $(d\iota_{xy})(\eta_{xy}^{-1}(T\R_+^{\bar{A}(x,y)} \oplus 0))$ to $T \R_+^{\bar{A}(x,y)}$ is full rank). We write $\nu \iota_{xy}$ for the union of the normals to the images of $(d\iota_{xy})_p(\eta^{-1}_{xy}( 0 \oplus V^0(x,y)))$ over $p \in im \iota_{xy}$.

The
discussion of Section \ref{sec:linear-algebra-of-pontrjagin-thom} then gives a family of isometric isomorphisms of bundles over $T_{\CC'}(x,y)$,
\[ \kappa_{xy}: \underline{W^0(x,y)} \to \nu_v \iota_{xy}. \] 
 Using the Riemannian metrics on $V(x)$, by choosing a sufficiently small $\epsilon > 0$, we can find diffeomorphisms $\phi_{xy}$ from the $\epsilon$-neighboroods of the zero section of $\nu \iota_{xy}$ to neighborhoods of $im \iota_{xy}$. We define $T_{\CC'''}(x,y)$ to be the $\epsilon$-neighborhood of the zero section of $\underline{W^0(x,y)}$, and we define the final embedding $\iota_{xy}: T_{\CC'''}(x,y) \to V^0(x)$ to the composition $\phi_{xy} \circ \kappa_{xy}$. The discussion about products in Section \ref{sec:linear-algebra-of-pontrjagin-thom} shows precisely that these maps define an (equivariant) embedding of $\CC'''$. 
\end{proof}

\begin{lemma}
\label{lemma:extensions-of-embeddings-exist}
    Assume that $\CC$ is proper. After shrinking $\CC'''$, any embedding of $(\CC, \CC''')$ produced as in Proposition \ref{prop:embeddings-exist} admits an extension. 
\end{lemma}
\begin{proof}
    We take the shrinking to be some $\epsilon$-neighborhood of the zero-sets of the defining sections. After this,
    this follows from an equivariant double induction, as in Proposition \ref{prop:embeddings-exist}, of the $G$-equivariant Tietze extension theorem. We omit the details; the only \emph{crucial} fact is that the properness of the category is what allows us to use the Tietze exentsion theorem, which gives equivariant extensions of continuous maps from \emph{closed} sets, which we take to be closed neighborhoods of various zero sets of Kuranishi sections. 
\end{proof}

There are also several relative variants of Proposition \ref{prop:embeddings-exist} which are needed for the construction of continuation maps and for the geometric setup that allows us to construct the cyclotomic structure on symplectic cohomology. These are proven identically to Proposition \ref{prop:embeddings-exist}, but with a slightly modification of conditions $D$ and $D'$. We explain these statements and modifications below. 

\begin{proposition}
\label{prop:relative-variant-of-embedding-construction}
    Let $\CC$ be a $G$-equivariant virtually smooth framed flow category with with the framing parameterized by an $E_2^s$-parameterization. Suppose that the objects of $\CC$ can be partitioned into two disjoint $G$-invariant subsets  $\CC_1$ and $\CC_2$ such that 
    \begin{itemize}
    \item If $x \in \CC_1$ and $y \in \CC_2$ then either $x > y$ or they are incomparable; 
    \item We have that $\max_{y \in \CC_2} \bar{A}(y) < \min_{x \in \CC_1} \bar{A}(x). $
    \end{itemize}
    Write $\CC_1$ and $\CC_2$ for the full flow subcategories of $\CC$ on the relevant objects, and suppose that (with the framings and embeddings of framings restricted from $\CC$) there are embeddings and extensions of embeddings $\bar{\mathbf{E}}_i$ and $\mathbf{E}_i$ for $\CC_i$ for $i=1,2$, respectively. Then there exists an $E^s$-parameterization $\mathcal{E}$, an embedding $\mathbf{E}$ of $\CC'_{\mathcal{E}}$  and an extension $\bar{E}$ of $\mathbf{E}$ such that the restriction of $\bar{\mathbf{E}}$ and $\mathbf{E}$ to $\CC_{\mathcal{E}, i} \subset \CC_{\mathcal{E}}$, the full subcategories on the objects of $\CC_i$, agrees with the stabililized embedding $\mathbf{E}^{\mathcal{E}|_{\CC_{\mathcal{E}, i}}}$ and the stabilized extension $\bar{\mathbf{E}}^{\mathcal{E}|_{\CC_{\mathcal{E}, i}}}$ for $i=1, 2$, respectively.  
\end{proposition}
\begin{proof}
This follows identically to the proofs of Proposition \ref{prop:embeddings-exist} and Lemma \ref{lemma:extensions-of-embeddings-exist}, with a few key modifications. Write $\mathcal{E}_{\CC'}$ for the $E^s$-parameterization underlying the parameeterization of the framing of $\CC'$. 

When defining $\mathbf{E}$, we first note that we can start the induction with $\mathcal{P} = \CC_2$ and $\mathcal{E}_\mathcal{P} = \mathcal{E}_{\CC'}|_{\mathcal{P}}$. Note that all $E$-parameterizations are built by inductively stabilizing the zero $E$-parameterization upwards along the poset $\CC$. Then, whenever we begin the outer loop on some $x \in \CC_1$, before the first run of the inner loop we set $V^*(x)$ to be a copy of the vector space that we would stabilize by to build $\mathcal{E}_{\CC'}$ from $\mathcal{E}_{\CC'}|_{\mathcal{P}}$.  With this choice, there will be canonical copy of $V^0(x,y)$ for any $y \in \CC_1$ in $\mathcal{E}_{\mathcal{P}; \{x\}}$, the stabilization of $\mathcal{E}_{\mathcal{P}}$ by this vector space $V^*(x)$ at $x$. The condition in the lemma statement on $\CC_1$ and $\CC_2$ then means that we can start the second loop at $\mathcal{Q} = \mathcal{P}_{\bar{x}} \cap \CC_1$, since we just choose the embeddings to be the composition of the original embeddings to $\R_+^{\bar{A}(x,y)} \times V^0(x,y)$ with the inclusions of the latter into $\R_+^{\bar{A}(x,y)} \times V^0(x,y) \times  V^0_{\mathcal{E}_{\mathcal{P};\{x\}}}(x)$ arising from this initial choice of $V^*(x)$. We then continue identically with the rest of the induction.

After the double induction has terminated, we use now stabilize the category and get trivializations of normal bundles as before; however we will choose $\epsilon$ in a point-dependent way. Two conditions characterize our choice of $\epsilon$:
\begin{itemize}
    \item There is no restriction on the size $\epsilon$ of the valid neighborhood for $x>y$ in $\CC_1$. In fact we choose $\epsilon=\infty$ for those pairs. With this choice, \emph{the nondisjointness condition for embeddings is still satisfied} for all $x \in \CC_1$ and all $a$ such that $a > \max_{y \in \CC_2} \bar{A}(y)$. 
    \item With this choice, for $x \in \CC_1$ and $y \in \CC_2$ we can simply shrink $\epsilon$ as we move away from the boundary such that the nondisjointness condition is still satisfied. 
\end{itemize}
Finally, we note that for $x > y$ in $\CC_1$, the embedding of $\CC_{\mathcal{E}}$ produced by the above method \emph{agrees} with the stabilized embedding, via direct linear algebra computation. 

To produce $\bar{\mathbf{E}}$, we argue as before but choose the  closed sets to which we apply the Tietze extension theorem to be the unions of the boundaries of $\partial (\R_+^{\bar{A}(x,y)} \times V^0(x,y))_*$ with closed neighborhoods of the zero locus.

\end{proof}

\begin{remark}
\label{rk:cofibrancy}
    Analogs of Proposition \ref{prop:relative-variant-of-embedding-construction} without the constraint of $\bar{A}$ of that proposition can be produced by a \emph{restratification} process, where one stabilizes $\CC$, relabels the boundary strata so that the conditions of Proposition \ref{prop:relative-variant-of-embedding-construction} hold, and proceeds onwards with the proof. We do not pursue this (technically useful but not immediately necessary) elaboration in this work. The constraint on $\bar{A}$ can be seen as a partial `cofibrancy' condition, of the form that `different boundary components are labeled by different strata'.  A strong version of this condition would be that \emph{the map $\bar{A}: Ob(\CC) \to \Z$ is injective}, which arises automatically for the vitually-smooth flow categories produced by the method of Section \ref{sec:global-charts} for \emph{action-nondegenerate Hamiltonians}, i.e Hamiltonians for which all actions are distinct. 
\end{remark}

\begin{remark}
    \label{eq:differences-between-embedded-framing-and-a-virtually-smooth-framing}
    Given an embedded virtually smooth flow category $\CC'_{\mathcal{E}}$ (where $\mathcal{E}$ is the $E^s$-parameterization arising in Proposition \ref{prop:embeddings-exist}), one can use the linearizations of the embedding maps to define framings of the thickenings of $\CC'_{\mathcal{E}}$. These framings will not agree with the original framings of $\CC'$, nor will they agree with the framings of $\CC'_{\mathcal{E}}$ induced by taking the framing of $\CC'$ and then subsequently stabilizing it by $\mathcal{E}$. However, the final Floer homotopy type is still controlled by the original framings of $\CC'$ (see Remark \ref{rk:how-do-framings-of-thickenings-play-a-role}), sufficiently so that in the geometric case  we can prove a comparison with integral Floer homology in Section \ref{sec:comparison-with-floer-homology} (for which of course getting the framings correct is necessary, as the framings control the signs of the corresponding maps in a certain spectral sequence). More generally, it should be possible to prove that after a further sufficient stabilization, the framings produced by using the differentials of the embedding maps are deformation-equivalent in a suitable sense to the framings produced by taking the original framings of $\CC'$ and stabilizing.  
\end{remark}
\begin{proposition}
\label{prop:cyclotomic-variant-of-embedding-construction}
Let $\CC'$ be a virtual smoothing of a $G$-equivariant flow category $\CC$ which has a semi-freely parameterized embedded framing. Suppose that for some normal subgroup $H \subset G$,  $(\CC')^H$ (thought of as a $G/H$-equivariant virtual smoothing of $\CC^H$, with its induced framing and parameterization) has an ($G/H$-equivariant) embedding $\mathbf{E}_H$ and an extension $\overline{\mathbf{E}}_H$.  Then there is an $E^s$-parameterization $\mathcal{E}$ such that $\CC'_{\mathcal{E}}$ admits an embedding $\mathbf{E}$ with extension $\overline{\mathbf{E}}$ such that their restrictions to $(\CC'_{\mathcal{E}})^{H} = ((\CC')^H)_{\mathcal{E}^H}$  are $\mathbf{E}^{\mathcal{E}^H}$ and $\bar{\mathbf{E}}^{\mathcal{E}^H}$ respectively. 
\end{proposition}
\begin{proof}
    This follows identically to the proofs of Proposition \ref{prop:embeddings-exist} and Lemma \ref{lemma:extensions-of-embeddings-exist}, with a few key modifications. As in the proof of Proposition \ref{prop:relative-variant-of-embedding-construction}, before running the inner loop, at each step of the outer loop we initialize $V^*(x)$ to be the stabilization which builds up $\mathcal{E}|_{(\CC')^H}$ from zero, so that there is a copy of $V^0_{(\CC')^H}(x,y)$ in $V_{\mathcal{E}_{\mathcal{P}; \{x\}}}(x)$. Using this we require that $\iota_{xy}$ restrict to compositions of the embeddings $\iota_{xy}^H$ of $(\CC')^H$ with the inclusion $V^0_{(\CC')^H}(x,y) \subset V_{\mathcal{E}_{\mathcal{P}; \{x\}}}(x)$; this is possible by the inductive conditions and because $\iota_{xy}^H$ is already $G_{xy}/H_{xy}$-equivariant.  After the induction terminates, the one concludes similarly to the proof of Proposition \ref{prop:relative-variant-of-embedding-construction}, choosing $\epsilon(x)$ for $x \in T_{\CC'}(x,y)^H$ to be infinity and having it decrease rapidly as one exits this locus, and similarly requiring that choosing the closed subsets to which the apply the Tietze extension theorem to be unions of closed neighborhoods of the zero-sections with $\R^c_+ \times V^0(x,y)^{H_{xy}}$.
\end{proof}
\begin{remark}
    Without the normality condition on $H$ and the $G/H$-equivariance of the embedding and extension of $(\CC')^H$, it is not completely clear to the author if the corresponding analog of the proposition above holds, primarily because there may be a condition for extending a general map $M^H \to V^H$ to a map $M \to V$ for $M$ a $G$-manifold and $V$ a $G$ $G$-representation. 
\end{remark}

\subsection{Restratfications.}
\begin{definition}
    A \emph{restratification} of a flow category $\CC_1$ with an integral action $\bar{A}_1$ is a flow category $\CC_2$ with an integral action $\bar{A}_2$ such that $Ob(\CC_1) = Ob(\CC_2)$, such that the gradings (if defined) of the two flow categories agree,  such that there are embeddings $\rho_{xy}: [\bar{A}_1(x,y)] \subset [\bar{A}_2(x,y)]$ for each $(x,y) \in \CC_2$ and such that $\CC_2(x,y)$ is $\CC_1(x,y)$ as a space, with the stratification on $\CC_2(x,y)$ determined by 
    \[
    \begin{gathered}
    \CC_2(x,y)([\bar{A}_2(x,y)] \setminus j) = \emptyset \text{ if }j \notin Im \rho_{xy},\\ 
    \text{ and otherwise }\CC_2(x,y)([\bar{A}_2(x,y)] \setminus j) = \CC_1(x,y)([\bar{A}_1(x,y)] \setminus \rho_{xy}^{-1}(j)).
    \end{gathered}\]
    with such that the composition maps in $\CC_2$ agreeing with the composition maps in $\CC_1$ if we think of these categories as enriched in spaces. 
\end{definition}

\begin{definition}
\label{def:restratification-global-data}
Suppose $\CC_1$ has an integral action. We can specify a restratification with the following construction. 

Given a flow category $\CC_1$, we can define a restratification $\CC_2$ associated to an order-preserving inclusion of sets 
\[ \rho: S_{\bar{A}} \to \{0, \ldots, \ell\} =: S_{\bar{A}_2} \]
for some integer $\ell > \max_{x \in \CC} \bar{A}(x)-1$. This is defined by setting 
\[ \bar{A}_2(x) = \rho(\bar{A}_1(x))-1, \]

We call the integer $\ell - |S_{\bar{A}}|$ the \emph{defect} of the restratification.
\end{definition}

\begin{lemma}
\label{lemma:restratifications-exist}
    If $\CC_1$ has a virtual smoothing $\CC'_1$ and $\CC_2$ is a restratification of $\CC_1$, then there is associated virtual smoothing $\CC'_2$ of $\CC_2$ for which, writing $\CC'_i(x,y) = (T_i(x,y), V_i(x,y), \sigma_i(x,y))$, we have 
    \[T_2(x,y) = (-1,1)^{[\bar{A}_2(x,y)] \setminus \rho_{xy}([\bar{A}_1(x,y)])} \times T_1(x,y) \]
    \[V_2(x,y) = \underline{\R^{[\bar{A}_2(x,y)  \setminus \rho_{xy}([\bar{A}_1(x,y)])}} \oplus \pi_{T_1(x,y)}^*V_1(x,y),  \]
  \[ \sigma_2(x,y)(a,b) = (id, \sigma_1(b)) \]
    where the stratification on $T_2(x,y)$ is 
    \[ T_2(x,y)([\bar{A}_2(x,y)] \setminus j) = (-1,1)^{[\bar{A}_2(x,y)] \setminus \rho_{xy}([\bar{A}_1(x,y)])}\times T_1(x,y)[\bar{A}_1(x,y)] \setminus \rho_{xy}^{-1}(j))\]
    if $j \in Im \rho_{xy}$, and otherwise $T_2(x,y)([\bar{A}_2(x,y)] \setminus j) = \emptyset$. The composition is defined using product of the composition on $\CC'_1$ and maps which use the embeddings 
    \[ \rho_{xyz}: ([\bar{A}_2(x,y)] \setminus Im \rho_{x,y}) \cup ([\bar{A}_2(y,z)] \setminus Im \rho_{y,z}) \to ([\bar{A}_2(x,z)] \setminus Im \rho_{x,z})\]
    associated to the maps of sets defining composition in $\CC_2$ to compose the $(-1,1)^{[\bar{A}_2(x_i,x_j)] \setminus \rho_{xy}([\bar{A}_1(x_i,x_j)])}$ factors and the $\R^{[\bar{A}_2(x_i,x_j)] \setminus \rho_{xy}([\bar{A}_1(x_i,x_j)])}$ factors.

     A framing of $\CC'_1$ defines a framing of $\CC'_2$, as specified under \eqref{eq:framings-for-stabilization-1} and \eqref{eq:framings-for-stabilization-2} below. Similarly, an embedding of the framing of $\CC'_1$, together with a choice of integer $M$, the \emph{restratification embedding shift}, satisfying
     \begin{equation}
         \label{eq:restratification-shift-choice-condition}
         0 \leq M \leq \rho(0)
     \end{equation} 
     defines an embedding of the framing of $\CC'_2$. If the parameterization of the framing of $\CC'_1$ was $E^s_2$ and the restratification was simple then the same holds for $\CC'_2$. Finally, embeddings of $\CC'_1$ extend to embeddings of $\CC'_2$ canonically, as do extensions of embeddings.

\end{lemma}
\begin{proof}
    It is clear that the above definition defines a virtual smoothing of $\CC$. We need to specify the framings. We set 
    \begin{equation}
        \label{eq:parameterizations-after-restratification}
         V^0_{\CC''}(x,y) = V^0_{\CC'}(x,y), \; \; V^1_{\CC''}(x,y) =\R^{[\bar{A}_2(x,y)  \setminus \rho_{xy}([\bar{A}_1(x,y)])} \oplus V^1_{\CC'}(x,y).
    \end{equation}
    We set the trivialization 
    \begin{equation}
        \label{eq:framings-for-stabilization-1}
        \eta^2_{xy}: TT_2(x,y) = T(-1,1)^{[\bar{A}_2(x,y)] \setminus \rho_{xy}([\bar{A}_1(x,y)])} \times TT_1(x,y) \to \underline{T\R}_+^{\bar{A}_2(x,y)} \times V^0_{\CC''}(x,y) 
    \end{equation} 
    (where we use the natural reordering isomorphism on the right) to be $\eta^2_{xy} = id \times \eta_{xy}$, where we use that $T(0,1) = \underline{\R}$ is canonically trivialized.

    We set the trivialization 
    \begin{equation}
        \label{eq:framings-for-stabilization-2}
        (\eta')^2_{xy}: V_2(x,y) = \underline{\R^{[\bar{A}_2(x,y)  \setminus \rho_{xy}([\bar{A}_1(x,y)])}} \oplus \pi_{T_1(x,y)}^*V_1(x,y) \to \R^{[\bar{A}_2(x,y)  \setminus \rho_{xy}([\bar{A}_1(x,y)])} \oplus V^1_{\CC'}(x,y) = V^1_{\CC''}(x,y)
    \end{equation}
    to be the $(\eta')^2_{xy} = id \times \eta_{xy}$

    We need to specify the parameterizations, and make sure that this new framing is embedded. We set 
    \[ V^0_{\CC''}(x) = V^0_{\CC'}(x), \; \; V^1_{\CC''}(x) = \R^{\bar{A}_2(x) - \bar{A}_1(x)-M} \times V^1_{\CC'}(x).\]
    The condition \eqref{eq:restratification-shift-choice-condition} ensures precisely that $\bar{A}_2(x) - \bar{A}_1(x)-M>0$ for all $x$. 
    We set the embeddings 
    \[ (\bar{\alpha}^2)^0_{xy}: V^0_{\CC''}(x,y) \oplus V^0_{\CC''}(y) \to V^0_{\CC''}(x)\]
    to be $\bar{\alpha}^2_{xy} = \bar{\alpha}_{xy}$, and 
    \[ (\bar{\alpha}^2)^1_{xy}: V^1_{\CC''}(x,y) \oplus V^1_{\CC''}(y)  = \R^{[\bar{A}_2(x,y)  \setminus \rho_{xy}([\bar{A}_1(x,y)])} \times V^1_{\CC'}(x,y) \times \R^{\bar{A}_2(y) - \bar{A}_1(y)-M} \times V^1_{\CC'}(y) \to \R^{\bar{A}_2(x) - \bar{A}_1(x)-M} \times V^1_{\CC''}(x)\]
    by noting that 
    \[ \bar{A}_2(x,y) - \bar{A}_1(x,y) = (\bar{A}_2(x) - \bar{A}_1(x)+M) - (\bar{A}_2(y) - \bar{A}_1(y)-M)\]
    and simply moving switching the third factor of the domain with the second and then applying $\bar{\alpha}^1_{xy}$ to the latter two factors. 
    Writing $(\mathcal{F}_0', \mathcal{F}_1')$ for the $F'_2$-parameterization of the framing of $\CC'_2$ and $(\mathcal{F}_0, \mathcal{F}_1)$ for the $F'_2$ parameterization of the framing of $\CC'_1$, we see that $\mathcal{F}'_0 = \mathcal{F}_0$. If $\mathcal{F}_1$ was an $F^s$ parameterization then $\mathcal{F}'_1$ is as well: the free parts are the same, while the canonical part of $\mathcal{F}'_1$ (with defining function $f'_1$, and the defining function of the canonical part of  $\mathcal{F}_1$ being $f_1$) by setting $f'_1(\rho(n)) = f_1(n)$ and otherwise $f'_1(n) = 1$ for $n \in S_{\bar{A}_2} \setminus \im \rho$.

    We extend embeddings of $\CC'_1$ to embeddings of $\CC'_2$ by setting 
    \[ \iota_{xy}: (-1,1)^{[\bar{A}_2(x,y)] \setminus \rho_{xy}([\bar{A}_1(x,y)])} \times T_1(x,y) \to  \R_+^{[\bar{A}_2(x,y)] \setminus \rho_{xy}([\bar{A}_1(x,y)])} \times \R_+^{[\bar{A}_1(x,y)]}\times V^0_{\CC'_1}(x,y)
    =\R_+^{\bar{A}_2(x,y)}\times V^0_{\CC'_1}(x,y)\]
    as $\iota_{xy}(a,p) = (a+(2,\ldots ,2), \iota_{xy}(p))$ where the $\iota_{xy}$ inside the parentheses is the original embedding map for $\CC'_1$, and $(2, \ldots, 2)$ is the vector containing only twos. We will choose, once and for all, a function 
    \begin{equation}
        \label{eq:stabilizing-function-for-restratification}
        \xi: \R_+ \to \R_*
    \end{equation}
    which extends continuously to $\R_+^*$ as a pointed map, and such that $\xi^{-1}(\infty)=0$ (and thus $\xi$ descends to a homeomorphism $\bar{\xi} = (\R_+)_*/\{0, \infty\} \to \R_*$). 
    Namely, we set
    \[ \xi|_{(1, \infty)}(p) = p-2, \xi(0) = \infty, \xi|_{0, 1}(p) = -1/p. \]
    We set the new extension of the embedding $e^2_{xy}$ be 
    
    \begin{equation}
        \label{eq:effect-of-restratification-on-extension-of-embedding}
        e^2_{xy} = ((\xi, \xi, \ldots, \xi) \wedge e_{xy}) 
    \end{equation}where we have used the decomposition \eqref{eq:parameterizations-after-restratification} of $V^1_{\CC''}(x,y)$ as well as \eqref{eq:canonical-smash-homeo}.
\end{proof}

\begin{lemma}
    Given a restratification $\CC_2$ of a flow category $\CC_1$ with associated embeddings $\rho^1_{xy}: [\bar{A}_1(x,y)] \to [\bar{A}_2(x,y)]$ and a restratification $\CC_3$ of  $\CC_2$ with associated embeddings $\rho^2_{xy}: [\bar{A}_1(x,y)] \to [\bar{A}_2(x,y)]$, $\CC_3$ is also a restratification of $\CC_1$ with associated embeddings $\rho^{12}_{xy} = \rho^2_{xy} \circ \rho^1_{xy}$. Moreover, the virtual smoothing of $\CC_3$ associated to $\rho^{12}_{xy}$ by Lemma \ref{lemma:restratifications-exist} and the virtual smoothing of $\CC_3$ built by first applying Lemma \ref{lemma:restratifications-exist} to $\rho^1_{xy}$ and subsequently to $\rho^2_{xy}$, are canonically isomorphic, and there is an isomorpism of the parameterizations of the framings and of the embeddings of the framings which takes one framing to the other, and takes the embeddings and extensions of embeddings from one to the other.
\end{lemma}
\begin{proof}(Sketch) We do not spell out all the details of this compositionality result. However, it follows straightforwardly from the canonical identification of sets 
\[ [\bar{A}_3(x,y)] \setminus \rho^{12}_{xy}([\bar{A}_1(x,y)]) = ([\bar{A}_3(x,y)] \setminus \rho^2_{xy}([\bar{A}_2(x,y)])) \cup ([\bar{A}_2(x,y)] \setminus \rho_{xy}([\bar{A}_1(x,y)])).\]
\end{proof}

\begin{proposition}
\label{prop:simple-restrat-effect-on-homotopy-type}
The effect of a restratification on the floer homotopy space is 
\begin{equation}
    \label{eq:effect-of-restrat-on-floer-homotopy-space}\{\CC'_2\} \simeq \Sigma^{\R^{S_{\bar{A}_2} \setminus \rho(S_{\bar{A}_1})}} \{\CC'_1\} \simeq \Sigma^\ell \{\CC'_1\}
\end{equation} 
\and 
\[ |\CC'_2| \simeq F_{\R^{\ell+M}}\Sigma^\ell |\CC'_1|.\]
\end{proposition}
\begin{proof}
    Write 
    \[S'_{< a} = \{a' \in S_{\bar{A}_2} \setminus \rho(S_{\bar{A}_1}): a' < a\},  S'_{> a} = \{a' \in S_{\bar{A}_2} \setminus \rho(S_{\bar{A}_1}): a' > a\}. \]
    Then for $a \in S_{\mathfrak{A}_1}$, under a restratification, the term corresponding to $z$ in the wedge summand \eqref{eq:floer-homotopy-type-new} changes by taking smash product with 
    \begin{equation}
    \label{eq:change-in-wedge-summand-under-restratification}
    (\R_+^{S'_{<\bar{A}_2(z)}} \times \R^{S'_{>\bar{A}_2(z)}})_*
    \end{equation}
    with the second factor in the product coming from the change in $V^1(x)^\perp$. For $z \in M(\CC)$ the boundary of the $\R_+$ factors is collapsed to the basepoint; thus via $\bar{\xi}$ (see \eqref{eq:stabilizing-function-for-restratification}) we can identify the effect of those factors as simply a smash product with $(\R^{S'_{<\bar{A}_2(z)}})_*$; thus we simply suspend the cells  of $\{\CC_1\}$ corresponding to $z \in M(\CC)$ by $\R^{S_{\bar{A}_2} \setminus \rho(S_{\bar{A}_1})}$, as desired. 
    
    For the other cells, we can use $\bar{\xi}$ to convert the $\R_+$ factors in \eqref{eq:change-in-wedge-summand-under-restratification} corresponding to $\{s \in S_{\bar{A}_2} \setminus : s < \min \rho(S_{\bar{A}_1})\}$, as before. For the other factors, we see that we are exactly using $\xi$ in $\eqref{eq:effect-of-restratification-on-extension-of-embedding}$ to define the action of the extension of the embedding, which is used to define take the quotient in \eqref{eq:floer-homotopy-type-new};  thus, identifying those factors with $\R$ factors via $\xi$ as well is consistent with the quotient in \eqref{eq:floer-homotopy-type-new}. This establishes \eqref{eq:effect-of-restrat-on-floer-homotopy-space}. The rest of the claim follows from computing the change in the shift space under the assignment \eqref{eq:parameterizations-after-restratification}.
\end{proof}

\section{Review of (Hamiltonian) Floer Homology}
\label{sec:hamiltonian-floer-review}

Hamiltonian Floer Homology, like other Floer theoretic invariants, arises formally from the study of the gradient flow of a closed $1$-form on an infinite-dimensional manifold. Specifically, let $(M^{2n}, \omega)$ be a symplectic manifold, such that $\omega$ has a chosen primitive $1$-form $\lambda$. To a $t$-dependent Hamiltonian function
\[ H = H_t: M \times S^1_t \to \R\]
Hamilton's equations 
\begin{equation}
\label{eq:hamiltons-equations}
dH = -i_{X_H} \omega
\end{equation}
associate a $t$-dependent Hamiltonian vector field $X_H = X_{H_t}$ with associated time $t$ flows given by $\phi^t_H$. The $1$-periodic trajectories of the flow of $X_H$ are fixed points of $\phi^1_H$, or equivalently the solutions to 
\begin{equation}
\label{eq:periodic-point-equations}
    x: S^1 \to M, \frac{dx}{dt} = X_{H_t}.
\end{equation}
Thinking of the periodic trajectories as elements of the free loop space $LM$, which behaves much like an infinite-dimensional manifold, one notices that \eqref{eq:periodic-point-equations} are formally the equations for the critical points of the symplectic action functional 
\begin{equation}
\label{eq:sympectic-action-functional}
    \mathcal{A}_H: LM \to \R, \mathcal{A}_H(x) = \int_x \lambda - x^*H dt. 
\end{equation}
Formally, the differential of the symplectic action functional takes the form 
\begin{equation}
    \label{eq:hamiltonian-1-form}
    d\mathcal{A}_H: TLM \to \R, d\mathcal{A}_H(x,v) = \int_0^1 \omega(\dot x, v) - dH(v) dt 
\end{equation}
where $x:S^1 \to M$ is a point of the free loop space $LM$ and $v \in C^\infty(x^*TM)$ is a tangent vector field along $x$, i.e. an element of $T_xLM$. Picking a $S^1_t$-dependent family $J_t$ of almost complex structures on $M$ which are all compatible with $\omega$ defines a $t$-dependent metric on $M$, and thus a metric on $LM$ via 
\begin{equation}
\label{eq:formal-inner-product}
    \langle v, w \rangle = \int_0^1 \omega_{x(t)}(v(t),J_t w(t)) dt
\end{equation}
where $x \in LM$ and $v, w \in T_xLM$. The downwards gradient flow trajectories of the symplectic action functional with respect to this metric are thus formally given by solutions to Floer's equations
\begin{equation}
\label{eq:floers-equation}
    u: Z \to M, \partial_s u + J_t (\partial_t u- X_{H_t}) = \partial_s u + J_t \partial_t u - \Grad H_t= 0
\end{equation}a
where $Z = \R_s \times S^1_t$ is thought of a holomorphic cylinder with complex structure $i$ such that $i \partial_s = \partial_t$.  

Let $Fix(H)$ denote the solutions to \eqref{eq:periodic-point-equations}; these are in canonical bijection with their starting points $x(0)$. The Hamiltonian $H_t$ is \emph{nondegenerate} if $(d\phi^1_H)_{x(0)}$ has no eigenvalues equal to $1$ for all $x \in Fix(H)$. Formally, the Hessian of $\mathcal{A}_H$ is given by 
\begin{equation}
    Hess(\mathcal{A}_H)(x; v, w) = \int_0^1 \langle J_t \Grad_t v + (\Grad_{v(t)} J_t) \partial_t u - \Grad_v \Grad H_t, w(t) \rangle dt
\end{equation}
where $x \in LM$ and $v, w \in T_xLM$. Just as in the context of finite-dimensional Morse theory, where Hessians of functions are self-adjoint with respect to the metric, the operator 
\begin{equation}
\label{eq:hessian-operator-form}
    \xi \mapsto  J_t \Grad_t v + (\Grad_{v(t)} J_t) \partial_t u - \Grad_v \Grad H_t
\end{equation}
is a formally self-adjoint operator with respect to the inner product \eqref{eq:formal-inner-product}. 

\begin{remark}
    The formal self-adjointness of the operator \eqref{eq:hessian-operator-form} is not obvious, but can be established with the same computation that establishes the self-adjointness of the Hessian $v \mapsto \Grad_v \Grad f$ for $f$ a Morse function on a finite-dimensional manifold. Namely, it follows from the torsion-freeness of the Levi-Civita connection, together with the identity 
    \begin{equation*}
        \begin{gathered}
            g(\Grad_X \Grad f, Y) - g(X, \Grad_Y \Grad f) = X g(\Grad f, Y) - g(\Grad f, \Grad_X Y) - Yg(\Grad f, X) + g(\Grad f, \Grad_Y X) 
            \\ = X(Y(f)) - Y(X(f)) - [X,Y]f = 0. 
        \end{gathered}
    \end{equation*} 
    for $g$ a metric and $X,Y$ vector fields. 
\end{remark}

An immediate observation from the formula \eqref{eq:hessian-operator-form} is that if $H$ is nondegenerate then $\mathcal{A}_H$ behaves formally like a Morse function. Making this assumption forces $Fix(H)$ to be finite if $M$ is compact (possibly with boundary) or if $H$ satisfies certain standard conditions near infinity which will be spelled out later. For any pair $x_\pm \in Fix(H)$, one then defines moduli spaces $\mathcal{M}(x_-,x_+)$ as the sets of solutions to \eqref{eq:floers-equation} satisfying the asymptotic conditions 
\begin{equation}
    \label{eq:asymptotic-conditions}
    \lim_{s \to \pm \infty} u(s, t) = x_\pm(t).
\end{equation}
Under standard conditions described below, the Gromov-Floer compactifications $\overline{\mathcal{M}}(x_-, x_+)$ of $\mathcal{M}(x_-, x_+)$ are compact Hausdorff spaces, and come with gluing maps 
\[ \cM(x_0, x_1) \times \ldots \times \cM(x_{k-1}, x_k) \to \cM(x_0, x_k)\]
for any $(x_0,\ldots, x_k) \in Fix(H)^{k+1}$. The specializations of the gluing maps to $k=2$ define the compositions of the  \emph{pre-flow category} $\CC(H, J)$ associated to $(H, J)$, which is a-priori just a topologically-enriched category, which has objects $Fix(H)$ and morphisms $\CC(H, J)(x_-, x_+) = \cM(x_-, x_+)$. 

\begin{remark}
    This is not yet a flow category in the sense of Definition \ref{def:flow-category} because we need to pick a convention for stratifying the spaces $\cM(x_-, x_+)$. When the Floer datum $(H, J)$ is regular, it is conventional to stratify by the difference of the Conley-Zehnder indices of the endpoints; however, this does not make sense for irregular Floer data since the index of a curve $u \in \cM(x_-, x_+)$ may be negative. 
\end{remark}

\paragraph{Conditions on symplectic manifolds and Floer data.}

\begin{definition}
    A \emph{Liouville domain} is a compact exact symplectic manifold with boundary $(M, \omega, \lambda)$ such that the Liouville vector field $X$ defined by $i_X \omega = \lambda$ points outward along $\partial M$. 
    
    Under this condition $\lambda|_{\partial M}$ a contact form on $\partial M$. Moreover, 
    The flow of the vector field $X$ defines a collar $[\epsilon, 1]_t \times \partial M$ around $\partial M$ sending $1$ to $\partial M$ and such that
    \[  \lambda|_{[\epsilon,1]_t \times \partial M} = e^t \lambda|_{\partial M},\]
    \[ X|_{[\epsilon,1]_t \times \partial M} = \partial/\partial t.\]
    A symplectic manifold with with contact boundary is a sympectic manifold with boundary such that an open neighborhood of its boundary is Liouville.

    A \emph{symplectically atoroidal} manifold $M$ is a symplectic manifold $(M, \omega)$ such that $\int_{T^2} u^*\omega = 0$ for every differentiable map $u: T^2 \to M$. Note that Liouville domains are necessarily symplectically atoroidal. 

    An \emph{admissible} symplectic manifold $M$ is a symplectic manifold that is symplectically atoroidal with contact boundary.
    \end{definition}
    \begin{definition}
    Let $M$ be an admissible symplectic manifold. 
    A pair $(H, J)$, where $H$ is an $S^1$-dependent Hamiltonian on $M$, and $J$ is an $S^1$-dependent almost complex structure taming $\omega$, is a \emph{Floer datum}. In the case when $M$ has contact boundary, we say that $(H, J)$ is a \emph{convex linear Floer datum} if 
    (after possibly shrinking the collar 
    \begin{itemize}
        \item $H|_{[\epsilon,1]_t \times \partial M} = ct$ for  constant $c>0$, called the \emph{slope} of $H$,  and 
        \item On the image of the collar, $J$ preserves $\ker \lambda$ and sends $X$ to the Reeb vector field of $\lambda|_{\partial M}$. 
    \end{itemize}
   A Floer datum is \emph{nondegenerate} if  $H^{\# k}$ is nondegenerate for all $k$, and (if $M$ is Liouville none of Reeb orbits of $\lambda|_{\partial M}$ has period an integral multiple of the slope of $H$. 

   An \emph{admissible} Floer datum $(H, J)$ is a nondegenerate Floer datum which is convex linear if $M$ is Liouville. 
\end{definition}

\begin{definition}
    Given a pair of admissible Floer data $(H_\pm, J_\pm)$ such that $H_\pm$ has slope $c_\pm$, an \emph{admissible continuation datum} (from $(H_+, J_+)$ to $(H_-, J_-)$)  is a family $(H_s, J_s)$ of admissible Floer data for $s \in \R$ such that
    \begin{itemize}
        \item The functions $H^s:\R \times S \times S^1 \to \R$ and $J^s: \R \times S^1 \times M \to \R \times S^1 \times End(TM)$ are smooth, and
        \item $(H^s, J^s)= (H_\pm, J_\pm)$ for $\pm s >>0$, and finally,
        \item (if $M$ is has contact boundary) slopes $c_s$ of $H_s$ are non-increasing, i.e. $c'_s \leq 0$. 
    \end{itemize}
\end{definition}

It is a standard result \cite{abouzaid-sh} that the Gromov-Floer bordifications of moduli spaces of Floer trajectories with Hamiltonian $H$ and complex structure $J$ are compact when $(H, J)$ an admissible Floer datum on an admissible symplectic manifold, as are the Gromov-Floer bordifications of continuation maps defined using $(H_s, J_s)$ when $(H_s, J_s)$ are an admissible continuation datum.

\begin{remark}
    In the case where $M$ is not exact, the action functional \eqref{eq:sympectic-action-functional} is not defined in general, but is only defined up to a constant. To make sense of the action functional, one chooses basepoints on each component of $LM$ and integrates the expression under the integral sign for $d\mathcal{A}_H$ \eqref{eq:hamiltonian-1-form} over a map from $S^1 \times I$ with one end sent to the basepoint loop and the other to the given loop $\gamma$. This gives a well-defined value for $\mathcal{A}_H(\gamma)$ due to the condition of being symplectically atoroidal. The rest of the discussion in this paper makes sense with this convention.
\end{remark}

\begin{definition}
\label{def:action-nondegenerate}
    A Hamiltonian $H: S^1 \times M \to \R$ is said to be \emph{action-nondegenerate} if the actions of all time-$1$ fixed points of its Hamiltonian flow are pairwise distinct. If $H = H_0^{\#k}$ then $H$ is said to be \emph{$C_k$-action-nondegenerate} if the sets of fixed points of $H$ with the same actions are exactly the $C_k$-orbits of $Fix(H)$. A Floer datum is ($C_k$-)action-nondegenerate if the underlying Hamiltonian is, and a continuation-map datum is ($C_k$-)action-nondegenerate if both the Hamiltonians $H_\pm$. 
\end{definition}

The following result is standard:
\begin{lemma}
\label{lemma:action-nondegeneracy-lemma}
    Action-nondegenerate Hamiltonians are $C^\infty$-dense in the space of all Hamiltonians. Similarly $C_k$-action-nondegenerate Hamiltonians are $C^\infty$-dense among all Hamiltonians $H$ of the form $H = H_0^{\#k}$.  $\blacksquare$
\end{lemma}

\subsection{Gradings in Hamiltonian Floer Homology}
\label{sec:grading-conventions}
We will now review our conventions for gradings in Hamiltonian Floer homology, and recall, via a standard analytic argument, why the vanishing of the polarization class (defined by Proposition \ref{prop:analytic-polarization-class}) implies that the gradings on Hamiltonian Floer Homology are absolute. 

First, recall that for a path of symplectic matrices $\gamma: [0,1] \to Sp(2n, \R)$ such that $\gamma(0) = id$ and $\gamma(1)$ does not have $1$ as an eigenvalue, Conley and Zehnder \cite{conley1984morse} defined an integer $CZ(\gamma) \in \Z$. Now, given a time-dependent Hamiltonian $H_t$ and a \emph{capped periodic orbit}, that is, a map $\tilde{x}: D^2 \to M$ such that  $x=u|_{\partial D^2}: S^1 \to M$ (where $S^1$ is parameterized  \emph{counter-clockwise}) is a $1$-periodic orbit of the Hamiltonian vector field of $H_t$, one can first use $\tilde{x}$ to pick a trivialization $\xi: u^*TM \to \R^{2n}$ as a symplectic vector bundle, and subsequently use $\xi$ and the linearized Hamiltonian flow of $H_t$ to define the Conley-Zehnder index of $u$ via $CZ(u) = CZ(t \mapsto \xi \circ d\phi^t_{x(0)} \circ \xi^{-1})$, where $\phi^t$ is the time-$t$ flow of $X_{H_t}$. The quantity $CZ(u)$ is always independent of the choice of trivialization $\xi$, and is unchanged upon deforming the capping $u$ while keeping the orbit $x$ fixed. Given two different cappings $\tilde{x}_1$ and $\tilde{x}_2$ of $x$, one can glue $\tilde{x}_1$ with the `reverse' $\tilde{x}_2'$ of $\tilde{x}_2$ to get a map $\tilde{x}_{S^2}: S^2 \to M$ sending the equator of $S^2$ to $x$; one can then establish a formula that $CZ(\tilde{x}_1) - CZ(\tilde{x}_2) = 2c_1(\tilde{x}_{S^2}^*TM)$; thus, if $c_1(TM)=0$ then the Conley-Zehnder index only depends on $x$, and can be written $CZ(x)$. 

Now, given a Floer trajectory $u$ satisfying \eqref{eq:floers-equation} with $u(s, t) \to x_\pm(t)$ as $s \to \pm \infty$, one an establish the formula 
\begin{equation}
    \label{eq:cz-index-formula}
    CZ(\tilde{x}_-) + \ind D_u = CZ(\tilde{x}_+)
\end{equation}
where $D_u$ is the linearized operator associated to $u$ \cite{robbin1995spectral} and one constructs the cap $\tilde{x}_+$ by concatenating $\tilde{x}_-$ with $u$ at the $s=-\infty$ end of $u$. When $c_1(M)=0$, the cappings can be dropped from the notation. Also, in the conventions of Conley and Zenhnder, one that for a $1$-periodic orbit corresponding to a nondegenerate fixed point of a time-independent Hamiltonian $H$, one has that 
\[ CZ(x) = \ind x - n\]
where $\ind x$ is the Morse index of $x$ as a critical point of $H$. Looking at the Floer equation in these coventions \eqref{eq:floers-equation}, this agrees with the formula \eqref{eq:cz-index-formula}: a $t$-independent solution to Floer's equation the pullback of a solution $\nu: \R \to M$ of the upwards gradient flow of $H$, and as such, the index of a critical point at $s = + \infty$ will be larger than that of the index of a critical point at $s = - \infty$. However, this is slightly inconvenient for our purposes, as we are trying to study the \emph{negative} gradient flow of the action functional \eqref{eq:sympectic-action-functional}, and so we want the index to \emph{decrease} from $s = - \infty$ to $s = + \infty$. Thus, we will Hamiltonian Floer homology (and Hamiltonian Floer cohomology) by the \emph{negative} of the Conley-Zehnder index:
\[ \widetilde{CZ}(\gamma) = -CZ(\gamma).\]
Thus, when $(H, J)$ are convex nondegenerate regular Floer data and $c_1(M)=0$, we will define Floer homology and cohomology groups
\[ HF^{cont}_*(H, J) = H_*(CF^{cont}_\bullet(H, J)), \;\;CF^{cont}_k(H, J) = \bigoplus_{\substack{x \in Fix(H) \\ \tilde{CZ}(x) = k}} \Z x,\;\; dx = \sum_{\substack{x \in Fix(H) \\ \widetilde{CZ}(y) = {k-1}}} \#\mathcal{M}(x,y)y, \]
\[ HF_{cont}^*(H, J) = H^*(CF^\bullet(H, J)), CF_{cont}^k(H, J) = \bigoplus_{\substack{x \in Fix(H) \\ \tilde{CZ}(x) = k}} \Z x, dx = \sum_{\substack{x \in Fix(H) \\ \widetilde{CZ}(y) = {k+1}}} \#\mathcal{M}(y,x)y, \]

The counts $\#\mathcal{M}(y,x)$ are the number of elements of this moduli space counted with sign, where the signs are defined by the coherent orientation procedure of \cite{floer1993coherent}. We will now briefly review another perspective on this procedure, as well as for the formula \eqref{eq:cz-index-formula}; a reference for this perspective in the setting of Lagrangian Floer homology can be found in \cite{seidel2008fukaya}, and for Hamiltonian Floer homology in \cite{abouzaid-sh}. To a cap for an orbit, we can associate the \emph{orientation operator} $D_u$ as follows. Thinking of $D^2$ as a Riemann surface via a subset of $\C$, there is a biholomorphism $\lambda$, well defined up to translation, of $D^2 \setminus 0 \to Z$ sending the outer boundary to the boundary at $s = + \infty$, and sending $1$ to the point at $s=+\infty$ and $t=0$. One can then define the Fredholm operator 
\[ D_{\tilde{x}}: L^{1, 2}(\tilde{x}^*TM) \to L^{2}(\tilde{x}^*TM)\]
by requiring that under the above biholomorphism, $D_{\tilde{x}}$ takes the form of the standard linearized Floer operator \eqref{eq:linearized-floer-equation}. Similarly, to a \emph{negative cap}, which is a map $\tilde{x}' :D^2 \to M$ such that $\tilde{x}|_{S^1}(t) = x(-t)$,  we can associated a corresponding operator operator $D_{\tilde{x}}$ which takes the form of the standard linearized operator under the biholomorphism $z \mapsto \lambda(1/z)$. This new operator will satisfy
\[ \ind D_{\tilde{x}'} = - \ind D_{\tilde{x}}.\]
In fact, one has that $CZ(u) = \ind D_u$, and the index formula \eqref{eq:cz-index-formula} is a consequence of linear gluing of operators 
\begin{equation}
\label{eq:gluing-of-operators}
D_{\tilde{x}_-}\# D_u = D_{\tilde{x}_+}
\end{equation} 
Similarly, the corresponding index formula for $\widetilde{CZ}$, namely 
\[ \widetilde{CZ}(\tilde{x}_-) =  \ind D_u + \widetilde{CZ}(\tilde{x}_+)\]
follows from a corresponding linear gluing argument using $D_{\tilde{u}}$. Now, given a Fredholm operator D, there is an determinant line 
\[ \det D = \ker D \tensor \coker^{\vee} D\]
and the orientation line $o(D)$, which is the rank $1$-free $\Z$-module such that the two generators of $o(D)$ are the two orientations of $det D$. 
Gluing formulae like \eqref{eq:gluing-of-operators} define isomorphisms 
\[ o(D_{\tilde{x}_-}) \tensor o(D_u) = o(D_{\tilde{x}_+});\]
thus, when $u$ has Fredholm index $1$ and $D_u$ is surjective, $o(D_u)$ is canonically isomorphic to $\Z$, and we get isomorphisms $o(u): o(D_{\tilde{x}_-}) \to o(D_{\tilde{x}_+})$ associated to elements of the moduli spaces $\mathcal{M}(y,x)$ used the definitions of the differentials for the Floer homology complexes. We turn the isomorphisms $o(u)$ into signs for Floer trajectories by choosing \emph{orientation data} $\mathfrak{o}$ given by a choice of cap $\tilde{x}$ for every fixed point $x$ choice of isomorphism $o(D_{\tilde{x}}) \simeq \Z$ for each such cap; by a gluing argument these choices induce such trivializations for \emph{all} caps of all fixed points, and thus one gets a sign for every index $1$ Floer trajectory between contractible fixed points. The same procedure is used to assign signs in the Floer cohomology complexes, but replacing $D_{\tilde{x}}$ with $D_{\tilde{x}'}$ everywhere and using the corresponding gluing formulae. The differentials on the Floer complexes above depend on this choice, but the final complexes are independent of $\mathfrak{o}$ up to isomorphism of complexes.

Let us now explain why this latter discussion can be made sense of, even in the setting of periodic orbits which are not contractible. Suppose that the polarization class 
\[ \rho \in KO^1(M)\]
defined by Proposition \ref{prop:analytic-polarization-class} vanishes. Now, writing $LM$ for a subspace of the free loop space of $M$ consisting of sufficiently regular loops, the discussion of Section \ref{sec:polarization-classes} explains that, writing 
\[ \mathcal{H} = L^2(S^1, \C^n), \mathcal{H}_1= W^{1,2}(S^1, \C^n),\]
we can trivialize the bundles $\mathcal{T}^{1,2}(LM)$ and $\mathcal{T}^{0,2}(LM)$ over $LM$ which are the corresponding Fredholm completions of the tangent bundle $TLM$ to trivial bundles $\underline{\mathcal{H}_1}$ and $\underline{\mathcal{H}}$, respectively. Moreover, in these trivializations, we have operators $A_\gamma \in S(\mathcal{H}, \mathcal{H}_1)$ corresponding to the Hessians \eqref{eq:hessian-operator-form}, and the operators  $D_{\tilde{u}}$, under the biholomorphism above, take the form 
\[ D_{\tilde{u}} = \partial_s + A_{u_s}\]
where $A_{u_s} \in S(\mathcal{H}, \mathcal{H}_1)$ is an unbounded self-adjoint fredholm operator on $\mathcal{H}$ of index zero defined on $\mathcal{H}_1$ (see \eqref{eq:hessians-after-trivialization}) (here $u_s\in LM$ is the restriction of $u$ to the circle around the origin with given $s$ coordinate under the biholomorphism $\lambda$). Thus, besides the trivializations of $\mathcal{T}^{i, 2}(LM)$, the operators $A_s$ will depend on the cap $\tilde{x}$; however, the operator $A_{+ \infty} = \lim_{s \to \infty} A_s \in S_*(\mathcal{H}, \mathcal{H}_1) = A_x$ will only depend on $x$. Set 
\[ A_0 = i \partial_t + \text{diag}(1, -1, 1, -1, \ldots).\]
The fact that $\rho = 0$ means that we can find  operators $\kappa(A_\gamma, t) \in \mathcal{S}(\mathcal{H}, \mathcal{H}_1), t\in[0,1]$ depending continuously on $\gamma$ and on $t$, such that $\kappa(A_\gamma, 1) = A_0.$ Now define a variant of $D_{\tilde{x}}$ which only depends on $x$ and on $\kappa$, namely, 
\[ D^o_x: L^{1, 2}(\R_s, \mathcal{H}) \cap L^2(\R_s, \mathcal{H}_1) \to L^2(\R_s, \mathcal{H})\]
\[D^o_x = \partial_s + \kappa(A_x, s'(-s))\]
where $s': \R \to (0,1)$ is the diffeomorphism \eqref{eq:make-r-compact}. Similarly, we can define $D^{o'}_x = \partial_s + \kappa(A_x, s'(s))$ as a variant of $D_{\tilde{x}'}$ that is manifestly independent of the capping. These operators will satisfy the analogs of the gluing rules \eqref{eq:gluing-of-operators}, and since $A_0$ is the operator one would assign to a hyperbolic fixed point of a Hamiltonian diffeomorphism, the proof of the formula \eqref{eq:cz-index-formula}
shows that 
\[ \ind D^o_x = CZ(x), \;\ind D^{o'}_x = \widetilde{CZ}(x)\]
if $x$ is a Hamiltonian orbit admitting a cap. Thus we define in general $\widetilde{CZ}(x)$ to be the index of the operator above, then we can define the Floer homology groups and complexes complexes generated by all periodic orbits as 
\[ HF_*(H, J) = H_*(CF_\bullet(H, J)),\;\; CF_k(H, J) = \bigoplus_{\substack{x \in Fix(H) \\ \widetilde{CZ}(x) = k}} \Z x, \;\;dx = \sum_{\substack{x \in Fix(H) \\ \widetilde{CZ}(y) = {k-1}}} \#\mathcal{M}(x,y)y, \]
\[ HF^*(H, J) = H^*(CF^\bullet(H, J)),\;\; CF^k(H, J) = \bigoplus_{\substack{x \in Fix(H) \\ \widetilde{CZ}(x) = k}} \Z x,\;\; dx = \sum_{\substack{x \in Fix(H) \\ \widetilde{CZ}(y) = {k+1}}} \#\mathcal{M}(y,x)y.\]
The signs above are defined by replacing $o(D_{\tilde{x}})$ and $o(D_{\tilde{x}})$ by $o(D^o_x)$ and $o(D^{o'}_x)$ in the definition of orientation data, and otherwise following the procedure outlined earlier.

\section{A simple construction of global charts in the exact case}
\label{sec:global-charts}

In the next section, we will construct a series of virtual smoothings of the pre-flow categories $\CC(H, J)$ associated to convex nondegenerate Floer data $(H, J)$, and the related pre-flow categories $\CC(H_s, J_s)$ associated to convex nondegenerate continuation data $(H_s, J_s)$. Moreover, when $(H, J)$ and $(H_s, J_s)$ consist of \emph{iterated} data (i.e they are invariant under the action of $C_k \subset S^1$), these smoothings will be equivariant and can be taken to satisfy additional compatibility conditions. Finally, in Section \ref{sec:producing-framings}, we will construct equivariant framings of all of these virtually smooth flow categories. 

Now, the existence of virtual smoothings of the pre-flow categories $\CC(H, J)$ is now new. In fact, these were constructed in \cite{bai2022arnold} and in \cite{rezchikov2022integral}. However, those constructions have the following two deficiencies:
\begin{itemize}
    \item Both of these constructions are based on the construction in \cite{AMS}, which relies heavily on  the space of genus zero curves in $\mathbb{P}^d$, $\overline{\mathcal{M}}_{0, n}(\mathbb{P}^d)$, which plays the role of a parameter space for domains (just as $\overline{Conf}_d$ does below). Now,  $\overline{\mathcal{M}}_{0, n}(\mathbb{P}^d)$ is a complex manifold, and this makes the resulting virtual smoothings into stably almost complex virtual manifolds (or virtual orbifolds, in the non-exact setting of \cite{AMS, bai2022arnold,rezchikov2022integral}). However, the tangent bundle of $\overline{\mathcal{M}}_{0, n}(\mathbb{P}^d)$ is \emph{nontrivial}; as such, it is more difficult to use this construction to produce \emph{framings}.
    \item The construction involves taking quotients by various compact Lie groups, and the resulting composition maps in the virtual smoothing are comparatively complicated, as one is forced to consider certain nontrivial group-theoretic relationships between $\prod_{i=1}^n \overline{\mathcal{M}}_{0, n_i}(\mathbb{P}^{d_i})$ and a certain subspace of $\overline{\mathcal{M}}_{0, \sum_i n_i- 2(n-1)}(\mathbb{P}^{\sum_i d_i})$. As such, even if one resolves the issue discussed previously regarding framings of the parameter spaces, constructing \emph{equivariant framings} of the resulting thickenings would involve a fair amount of careful equivariant index theory involving various Lie groups. The reader may compare this with the construction in Section \ref{sec:producing-framings}, which is relatively straightforward, but already must be written carefully so as to avoid potential complications associated to induction up the flow category. (Of course, other constructions which \emph{do} run into these issues are not written in this paper, but the author considered many of them prior to  settling on the construction of Section \ref{sec:producing-framings}.) Part of the purpose of this paper is to clean up the analytic construction of framings in the setting of Hamiltonian Floer homology, which has previously been treated before in ways that are opaque to the author \cite{cohen2007floer, large2021spectral}; the author was unsure based on these previous works if the construction of equivariant framings could be done without issue.
\end{itemize}

One might think that one does not need global charts in the setting of this paper, since we are in the symplectically aspherical case; however, because we are working with equivariant Floer data, we may not be able to achieve transversality while preserving equivariance.  Since solutions to Floer's equation have $\mathbb{R}$-symmetry,  to choose the perturbations of Floer's equation consistently, we must find a way to get rid of this symmetry. One classic method is to stabilize the domains by looking at their intersections with certain hypersurfaces \cite{cieliebak-mohnke, pardon2016algebraic}. However, choices of local hypersurfaces cannot produce global charts, while choosing global hypersurfaces transverse to all curves may or may not exist. However, because we are working in the exact setting, we have a canonical function which decreases monotonically along each Floer trajectory: the \emph{action}. We can thus kill the $\mathbb{R}$-symmetry of Floer trajectories by requiring that distinguished $s$-coordinates of the Floer trajectories have particular values of the action.

\subsection{Thickenings for Hamiltonian Floer Homology}
\label{sec:thickenings-for-hamiltonian-floer-homology}
In this section, we explain how to define a convenient Kuranishi chart on an \emph{uncompactified} moduli space of Floer trajectories. 

Set
\[Conf_r' = \{ (x_1, \ldots, x_n) \in \mathbb{R}^n \, | \, x_1 < x_2 < \ldots < x_r \}. \]
There is an action of $\mathbb{R}$ on $Conf'_n(\mathbb{R})$ by simultaneous translation of the $x_i$. We define $Conf_r(\mathbb{R}) = Conf'_r(\mathbb{R})/\mathbb{R}$.

The space $Conf_r(\mathbb{R})$ parameterizes domains for the Floer equation equipped with a sequence of labeled $s$-coordinates. In particular, there is a universal curve 
\[ \pi: Conf_r(\mathbb{R})^Z \to Conf_r(\mathbb{R}) \]
with fibers isomorphic to $Z=\R_s \times S^1_t$ in a way that is canonical to $s$-translation. (More naturally, there is a universal curve 
\[ Conf_r(\mathbb{R})^{\R} \to Conf_r(\mathbb{R})\]
with each fiber a copy of $\R$ with $r$ marked points;  the space $Conf_r(\mathbb{R})^{\R}$ is the quotient of $Conf_r(\mathbb{R})^Z$ by a fiberwise $S^1$ action which rotates the $S^1$ coordinate in $Z$.)
One constructs the universal curve via 
\[Conf_r(\mathbb{R})^Z = (Z \times Conf_r'(\mathbb{R}))/\R \]
where the quotient is by the $\R$ action which simultaneously translates $Z$ and $(x_1, \ldots, x_r)$. Given $a \in Conf_r(\mathbb{R})$, we write $Conf_r(\mathbb{R})^Z_a$ for $\pi^{-1}(a)$.

We now introduce some additional helpful notation. Since a totally ordered finite set is determined up to unique isomorphism by its cardinality, given another totally ordered set $S$ such that $|S|=r$, we can write $Conf_S$ for $Conf_r$ but with the markings $x_j$ thought of as labeled by points of $S$ in an \emph{order-reversing} fashion; thus for $u \in Conf_S$, the fiber $\pi^{-1}(u)$ is equipped with a collection of marked $s$-coordinates $x_{\sigma}$ for $\sigma \in S$. We will similarly use notation $Conf_S^Z$, etc, to denote the corresponding spaces with $S$ replaced by $|S|$ but with the marked coordinates labeled by elements of $S$. The order reversing convention is because the indices of the points are increasing from left to right; but we will shortly be labeling the points with \emph{actions}, which are \emph{decreasing} from left to right.

Given a Floer trajectory $u: Z \to M$ and a sequence of real numbers $\mathcal{A}_1>  \ldots> \mathcal{A}_r$ which are not actions of any Hamiltonian orbit, one gets a corresponding point $S(u) \in Conf_r(\mathbb{R})$ by
\[ S(u) = (x_1, \ldots, x_n) \text{ where } \bar{A}(u_{x_i}) = \mathcal{A}_i. \]
Translating $u$ in the $s$ direction defines the same point in $Conf_r(\mathbb{R})$, so the above prescripttion defines a continuous, well-defined map  $\mathcal{M}(x,y) \to Conf_r(\mathbb{R})$ for every pair $(x, y)$ of Hamiltonian trajectories of $H$ for which $\mathcal{A}(x) > \mathcal{A}_1$ and $\mathcal{A}_r > \mathcal{A}(y)$. 

Write
\[ C^\infty(\Omega^{0,1}_Z \tensor TM)' \subset C^\infty(Z \times M, \Omega^{0,1}_Z \tensor TM) \]
for the linear subspace of sections which are zero on $(((-\infty, T) \cup (T, \infty))\times S^1) \times M$ for some $T$, and are also zero on $Z \times U$ for some neighborhood $U$ of the images of all time-1 orbits of $H$. 

Let us choose an orthogonal $S^1$-representation $V$ and an $S^1 \times \R$-equivariant map
\[ \tilde{\lambda}_V: V \times Conf'_r(\mathbb{R}) \to C^\infty(\Omega^1 Z \tensor TM)'\]
Here, $S^1$ acts trivially on $Conf'_r(\mathbb{R})$ and acts on the via pushforward of sections under the action of $S^1$ on $Z$ by rotation. The other factor $\R$ acts on the codomain by pushforward of sections under $s$-translation in the $Z$ factor, and $\R$ acts trivially on $V$. Moreover, we also require that $\tilde{\lambda}_V$ a linear map on each fiber of $V \times Conf'_r(\mathbb{R})$ over $Conf'_r(\mathbb{R})$. 

Equivalently, this is the data of an $S^1$-equivariant bundle map over $Conf_r(\mathbb{R})^Z\times M$ of the form
\begin{equation}
    \label{eq:perturbation-data-1}
    \lambda_V: \underline{V} \to \Omega^{0,1}_{Conf_r(\mathbb{R})^Z/Conf_r(\mathbb{R})} \tensor TM,
\end{equation} 
where $\Omega^{0,1}_{Conf_r(\mathbb{R})^Z/Conf_r(\mathbb{R})}$ denotes the bundle of fiberwise $(0,1)$-forms, and where for every $a \in Conf_r$, the restriction of  $\lambda_V$ to $Conf_r(\mathbb{R})^Z_a \times M$ is zero at the ends and is zero on $Conf_r(\mathbb{R})^Z_a\times U$ for some neighborhoods $U$ of the orbits of $H$.

Then we can define the moduli space
\begin{equation}
    \label{eq:one-smoothed-moduli-space}
    T(x_-, x_+) = \left\{ \begin{array}{c|c} 
    (u,a,v) : v \in V, a \in Conf_r(\mathbb{R}), &  \partial_s u(s,t) + J_t \partial_t u(s,t) -\Grad H(u(s,t)) 
 = \lambda_V(((s,t), a), v), \\
  u: Conf_r(\mathbb{R})^Z_a \to M  & \lim_{s \to \pm \infty} u(s, \cdot) = x_\pm(\cdot). 
 \end{array}\right\}.
\end{equation} 

In the next sections we will define all moduli spaces as families of maps from fibers of universal curves over certain parameter spaces of domains. However, it may be helpful for the reader's understanding to note here that the above moduli space can be equivalently written as 

\begin{equation}
    \label{eq:one-smoothed-moduli-space}
    T(x_-, x_+) = \left\{ \begin{array}{c|c} 
    (u,a',v) : v \in V, a' \in Conf'_r(\mathbb{R}), &  \partial_s u(s,t) + J_t \partial_t u(s,t) -\Grad H(u(s,t)) 
 = \widetilde{\lambda}_V(v, a'), \\
  u: Z \to M  & \lim_{s \to \pm \infty} u(s, \cdot) = x_\pm(\cdot). 
 \end{array}\right\}/\mathbb{R}
\end{equation} 
where the $\R$ action translates $u$ and $a'$ simultaneously.

This moduli space acquires a $C_k$ action if $H_t, J_t$ are $C_k$ equivariant. Moreover for any choice of $(H_t, J_t)$, standard transversality theory shows that there is a choice of $V$ and $\lambda_V$ such that the above moduli space is a smooth $C_k$-manifold.

There is a natural section $\sigma: T(x_-, x_+) \to \underline{V \times \R^r}$ given by 
\[ \sigma(v, a,v) = (v, \mathcal{A}_1 - \mathcal{A}(u|_{x_1}), \ldots \mathcal{A}_r - \mathcal{A}(u|_{x_r})) \]
where $a = (x_1, \ldots, x_r)$. The zero-section of $\sigma$ is exactly the moduli space of Floer trajectories, because if $v = 0$ then $u$ solves the Floer equation, and for every Floer trajectory there is a unique collection of points on $\mathbb{R}$ (up to simultaneous translation of $u$ and of the points) such that the actions at the $s$-coordinates given by the points agree with the actions $\mathcal{A}_1, \ldots, \mathcal{A}_r$. 

In order to build smoothings of $\CC(H, J)$, we must find compatible ways of choosing $\lambda_V$ over the moduli spaces $Conf_r(\mathbb{R})$ such that products of moduli spaces \eqref{eq:one-smoothed-moduli-space} glue together to  topological $\langle k \rangle$-manifolds, and subsequently find smoothings of stabilizations of these manifolds, as in \cite{AMS, rezchikov2022integral, bai2022arnold}. 

\subsection{Parameter spaces and integralization data for global charts}
\label{sec:parameter-spaces-for-global-charts}
The real locus of genus zero Deligne-Mumford space $\overline{\M}_{0, n}(\R)$ is widely used as a parameter space of domains in various pseudoholomorphic curve constructions \cite{fukaya1997zero, seidel2008fukaya}. In this section, we define a certain family of moduli spaces built out of the spaces $\overline{\M}_{0, n}(\R)$ which we will subsequently use to compatibly choose perturbation data.

\begin{figure}[h!]
    \centering
    \includegraphics[width=\textwidth]{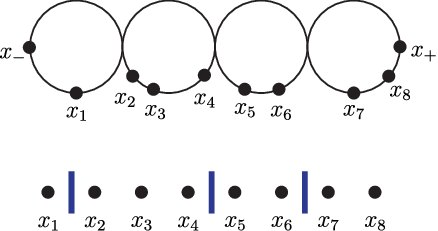}
    \caption{\emph{Parameter spaces used for global charts, and their stratifications.} The moduli space of marked genus zero real curves $\overline{M}_{0, r+2}(\R)$ contains a subset $\overline{Conf}_r$ which is a natural partial compactification of $Conf_r$. A typical curve in this subset is depicted on top. Each component of this curve is, after removing $x_\pm$ and the nodes, biholomorphic to a strip in a way unique up to $s$-translation, with the marked points giving marked $s$ coordinates. Replacing these strips by cylinders gives moduli spaces of domains of broken Floer trajectories with marked $s$-coordinates, which are used in the stabilization procedure of this section.  (Bottom) One stratifies $\overline{Conf_{r}}$ as a $\langle r-1\rangle$ manifold by identifying $[r-1]$ with the spaces in between $r$ elements in a row, and placing a marker (i.e. removing an element from $[r-1])$ for every node/break-point of the corresponding curve.}
    \label{fig:curve-moduli-space}
\end{figure}

The connected components of $\overline{\M}_{0, n}(\R)$ are labeled by cyclic orderings of the marked points $x_1, \ldots, x_n$. Let us call the labeled points of $\overline{\M}_{0, r+2}(\R)$ as $x_-, x_1, \ldots, x_r, x_+$, and 
let $\overline{Conf}_r \subset \overline{\M}_{0, r+2}(\R)$ be the subset consisting of marked bordered nodal Riemann surfaces for which the cyclic ordering of points on the boundary of the curve is exactly $(x_-, x_1, \ldots, x_r, x_+)$, and such that every vertex of the labeled tree associated to $a$ has at most two neighbors, and such that the extremal vertices of the tree (of which there are at most two) contain $x_-$ and $x_+$, respectively (Figure \ref{fig:curve-moduli-space}).

It is well known that $\overline{\M}_{0, r+2}(\R)$ is a manifold with corners, and it is straightforward to see that $\overline{Conf}_r(\R)$ is thus also a manifold with corners, and in fact a smooth $\langle r-1 \rangle$-manifold. Moreover, the locus of curves with exactly one component in $\overline{Conf}_r$, i.e. the top-dimensional open stratum, is in bijection with $Conf_r$. Now $\overline{\mathcal{M}}_{0, r+2}(\R)$ admits a map from the corresponding universal curve $\overline{\mathcal{M}}_{0, r+2}(\R)^{\bar{D}}$, the fibers of which can be identified with trees of oriented disks with certain marked points on their boundaries.   we define the universal curve $Conf_r^{\bar{D}}$ to be the restriction of $\overline{\mathcal{M}}_{0, r+2}(\R)^{\bar{D}}$ to $Conf_r$. A typical element fiber of $Conf_{r}^{\bar{D}}$ is depicted in Figure \ref{fig:curve-moduli-space}. Every irreducible component of such a fiber is a disk with either zero, one, or two nodes, and if it has $w$ nodes then the marked points it carries intersect the subset $\{x_+, x_+\}$ in a set of cardinality $2-w$.  

Let us define $\overline{\mathcal{M}}_{0, r+2}(\R)^{D}$ by removing from every fiber the nodes and the marked points $x_{\pm}$. Then every component of the a fiber of $\overline{Conf}_r^{D}$ is biholomorphic to $\R \times [0,1]$ in a manner unique up to $s$-translation such that the marked points on the fiber are sent to the $\R \times \{0\}$ boundary. This defines a relation on $\overline{Conf}_r^{D}$ by saying that two points are equivalent if they lie in the same fiber over $\overline{Conf}_r$, on the same component of the fiber, and if under any (or equivalently all) biholomorphisms as above they have the same $s$ coordinate. The quotient of $\overline{Conf}_r^{D}$ by this relation is $\overline{Conf}_r^{\R}$; its fibers are disjoint oriented copies of $\R$ each with a coordinate $s$ defined up to $s$ translation, and with marked $s$ coordinates defined up to the same $s$ translation. Similarly, taking the quotient of $\overline{Conf}_r^{\bar{D}}$ by the same relation as before defines  $\overline{Conf}_r^{\bar{\R}}$, which is a compactification of $\overline{Conf}_r^{\R}$. We set 
\[ \overline{Conf}_r^{Z} = S^1 \times \overline{Conf}_r^{\R}, \overline{Conf}_r^{\bar{Z}} = \overline{Conf}_r^{\bar{\R}};\]
these will be universal curves over $\overline{Conf}_r$ as well, with the former parameterizing domains of Floer trajectories with auxilary marked $s$-coordiantes on each component, well defined up to $s$-translation; and the latter parameterizing the same but with the condition that adjacent Floer trajectories must limit at corresponding ends to the same Hamiltonian orbit.

\newcommand{\auxactions}{{\mathfrak{A}}}
\newcommand{\paramspace}{{\overline{M}}}

\paragraph{Integralization data and relations between parameter spaces.}
Fix a convex linear floer datum $(H, J)$. 
In this section, write $\mathcal{A}$ for a fixed function which can be written as 
\[ \mathcal{A} = c\mathcal{A}_H + d \text{ for some } c> 0, d \in \R \]
where $\mathcal{A}_H$ is the action \eqref{eq:sympectic-action-functional} associated the Hamiltonian $H$. We will refer to $\mathcal{A}$ as the \emph{adapted action}. 
Suppose now that we have chosen a finite set of action values $\auxactions \subset \R$ which are not adapted actions of any Hamiltonian orbit of $H$. We assume that 
\begin{equation}
\label{eq:action-condition}
    \text{For every }x \in Fix(H)\text{ there is some }A \in \auxactions \text{ such that for all } y \in \CC, \mathcal{A}(x) > \mathcal{A}(y), \text{ we have } \mathcal{A}(x) > A > \mathcal{A}(y). 
\end{equation}

Given such a choice, the pre-flow category $\CC(H, J)$, which is a priori just a topologically enriched category, becomes a flow category with integral action defined below. We set the integral action to be
\[\bar{A}(x) = \#\mathfrak{A}(x)-1, \text{ where } \mathfrak{A}(x) =\{A \in \auxactions : \mathcal{A}(x) > A\}, \]

For any $y \in \CC(H, J)$ which is not a maximum of $\CC$, we will write $A_y^+$ for the smallest element of $\mathfrak{A}$ that is greater than $\mathcal{A}(y)$.

As usual we will write write $\bar{A}(x,y) = \bar{A}(x) - \bar{A}(y) -1$; then $\bar{A}(x,y)$ is the cardinality of the set
\[ S_{\CC(x,y)} =\mathfrak{A}'(x,y) := \mathfrak{A}(x,y) \setminus A_y^+, \text{ where }\mathfrak{A}(x,y) = \{A \in \mathfrak{A}: \mathcal{A}(x) > A > \mathcal{A}(y).\}\]
This convention specifies the lifts of the Floer moduli spaces $\cM(x,y)$ to $\langle k \rangle$-spaces uniquely. 

 Write 
\[F(x,y) = \{ z \in Fix(H): \mathcal{A}(x) > \mathcal{A}(z) > \mathcal{A}(y).\} \]

For every pair $x>y$ of objects of $\CC(x,y)$, we define 
\begin{equation}
    \paramspace(x,y) = \paramspace(\bar{A}(x), \bar{A}(y)) = \overline{Conf}_{\mathfrak{A}(x,y)}. 
\end{equation}
The notation means that we have multiple names for the same manifold $\paramspace(\bar{A}(x), \bar{A}(y))$ (which indeed manifestly only depends on $(\bar{A}(x), \bar{A}(y))$), but if $(n_1, m_1) \neq (n_2, m_2)$ then we do not think of $\paramspace(n_1, m_1)$ as equal to $\paramspace(n_2, m_2)$. 
In any case, this is a $\langle \bar{A}'(x,y)\rangle$-manifold: thinking of the set as labeling the \emph{spaces between} the marked actions of $\mathfrak{A}(x,y)$, we say that $a \in \paramspace(x,y)$ is the stratum corresponding to $\bar{A}' (x,y) \setminus S$ if the break-points between the components of $\pi^{-1}(a)$ lie in the spaces corresponding to the elements of $S$ (Figure \ref{fig:curve-moduli-space}).

We define the universal curve $\overline{M}(x,y)^{\overline{Z}}$  over $\overline{M}(x,y)$ to be the corresponding universal curve over $\overline{Conf}_{\mathfrak{A}(x,y)}(\mathbb{R})$, and similarly for all other universal curves 
$\overline{M}(x,y)^{Z}$, $\overline{M}(x,y)^{\overline{\R}}$, and $\overline{M}(x,y)^{\R}$.

We have an $F^s$-parameterization $\mathcal{F}^1_{can}$ associated to an $E^s$-parameterization $\mathcal{E}^1_{can}$ given by 
\begin{equation}
    \label{eq:canonical-obstruction-f-parameterization}
    V^1_{can}(x,y) = \R^{\mathfrak{A}(x,y)},  V^1_{can}(x) = \R^{\mathfrak{A}(x)}.
\end{equation}
\begin{remark}
    The coordinates of $\R^{\mathfrak{A}(x,y)}$ will be used below to remove the degrees of freedom introduced by the marked points via an obstruction section. There is no analog of $V^0_{can}(x,y)$, which would be used to frame the parameter space $\paramspace(x,y)$ themselves; instead, they will be framed via $T\R_+^{\bar{A}'(x,y)}$ in Section \ref{sec:producing-framings} below. Thus there is no $V^0_{can}$ because we separate out the `corner' framings from the remaining data of the framing of the thickening in Definition \ref{def:framing-of-flow-category}.
\end{remark}

\begin{remark}
\label{rk:can-incorporate-shift-in-canonical-param}
    Let $\mathfrak{A}_{lower} = \cap_{x \in \CC} \mathfrak{A}(x).$ Then for any subset $\Delta \mathfrak{A}_{lower} \subset \mathfrak{A}_{lower}$, we can replace $\R^{\mathfrak{A}(x)}$ in the defintion of $V^1_{can}$ above \eqref{eq:canonical-obstruction-f-parameterization} with $\R^{\mathfrak{A}(x) \setminus \Delta \mathfrak{A}_{lower}}$ and this would still give rise to the same $F^s$-parameterization $\{V^1_{can}(x,y)\}$. Thus the rest of the construction of the section would go through as before, but eventually when constructs the Floer homotopy type one would end up having a homotopy type that is de-suspended by $\R^{\Delta \mathfrak{A}_{lower}}$. This variant of the construction arises when one studies extensions of integralization data (as in Lemma \ref{lemma:extend-perturbation-data-across-restratifications}), e.g. when one thinks of the integralization datum for a continuation map $\CC(H_s, J_s)$ (see Section \ref{sec:continuation-maps-geometry}) as restricting to an integralization datum for $\CC(H_-, J_-)$, in which case $\mathfrak{A}_{lower}$ is large and incorporates the extra element $A_c$ which does not have a corresponding coordinate of $V^1_{can}$ (see Remark \ref{rk:shift-in-integralization-due-to-continuation}). See also the restratification embedding shifts of Proposition \ref{lemma:restratifications-exist}.

\end{remark}
For every $z \in F(x,y)$, there are natural inclusions $\iota_{xzy}: \overline{M}(x,z) \times \overline{M}(y,z) \to \overline{M}(x,z).$
There are natural associative inclusion maps 
\begin{equation}
    \iota_{xzy}: \paramspace(x,y) \times \paramspace(y,z) \to \paramspace(x,z).
\end{equation}
Moreover every $\chi_1 \in \overline{M}(x,z)$ and $\chi_2 \in \overline{M}(z,y)$, there is an identification
\begin{equation}
    \label{eq:identification-of-universal-curve-domains}
    \iota'_{xzy}: \overline{M}(x,z)^Z_{\chi_1} \cup \overline{M}(z,y)^Z_{\chi_2} \to \overline{M}(x,y)^Z_{\iota_{xzy}(\chi_1, \chi_2)}
\end{equation}
which covers the above inclusions, and which identify the fibers of the latter with the unions of the fibers of the former. For any $x_- > x_1 > \ldots > x_r > x_+$ in $\CC(H, J)$, the  maps $\iota_{xzy}$ define maps 
\[\paramspace(x_-, x_1) \times \ldots \times \paramspace(x_r, x_+)\to  \paramspace(x_-, x_+)   \]
by composition, and any such composition defines the same map. There are corresponding well-defined identifications of fibers of universal curves like \eqref{eq:identification-of-universal-curve-domains} defined by composing the maps \eqref{eq:identification-of-universal-curve-domains}.

Let us now make the following 
\begin{equation}
\label{eq:assume-action-nondegeneracy}
\parbox{0.8\textwidth}{ \textbf{Assumption:} $H$ is action-nondegenerate, or more generally, if $H= H_0^{\#k}$ and we are working in the equivariant setting then $H$ is $C_k$-action-nondegenerate (which is a weaker condition).}
\end{equation}
\begin{definition}
\label{def:integralization-data}
A choice of \emph{integralization data} $\bar{\mathfrak{A}}$ for $\CC(H, J)$ is a tuple $(\mathcal{A}, \mathfrak{A})$ such that $\mathcal{A}$ agrees with the action of the Hamiltonian $H$ up to composition with a  affine-linear orientation-preserving function and $\mathfrak{A}$ satisfies \eqref{eq:action-condition}. 
\end{definition}

\subsection{Compatible systems of perturbation data}
\label{sec:perturbation-data}

\begin{definition}
\label{def:system-of-perturbation-data}
    A \emph{compatible system of perturbation data} for $\CC(H, J)$ is a tuple 
\begin{equation}
\label{eq:system-of-perturbation-data}
    \mathfrak{P} = (\bar{\mathfrak{A}}, \mathcal{V}^1 = (\{V^1(x,y)\}, \{\tau^1_{xyz}\}), \{\lambda_{V^1(x,y)}\})
\end{equation}
where $\bar{\mathfrak{A}} = (\mathcal{A}, \mathfrak{A})$ is a choice of integralization data, $\mathcal{V}^1$ is an $F^s$-parameterization arising from an $E^s$-parameterization $\mathcal{E}^1$, and for every such $(x,y)$, the \emph{perturbation function} 
\begin{equation}
    \label{eq:perturbation-data-2}
    \lambda_{V^1(x,y)}: \underline{V^1(x,y)} \to \Omega^{0,1}_{\overline{M}(x,y)^Z/\overline{M}(x,y)} \tensor TM,
\end{equation} 
is a bundle map over $\overline{M}(x,y)^Z \times M$ which is zero when restricted to $U \times M$ for some neighborhood $U$ of the (removed) nodes of the fibers of the universal curve, and is also equal to zero when restricted to $\overline{M}(x,y)^Z \times V$ for some neighborhood $V$ of the images of all Hamiltonian orbits of $H$. 

Moreover, for every triple $x > z > y$ objects of $\CC(H, J)$, and every pair $\chi_1 \in \overline{\M}(x,z)$, $\chi_2 \in \overline{\M}(z,y)$, we have 
\[ (\iota'_{xzy})^*\lambda_{V(x,y)}( \cdot, \tau_{xyz}(v_1, v_2)) = \lambda_{V(x,z)}( \cdot, v_1) + \lambda_{V^1(z,y)}(\cdot, v_2) \]
where on the right hand side we take the extension by zero of the $TM$-valued $1$-form across the other component of the domain of $i'_{xzy}$ \eqref{eq:identification-of-universal-curve-domains}. 
\end{definition}

Given a compatible system of  perturbation data, for every $x_- > x_+$ in $\CC(H, J)$, we define 
\begin{equation}
\label{eq:thickenings-for-floer-trajectories}
    T(x_-, x_+) = \left\{ \begin{array}{c|r} 
    (u,a,v) : v \in V^1(x_-, x_+), &  \partial_s u(s,t) + J_t \partial_t u(s,t) -\Grad H(u(s,t))= \lambda_V(((s,t), a), v),\\
 a \in \paramspace(x,y), u: \paramspace(x,y)^Z_a \to M &   \lim_{s \to -\infty} u_1(s, \cdot) = x_-, \lim_{s \to \infty} u_r(s, \cdot) = x_+\\ 
   u = u_1 \cup \ldots \cup u_r  & \lim_{s \to + \infty} u_j(s, \cdot) = \lim_{s \to - \infty} u_{j+1}(s, \cdot), j=1, \ldots, r-1
 \end{array}\right\}.
\end{equation}
Above  $u = u_1 \cup \ldots \cup u_r$ is the decomposition of $u$ as a union of maps from the connected components of the domain, each one of which is canonically isomorphic to $Z$ up to a global $s$-translation. There is a natural map 
\[ \pi: T(x_-, x_+) \to \paramspace(x_-, x_+), \pi(u,a,v) = a.\]

Now, for every sequence $x_- > x_1 > \ldots > x_r > x_+$ in $\CC(H, J)$, there is an inclusion 
\begin{equation}
    \label{eq:thickening-boundary-inclusions}T(x_-, x_1) \times \ldots \times T(x_r, x_+)\subset T(x_-, x_+)
\end{equation}
in which one concatentates the domains of the maps $(u_1, \ldots, u_{r+1})$ and applies the (well defined) isometry 
\[ V^1(x_-, x_1) \oplus \cdots \oplus V^1(x_r, x_+) \to V^1(x_-, x_+) \]
defined by the the $F'$-parameterization $\mathcal{V}^1$ to the corresponding vectors $(v_1, \ldots, v_{r+1})$. Moreover the inclusions \eqref{eq:thickening-boundary-inclusions} fit into commutative diagrams 
\begin{equation}
\label{eq:thickening-composition-maps}
    \begin{tikzcd}
        T(x_-, x_1) \times \ldots \times T(x_r, x_+)\ar[r] \ar[d] &  T(x_-, x_+) \ar[d]\\
        \paramspace(x_-, x_1) \times \ldots \times \paramspace(x_r, x_+)\ar[r] &  \paramspace(x_-, x_+)
    \end{tikzcd}
\end{equation}
where the vertical arrows are products of projections, and the lower horizontal arrows defined by composing the maps $\iota_{xzy}$.

The space $T(x_-, x_+)$  carries a vector bundle $\underline{V^1_T(x_-, x_+)}$ with fiber 
\begin{equation}
    \label{eq:obstruction-bundle-for-my-stabilization}
    V^1_T(x_-, x_+) = V^1(x_-, x_+) \oplus V^1_{can}(x,y).
\end{equation} 
These vector spaces are the vector spaces associated to the $F^s$-parameterization $\mathcal{V}^1_T = \mathcal{V}^1 \oplus \mathcal{V}^1_{can}$, which is embedded into the $E^s$-parameterization $\mathcal{E}^1 \oplus \mathcal{E}^1_{can}$. 
There is a section $\sigma(x_-, x_+)$ of this vector bundle defined as follows. We define the function 
\[ \delta A(x,y): T(x,y) \to \R^{\mathfrak{A}(x,y)}, \delta A(x,y) (u, a, v) = \left(\mathcal{A}(u|_{s_A}) - A\right)_{A \in \mathfrak{A}(x,y)} \]

We then define 
\begin{equation}
\label{eq:obstruction-section-floer-trajectories}
    \sigma(x_-, x_+)(u, a, v) = (v, \delta A(x_-, x_+)(u, a,v))
\end{equation}

\begin{figure}[h!]
    \centering
    \includegraphics[width=\textwidth]{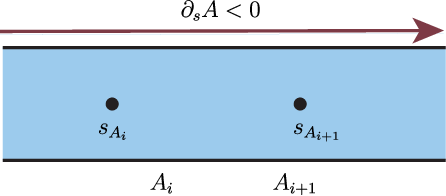}
    \caption{\emph{The stabilization procedure.} The action (and thus the modified action) is strictly decreasing along a Floer trajectory. Thus, given a (broken) Floer trajectory, there a unique $s$-coordinate at which the corresponding loop has a (modified) action $A_i$ in between the modified actions of the endpoints. This defines the embedding of spaces of Floer trajectories into the thickenings \eqref{eq:thickenings-for-floer-trajectories} of the global charts of this section; in general the obstruction sections $\delta A$ measure the difference in modified action between the value at $s_{A_i}$ and the value $A_i$ itself. Note that if one wished to add on additional marked modified actions $A_i$, then (for small deviations from zero) there would be a unique location to put $s_{A_i}$ such that $\mathcal{A}(u|_{s_{A_i}})$ differs from $A_i$ by a given value; this defines the correspondence between restratifications of $\CC'(H, J)$ and extensions of integralization data used in Lemma \ref{lemma:extend-perturbation-data-across-restratifications} and Proposition \ref{prop:cyclotomic-compatibility-all-perturbation-data-choices}.}
    \label{fig:description-of-marking-fixing-condition}
\end{figure}

It is clear that $\sigma(x_-,x_+)^{-1}(0)$ is exactly the space of Floer trajectories $\CC(H, J)(x_-, x_+)$ (Figure \ref{fig:description-of-marking-fixing-condition}). Thus, the tuples $(T(x_-, x_+), V^1_T(x_-, x_+), \sigma(x_-, x_+))$ will serve as our Kuranishi thickenings of $\CC(H, J)(x_-, x_+)$, once we choose appropriate compatible systems of perturbation data and subsequently find smoothings of $T(x_-, x_+)$. 

\paragraph{Equivariant perturbation data. }

Suppose now that $(H, J)$ are $C_k$-invariant, so that $\CC(H, J)$ is a $C_k$-equivariant pre-flow category. To define the spaces $T(x,y)$ we needed to  perturbation data in an equivariant way.
Now, there is a trivial $C_k$-action on $\paramspace(x,y)$. 
This action lifts to an action on $\paramspace(x,y)^Z$ and the other universal curves by also rotating the fibers of the universal curves by $C_k$. 

A $C_k$-equivariant compatible system of perturbation data $(\bar{\mathfrak{A}}, \mathcal{V}^1 = (\{V^1(x,y)\}, \{\tau^1_{xyz}\}), \{\lambda_{V^1(x,y)}\})$ is a compatible system of perturbation data such that $\mathcal{V}^1$ is a $C_k$-equivariant $F$-parameterization and such that the maps \eqref{eq:perturbation-data-2} are $C_k$-equivariant, where the action of $C_k$ on $\paramspace(x,y)^Z$ given by rotating the fibers of all universal curves is the action described above, and the action on $\underline{V^1(x,y)}$ is the natural equivariant structure induced from this action and the $C_k$-representation $V^1(x,y)$.  

When $H$ and $J$ are $C_k$-invariant and the spaces $T(x,y)$ are defined using a $C_k$-equivariant compatible system of perturbation data, the spaces $T(x,y)$ acquire a $C_k$-action the product of the $C_k$ action on $V^1(x,y)$ and the action on $\paramspace(x,y)^Z$. The inclusions \eqref{eq:thickening-boundary-inclusions} are $C_k$ equivariant. 

Letting $C_k$ act on $\underline{V^1_T(x,y)}$ by the action induced from the $C_k$-equivariant structure of $\mathcal{V}^1_T$ covering the $C_k$-action on the spaces $\bigsqcup_{x_-, x_+} T(x_-, x_+)$, the sections $\sigma(x_-, x_+)$ become $C_k$-equivariant in the sense that $\sigma(x_-, x_+)(gy) = g\sigma(x_-, x_+)(y)$ for all $g \in C_k$ and all $y \in T(x_-, x_+)$ (where the action of $g$ takes $T(x_-, x_+)$ to $T(gx_-, gx_+)$.)

\paragraph{Existence of regular systems of perturbation data.}

Given $u \in \CC(H, J)(x_-, x_+)$, there is a unique $S(u) \in M(x_-, x_+)$ given by the composition of the above identification $\CC(H, J)(x_-, x_+) \subset T(x_-, x_+)$ with the projection to $\paramspace(x,y)$. 

\begin{proposition}
\label{prop:compatible-perturbation-data-exist}
There exists a compatible system of perturbation data $\mathfrak{P}$ such that for every $u \in \CC(H, J)(x_-, x_+)$,
\begin{equation}
\label{eq:stabilized-linearized-operator}
    D^\mathfrak{P}_u: W^{k, 2}(Z, u^*TM) \oplus V^1(x_-, x_+) \to W^{k-1, 2}(Z, \Omega^{0,1} \tensor u^*TM)
\end{equation}
\[ D^\mathfrak{P}_u(\xi, v) = D_u(\xi) - \lambda_{V^1(x_-, x_+)}|_{\paramspace(x,y)^Z_{S(u)}}(v)\]
is surjective for all $k >1$. We call such a compatible system of perturbation data a \emph{regular} system.

If $(H, J)$ are $C_k$-invariant then this compatible system of perturbation data can be chosen to be $C_k$-equivariant. 

\end{proposition}
\begin{proof}
    This follows from an equivariant double induction argument similar to that explained in the proof of Proposition \ref{prop:embeddings-exist}. The outer loop is on downward-closed $G$-sub-posets $\mathcal{P} \subset \CC(H, J)$. The inductive hypothesis is that there is a chosen free $E$-parameterization with $\mathcal{E}_{\mathcal{P}} = \{V_{\mathcal{P}}(x)\}$ of $\CC(H, J)$ with associated $F$-parameterization $\mathcal{F}_{\mathcal{P}}= \{V_{\mathcal{P}}(x,y)\}$   together with a compatible regular system of perturbation data $\mathfrak{P}_{\mathcal{P}} = (\bar{\mathfrak{A}}|_{\mathcal{P}}, \mathcal{F}_{\mathcal{P}}|_{\mathcal{P}}, \{\lambda_{V_{\mathcal{P}}}\})$ for the the full flow subcategory on the objects of $\mathcal{P}$, such that $\mathcal{E}_{\mathcal{P}}$ can be built from the zero $E$-parameterization by stabilizing only at elements of $\mathcal{P}$.  To proceed with an iteration of the outer loop we choose $x \in \CC(H, J) \setminus \mathcal{P}$ as in Proposition \ref{prop:embeddings-exist} and introduce $\mathcal{P}_x$, $\mathcal{P}_{\bar{x}},$ and $\mathcal{P}'$ as in that setting. The inner loop is over upwards-closed $G_x$-sub-posets $\mathcal{Q} \subset \mathcal{P}_x$ containing $x$, as in Proposition \ref{prop:embeddings-exist}. To make one step of the inner loop we will find $y \in \mathcal{P}_x \setminus \mathcal{Q}$ as in Proposition \ref{prop:embeddings-exist} and set $\mathcal{Q}' = \mathcal{Q} \cup G_x y$, $\mathcal{Q}_y = \mathcal{Q} \cup \{y\}$ as in that setting. Initially we set $\mathcal{E}_{\mathcal{P}_x; \{x\}} = \mathcal{E}_{\mathcal{P}}$ and over each iteration of the inner loop we will produce $\mathcal{E}_{\mathcal{P}_x; \mathcal{Q}'}$ by stabilizing $\mathcal{E}_{\mathcal{P}_x; \mathcal{Q}}$ at $x$ by a $G_x$-representation $V'_{xy}$. The inductive condition for the inner loop is that the $E$-parameterization $\mathcal{E}_{\mathcal{P}_x, \mathcal{Q}}$ is part of a compatible regular system of perturbation data $\mathfrak{P}_{\mathcal{P}_x; \mathcal{Q}} = (\bar{\mathfrak{A}}|_{\mathcal{Q}}, \mathcal{E}_{\mathcal{P}_x, \mathcal{Q}}|_\mathcal{Q}, \{\lambda_{V_{\mathcal{P}_x, \mathcal{Q}}(x,y)}\})$ over the full flow subcategory on the objects of $\mathcal{Q}$ such that $\mathfrak{P}_{\mathcal{P}_x; \mathcal{Q}}|_{\mathcal{P} \cap \mathcal{Q}} = \mathfrak{P}_{\mathcal{P}}|_{\mathcal{P} \cap \mathcal{Q}}$.

    When running the inner loop at $y$, we do nothing if $y \notin \mathcal{P}_{\bar{x}}$. Otherwise, writing $V_{\mathcal{P}, \mathcal{Q}}(x,y)$ for the vector space associated by $\mathcal{E}_{\mathcal{P}; \mathcal{Q}}$ to $x>y$, and $V_{\mathcal{P}; \mathcal{Q}'}(x,y) = V_{\mathcal{P}, \mathcal{Q}}(x,y) \oplus V'_{xy}$ for the vector space that will be associated to $x>y$ by $\mathfrak{P}_{\mathcal{P}_x'; \mathcal{Q}}$ after the end of this iteration of the inner loop. We now note that the inductive condition for $\mathcal{Q}'$  fixes the values of $(\lambda_{V_{\mathcal{P}; \mathcal{Q}'}(x,y)})|_{\underline{V_{\mathcal{P}, \mathcal{Q}}(x,y)}}$ when restricted to the product of an open subset $U_{xy}$ of $\partial \cM(x,y)^Z$ with $M$. To see this one first notes that if $x > x_1 > y$ and $x > x_2 > y$ in $\CC(H, J)$, then either $\bar{A}(x_1) = \bar{A}(x_2)$, in which case the images of $\paramspace(x,x_1) \times \paramspace(x_1, y)$ and $\paramspace(x, x_2) \times \paramspace(x_2, y)$ agree in $\paramspace(x,y)$, or $\bar{A}(x_1) \neq \bar{A}(x_2)$, in which case the intersection of the images is the image of $\paramspace(x,x_1) \times \paramspace(x_1, x_2) \times \paramspace(x_2, y)$ in the corresponding codimension $2$ stratum. In the first case, we must (by Assumption \eqref{eq:assume-action-nondegeneracy}) have that $x_1 = gx_2$ for some $g \in C_k$, in which case one sees that the perturbation function defined on this set via its identification with $\paramspace(x,x_1) \times \paramspace(x_1, y)$ and via its identification with $\paramspace(x,x_2) \times \paramspace(x_2, y)$ agree by the fact that the previous choices satisfied the equivariance properties required. In the second case, the perturbation function defined over this set will be well defined just by associativity of the maps associated to the $F$-parameterization corresponding to $\mathcal{E}_{\mathcal{P}; \mathcal{Q}}$ and the compatibility conditions for a system of perturbation data.
    By a simple induction on the face poset of $\cM(x,y)^Z$, one extends the $G_{xy}$ equivariant bundle map $\lambda_{V_{\mathcal{P}; \mathcal{Q}'}(x,y)}|_{V_{\mathcal{P}, \mathcal{Q}}(x,y) \times U_{xy} \times M}$ to a $G_{xy}$-equivariant bundle map 
    $\lambda_{V_{\mathcal{P}, \mathcal{Q}}(x,y)}$ defined on $\underline{V_{\mathcal{P}, \mathcal{Q}}(x,y)}$. 
    
    By the gluing theorem, if we do not perform an additional stabilization of the current $E$-parameterization at this step (and so $V'_{xy} = 0$) then the operators $D^{\mathfrak{P}_{\mathcal{P}_x, \mathcal{Q}'}}_u$ would be surjective for all $u \in \CC(H, J)(x,y)$ except for those in a compact subset of $\pi^{-1}(M(x,y))$.  We thus stabilize sufficiently at $x$ such that $V'_{xy}$ is a sufficiently large $G_{xy}$-representation such that one can set the $G_{xy}$-equivariant map $\lambda_{V_{\mathcal{P}; \mathcal{Q}'}(x,y)})|_{V'_{xy}}$ to be a map which zero in an open neighborhood of $\partial \cM(x,y)^Z \times M$ and such that  $D^{\mathfrak{P}_{\mathcal{P}_x, \mathcal{Q}'}}_u$ is surjective for all $u \in \CC(H, J)(x,y)$ for every $u$ in this comapct subset. Of course, we have $\lambda_{V_{\mathcal{P}; \mathcal{Q}'}(x,y)})|_{V_{\mathcal{P}, \mathcal{Q}}(x,y)} = \lambda_{V_{\mathcal{P}, \mathcal{Q}}(x,y)}$.
    
    This choice of $V'_{xy}$ and $\lambda_{V_{\mathcal{P}_x; \mathcal{Q}'}(x,y)})$ fixes $\mathfrak{P}_{\mathcal{P}_x; \mathcal{Q}'}$ by the equivariance and compatibility conditions of a compatible regular system of perturbations together with the inductive hypothesis, once we also require that $\lambda_{V_{\mathcal{P}_x; \mathcal{Q}'}(x,y')})|_{\underline{V'_{xy}}} = 0$ for all $y' \notin G_x y$ (this is compatible with the previous choices, since $\lambda_{V_{\mathcal{P}; \mathcal{Q}'}(x,y)})|_{V'_{xy}}$ is zero near the boundary of $\cM(x,y)^Z \times M$).  Once we have reached $\mathcal{Q} = \mathcal{P}_x$; one notices that one has produced the $E$-parameterization underlying $\mathfrak{P}_{\mathcal{P}_x; \mathcal{P}_x}$ from that underlying $\mathfrak{P}_{\mathcal{P}_x; {x}}$ by stabilizing at $x$ by some $G_{x}$-representation $V'_x$; one thus sets the $E$-parameterization underlying $\mathfrak{P}_{\mathcal{P}'}$ to be that produced by stabilizing the $E$-parameterization underlying $\mathfrak{P}_{\mathcal{P}}$ at $x$ by $V'_x$ (see Lemma \ref{lemma:how-to-stabilize-equivariantly}) and subsequently fixing $\mathfrak{P}_{\mathcal{P'}}$ by the inductive condition, equivariance, and the condition that that $\mathfrak{P}_{\mathcal{P}'}|_{\mathcal{P}} = \mathfrak{P}|_{\mathcal{P}}$ and $\mathfrak{P}_{\mathcal{P}'}|_{\mathcal{P}_x} = \mathfrak{P}_{\mathcal{P}_x, \mathcal{P}_x}$.  Upon termination of the outer loop we set $\mathfrak{P} = \mathfrak{P}_{\CC(H, J)}$. 
\end{proof}

\subsection{Topological $\langle k \rangle$-manifolds and Microbundles}
In fact, we will find smoothings of $T(x_-, x_+)$ following the strategy of \cite{AMS, rezchikov2022integral, bai2022arnold}. We briefly review the relevant notions.

\begin{definition}
A microbundle of rank $n$ over a topological space $X$ is a diagram $\{ X \xrightarrow{s} E \xrightarrow{p} X \}$ where $E$ is a topological space, and $p$ and $s$ are continuous maps (called the \emph{projection} and the \emph{section}, respectively) such that $p \circ s = id_X$, and moreover for each $x \in X$, there are neighborhoods $x \in U \subset X$ and $s(x) \in V \subset E$, as well as a homeomorphism $h: U \times \R^n \to V$ with $p \circ h = pr_U$, and $h |_{U \times \{0\}} = s$. The space $E$ is the \emph{total space} of the microbundle. When the context is clear, the we will abbreviate the microbundle as $E$.

A morphism of microbundles is an equivalence class of continuous map between the total spaces which commute with the projections and the sections, respectively. Two such maps are equivalent if they agree in a neighborhood of the image of the section. 
\end{definition}

Given a vector bundle $\pi: V \to X$, the \emph{associated microbundle} is $\{X \xrightarrow{0} V \xrightarrow{\pi} X$, which we will also denote by $V$. Given a microbundle $E$, a vector bundle $V$ together with an isomorphism of the associated microbundle $V \to E$ is called a \emph{vector bundle lift} of $E$. 

Given a pair of microbundles $E_1, E_2$ over $X$, the direct sum $E_1 \oplus E_2$ is the microbundle 
\begin{equation}
\label{eq:microbundle-direct-sum}
    X \xrightarrow{s_1 \times s_2 \circ \Delta} (p_1 \times p_2)^{-1}(\Delta(X))  \xrightarrow{p_1 \times p_2} X. 
\end{equation}
where the total space is a subspace of $E_1 \times E_2$.

Given a map of spaces $f: X \to Y$ and a microbundle $E$ over $Y$, the pullback microbundle $f^*E$ is given by 
\begin{equation}
    X \xrightarrow{x \mapsto (s(f(x), x)} f^*E = \{(e, x) \in E \times X | p(e) = f(x) \} \xrightarrow{(e, x) \mapsto x} X.
\end{equation}

These notions admit all admit $G$-equivariant versions.

\begin{definition}
Let $X$ be a topological manifold. The \emph{tangent microbundle} $T_\mu X$ (also denoted $\tau X$) is the microbundle $X \xrightarrow{\Delta} X \times X \xrightarrow{p_1} X$, where $\Delta$ is the diagonal map and $p_1$ is projection to the first factor. If $G$ is a topological $G$-manifold $T_\mu X$ is a $G$-microbundle via the diagonal action on $X \times X$.
\end{definition}

The easiest way to define the tangent microbundle of a $\langle k \rangle$-manifold $M$ is to define it via the tangent microbundle of its \emph{double} \cite[Section 2.2.3]{rezchikov2022integral}
\begin{equation}
    \label{eq:double}
        \mathcal{D}(M) := \bigsqcup_{g \in (\Z/2)^k} g M/ \sim
    \end{equation}
where we glue the $2^k$ copies of $M$ in the natural way along their strata, such that $\mathcal{D}(M)$ acquires a $(\Z/2)^k$, comm, with $M$ the closure of a fundamental domain for this $(\Z/2)^k$ action. 

\begin{definition}
Given a topological $G-\langle k \rangle$-manifold $M$, its tangent microbundle $T_\mu M$ is the pullback of $T_\mu \mathcal{D}(M)$ under the canonical inclusion $M \to \mathcal{D}(M)$. A \emph{vector bundle lift} of $T_\mu M$ is defined to be a $G$-equivariant vector bundle $V \to M$ together with a $G$-equivariant isomorphism $\eta: V_\mu \to T_\mu M$, where $V_\mu$ is the microbundle associated to $M$, such that the $G$-equivariant map $\mathcal{D}\eta: \mathcal{D}(V_\mu) \to \mathcal{D}(T_\mu M)$ can be enhanced to a $G \times (\Z/2)^k$-equivariant map where the $G \times (\Z/2)^k$-equivariant structure on $\mathcal{D}(V_\mu) = (\mathcal{D}(V))_\mu$ arises by restriction from an enhancement of the $G$-equivariant structure of the vector bundle $\mathcal{D}V \to \mathcal{D}M$ to a $G \times (\Z/2)^n$-equivariant  vector bundle structure.
\end{definition}

The local geometry of $\langle k \rangle$-manifolds immediately shows:
\begin{lemma}
\label{lemma:decomposition-of-vector-bundle-lift}
    Let $V \to \tau M$ be a vector bundle lift of the tangent microbundle of a $G-\langle k \rangle$-manifold $M$. Then there are decompositions 
    \[ V|_{X(S)} = V(S) \oplus \underline{\R}^{[k]  \setminus S}\]
    into $G$-subbundles where the $G$-action on the second factor is trivial. 
\end{lemma}
\begin{proof}
  This follows from the fact that locally near $x$ in $X(S)$ the fixed loci under various subgroups of $(\Z/2)^k$ have certain dimensions (one loses a dimension for each $\Z/2$-factor associated to each element of $[k] \setminus S$) together with the invariance of dimension under homeomorphism.
\end{proof}

The notion of topological submersion is defined in \cite{AMS} and \cite{rezchikov2022integral}.

The \emph{vertical tangent microbundle} $T^{vt}_\mu(\pi)$ of $\pi$ is the microbundle 
\begin{equation}
    M \xrightarrow{\Delta} M \times_B M \xrightarrow{\pi \times \pi} M
\end{equation}
where $M \times_B M$ is the fiber product of $\pi$ with itself, $\Delta$ is the diagonal, and $\pi \times \pi$ is the canonical map out of the fiber product. 

We know review two propositions from \cite{rezchikov2022integral}, which are immediately inherited from analogous propositions in \cite{AMS} via the functorial properties of the doubling construction:
\begin{proposition}
\label{prop:submersion-microbundle-splitting}
    There is a map $P: T_\mu M \to T^{vt}_\mu M$ such that $P \oplus \tau: T_\mu M \to T^{vt}_\mu M \oplus \pi^* T_\mu B$ is an isomorphism. 
\end{proposition}

Below, we freely use the notion, introduced in \cite{AMS}, of a fiberwise smooth $C^1_{loc}$-structure on a topolgoical submersion. In particular, if a topological submersion $\pi: M \to B$ has a fiberwise-smooth $C^1_{loc}$ structure, then there is an associated vertical tangent bundle $T^{vt} M$, with fiber over $p$ equal to $T_p (\pi^{-1}(\pi(p))$.

\begin{proposition}
\label{prop:vector-bundle-lifts-of-c1-loc-submersions}
    If a topological submersion $\pi: M \to B$ has a fiberwise-smooth $C^1_{loc}$ structure and $B$ is a smooth manifold, then $T^{vt} M$ is a vector bundle lift of the microbundle  $T^{vt}_\mu M$. As such, by Proposition \ref{prop:submersion-microbundle-splitting}, the tangent microbundle of $\pi$ has a vector bundle lift, namely $T^{vt}M \times \pi^{*} TB$. 
\end{proposition}
\subsection{Gluing theorem}
Now, once we have chosen a compatible system of surjective perturbation data, the methods of \cite{pardon2016algebraic}, as in \cite{AMS}, establish that the spaces $T(x_-, x_+)$ are topological $\langle k \rangle$-manifolds. 

\begin{proposition}
\label{prop:pardon-gluing}
    Fix $(u, a,0) \in T(x_-, x_+)$, with $u = u_1 \cup \ldots \cup u_r$, and the associated sequence of Hamiltonian orbits $(x_-, x_1, \ldots, x_{r-1}, x_+)$. Then $a$ is the image of a point $(a_1, \ldots, a_r)$ in $\cM(x_-, x_1) \times \ldots \times \cM(x_{r-1}, x_+)$, and $(u, a, v)$ is the image of $(\{(u_i, a_i, v_i)\}) \in T(x_-, x_1) \times \ldots \times T(x_{r-1}, x_+)$
    \begin{itemize}
        \item A neighborhood $U_1$ of $(a_1, \ldots, a_r)$, and
        \item A neighborhood $U_2$ of 
            \[\pi^{-1}(a_1) \times \ldots \times \pi^{-1}(a_r) \subset T(x_-, x_1) \times \ldots \times T(x_{r-1}, x_+),\] and 
        \item an $\epsilon > 0$, such that 
        \item There is a map 
        \[[0, \epsilon)^{r-1} \times U_1 \times U^2 \to T(x_-, x_+)\]
        which is a homeomorphism onto its image, sends $0 \times a \times (\{(u_i, a_i, v_i)\})$ to $(u,a,v)$, and 
        \item Covers a stratum-preserving homeomorphism 
        \[ [0, \epsilon)^{r-1} \times U_1 \to \cM(x_-, x_+)\]
        upon applying projection to $[0, \epsilon)^{r-1} \times U_1$ onto the domain, and $\pi$ onto the codomain. (This latter homeomorphism is the product of a topological-$\langle k \rangle$-manifold chart on the $\overline{Conf}_{k_i}$ factors, and the standard isomorphism which concatenates the $\R^{\bar{A}(x_i) - \bar{A}(x_{i+1})}$ factors.)
    \end{itemize}
    Thus, $T(x_-, x_+)$ is a topological $\langle k \rangle$-manifold. Moreover, the map $\pi$ has a fiberwise $C^1_{loc}$-structure. 
\end{proposition}
\begin{proof}
    This follows immediately from Pardon's gluing theorems \cite{pardon2016algebraic}, as the corresponding result does in \cite{AMS}. The only difference in the analysis is that the hypersurface constraints of \cite{pardon2016algebraic} are not present here, but this only simplifies the required estimates. As in \cite[Corollary 6.29]{AMS}, the fiberwise $C^1_{loc}$ structure follows from \cite[Proposition C.11.1]{pardon2016algebraic}. 
\end{proof}

We note that the vertical tangent bundle is explicit:
\begin{lemma}
\label{lemma:vertical-tangent-bundle-formula}
    The vertical tangent bundle of $\pi: T(x_-, x_+) \to \cM(x_-, x_+)$ has fiber at $(u, a, v)$, with $u= u_1 \cup \ldots \cup u_r$, $u_j \in T(x_{j-1}, x_j)$ canonically isomorphic to 
    \[\oplus_{j=1}^r\ker D_{u_j}^{V(x_{j-1}, x_j)},\] where we set $x_0 = x_-$ and $x_r = x_+$. 
\end{lemma} 

\subsection{Smoothing theorem}
\begin{remark}
    The notion of a virtually smooth flow category works fine in the topological setting when there are no stabilizing vector bundles (so that no connections are needed.) 
\end{remark}

\begin{proposition}
Let $G$ be a finite group. 

    Let $X$ be a $G-\langle k \rangle$-manifold, and suppose that \emph{$X$ has a boundary-smooth structure}, namely that for each $X(S)$, $S \subset [k]$, $S\neq [k]$ is equipped with a smooth structure such that if $T \subset S$ then the restriction of the smooth structure on $X(S)$ to $X(T)$ agrees with the smooth structure on $X(T)$. Then there is a $G$-representation $R(\sigma)$ and a smooth structure on $X \times R(\sigma)$ which restricts to the induced smooth structures on each $X(S) \subset R(\sigma)$. 
\end{proposition}

\begin{definition}
    \label{def:collars}
The \emph{collaring} $Coll(X)$ of a $G-\langle k \rangle$-manifold  $X$ defined as 
    \begin{equation}
    \label{eq:collar}
        Coll(X) = \bigcup_{S \in 2^{[k]}}X(S) \times[0,1]^{[k] \setminus S}/\sim
    \end{equation}
with the equivalence relation is generated by the relations that if $(x, t) \in X(S) \times [0,1]^{[k]\setminus S}$ and $(y, u) \in M(T) \times [0,1]^{[k \setminus T]}$, then $(x, t) ~(y, u)$ if $S <T$ and $x=y$, and if $t|_{([n] \setminus S) \setminus ([n] \setminus T)} = 1$. This is naturally a $\langle k \rangle$-manifold, and in the setting where the boundary strata of $X$ have compatible smooth structures, the $\langle k\rangle$-manifold $Coll(X) \setminus X([k])$ is naturally a smooth $\langle k \rangle$-manifold (with a nonempty top open stratum).

A \emph{collar} on a topological $G-\langle k \rangle$-manifold $X$ is a stratum-preserving homeomorphism $Coll(X) \to X$. If $X$ has a boundary-smooth structure, a collar is boundary-smooth if it induces a diffeomorphism $Coll(X)(S) \to X(S)$ for all $S \neq [k]$. 
\end{definition}

\begin{figure}[h!]
    \centering
    \includegraphics[width=\textwidth]{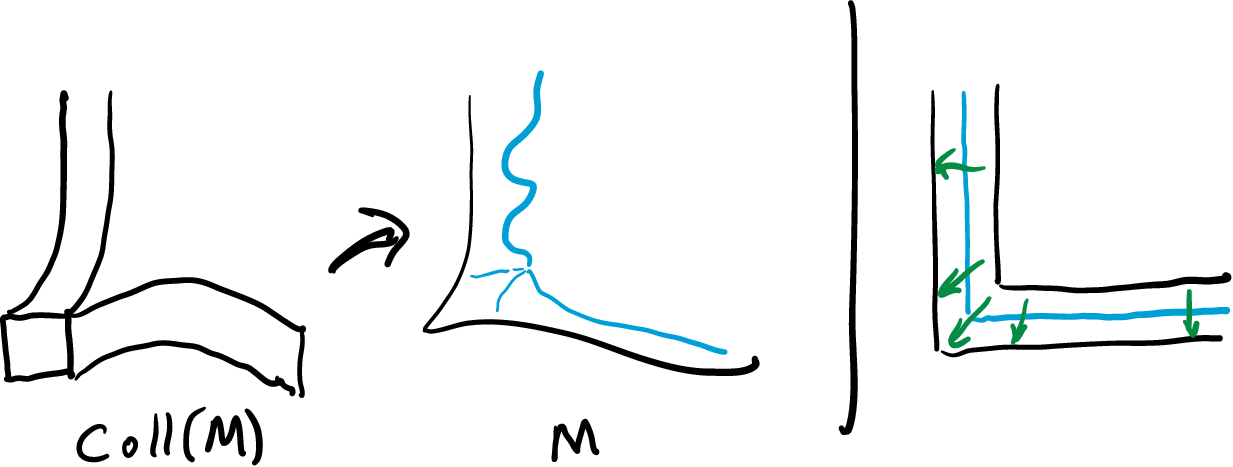}
    \caption{A collar of a manifold is a homeomorphism from $Coll(M)$ to $M$ (left). We can ask for this collaring to be compatible with various smooth structures, e.g. smooth structures along the boundary. A boundary-smooth structure on the boundary of induces a smooth structure on a neighborhood of the boundary via Proposition \ref{prop:corner-collars}. The collar $Coll(M) \setminus M$ admits a retraction onto its own boundary (right); if the manifold vector bundle lift of its tangent bundle which is smooth the boundary strata, then this lift can be deformed to one that is smooth on the image of the collar using this retraction (Lemma \ref{lemma:extend-smooth-vector-bundle-lift-to-one-which-agrees-on-a-collar}). One can then use the (equivariant) stable smoothing theorem of Lashof \cite{lashof2006stable} to extend the smooth structure on the image of the collar to the interior after a stabilization of the manifold (Proposition \ref{prop:relative-smoothing}).}
    \label{fig:corner-collars}
\end{figure}

\begin{definition}
    A \emph{nice} submanifold of a $G-\langle k \rangle$-manifold $X$ is a subset $Y \subset X$  preserved by the $G$-action such that for every $y \in Y$ there exists $\langle k \rangle$-manifold chart $\psi: \R^{[k] \setminus S}_+ \times \R^m \to \to X$ such that the image of this chart intersects $Y$ in $\psi(\R^{[k] \setminus S}_+ \times (\R^\ell  \to \oplus 0))$. 
\end{definition}

\begin{definition}
\label{def:elegant-metric}
An \emph{elegant} metric on a smooth $G -\langle k \rangle$-manifold $Coll(X)$ such that there exists a smooth collar $\psi: Coll(X) \to X$ with the metric on $\psi(X(S) \times [0,1]^{[k] \setminus S})$ agreeing with the product of the metric on the submanifold $\psi(X(S))$ with the pushforward under $\psi$  of the standard metric flat metric on the cube (of some, or equivalently any, side length size $\epsilon$). We call such a collar an \emph{elegant} collar.
\end{definition}

\begin{proposition}
\label{prop:corner-collars}
    A topological $G$-$\langle k \rangle$-manifold $X$ with a boundary-smooth structure has a boundary-smooth collar. If $X$ has a nice $G-\langle k \rangle$-submanifold $Y$ which is equipped with a smooth structure and an elegant metric, then the collar can be chosen so that its restriction to $Coll(Y)$ is an elegant collar onto $Y$. 
\end{proposition}
This proposition is a simple adaptation of the argument of \cite{Connelly1971}, and a proof will be given below. A boundary-smooth collar defines a smooth structure on an open neighborhood of the boundary of $X$, via the image of the smooth structure on $Coll(X) \setminus X$. 

\begin{lemma}
\label{lemma:extend-smooth-vector-bundle-lift-to-one-which-agrees-on-a-collar}
    Suppose that a topological $G-\langle k \rangle$-manifold $X$ with a boundary-smooth structure also has a vector bundle lift $V \to \tau X$ of its tangent microbundle $\tau X$ which is smooth on each boundary stratum in the sense that for each $S \subsetneq [k]$ and each point $x$ in the interior of $X(S)$, the composition of the exponential map of $T_xX(S) \to X(S)$ with the inverse of the microbundle lift lies in the subbundle $V(S)$ of Lemma \ref{lemma:decomposition-of-vector-bundle-lift}, and such that this composition is smooth. Then there is another vector bundle lift $(V, V_\mu \to \tau X)$ which is smooth on each boundary stratum, such that it is also smooth on the image of a boundary-smooth collar on $X$ (with the induced smooth structure on the interior of the collar). Moreover this vector bundle lift is homotopic to the original vector bundle lift, where the homotopy is constant over $\partial X$, and for which the decompositions of Lemma \ref{lemma:decomposition-of-vector-bundle-lift} agree for these two vector bundle lifts. If $Y$ is a nice $G-\langle k \rangle$ submanifold of $X$ which has a smooth structure and an elegant metric, with $V|_Y$ having a factor of the form $TY$ such that the vector bundle lift restricted to $TY$ is the canonical vector bundle lift given by the exponential map with respect to the metric, then this homotopy of the vector bundle lifts can be chosen to be smooth on $Y$, with all vector bundle lifts $TY \to \tau Y$ being given by exponential maps with respect to some elegant metric on $Y$. 
\end{lemma}
\begin{proof}(Sketch.)
Take a boundary-smooth collar  $\xi: Coll(X) \to X$, which restricts to an elegant collar on $Y$ if $Y$ exists.  We have that $\xi^*\tau X \simeq \tau Coll(X \setminus X$. Now, one can construct a retraction $r: Coll(X) \setminus X \to \partial(Coll(X) \setminus X)$ as indicated in Figure \ref{fig:corner-collars} (right) with the property that $r^* \tau Coll(X)\setminus X \simeq \tau Coll(X)\setminus X$. Thus, $(r^*\xi^*V, r^*\xi^*(V_\mu \to \tau X)$ defines a vector bundle lift of $r^*\xi^* \tau X \simeq \tau (Coll(X) \setminus X)$, which can moreover be seen to be smooth on is smooth on $Coll(X) \setminus X$. Since $\xi|_{Coll(X) \setminus X}$ is a homeomorphism onto its image, we can use this to push forward $r^*\xi^*V$ and this latter vector bundle lift to a vector bundle lift over $\xi(Coll(X) \setminus X)$ which is smooth over its image in the sense desired in the Lemma. But in fact $\xi(Coll(X) \setminus X) = \xi(X)$ is homeomorphic to $X$, so we can extend the above vector bundle lift by the vector bundle lift given by pushing forward the vector bundle lift given in the problem under $X \subset Coll(X) \setminus X \xrightarrow{\xi} X$. But now we can \emph{shrink the collar} (again, see Figure \ref{fig:corner-collars}, right), producing the desired homotopy of vector bundle lifts. In the case where there is a submanifold $Y$, it is straightforward to verify that this new vector bundle lift and homotopy are still smooth on $Y$.
\end{proof}

\begin{lemma}
\label{lemma:patch-together-homotopies-of-lifts}
    Given a $G-\langle k \rangle$-manifold $X$ with a vector bundle $V$ and a vector bundle lift $h: V_\mu \to \tau X$ of its tangent microbundle. Write $V(S)$ for the subbundle of $V|_{X(S)}$ defined as in Lemma \ref{lemma:decomposition-of-vector-bundle-lift}). Suppose we are given homotopies of vector bundle lifts $h^t_j: V([k] \setminus \{j\})_\mu \to \tau X_j$, $t \in [0,1]$, $j \in [k]$, which are compatible in the sense that $h^t_j|_{V([k] \setminus \{j, \ell\})} = h^t_\ell|_{V([k] \setminus \{j, \ell\})}$ for all $t$ and all $j, \ell \in [k]$, and such that $h^0_j = h|_{V([k] \setminus j)}$. Then there exists a homotopy of vector bundle lifts $h_t: V_\mu \to \tau X$ which restricts to all of the $h^t_j$ on $V([k] \setminus \{j\})_\mu$. Moreover, if $Y$ is a nice $G-\langle k \rangle$ submanifold of $X$ which has a smooth structure and an elegant metric, with $V|_Y$ having a factor of the form $TY$ such that all the $h^t_j|_{TY_j}$ are canonical lifts  associated by the exponential map to some elegant metrics on $Y_j$, then this homotopy of the vector bundle lifts can be chosen to be smooth on $Y$, with all vector bundle lifts $TY \to \tau Y$ being given by exponential maps with respect to some elegant metric on $Y$. 
\end{lemma}
\begin{proof}(Sketch.) Let us first assume that $Y$ is not present.
We first note that the fibers of the retraction $r$ of the previous lemma (see Figure \ref{fig:corner-collars}, right) are intervals; in particular, we have a homeomorphism 
\begin{equation}
    \label{eq:corner-collars-give-non-corner-collars}
    Coll(X)\setminus (X \setminus \partial X) \simeq \partial X \times [0,1].
\end{equation} Using this homeomorphism we can extend the decompositions of Lemma \eqref{lemma:extend-smooth-vector-bundle-lift-to-one-which-agrees-on-a-collar} to the interior of the image of the collar.  We will first construct a microbundle lift $h^1:V \to \tau X$ such that $h^1|_{V([k] \setminus j)|_{X_j}} = h^1_j$; the rest of the homotopy $h^t$ is constructed by `shrinking collar' as before. Now, there is is a subspace $\tilde{V} \subset V$ which fiber over $x \in X(S)$ given by $V(S) \times \partial \R^{[k] \setminus S}_+$ under the decomposition of Lemma \eqref{lemma:extend-smooth-vector-bundle-lift-to-one-which-agrees-on-a-collar}, and the  homotopies $h^t_j$ piece together to give maps $\tilde{h}^t: \tilde{V} \to \partial X$ which on each fiber are a homeomorphism onto their image.  We now find a subbundle $\underline{\R}_\tau \subset V|_{im Coll(X) \setminus X}$  with coordinate $\tau$ which is `perpendicular' to the subsets $0 \oplus \partial \R_+^{[k] \setminus S}$ in the decompositions of Lemma \ref{lemma:extend-smooth-vector-bundle-lift-to-one-which-agrees-on-a-collar} in the sense that there are splittings of the associated microbundle to $V$ as 
\begin{equation}
    \label{eq:microbundle-decompositions}
    V_\mu \simeq V(S) \oplus (Op(0) \subset \underline{\partial \R_+^{[k] \setminus S}}) \oplus \underline{\R}
\end{equation} 
over a neighboorhood of $r^{-1}(X(S) \times [0,1]^S)$ (where this uses the extension of the decompositions of Lemma \ref{lemma:extend-smooth-vector-bundle-lift-to-one-which-agrees-on-a-collar} to this neighborhood) which restrict to one another as $S$ varies. More geometrically if $v$ is the inwards pointing vector in $\underline{\R}$ the sets $V^S \oplus (\partial \R_+^{[k] \setminus S}]_+ \epsilon v$ foliate a neighborhood of the zero section in $v$. One then sets $h^1$ restricted to the restriction of $\tilde{V}'$ to the image of $X$ under the collar to be $h$; the idea of the construction this construction is to use this trivial subbundle to glue the homotopies $\tilde{h}^t$ to an extension of such a partially defined $h^1$ to a vector bundle lift which ends up restricting to $\tilde{h}^1$ on $\tilde{V}|_{\partial X}$ as desired. Specifically, we define the set $Q$ to be the parallelogram in $[0,1]_{\tilde{t}} \times [-1,1]_\tau$ given $0 \leq \tau + \tilde{t} \leq 1$. We then set, under the decompositions \eqref{eq:microbundle-decompositions}, writing $x \in Coll(X) \setminus (X \setminus \partial X))$ as $(x', \tilde{t})$ under the homeomorphism \eqref{eq:corner-collars-give-non-corner-collars} to be 
to be 
\[h^1(\tilde{t}, x'; v, \rho+\tau) = (\tilde{t}, \tilde{h}^{\tilde{t}}(x'; v, \rho)) \]
whenever $(\tilde{t}, \tau) \in Q$, whenever $(v, \rho, \tau)$ is the element of $V_\mu$ viewed under the decomposition \eqref{eq:microbundle-decompositions} (given by the foliation arising by translating the boundary of the positive orthant).  This defines $h^1$ everywhere except for $V|_{\partial X} \setminus \tilde{V}'$ and uniquely characterizes $h^1$. The rest of the deformation $h^t$ is constructed as in the previous lemma by `shrinking the collar' on as one proceeds in the homotopy $\tilde{h}^t$.  The modification when $Y$ is present follows as earlier.$\blacksquare$
\end{proof}

\begin{proposition}
\label{prop:relative-smoothing}
    Given a topological $G$-manifold $X$, with a lift of its tangent microbundle, together with smooth structure on an open set $U \subset X$ such that the composition of the exponential map $TU \to X \times X$ with the inverse to the lift of the tangent microbundle is smooth near the zero section, and a pre-compact open subset $K \subset X$, there exists a $G$-representation $V$ and a smooth structure on $V \times (U \cup K)$ which restricts to the smooth structure on $V \times U$ induced from the smooth structure on $U$. Moreover, if there is a locally flat submanifold $N$ of $K$ that has a smooth structure, the smooth structure on $V \times (U \cup K)$ can be chosen to restrict to the induced smooth structure on $V \times N$. 
\end{proposition}
\begin{proof}
    This follows from the equivariant stable smoothing  theorem \cite{lashof2006stable}, or rather, from its relative version, which is spelled out in \cite{bai2022arnold}. 
\end{proof}

Combining the above propositions, we immediately have the following result. 
\begin{proposition}
\label{prop:smoothing-a-k-manifold}
Given a topological $\langle k \rangle$-manifold $X$ with a boundary-smooth structure, and a precompact open subset $X' \subset X$, write $X''$ for the union of $X'$ with the image of a boundary-smooth collar on $X$.  Then there is a $G$-representation $V$ and a smooth structure on $X'' \times V$ which restricts, for each $S \subset [k]$, $S \neq [k]$, to the smooth structures on $X''(S) \times V$ induced from the smooth structures on $X(S)$.$\blacksquare$
\end{proposition}

We will repeatedly use Proposition \ref{prop:smoothing-a-k-manifold} to find smoothings of flow categories:
\begin{proposition}
\label{prop:smoothing-a-flow-category}
There exists an equivariant $E$-parameterization $\mathcal{E}$ with associated $F$-parameterization $\mathcal{F}$ such that the stabilization $\CC'(H, J)_{\mathcal{F}}$ of the (virtual) \emph{topological smoothing} $\CC'(H, J)$ of $\CC(H, J)$ has a lift to a (virtual) \emph{smoothing} of $\CC(H, J)$. 
\end{proposition}

\begin{proposition}
\label{prop:smooth-flow-category-floer-trajectories}
    There exists a regular system of perturbation data $\mathfrak{P}$ for $\CC(H, J)$ such that setting
    \[ \CC'(H, J)_{\mathfrak{P}}(x_-, x_+) = (T(x_-, x_+), V^1_T(x_-, x_+), \sigma(x_-, x_+))\]
    with these quantities as defined as in Equations \eqref{eq:thickenings-for-floer-trajectories}, \eqref{eq:obstruction-bundle-for-my-stabilization}, \eqref{eq:obstruction-section-floer-trajectories}, and compositions defined via \eqref{eq:thickening-composition-maps}, 
    we have that $\CC'(H, J)_{\mathfrak{P}}$ is a topologically smooth flow category that can be given a smooth structure. 
\end{proposition}
\begin{proof}
    This follows from Proposition \ref{prop:smoothing-a-flow-category} and Definition \ref{def:trivial-stabilization-of-perturbation-data}. 
\end{proof}

\begin{remark}
    If $X$ a $G-\langle k \rangle$-manifold with a boundary-smooth structure, then by Proposition \ref{prop:relative-smoothing},  $Coll(X) \times V$ has a smooth structure which restricts to the smooth structures on $Coll(X)(S) \times V$ for $S \subset [k]$, $S \neq [k]$. By repeatedly using the collaring construction, we could prove an analog of Proposition \ref{prop:smoothing-a-flow-category} without Proposition \ref{prop:corner-collars}, and proceed with the rest of the paper. However, this would involve collaring $T_{\CC'(H,J)}(x,y)$ many times, which makes the statement of the analog of Proposition \ref{prop:smoothing-a-flow-category} somewhat inelegant, and one would be forced to prove a number of invariance results associated to the effect of collaring a flow category on downstream constructions, e.g. on the Floer homotopy type, would would add a similar amount of complexity. The strategy of repeatedly collaring flow categories is taken in previous works \cite{abouzaid2021arnold, rezchikov2022integral, bai2022arnold} and we avoid this once and for all with the propositions of this section.
\end{remark}

\begin{lemma}
\label{lemma:extend-corner-chart-to-chart-smooth}
    Let $Y$ be a smooth $\langle k \rangle$-manifold, and let $y \in Y(S)$ lie in the interior of $Y(S)$ and $S \subset [k]$, $S \neq [k]$. A smooth $\langle k \rangle$-manifold boundary-chart at $y$ is a map 
    \[\psi: \R^{[k] \setminus S}_+ \setminus (\R^{[k] \setminus S}_+)^\circ \times \R^m \to Y \]
    where $(\R^{[k] \setminus S}_+)^\circ$ denotes the interior of $\R^{[k] \setminus S}_+$, satisfying the following properties:
    \begin{itemize}
        \item For $y \in V$, we have that $S(y) \geq S(\psi(0,0))$, and moreover
    	\item For $[k] \supset T \supset S$ with $T \neq S$, have that $\psi^{-1}(X(T)) \subset \R^{[k] \setminus S}_+(r_S(T)) \times \R^{m}$.
    \end{itemize}
    Then the restriction of $\psi$ to some neighborhood of $0$ is the restriction of a smooth $\langle k \rangle$-manifold chart at $y$ to $\R^{[k] \setminus S}_+ \setminus (\R^{[k] \setminus S}_+)^\circ \times \R^m$. 

    Moreover, if $Y$ is a $G -\langle k \rangle$-manifold where $G$ is finite, then if $\psi$ is $Stab(x)$-equivariant for a $Stab(x)$ on the domain which is the product of the trivial action on $\R^{[k] \setminus S}_+ \setminus (\R^{[k] \setminus S}_+)^\circ$ and a faithful linear action on $\R^m$, then the extending smooth-manifold chart can be taken to be $Stab(x)$-equivariant as well, with domain  neighborhood of zero in $\R^{[k] \setminus S}_+ \setminus (\R^{[k] \setminus S}_+)^\circ \times \R^m$ with the trivial $Stab(x)$ on the first factor and the same linear action on the second factor.
\end{lemma}
\begin{proof}
Write $\partial \R^{[k] \setminus S}_+ = \R^{[k] \setminus S}_+ \setminus (\R^{[k] \setminus S}_+)^\circ.$

    Let $\phi: \R^{[k] \setminus S}_+ \times \R^m \sup U \to Y$ be a smooth $\langle k \rangle$-manifold chart at $y$. After shrinking $U$, there is a map 
    \[ \psi^{-1} \circ \phi: \partial \R^{[k] \setminus S}_+ \times \R^m \to V \supset \partial \R^{[k] \setminus S}_+ \times \R^m\]
    for some neighborhood $V$ of $U$. The components of this map define functions 
    \[ f_i:  \partial \R^{[k] \setminus S}_+ \times \R^m \cap U  \to V \to \R_+, i \in \R^{[k] \setminus S}\]
    as well as a function 
    \[ g: \partial \R^{[k] \setminus S}_+ \times \R^m \to \R^m\]
    such that both $g$ and $f_j$ are smooth when restricted to each stratum $\partial \R^{[k] \setminus S}_+([k] \setminus j) \times \R^m \cap U$ for $j \in [k] \setminus S$, and that $f_j^{-1}(0) = \R^{[k] \setminus S}_+([k] \setminus j) \times \R^m $, with $df_j(x)$ nonzero when $x \in \R^{[k] \setminus S}_+(T) \times \R^m$ for $T \subsetneq [k] \setminus j$. Write $g_y: \R^m \to \R^m$ for the restriction of $g$ to $\{y\} \times \R^m$. Then we have that $dg_0 \neq 0$. Thus, by the inverse function theorem, after possibly shrinking $U$, if we extend $g$ to a smooth function $U \to \R^m$ and $f_i$ to smooth functions $U \to \R_{\geq 0}$ such that $f_j^{-1}(0) = (\R^{[k] \setminus S}_+([k] \setminus j) \times \R^m ) \cap U$, with $df_j|_{f_j^{-1}(0) \neq 0}$ then the inverse function theorem lets us extend a restriction of $\psi$ to neighborhood of $0$ to a smooth $\langle k \rangle$-chart around $0$. 

    The statement that $g$ can be equivariantly extended to a smooth function is standard, see for example \cite{rezchikov2022integral}. To extend each $f_j$ to smooth functions with the desired properties, we note that we can write $f_j = x_j \tilde{f}_j$ for $x_j$ the coordinate which is zero on $\R^{[k] \setminus S}_+([k] \setminus j) \times \R^m$ and $\tilde{f}_j \neq 0$ on a neighborhood of $\R^{[k] \setminus S}_+([k] \setminus j) \times \R^m$. Thus after extending $\tilde{f}_j$ to $U$ and defining $f_j = x_j \tilde{f}_j$ on $U$, after possibly shrinking $U$ we will have an extension of $f_j$ with the desired properties.
\end{proof}

\begin{lemma}
\label{lemma:boundary-smooth-corner-charts-exist}
    Given a topological $\langle k \rangle$-manifold $X$ with a boundary-smooth structure and a point $x \in X(S)$ with $x$ lying in the interior of $X(S)$, there exists a topological $\langle k \rangle$-manifold chart $\psi: \R^{[k] \setminus S}_+ \times \R^{m} \to X$, $m = \dim X - k + |S|$, such that $\psi$ is smooth restricted to any boundary hypersurface of its domain. Moreover, this chart can be taken to be $Stab(x)$ equivariant with the action on the domain being trivial on the first factor and linear and faithful on the second factor. 
\end{lemma}

\begin{proof}
    Repeated application of Lemma \ref{lemma:extend-corner-chart-to-chart-smooth} allows us to construct maps 
    \[ \psi: \partial \R^{[k] \setminus S}_+  \times \R^m \to X \]
    satisfying the following properties:
    \begin{itemize}
    \item $\psi((0,0)) = x$;
        \item For $y \in V$, we have that $S(y) \geq S(\psi((0,0))$, and moreover
    	\item For $[k] \supset T \supset S$ with $T \neq S$, have that $\psi^{-1}(X(T)) \subset \R^{[k] \setminus S}_+(r_S(T)) \times \R^{m}$ and finally 
     \item $\R^{[k] \setminus S}_P([k] \setminus j) \times \R^m \to X([k] \setminus j)$ is a smooth $\langle k \rangle$-manifold chart at $x$, for each $j \in [k] \setminus S$. 
    \end{itemize}
    
    Given the above data, the lemma follows from an application of the tubular neighborhood theorem for topological manifolds with boundary: using the fact that topological $\langle k \rangle$-manifolds are topological manifolds with boundary, one can produce a global topological collar on $X$, and a use a choice of topological collar on $\R^{[k] \setminus S}_+ \times \R^{m}$ thought of as a topological manifold with boundary to define an extension of $\psi$ to an open neighborhood of $(0,0) \in   \R^{[k] \setminus S}_+  \times \R^m$ which is a topological $\langle k \rangle$-manifold chart at $x$.
\end{proof}

\begin{proposition}
\label{prop:one-collar-coordinate}
    Let $X$ be a $G-\langle k \rangle$-manifold with a boundary-smooth structure. Choose a boundary-hypersurface $X_j \subset X$, $j \in [k]$. There exists a map
    \[ \phi_j: X_j \times [0,1]\cup_{X_j} X  \to X \]
    (where we identify $X_j \times \{1\}$ on the left with $X_j \subset X$ on the right) which is a homeomorphism onto its image, such that 
    \begin{itemize}
        \item $\phi_j|_{X_j(S) \times \{0\}}$ is a diffeomorphism onto $X(S)$ for $S \subset [k] \setminus j$;
        \item For $k \neq j$, 
        \[\phi_j|_k: X_j \cap X_k \times[0,1] \cup_{X_j \cap X_k} X_k \to X_k\]
        is a diffeomorphism onto its (open) image, and 
        \item $\phi_j$ is $G$-equivariant.
    \end{itemize}
    Moreover, if $X$ has a nice $G-\langle k \rangle$-submanifold with a smooth structure an an elegant metric, then $\phi$ can be chosen to be elegant along $Y \cap X_j$, i.e $\phi_j(Y_j \times [0,1] \cup_{Y_j}Y) \subset Y$, and the metric on the image of $\phi_j:Y_j \times [0,1]$ is the product of the induced metric on $\phi_j(Y_j)$ with the metric on an interval.
\end{proposition}

\begin{proof}[Proof of Proposition \ref{prop:corner-collars}]
This follows by repeatedly applying Proposition \ref{prop:one-collar-coordinate}. For $r = 1, 2, \ldots, k$, define 
\[ \mathcal{S}_r \subset 2^{[k]}, S_r = 2^{[k]\setminus 1} \cup 2^{[k] \setminus 2} \cup \ldots \cup 2^{[k] \setminus r}.\]
Define
\[ Coll_r(X) \subset Coll(X)\]
\[ Coll_k(X) = \cup_{T \in \mathcal{S}_r} X(T) \times[0,1]^{[k] \setminus T}/\sim.\]
Note that $Coll_1(X) =  Coll(X_1) \times [0,1].$ Let us imagine first that $Y$ is not present. Then,
applying Proposition \ref{prop:one-collar-coordinate} to $j=1$ and composing $\phi_1$ with the collar $Coll(X_1) \to X_1$ defines a map $f_1: Coll_1(X) \to X$ which is stratum preserving and is a diffeomorphism onto its image when restricted to each boundary stratum. This is the first step of the induction. For $1 < r \leq k$, assume that we have defined a map $f_{r-1}: Coll_{r-1}(X) \to X$ which is a diffeomorphism onto its image when restricted to each boundary stratum. Invoking Proposition \ref{prop:one-collar-coordinate} for $j=r$ we get a map $\phi_r$ as in the proposition. This map in turn defines $\bar{\phi}_r$ via the composition $X \to X \times [0,1] \cup X \xrightarrow{\phi}_r X$, where the first arrow is the inclusion into the second term of the union. Now note that $Coll_r(X) = Coll_{r-1}(X) \cup Coll_{r-1}(X)_r \times [0,1]$. We define $f_r: Coll_r(X) \to X$ to be $\bar{\phi}_r \circ f_{r-1}$ on $Coll_{r-1}(X)$, and $f_r|_{Coll_{r-1}(X)_r \times [0,1]}$ to be $\phi_r(f_{r-1} \times id)$. The variant where $Y$ is present is straightforward via the corresponding variant of Proposition \ref{prop:one-collar-coordinate} (choosing $\phi_j$ such that $\phi_j(X_j \times [0,1])$ lies  in a sufficiently small neighborhood of $X_j$).  
\end{proof}

\begin{proof}[Proof of Proposition \ref{prop:one-collar-coordinate}]
This follows from the argument of \cite{Connelly1971}. We use the notation of that paper in this proof. We note first that $X_r$ is locally collared in $X$ for any $1 \leq r \leq k$. Moreover, our spaces are paracompact Hausdorff topological manifolds with boundary, so every open cover admits a star-finite open refinement \cite{gauld1974topological}. We choose the open sets $U_1, \ldots, U_n$ as in the proof \cite{Connelly1971} so that $\overline{V} \subset U_1$  and such that $h_2, \ldots, h_n$ arise from local charts produced as in Lemma \ref{lemma:boundary-smooth-corner-charts-exist}. Moreover we require that $G$ permutes the neighborhoods $U_2, \ldots, U_n$ and acts freely on this set, and moreover such that if $g U_a = U_b$ then $h_b \circ g = h_a$. We then find a star-finite refinement $V_1, \ldots, V_n$ such that $\overline{V}_i \subset U_i$ and such that $\overline{V} \subset V_1$. The maps $f_1$ and $g_1$ as in \cite{Connelly1971} is defined using $\tilde{\phi}_j$. For each inductive step, we require that $\phi_i$ is $Stab(U_i)$-equivariant and such that $\phi_i|_{H_ig_{i-1}(X_k)}$ is a diffeomorphism for each $k \neq j$; this can be achieved by using the smooth equivariant Tietze theorem to define $\lambda_i$ on $\overline{U}_i \cap \partial X$ first and then extending using the equivariant Tietze extension theorem. After defining $\phi_i$ one defines $\phi_b$ for all $b$ such that $U_b = gU_i$ by $G$-equivariance, and then proceeds onwards with the induction choosing the first $U_\ell$ such that $\phi_\ell$ has not yet been defined. The variant with $Y$ present just involves using the local collarings to be elegant near $Y$.
\end{proof}

\begin{proof}[Proof of Proposition \ref{prop:smoothing-a-flow-category}]
This follows from the equivariant double induction as in Propositions \ref{prop:embeddings-exist} or \ref{prop:compatible-perturbation-data-exist}. The outer loop of the induction is on downwards closed $G$-subposets $\mathcal{P} \subset \CC(H, J)$. The inductive hypothesis is that there is an $E$-parameterization $\mathcal{E}_\mathcal{P}$ with associated $F$-parameterization $\mathcal{F}_\mathcal{P}$ such that we have chosen a smooth lift of the full subcategory $\CC'(H, J)_{\mathcal{F}_\mathcal{P}; \mathcal{P}}$ of $\CC'(H, J)_{\mathcal{F}_\mathcal{P}}$ on the objects of $\mathcal{P}$. To expand $\mathcal{P}$ we to $\mathcal{P}_x$ choose a minimal $x \in \CC(H, J) \setminus \mathcal{P}$ as usual, and introduce $\mathcal{P}_x, \mathcal{P}_{\bar{x}},$ and $\mathcal{P}'$ as usual. Clearly $\CC'(H, J)_{\mathcal{F}_\mathcal{P}; \mathcal{P}_x}$ is a $G_x$-equivariant topologically smooth flow category. Over one run of the inner loop, we will define an $E$-parameterization $\mathcal{E}_x$ which is the stabilization along the $G$-orbit of $x$ of the zero $E$-parameterization  (Lemma \ref{lemma:how-to-stabilize-equivariantly}). Write $\mathcal{F}_x$ for the associated $F$-parameterization; note that for $w, w' \in \mathcal{P}'$, we have that the vector space associated by $\mathcal{F}_x$ to the pair $(w, w')$ is nonzero exactly if $w \in Gx$, and in the latter case is exactly the vector space associated by $\mathcal{E}_x$ to $w$.  We then set $\mathcal{E}_{\mathcal{P}''} = \mathcal{E} \oplus \mathcal{E}_x$, $\mathcal{F}_{\mathcal{P}''} = \mathcal{F} \oplus \mathcal{F}_x$. One run of the inner loop will also define a smooth lift of $\CC'(H, J)_{\mathcal{F}_{\mathcal{P}'}; \mathcal{P}'}$ which agrees with the previously chosen smooth lift of $\CC'(H, J)_{\mathcal{F}_{\mathcal{P}''}; \mathcal{P}} = \CC'(H, J)_{\mathcal{F}_{\mathcal{P}}; \mathcal{P}}$. This smooth lift of $\CC'(H, J)_{\mathcal{F}_{\mathcal{P}'}}$ is defined by transport of structure (i.e. by the $G/G_x$-action) via a choice of smooth lift of the $G_x$-equivariant topologically smooth flow category $\CC'(H, J)_{\mathcal{F}_{\mathcal{P}'}; \mathcal{P}'_x}$, which in turn only depends on the $G_x$-representation $\tilde{V}(x)$ that $\mathcal{E}_x$ assigns to $x$, which itself defines $\mathcal{E}_x$. 

To choose this latter smooth lift and to define the representation $\tilde{V}(x)$, we induct downwards on upwards-closed $G_x$-subposets $\mathcal{Q} \subset \mathcal{P}'$. At each step we will enlarge $\tilde{V}(x)$; for shorthand, write $\mathcal{F}_{\mathcal{P}'; \tilde{V}(x)}$ for $\mathcal{F}_\mathcal{P'}$ defined via this choice of $\tilde{V}(x)$. The inductive assumption is that we have chosen a smooth lift of $\CC'(H, J)_{\mathcal{F}_{\mathcal{P}'; \mathcal{Q}}}$  which agrees with the smooth lifts chosen in previous iterations of the outer loop. We will write $\mathcal{F}_{\mathcal{P}'; \mathcal{Q}}$ for the corresponding $F$-parameterization with obtained from $\mathcal{F}_{\mathcal{P}}$  by stabilizing by $\tilde{V}(x)$ (with its value at the end of the step corresponding to $\mathcal{Q}$) at $x$. We will enlarge $\mathcal{Q}$ to $\mathcal{Q}' = \mathcal{Q} \cup \{G_xy\}$ with $y$ maximal in $\mathcal{P}_x \setminus \mathcal{Q}$ as usual. At the end of the step let us denote the new value of $\tilde{V}(x)$ by $\tilde{V}(x)'$; we will simply set $\tilde{V}(x)' = \tilde{V}(x) \oplus R$ where $R$ is some sufficiently large $G_x$-representation. To choose $R$, we note that the thickening associated to $\CC'(H, J)_{\mathcal{F}_{\mathcal{P}'; \tilde{V}(x)}}(x, y)$  is a topological $\langle k \rangle$-manifold with a boundary-smooth structure. 

We now wish to apply Proposition \ref{prop:smoothing-a-k-manifold} to give a stabilization of $\CC'(H, J)_{\mathcal{F}_{\mathcal{P}'; \mathcal{Q}}}(x, y)$ a smooth structure. We have that this thickening has a microbundle lift $h$ coming from Lemma \ref{lemma:vertical-tangent-bundle-formula}; the previous steps of this induction will define homotopies $h^j_t$ of the restrictions of this microbundle lift to the boundary strata of the thickening, satisfying the conditions of Lemma 
\ref{lemma:patch-together-homotopies-of-lifts}, and such that $h^j_1$ is smooth on the $j$-th boundary face of the thickening with respect to its smooth structure. Using Lemma \ref{lemma:patch-together-homotopies-of-lifts} we extend this homotopy to a homotopy $h_t$ of microbundle lifts with $h_0 = 0$, such that the new microbundle lift $h^1$ is boundary-smooth in that it satisfies the conditions of Lemma \ref{lemma:extend-smooth-vector-bundle-lift-to-one-which-agrees-on-a-collar}. We now use Proposition \ref{prop:smoothing-a-k-manifold}
  with $X'$ chosen to be a neighborhood of the zero locus of the Kuranishi section of $\CC'(H, J)_{\mathcal{F}_{\mathcal{P}'}}(x, y)$; this implicitly defined further homotopy of the vector bundle lift $h^1$ to a vector bundle lift which is smooth with respect to this smooth structure; moreover, this homotopy is constant on the boundary strata of the thickening. We use these homotopies in the later induction steps when utilizing Lemma \ref{lemma:patch-together-homotopies-of-lifts}   We now set $R$ to be the representation produced by this application Proposition \ref{prop:smoothing-a-k-manifold} applied to this thickening and this choice of $X'$.

  The result of Proposition \ref{prop:smoothing-a-k-manifold} then defines a $G_{xy}$-equivariant smoothing of the thickening of $\CC'(H, J)_{\mathcal{F}_{\mathcal{P}'; \mathcal{Q}'}}(x, y)$, which induces corresponding smoothings of $\CC'(H, J)_{\mathcal{F}_{\mathcal{P}'; \mathcal{Q}'}}(x, y')$ for $y' \in Gy$ via transport of structure (using the $G_x$ action). But this is exactly the data needed to fulfill the inductive hypothesis for $\mathcal{Q}' \supset \mathcal{Q}$. 

We have explained the inductive steps for both the inner and the outer inductions; each induction terminates by the finiteness of the number of objects of $\CC(H, J)$, and the base cases can be taken to be vacuous. Note that we can use Definition \ref{def:trivial-stabilization-of-perturbation-data} to interpret the stabilization of the virtual smoothing required as arising from a stabilization of the perturbation data.
\end{proof}

\subsection{Stabilization}
Given a surjective Fredholm operator $F: B_1 \to B_2$, between a pair of Banach spaces, and another map $A: W \to B_2$ where $W$ is a finite-dimensional inner product space , there is an associated operator
$F^W: V_1 \oplus W \to V_2$ by $F^W(v,w) = Dv - A w$. Then $\ker D^A$ is canonically isomorphic to $\ker D \oplus W$ as follows: there is an exact sequence 
\[ 0 \to \ker D \xrightarrow{i} \ker D^A \xrightarrow{p} W \to 0\]
where $i$ is induced by the inclusion $V_1 = V_1 \oplus 0 \subset V_1 \oplus W$, and $p$ is induced by the projection $V_1 \oplus W \to W$. In the case when all these spaces are inner product spaces with $i$ an isometry, the map $p$ has a \emph{canonical} section (by sending $w$ to the unique pair $(v, w)$ such that $Dv + A w = 0$ and such that $v$ is perpendicular to $\ker D$). This canonical section induces a canonical isometry $W \to (\ker D)^\perp \subset \ker D^A$ by using the polar decomposition. 

\begin{proposition}
\label{prop:stabilizing-a-regular-flow-category}
    Let 
    $(\mathfrak{P}_1 = (\bar{\mathfrak{A}}, \mathcal{V}^1, \{\lambda_{V^1(x,y)}\})$ be a regular equivariant system of perturbation data, and let $\mathfrak{P}_2 = (\bar{\mathfrak{A}}, \mathcal{V}^2, \{\lambda_{V^1_2(x,y)}\})$ be an equivariant system of perturbation data such that 
    \[ \mathcal{V}_2 = \mathcal{V}_1 \oplus \mathcal{F}\]
    where $\mathcal{F}$ is some $F'$-parameterization of $\CC(H, J)$, and such that 
    \[ (\lambda_{V^1_2(x,y)})|_{V^1(x,y)} = \lambda_{V^1(x,y)} \text{ for all } x > y \text{ in } \CC(H, J). \]
    Denote the corresponding topological flow categories by $\CC'_1(H, J)$ and $\CC'_2(H, J)$. Then a restriction of $\CC'_2(H, J)$ is the stabilization by $\mathcal{F}$ of $\CC'_1(H, J)$.
\end{proposition}

\begin{definition}
\label{def:stabilization-of-perturbation-data}
    With the notation of Proposition \ref{prop:stabilizing-a-regular-flow-category}, we say that $\mathfrak{P}_2$ is the stabilization of $\mathfrak{P}_1$ by the data $\mathcal{F}_{\lambda} = (\{V^1(x,y)^\perp\}, \{\tau_{xyz}\}), \{\lambda^1_2(x,y)|_{V^1(x,y)^\perp}\})$, which we may refer to as an \emph{$F'$-parameterization equipped with perturbation data.}
\end{definition}
\begin{definition}
\label{def:trivial-stabilization-of-perturbation-data}
    With the notation of Definition \ref{def:stabilization-of-perturbation-data}, we say that $\mathcal{F}_\lambda$ is  equipped with \emph{trivial} perturbation data when all the maps $\lambda^1_2(x,y)|_{V^1(x,y)^\perp}$ are identically zero. As such, a stabilization of $\mathcal{P}_1$ by $(\{V^1(x,y)^\perp\}, \{\tau_{xyz}\}), \{0\})$ is entirely determined by the underlying $F'$-parameterization $\mathcal{F} = (\{V^1(x,y)^\perp\}, \{\tau_{xyz}\}))$, and we will call $\mathcal{P}_2$ the stabilization of $\mathcal{P}_1$ by $\mathcal{F}$. One sees immediately from the definitions that in this setting, there is no restriction needed to get an isomorphism of $\CC'_2(H, J)$ with the stabilization of  $\CC'_1(H, J)$ by $\mathcal{F}$, and the proof of Proposition \ref{prop:stabilizing-a-regular-flow-category} is essentially trivial. 
\end{definition}
\begin{proof}[Proof of Proposition \ref{prop:stabilizing-a-regular-flow-category}]
  Write $\CC'_i(H, J)(x,y) = (T_i(x,y), \underline{V^1_i(x,y)}, \sigma_i(x,y))$ for $i=1,2$. There is an inclusion $T_1(x,y) \subset T_2(x,y)$ for every $x>y$ objects of $\CC(H, J)$. Moreover, there are fiberwise smooth structures for the topological submersions $\pi_i: T_i(x,y) \to \paramspace(x,y)$, and the $C^1_{loc}$ structure on these submersions shows that the union of the fiberwise normal bundles defines a bundle
\[NT_1(x,y) = \cup_{a \in \paramspace(x,y)}N(T_1(x,y))_a \]
over $T_1(x,y)$.

There is a continuous map $T_2(x_-, x_+) \to  V^1(x,y)^\perp$ given by the composition 
\[(u, a, (v_1, v_2)) \mapsto v_2. \]
One sees immediately that this map is smooth on each fiber, and its fiberwise differential defines a map
\[T_vT_2(x_-, x_+) \to V^1(x,y)^\perp\]
which is fiberwise linear.
The $C^1_{loc}$ charts show that this map is continuous. Moreover, by Lemma \ref{lemma:vertical-tangent-bundle-formula} and the surjectivity of $\mathcal{P}_1$ we see that this map is surjective when restricted to each fiber, of $T_vT_2(x_-, x_+)$, and has kernel $T_vT_1(x_-, x_+)$. As such, one has a bundle isomorphism $NT_1(x_-, x_+) \simeq \underline{V^1(x,y)}^\perp$; write $r_{x_-, x_+}$ for the inverse of this bundle isomorphism. Over the image of $T_1(x_-, z) \times T_1(z,x_+)$ in $T_1(x_-, x_+)$, there are canonical decompositions 
\begin{equation}
    \label{eq:vertical-tangent-bundle-decompositions}
    T_vT_i(x_-, z)\boxplus T_vT_i(z,x_+) \simeq T_vT_i(x_-, x_+) \text{ for } i =1,2,
\end{equation}
and the compatibility condition for perturbation data implies that $r_{x_-, x_+}\underline{\tau'_{x_-,z,x_+}}= r_{x_-, z} \oplus r_{z, x_+}$. Now, the bundles $T_vT_i(x_-, x_+)$ all carry inner products induced by sum of the $L^2$ inner product on rapidly decaying sections of $u^*TM$ over $Z$ and the $L^2$ inner products on $V^1(x_-, x_+)$ and $V^1_2(x_-, x_+)$. Thus the decompositions \eqref{eq:vertical-tangent-bundle-decompositions} are actually orthogonal. In particular, defining $r'_{x_-, x_+}$ to be the bundle isometry $NT_1(x_-, x_+) \simeq \underline{V^1(x,y)}^\perp$ produced by applying polar decomposition to $r_{x_-, x_+}$, we have that $r'_{x_-, x_+}\underline{\tau'_{x_-,z,x_+}}= r'_{x_-, z} \oplus r'_{z, x_+}$. 

The above data inner products on the $T_vT_i(x_-, x_+)$ defines Riemannian metrics on the (smooth) fibers of $T_-(x_-, x_+)$ over $\paramspace(x_-, x_+)$ (here the smooth structure comes from the description of the the thickenings as products of zero sets of Fredholm sections over Banach manifolds, \emph{not} via the previously described smooth structures!) such that the restriction of the Riemannian metrics on the fibers of $T_2(x_-, x_+)$ agree with those on the fibers of $T_1(x_-, x_+)$, and such that the induced Riemannian metrics on the normal bundles $NT_1(x_-, x_+)$ are those given by restricting the metrics on $T_vT_2(x_-, x_+)|_{T_1(x_-, x_+)}$ to the perpendicular to $T_vT_1(x_-, x_+)$. Moreover, the charts associated to the gluing construction show that these Riemannian metrics vary continuously from fiber to fiber. 

Putting all this together, we see that by applying the tubular neighborhood theorem to the inclusion $T_1(x_-, x_+) \subset T_2(x_-, x_+)$ fiberwise that the restriction of the fiberwise exponential map to an open neighborhood $U(x_-, x_+)$ of the zero section of $NT_1(x_-, x_+)$ defines a homeomorphism onto its image 
\[ r^{exp}_{x_-, x_+}: U(x_-, x_+) \to T_2(x_-, x_+)\]
which is a fiberwise local diffeomorphism. The compatibility conditions for the metrics imply that \[\iota_{x_-, z, x_+} (r^{exp}_{x_-,z}(a_1), r^{exp}_{z, x_+}(a_2)) = r^{exp}_{x_-, x_+}(\underline{\tau'_{xzy}}(a_1, a_2)) \]
for all $x_- > z > x_+$ in $\CC(H, J)$ and all $a_1 \in T(x_-, z)$ and $a_2 \in T(z, x_+)$ for which the above maps are defined. Thus, after inductively restricting the neighborhoods $U(x_-, x_+)$, these maps define an isomorphism from a restriction of the stabilization of $\CC'_1(H, J)$ to a restriction of $\CC'_2(H, J)$. 
\end{proof}

\subsection{Continuation map}
\label{sec:continuation-maps-geometry}

Let us suppose that we have chosen all the data of the previous paragraphs for a pair of convex nondegenerate Floer data $(H_-, J_-)$ and $(H_+, J_+)$. Given convex continuation data $(H_s, J_s)$, we can define continuation map moduli spaces 
\[ \mathcal{M}(H_s, J_s)(x_-, x_+) =\{ u: Z \to M : \partial_s u + J_s \partial_t u = \Grad H_s, \lim_{s \to \pm \infty} u_s = x_\pm\}\]
for $x_\pm \in Fix(H_\pm)$. The
the Gromov-Floer bordifications $\overline{\mathcal{M}}(H_s, J_s)(x_-, x_+)$ are compact. If $H_s$ and $J_s$ are $C_k$-equivariant then there is a loop rotation action of $C_k$ on the disjoint union of the above spaces over pairs $(x_-, x_+) \in Fix(H_-) \times Fix(H_+)$. We will explain how to choose `perturbation data' such that these spaces are the underlying spaces of a virtually smooth equivariant flow category which contains stabilizations of $\CC'(H_\pm, J_\pm)$ as flow subcategories. 

\paragraph{Integralization data for continuation maps. }
The first bit of data we chose to define $\CC'(H_\pm, J_\pm)$ were the auxiliary actions $\mathfrak{A}_\pm$ and the corresponding integralization data $\bar{\mathfrak{A}_\pm}$.

In the construction of Section \ref{sec:thickenings-for-hamiltonian-floer-homology}, we relied crucially on the property that the action functional is strictly decreasing along Floer trajectories in order to canonically place markings along the domains of Floer trajectories, thus embedding spaces of Floer trajectories into the thickenings constructed in that section. 
Because the action functional does not strictly decrease along continuation map trajectories, we will have to include, as part of the data of an integralization datum for continuation maps, an auxiliary function $c(s)$ which which will be used to artificially shift the action functional to a new function which strictly decreases along continuation trajectories. The following lemma bounds how rapidly the action functional can increase along continuation trajectories:

\begin{lemma}
\label{lemma:a-priori-bound-on-continuation-error}
    Let $(H_s, J_s)$ be a convex continuation datum. If $u \in \cM(H_s, J_s)(x_-, x_+)$ then, as functions of $s \in \R$, we have 

    \[ \frac{d}{ds} \mathcal{A}_{H_s}(u_s) < \max_{(x, t) \in M \times S^1} |\partial_s H_s(x,t)|  \]
    and in particular 
    \[ \mathcal{A}_{H_+}(x_+) - \mathcal{A}_{H_-}(x_-) < \int_{-\infty}^\infty \max_{(x, t) \in M \times S^1} |\partial_s H_s(x,t)|  \]
\end{lemma}
\begin{proof}
    This is essentially standard. When $H_s$ is independent of $s$, we have that $(d/ds)\mathcal{A}_{H_s}(u_s) < 0$ since $u_s$ is in this case the gradient flow of this action functional. The only new term that arises in the computation of $(d/ds)\mathcal{A}_{H_s}(u_s)$ when $H_s$ is not independent of $t$ is 
    \[ \int_0^1 - (\partial_s H_s)(u_s(t)) dt \leq \max_{(x, t)\in M \times S^1} |\partial_s H_s(x)|.\]
    Since $\partial_s u_s|_{s = s_0}$ agrees with the vector field along $u_{s_0}$ corresponding to the gradient flow of $\mathcal{A}_{H_{s_0}}$, at $s=s_0$, the remanining terms in $((d/ds)\mathcal{A}_{H_s}(u_s))|_{s = s_0}$ sum to a strictly negative quantity; thus the bound above proves the claims of the lemma. 
\end{proof}

\begin{definition}
\label{def:integralization-datum-continuation-map}
    An \emph{integralization datum} for $\CC(H_s, J_s)$ is a tuple $(\mathcal{A}^{cont}, \mathfrak{A}_{cont}, \{A_c\})$ such that
    \begin{enumerate}[(A)]
        \item There exist functions $c: \R \to \R$ and constants $C_\pm \in \R$ and $D > 0$ such that 
        \begin{itemize}
            \item $\mathcal{A}^{cont} = D \mathcal{A}_{H_s} + c(s)$ as functions $\R \times LX \to \R$;
            \item $c(s) = C_\pm$ for $\pm s>> 0 $;
    \item $\frac{d}{ds} c(s) < D\max_{(x,t) \in M \times S^1} |\partial_s H_s(x)|$. 
    \end{itemize}
    \item $\mathfrak{A}_{cont} \subset \R$ can be written as
    \begin{equation}
        \label{eq:decomposition-of-marked-actions-for-continuation-map}
        \mathfrak{A}_{cont} =(\mathfrak{A}_- + C_-) \cup (\mathfrak{A}_+ + C_+), 
    \end{equation}
    with
    $\max_{A \in \mathfrak{A}_+ + C_+} A < \min_{A \in \mathfrak{A}_- + C_-} A;$.
    \item In addition, we have that $(D \mathcal{A}_{H_\pm} + C_\pm, \mathfrak{A}_\pm + C_\pm,)$ is an integralization datum for $\CC(H_\pm, J_\pm)$. 
    \item Finally, for every $x_\pm \in Fix(H_\pm)$, we have that the set $\mathfrak{A}(x_-, x_+)$ contains at least two elements $A_c^-, A_c^+$ such that $A_c^+ > A_c > A^c_-$. 
    \end{enumerate}
\end{definition}

A crucial observation is that Lemma \ref{lemma:a-priori-bound-on-continuation-error} and the condition on $c(s)$ implies that for every pair $(x_-, x_+) \in Fix(H_-) \times Fix(H_+)$, and every continuation map solution $u \in \mathcal{M}(H_s, J_s)(x_-, x_+)$, we have that
    \begin{equation}
    \label{eq:modified-action-is-never-constant}
        \frac{d}{ds} \mathcal{A}^{cont}(s) < 0 \text{ for all } s. 
    \end{equation} 

\begin{lemma}
\label{lemma:choose-integralization-data-for-continuation-map}
    Given $\bar{\mathfrak{A}}_\pm$  arbitrary integralization data for $\CC(H_\pm, J_\pm)$, there exists an integralization datum $\bar{A}$ for $\CC(H_s, J_s)$ such that the integralization data defined in (C) of Definition \ref{def:integralization-datum-continuation-map} agree with $\bar{\mathfrak{A}}_\pm$.
\end{lemma}    

\begin{proof}
The set of $s \in \R$ where $\frac{d}{ds} \mathcal{A}_{H_s}(u_s)$ for $u \in \overline{\mathcal{M}}(H_s, J_s)(x_-, x_+)$ is all contained in the compact subset of $\R$ are where $\frac{d}{ds} H_s \neq 0$; thus, by compactness $\overline{\mathcal{M}}(H_s, J_s)(x_-, x_+)$ and finiteness of $Fix(H_\pm)$, by choosing $c(s)$ to be a sufficiently large rescaling of any function which has strictly negative derivative over that region will satisfy the first (and most difficult to satisfy) constraint. The remaining constraints together with the choice of $c$ characterize $\bar{\mathfrak{A}}$. 
\end{proof}

\paragraph{Shifting integralization data.}
Let $\bar{\mathfrak{A}} = (\mathcal{A}, \mathfrak{A})$ be an integralization datum for $\CC(H, J)$. For $D > 0$ and any constant $C$, define 
\[ D \bar{\mathfrak{A}} + C = (D \mathcal{A} + C, D \mathfrak{A}+C).\]
then this is also an integralization datum for $\CC(H, J)$ which defines the same stratification of the underlying pre-flow category. Similarly, if $\mathfrak{P}$ is a regular compatible system of perturbation data for $\CC(H< J)$ with underlying integralization datum $\bar{\mathfrak{A}}$, then, defining $D\mathfrak{P}+C$ to be tuple produced from $\mathfrak{P}$ by replacing $\bar{\mathfrak{A}}$ by $D \bar{\mathfrak{A}} + C$ we see that this is also a regular system of perturbation data, and the resulting categories $\CC'(H, J)$ produced from $\mathfrak{P}$ and from $D \bar{\mathfrak{A}} + C$ are isomorphic after rescaling the codomains of the functions $\delta A$ by $D$. 

\paragraph{Parameter spaces and stratifications.}

To define the parameter spaces for our virtual smoothings of continuation map moduli spaces, we must introduce a variant of $\overline{Conf}_{r}$.

Let $\bar{\mathfrak{A}}$ be an integralization datum for $\CC(H_s, J_s)$; write 
\[ \mathfrak{A} = \mathfrak{A}_- \sqcup \mathfrak{A}_+\]
where this is the decomposition \eqref{eq:decomposition-of-marked-actions-for-continuation-map}. Choose an $A_c$ such that $\min \mathfrak{A}_- > A_c > \max \mathfrak{A}_+$; this is possible by (B) of Definition \ref{def:integralization-datum-continuation-map}. Write $\underline{\mathfrak{A}} = \mathfrak{A} \cup A_c$, with the induced ordering. Write $\underline{\mathfrak{A}}(x) = \{A \in \underline{\mathfrak{A}} : A < \mathcal{A}(x)\}$, where $\mathcal{A}$ is the augmented action underlying $\bar{\mathfrak{A}}$. We give the pre-flow category $\CC(H_s, J_s)$ the integral action 
\[ \bar{A}(x) = \# \underline{\mathfrak{A}}(x)-1.\]
This is an integral action if we set 
\[  S_{\CC(x,y)} =\underline{\mathfrak{A}}'(x,y) := \underline{\mathfrak{A}}(x,y) \setminus A_y^+, \text{ where }\underline{\mathfrak{A}}(x,y) = \{A \in \underline{\mathfrak{A}}: \mathcal{A}(x) > A > \mathcal{A}(y)\},\]
and now $A_y^+$ is the smallest element of $\underline{\mathfrak{A}} = \mathfrak{A} \cup A_c$ which is greater than $\mathcal{A}(y)$. More concretely, we can write this as 
\[\bar{A}(x) = \bar{A}_+(x) \text{ for } x \in Fix(H_+), \]
\begin{equation}
\label{eq:shifted-integral-action-around-contiunuation-map}
    \bar{A}(x) = \bar{A}_-(x) + \# \mathfrak{A}_++1 \text{ for }x \in Fix(H_-).
\end{equation}

We will now define a open subset $\overline{Conf}^o_{\underline{\mathfrak{A}}(x,y)} \subset (\overline{Conf}_{\mathfrak{A}(x,y)})^\R$. An element of $(\overline{Conf}_{\mathfrak{A}(x,y)})^\R$ is an element $a$ of $(\overline{Conf}_{\mathfrak{A}(x,y)})$ together with an additional marked $s$ coordinate, which we will call $s_{A_c}$, on one of the components of the universal curve over $a$. The subset $\overline{Conf}^o_{\underline{\mathfrak{A}}(x,y)}$ consists of those elements such that $(a, s_{A_c})$ satisfying the property that $a$ lies in a stratum where there is a component of the universal curve containing marked points labeled by actions from both $\mathfrak{A}_+$ and from $\mathfrak{A}_-$, and $s_{A_c}$ is also is on that component.

\begin{figure}[h!]
    \centering
    \includegraphics[width=\textwidth]{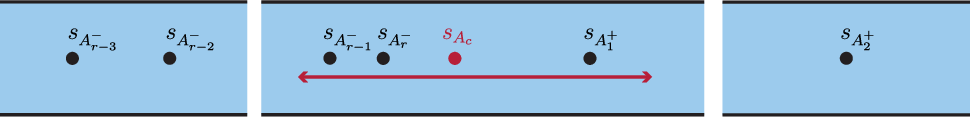}
    \caption{\emph{Parameter spaces of curves for continuation map.} When defining global charts for the flow categories corresponding to the continuation map, one must modify the parameter spaces in the way described above. The marked actions from $\mathfrak{A}_\pm$ are denoted $A^\pm_i$ above. One adds an extra marked $s$-coordinate (labeled $s_{A_c}$ which is used to `center' the corresponding component of the domain so as to write the continuation map equation; however, that marked point has no component of the obstruction section constraining it to lie on the modified action value $A_c$, and moreover it can move anywhere on the corresponding component of the domain, not being required to lie in between $s_{A^-_r}$ and $s_{A_1^+}$ or otherwise. (Indeed, there is no reason that the modified action value in the `middle' of a coninuation map trajectory would be below those two values). A precise description of this space is given in the text. The discussion for continuation map homotopies is similar. Because of condition $D$ of Definition \ref{def:integralization-datum-continuation-map}there are no breakings of broken continuation trajectories in regions with modified action between $A_r^-$ and $A_1^+$, one does not need to add corresponding strata to the moduli space. In the stratification, we think of the set organizing the stratification of this space as the spaces in bewteen the subsequent marked $s$-coordiantes points (even though $s_{A_c}$ may be out of order on a given curve), and one removes an element from the set for every breakpoint (one never removes the elements corresponding to the spaces to the left and right of $s_{A_c}$ due to Condition (D), thus dealing with the issue regarding the fact that $s_{A_c}$ may not lie in order).  }
    \label{fig:parameter-spaces-for-continuation-map}
\end{figure}

Notice that we may have that $s_{A_c} < s_{\min \mathfrak{A}_-}$ or $s_{A_c} > s_{\max \mathfrak{A}_+}$; these are precisely the constraints that would hold if we considered $\overline{Conf}_{\underline{\mathfrak{A}(x,y)}}$ instead, and are precisely the constraints that we wish to ignore in our construction. See also the discussion in Figure \ref{fig:parameter-spaces-for-continuation-map}, which also depicts the corresponding parameter spaces. 

We now define a stratification of $\overline{Conf}^o_{\underline{\mathfrak{A}}(x,y)}$ making it into a $\langle \underline{\mathfrak{A}}(x,y)'\rangle$-manifold. Recall as before that we identify $\underline{\mathfrak{A}}'(x,y)$ with the regions in between the elements of $\underline{\mathfrak{A}}(x,y)$; given an element $\tilde{a} \in \overline{Conf}^o_{\underline{\mathfrak{A}}(x,y)}$  we put it in the stratum $S(\tilde{a}) =\underline{\mathfrak{A}}'(x,y) \setminus S^c(\tilde{a})$, where $S^c(\tilde{a})$ contains the elements corresponding to the regions where there are break-points in between marked coordinates labeled by elements of $\mathfrak{A}(x,y)$. (See Figure \ref{fig:parameter-spaces-for-continuation-map}.)
Now, for $x>y$ $\CC(H_t, J_t)$ with $x \in Fix(H_-)$ and $y \in Fix(H_+)$, we define 
\[ \paramspace(x,y) = \paramspace(\bar{A}(x), \bar{A}(y)) = \overline{Conf}^o_{\underline{\mathfrak{A}}(x,y)}. \] 

We define the universal curves over over $\overline{M}(x,y)$ to be the pullback of the corresponding universal curves from the first factor. 

One has maps  $\iota_{xzy}$ covered by maps $\iota'_{xzy}$ identifying  universal curves as in \eqref{eq:identification-of-universal-curve-domains}.  Systems of perturbation data and their compatibility are defined as earlier but with the new notation. Specifically, given a point $p$ in $\paramspace(x,y)^Z$, let $\tilde{s}(p)=-\infty$ if $p$ lies on a component to the right of $s_{A_c}$, $\tilde{s}(p)=+\infty$ if $p$ lies on a component to the left of $s_{A_c}$, and otherwise $\tilde{s}(p)$ is the value of $s(p)$ if we apply an $s$-translation taking $s_{A_c}$ is to zero.  The definition of $T(x_-, x_+)$ for $x_\pm$ both in $Fix(H_-)$ or $Fix(H_+)$ are unchanged; however, if $x_\pm \in Fix(H_\pm)$, we define 
\begin{equation}
\label{eq:thickenings-for-floer-contiuation-homotopy}
    T(x_-, x_+) = \left\{ \begin{array}{c|r} 
    (u,a,v) : v \in V^1(x_-, x_+), &  \partial_s u(s,t) + J_{\tilde{s},t} \partial_t u(s,t) -\Grad H_{\tilde{s}}(u(s,t))= \lambda_V(((s,t), a), v),\\
 a \in \paramspace(x,y), u: \paramspace(x,y)^Z_a \to M &   \lim_{s \to -\infty} u_1(s, \cdot) = x_-, \lim_{s \to \infty} u_r(s, \cdot) = x_+\\ 
   u = u_1 \cup \ldots \cup u_r  & \lim_{s \to + \infty} u_j(s, \cdot) = \lim_{s \to - \infty} u_{j+1}(s, \cdot), j=1, \ldots, r-1
 \end{array}\right\}
\end{equation}
where we interpret $J_{\pm\infty}$ as $J_\pm$. 

Inclusions like \eqref{eq:thickening-boundary-inclusions} naturally exist with this definition, and there are vector bundles $\underline{V^1(x,y)}$ defined as in \eqref{eq:obstruction-bundle-for-my-stabilization} when $x$ and $y$ are fixed points of the same Hamiltonian, and if $x \in Fix(H_-)$ and $y \in Fix(H_+)$ is 
\begin{equation}
    \label{eq:obstruction-bundle-for-my-stabilization-continuation-map}
    V^1_T(x,y) = V^1(x,y) \oplus V^1_{can}(x,y), V^1_{can}(x,y) = \R^{\mathfrak{A}(x,y)}.
\end{equation} 

\begin{remark}
\label{rk:shift-in-integralization-due-to-continuation}
As before $V^1_{can}(x,y)$ comes from the $E^s$-parameterization $V^1_{can}(x)$ \eqref{eq:obstruction-bundle-for-my-stabilization-continuation-map}. However, note that $\mathfrak{A}(x,y)$ differs in general from $\mathfrak{A}(x,y)$ by the element $A_c$;  thus, we will \emph{not} have a component of obstruction bundle corresponding to the marked point $s_{A_c}$, which in turn corresponds to the fact that one does not quotient the space of solutions to the continuation equation by an $\R$-action. 
\end{remark}

The obstruction section, as in all cases, is 
\[\sigma(x,y): T(x,y) \to \underline{V^1_T(x,y)}, \sigma(x,y)(u,a,v) = (v, \delta A(u,a))  \]
where 
\[ \delta A(u,a) = (\mathcal{A}(u|_{s_A}) - A)_{A \in \mathfrak{A}(x,y)}.\]

As before we have maps like \eqref{eq:thickening-composition-maps}, and the property \eqref{eq:modified-action-is-never-constant} implies that $\sigma(x_-,x_+)^{-1}(0)$ is exactly the space of broken continuation map trajectories trajectories $\CC(H_s, J_s)(x_-, x_+)$. Thus we define a virtual topological flow category 
$\CC(H_s, J_s)(x,y) = (T(x_-, x_+), V^1_T(x_-, x_+), \sigma(x_-, x_+))$.

\begin{proposition}
\label{prop:perturbation-data-continuation-map-exist}
    Suppose that we have fixed compatible regular systems of perturbation data $\mathfrak{P}_\pm$ over $\CC(H_\pm, J_\pm)$. There exists a compatible regular system of perturbation data $\mathfrak{P}$ over $\CC(H_s, J_s)$ (with integralization data defined as above) such that the restriction of $\mathfrak{P}$ to $\CC(H_\pm, J_\pm)$ is a stabilization of $\mathfrak{P}_\pm$ by an $E^s$-parameterizaton equipped with trivial perturbation data. 
\end{proposition}
\begin{proof}
    Recall that the parameterizations underlying a compatible system of perturbation data are free. Thus, $\mathfrak{P}_{-}$ extends to a compatible system of perturbation data over $\CC(H_s, J_s)$ which restricts to $\CC(H_-, J_-)$ in a unique fashion given the constraing that the underlying $E^s$-parameterization of the extension (which we will also denote by $\mathfrak{P}_-$) assigns to any element $x \in \CC(H_+, J_+)$ the zero vector space. We can then run the argument in the proof of Proposition \ref{prop:compatible-perturbation-data-exist} twice, each time with the following modifications:
    \begin{itemize}
        \item The first time, the goal is to extend $\mathfrak{P}_+$ to a compatible system of perturbation data over $\CC(H_s, J_s)$ (also denoted by $\mathfrak{P}_+$) such that for $x > y$ in $\CC(H_-,J_-)$, $\lambda_{V_{\mathfrak{P}_+}(x,y)} = 0$. Write $\mathcal{E}^s_+$ for the $E$-parameterization underlying $\mathcal{P}_+$; this $E$-parameterization automatically extends to $\CC(H_s, J_s)$ (and will be denoted by the same symbol) by requiring that it is built from the zero parameterization only by stabilizing at objects of $\CC(H_+, J_+)$, and otherwise such that canonical part of the paramerization assigns the zero vector space to all markings in $\underline{\mathfrak{A}} \setminus \mathfrak{A}_+$. The goal of this bullet is accomplished by running the induction in Proposition \ref{prop:compatible-perturbation-data-exist} starting at $\mathcal{P} = \CC(H_+, J_+)$ with $\mathcal{E}_\mathcal{P} = \mathcal{E}_+$ and $\mathfrak{P}_{\mathcal{P}} = (\bar{\mathfrak{A}}, \mathcal{F}_\mathcal{P}|_{\mathcal{P}}, \{\lambda_{V_\mathcal{P}})$ (with notation as in the proof of Proposition \ref{prop:compatible-perturbation-data-exist}); the inductive condition is weakened to the condition that $\mathfrak{P}_{\mathcal{P}}$ restricts to $\mathfrak{P}_+$ along $\CC(H_+, J_+)$ and such that $\lambda_{V_{\mathfrak{P}_{\mathcal{P}}}(x,y)} = 0$ for $x, y \in \mathcal{P} \cap \CC(H_-, J_-)$. At each step of the inner loop if $y \in \CC(H_-, J_-)$ we set $\lambda_{V_{\mathcal{P}_x, \mathcal{Q}'}}(x,y) = 0$ and for other $y$ we simply do \emph{nothing} when exectuting the inner loop. The final value $\mathfrak{P}_{\CC(H_s, J_s)}$ defines the extension $\mathfrak{P}_+$. 
        \item At this point we have the compatible system of perturbation data  $\mathfrak{P}_\pm = \mathfrak{P}_- \oplus \mathfrak{P}_+$ (where direct sum is defined by taking direct sum of $F^s$-parameterizations and direct sum of perturbation data maps $\lambda$) which is regular when restricted to $\CC(H_\pm, J_\pm)$. Thus for any $x \in \CC(H_-, J_-)$ and $y \in \CC(H_+, J_+)$ there is a compact subset $K_{xy} \subset \CC(H_s, J_s)(x,y)$ defined by the property that $u \in K$ if $D^{\mathfrak{P}_+ \oplus \mathfrak{P}_-}_u$ is not surjective. We then rerun the entire double induction of Proposition \ref{prop:compatible-perturbation-data-exist} starting at $\mathcal{P}=0$ but with the inductive condition in the outer loop being that $\lambda_{V_{\mathfrak{P}_{\mathcal{P}}}}(x,y) =0$ for both $x$ and $y$ in $\CC(H_-, J_-)$ or in $\CC(H_+, J_+)$; and also that $D^{\mathfrak{P}_+ \oplus \mathfrak{P}_-}_u$ is surjective for $u \in K_{xy}$; it is straightforward to modify the inner loop to achieve this condition. After the double induction terminates we set $\mathfrak{P}_{\CC(H_s, J_s)} =: \mathfrak{P}'$.         
    \end{itemize}
    We finalize by setting $\mathfrak{P} = \mathfrak{P}_- \oplus \mathfrak{P}_+ \oplus \mathfrak{P}'$. 
\end{proof}

We will write $\CC'(H_\pm, J_\pm)$ for the topologically-smooth flow categories associated to $\mathfrak{P}_\pm$. We have the following proposition showing the smoothings we choose for either of these are compatible with a smoothing for the continuation-map flow category: 
\begin{proposition}
\label{prop:continuation-map-compatible-smoothings}
 The category  $\CC'(H_s, J_s)$ defined using $\mathfrak{P}$ produced by Proposition \ref{prop:perturbation-data-continuation-map-exist},  is topologically smooth. 
 
Moreover, let $\mathcal{E}_\pm$ be $E$-parameterizations such that $\CC'(H_\pm, J_\pm)_{\mathcal{E}_\pm}$ have been given smooth structures. Then there exists an $E^s$-parameterization $\mathcal{E}$ such that such that $\mathcal{E}|_{\CC(H_\pm, J_\pm)}$ is a semi-free stabilization of $\mathcal{E}_\pm$ by some $\mathcal{E}_\pm'$, and such that $\CC'(H_s, J_s)_{\mathcal{E}}$ has a smooth structure 
with with $\CC'(H_s, J_s)_{\mathcal{E}}|_{\CC(H_\pm, J_\pm)} = (\CC'(H_\pm, J_\pm)_{\mathcal{E}_\pm})_{\mathcal{E}_\pm'}$ as virtually smooth flow categories. 
\end{proposition}
\begin{proof}
    The fact that $\CC'(H_s, J_s)$ is topoloogically smooth follows from the corresponding analog of Proposition \ref{prop:pardon-gluing}, which is proven by the same arguments from \cite{pardon2016algebraic}. By the conditions that $\mathfrak{P}$ satisfies we have that $\CC'(H_s, J_s)|_{\CC(H_\pm, J_\pm})$ is indeed a stabilization of $\CC'(H_\pm, J_\pm)$ by some $E^s$-parameterization $\mathcal{E}_\pm^0$. Extend each of $\mathcal{E}_\pm$ to $\CC(H_s, J_s)$ by requiring that the extension is built from the zero parameterization by stabilizing only at elements of $\CC(H_\pm, J_\pm)$. Write 
    \[ \CC''(H_s, J_s) =  \CC'(H_s, J_s)_{\mathcal{E}_- \oplus \mathcal{E}_+}|_{\CC(H_\pm, J_\pm)};\]
    we have that 
    \[ \CC''(H_s, J_s) = \CC'(H_\pm, J_\pm)_{\mathcal{E}_\pm^0 \oplus \mathcal{E}_- \oplus \mathcal{E}_+} = (\CC'(H_\pm, J_\pm)_{\mathcal{E}_\pm})_{\mathcal{E}_\mp \oplus \mathcal{E}^0_\pm}\]
    and as such the full subcategories $\CC''(H_s, J_s)_{\CC(H_\pm, J_\pm)}$ are equipped with smooth structures induced from the smooth structures given on $\CC'(H_\pm, J_\pm)_{\mathcal{E}_\pm}$. We now rerun the equivariant double induction in the proof of Proposition \ref{prop:smoothing-a-flow-category} on $\CC''$, with the proviso that for pairs $(x, y)$ (with notation as in the proof of that proposition) for which both $x$ and $y$ are in $\CC(H_-, J_-)$ or in $\CC(H_+, J_+)$, we use the smooth structure on $T(x,y)$ induced by its presentation as the product of $T_{\CC'(H_s, J_s)}(x,y)$ with a vector space. The conclusion of that double induction outputs an $E^s$-parameterization $\mathcal{E}_{\CC(H_s, J_s)}$, and we can set  $\mathcal{E} = \mathcal{E}_{\CC(H_s, J_s)} \oplus \mathcal{E}_- \oplus \mathcal{E}_+$. \end{proof}

We have now a sufficient number of results to establish the existence of continuation maps between Floer homotopy types associated to Hamiltonian Floer homologies. We combine the relevant results in Section \ref{sec:geometric-results-review}.

\subsection{Cyclotomic compatibility}
In previous sections we have described how to produce $C_k$-equivariant virtually smooth flow categories  $\CC'(H, J)$ and $\CC'(H_s, J_s)$ when $(H, J)$ or $(H_s, J_s)$ are $C_k$-invariant. This latter condition of course means e.g. that $(H, J) = (H^{\# k} , J^{\# k})$. In this section we explain how to pick perturbation data and smoothings such that the virual smoothings $\CC'(H^{\#k}, J^{\#k })$ and $\CC'(H^{\#k}_s, J^{\#k}_s)$ are all compatible \emph{as $k$ varies}. We will only work out the case of $\CC'(H^{\#k}_s, J^{\#k}_s)$ since it contains the case of $\CC'(H^{\#k}, J^{\#k })$ as a subproblem.

Specifically, we prove the following 
\begin{proposition}
\label{prop:cyclotomic-compatibility-all-perturbation-data-choices}
    Let $\mathfrak{P} = (\bar{\mathfrak{A}}, \mathfrak{F}, \{\lambda_{V_{\mathfrak{P}}(x,y)}\})$ be a compatible regular system of perturbation data for $\CC(H_s, J_s)$ such that $\CC'(H_s, J_s)$ (defined using $\mathfrak{P}$ has a smooth structure; if $(H_s, J_s)$ are $C_\ell$-invariant we assume that all the associated data are $C_\ell$-equivariant as well. Suppose also that $\bar{\mathfrak{A}}$ satisfies the conditions in Lemma \ref{lemma:cyclotomic-compatibility-for-integralization-data-continuation-map} below. Fix $k$ a natural number greater than one. There exists a compatible $C_{k\ell}$-equivariant regular system of perturbation data $\mathfrak{P}_k = (\bar{\mathfrak{A}}_k, \mathfrak{F}_k, \{\lambda_{V_{\mathfrak{P}_k}(x,y)}\})$ for $\CC(H^{\#k}_s, J^{\#k}_s)$ such that
    \begin{itemize}
        \item The flow category $\CC(H^{\#k}_s, J^{\#k}_s)^{C_k}$ (with the stratification induced from $\bar{\mathfrak{A}}_k$ is a restratification of $\CC(H_s, J_s)$;
        \item The topologically smooth category $\CC'(H^{\#k}_s, J^{\#k}_s)$ defined using $\mathfrak{P}_k$ has a $C_{k\ell}$ smooth structure, which 
        \item Restricts to a smooth structure on the virtually smooth  $C_{k\ell}/C_k$-equivariant flow category $\CC'(H^{\#k}_s, J^{\#k}_s)^{C_k}$, which in turn is isomorphic to a  stabilization ( by an $E^s$-parameterization) of a restratification $\CC'(H_s, J_s)$ as a $C_\ell = C_{k\ell}/C_k$-equivariant virtually smooth flow category. 
    \end{itemize}
\end{proposition}

We prove this proposition by combining several lemmas proven below. 
\begin{proof}[Proof of Proposition \ref{prop:cyclotomic-compatibility-all-perturbation-data-choices}]
Combine Lemmas \ref{lemma:cyclotomic-compatibility-for-integralization-data-continuation-map}, \ref{lemma:change-of-perturbation-datum-and-restratification}, \ref{lemma:change-of-perturbation-datum-and-restratification-exists}, and \ref{lemma:smoothings-cyclotomic-compatibility}. 
\end{proof}

\paragraph{Shifting constants and a-priori bounds on actions of periodic points.}

\begin{lemma}
\label{lemma:apriori-bound-on-action-values}
    Given a an admissible symplectic manifold $(M, \omega)$ and a Hamiltonian function $H: M \times S^1 \to \R$, there is a constant $C_H$ depending only on $M$ and $H$ such that $\mathcal{A}_{H^{\#k}}(Fix(H^{\#k})) \subset [kC_H, -kC_H]$. 
\end{lemma}
\begin{proof}
    Pick a compatible almost complex structure $J$ on $M$. The action functional is 
    \[\mathcal{A}(\gamma) = \int_0^1 \lambda(\gamma') - H(\gamma). \]
    Recall that we have that $\gamma' = - \Grad H$; since $M$ is compact we have the bound 
    \[ \mathcal{A}(\gamma) \leq \|\lambda\|_{C^0} \|H\|_{C^1}\} + \|H\|_{C^0} \leq C_1 \|H\|_{C^1};\]
    for some constant $C_1$, where we write $\|H\|_{C^1} = \sup_{t \in S^1} \|H_t\|_{C^1}$ and similarly for $\|H\|_{C^0}$. Notice that $C_1$ is independent of $H$; the proposition then follows from the fact that $\|H^{\#k}\|_{C^1} = k\|H\|_{C^1}$.
\end{proof}

\paragraph{Compatibly choosing integralization data.}
In order to proceed we will need to change the integralization data for $\CC(H_s, J_s)$. The integralization data will change in a comparatively minor way: it will be \emph{extended}.

\begin{definition}
     Let $\bar{\mathfrak{A}} = (\mathcal{A}, \mathfrak{A})$. and $\bar{\mathfrak{A}}' = (\mathcal{A}', \mathfrak{A}')$ be a pair of integralization data for $\CC(H, J)$ or for $\CC(H_s, J_s)$ (in the latter case we add on the extra quantities $\{A_c\}$ and $\{A'_c\}$ to the respective tuples to the tuple. We say that $\bar{\mathfrak{A}}'$ is an \emph{extension} of $\bar{\mathfrak{A}}$ if 
    \[ \mathfrak{A}' = \mathfrak{A} \sqcup \mathfrak{A}_{new},\]
    and in the case that we are discussing $\CC(H_s, J_s)$ we also have that $A_c = A'_c$. 
\end{definition}

On the way to proving Lemma \ref{prop:cyclotomic-compatibility-all-perturbation-data-choices}, we will prove the following helpful Lemma as well: 

\begin{lemma}
\label{lemma:extend-perturbation-data-across-restratifications}
    Let $\mathfrak{P}$ be a perturbation datum for $\CC'(H_s, J_s)$ with underlying integralization datum $\bar{\mathfrak{A}}$. Let $\bar{\mathfrak{A}}'$ be an extension of $\bar{\mathfrak{A}}$.  
    Then there is a perturbation datum $\mathfrak{P}'$ for $\CC(H_s, J_s)$ with underlying integralization datum $\bar{\mathfrak{A}}'$ such that the virtual topological smoothing $\CC'(H_s, J_s)_{\mathfrak{P}'}$ associated to $\mathfrak{P}$ has agrees with a restriction of a  restratification of $\CC'(H_s, J_s)_{\mathfrak{P}}$, the topological smoothing associated to $\mathfrak{P}$. The same holds with $\CC(H_s, J_s)$ replaced by $\CC(H, J)$.
\end{lemma}

\begin{remark}
As before, there is an action of orientation-preserving affine rescalings of $\R$ on integralization data for $\CC(H_s, J_s)$ namely, $C\bar{\mathfrak{A}} + D = (C \mathcal{A}+D, C \mathfrak{A}+D, \{CA_c + D\})$, and given a perturbation datum $\mathfrak{P}$ with integralization data $\bar{\mathfrak{A}}$, there is a corresponding perturbation datum $C \mathfrak{P} + D$ with underlying integralization datum $C\bar{\mathfrak{A}}+D$ defining a virtual smoothing which is canonically isomorphic to that defined via $\mathfrak{P}$. 
\end{remark}

\begin{remark}
In Proposition \ref{prop:cyclotomic-compatibility-all-perturbation-data-choices}, the integralization datum $\bar{\mathfrak{A}}_k$ will satisfy the property that if we think of it as an integralization datum for $\CC(H_s^{\# k}, J_s^{\#k})^{C_k}$, it will have the form $\bar{\mathfrak{A}}_k = k \bar{\mathfrak{A}}'$ where $\bar{\mathfrak{A}}'$ is an extension of $\bar{\mathfrak{A}}$ under the identification $\CC(H_s^{\# k}, J_s^{\#k})^{C_k} \simeq \CC(H_s, J_s)$. The lemma below explains how to construct $\bar{\mathfrak{A}}_k$ from $\bar{\mathfrak{A}}$.
\end{remark}

\begin{remark}
    Note that $\mathcal{A}_{H^{\#k}}(\gamma^k)$ for a $C_k$-invariant loop $\gamma^k(t) = \gamma(kt)$ agrees with $k\mathcal{A}_H(\gamma)$. Thus, implicitly in the remark above in the Lemma below, we are using the fact that if we identify $\CC(H_s^{\# k}, J_s^{\#k})^{C_k}$ with $\CC(H_s, J_s)$, then the integralization datum $\bar{\mathfrak{A}}$ for $\CC(H_s, J_s)$ corresponds to the integralization datum $k \bar{\mathfrak{A}}$ for $\CC(H_s^{\# k}, J_s^{\#k})^{C_k}$.
\end{remark}

\begin{lemma}
\label{lemma:cyclotomic-compatibility-for-integralization-data-continuation-map}
    Let $\bar{\mathfrak{A}} = (\mathcal{A}, \mathfrak{A}, \{A_c\})$ be an integralization datum for $\CC(H_s, J_s)$. 
    Write $\mathcal{A} = D\mathcal{A}_{H_s} + c(s)$ with $c(s) \to C_\pm$ as $s \to \pm \infty$. Suppose that 
    \begin{equation}
        \label{eq:total-decrease-condition}
        C_- - C_+ > DC_{H_-} + DC_{H_+} + \int_{-\infty}^\infty D\max_{(x,t) \in M \times S^1} |\partial_s H_s|ds
    \end{equation} 
    and 
    \begin{equation}
        A_c \in (C_- -DC_{H_-}, C_+ +DC_{H_+}). 
    \end{equation}
    where $C_{H_\pm}$ are the constants defined in Lemma \ref{lemma:apriori-bound-on-action-values}.
    Identify the objects of $\CC(H_s, J_s)$ with a subset of the objects of $\CC(H_s^{\#k}, J_s^{\#k})$ via iteration. 
    Then there exists a set $\mathfrak{A}^k_n\subset \R \setminus k \mathfrak{A}$ such that
    \[ \bar{\mathfrak{A}}_k = (D\mathcal{A}_{H_s^{\#k}} + kc(s), k \mathfrak{A} \cup \mathfrak{A}^k_n, \{kA_c\}) \]
    is an integralization datum for $\CC(H_s^{\#k}, J_s^{\#k})$. Moreover, viewing $\bar{\mathfrak{A}}_k$ as an integralization datum for $\CC(H_s^{\#k}, J_s^{\#k})^{C_k} \simeq \CC(H_s, J_s)$, the  integral action on $\CC(H_s, J_s)$ induced via $\bar{\mathfrak{A}}_k$ is a restratification of the  integral action induced on $\CC(H_s, J_s)$ via $\bar{\mathfrak{A}}$.
\end{lemma}
\begin{proof}

    Note first that $H_s^{\#k}(x,t) = k H_s(x, kt)$, so 
    \[ \max_{(x,y) \in M \times S^1} |\partial_s H_s^{\#k}(x,t)| = k\max_{(x,y) \in M \times S^1} |\partial_s H_s(x,t)|. \]
    Therefore, $kc(s)$ satisfies Condition $A$ from Definition \ref{def:integralization-datum-continuation-map} for $H_s^{\#k}$.  In order for the set $\mathfrak{A}_{cont, new}$ to exist, it is sufficient to have that 
    \[ \min_{x \in \CC(H_-, J_-)} D \mathcal{A}_{H_-^{\#k}}(x) + k C_- > \max_{x \in \CC(H_+, J_+)} \mathcal{A}_{H_-^{\#k}}(x) + k C_+; \]
    this is guaranteed by the condition \eqref{eq:total-decrease-condition} as well as Lemma \ref{lemma:apriori-bound-on-action-values}. The last statement is elementary.
\end{proof}

\paragraph{Compatibly choosing perturbation data.}
We will now explain how to use the integralization data $\bar{\mathfrak{A}}_k$ produced by Lemma \ref{lemma:cyclotomic-compatibility-for-integralization-data-continuation-map} to extend perturbation data $\mathfrak{P}$ for $\CC(H_s, J_s)$ to perturbation data $\mathfrak{P}'$ for $\CC(H_s^{\#k}, J_s{\# k})$ which are compatible with $\mathfrak{P}$ in an appropriate fashion. Before we can explain the compatibility conditions, we must first explain some geometric statements about the change in the parameter spaces $\paramspace(x,y)$ associated to $\CC(H_k, J_k)^{C_k}$ induced by the extension $k \bar{\mathfrak{A}} \to \bar{\mathfrak{A}}_k$. This process is geometrically rather elementary, as one can see by ``adding new marked points'' to Figure \ref{fig:description-of-marking-fixing-condition}.

Write $\paramspace(x,y)$ for the parameter spaces associated by $\bar{\mathfrak{A}}$ to $\CC(H_s, J_s)$, and write $\paramspace_k(x,y)$ for the parameter spaces associated by $\bar{\mathfrak{A}}_k$ to $\CC(H_s^{\#k}, J_s^{\#k})$. Given an object $x \in \CC(H_s, J_s)$, write $x^k$ for the corresponding object of $\CC(H_s^{\#k}, J_s^{\#k})$. Then thinking of $k \bar{\mathfrak{A}}$ as an integralization datum for $\CC(H_s^{\#k}, J_s^{\#k})^k$, this data assigns parameter spaces to $(x^k, y^k)$ which are canonically isomorphic to $\paramspace(x,y)$. 

Recall that
\[\paramspace(x,y) = \overline{Conf}^c_{\underline{\mathfrak{A}}(x,y)}, \paramspace_k(x,y) = \overline{Conf}^c_{\underline{\mathfrak{A}_k}(x,y)}. \]

There open submanifold $U(x^k, y^k) \subset \overline{Conf}^c_{\underline{\mathfrak{A}_k}(x^k,y^k)}$ corresponding to those curves for which if we removed the marked points labeled by elements of  $\mathfrak{A}_k \setminus k \mathfrak{A}$, we would get a stable domain. Thus there is a projection 
\begin{equation}
    \label{eq:forgetting-markings-added-under-extension-of-integralization-data}
    m: U(x^k, y^k) \to \paramspace(x,y).
\end{equation}
We have that 
\begin{equation}
    \label{eq:preservation-by-domain-composition-equation}
    \iota_{x^k z^k y^k}(U(x^k, z^k) \times U(z^k, y^k)) = U(x^k, y^k).
\end{equation}
Write $U(x^k, y^k)^Z$ for the preimage of $U(x^k, y^k)$ in the total space of of the universal curve to $\paramspace(x^k, y^k)^Z$.

Note that the definition of a compatible system of perturbation data (Definition \ref{def:system-of-perturbation-data}) continues to make sense if the maps $\lambda_{V^1(x,y)}$ (in the notation of that definition) are only partially defined over the $U(x^k, y^k)$, i.e. 

Writing $\mathfrak{P} = (\bar{\mathfrak{A}}, \mathcal{F}, \{\lambda_{V_\mathcal{F}(x,y)}\}$ for a compatible regular system of perturbation data for $\CC(H_s, J_s)$, we can define the data
\begin{equation}
\label{eq:induced-system-of-perturbation-data}
    \mathfrak{P}'_0 = (\bar{\mathfrak{A}}_k, \mathcal{F}_0, \{\lambda_{V_\mathcal{F}(x^k,y^k)}\}_{x, y \in Fix(H)}), \lambda_{V_\mathcal{F}(x^k,y^k)}(a,v) = \lambda_{V_\mathcal{F}(x,y)}(m(a), v).
\end{equation}

The above data is similar to a compatible system of perturbation data, but instead is only a partially-defined such system. We explain the notation. 
The $F^s$-parameterization $\mathcal{F}_0$ is the $F^s$-parameterization obtained from $\mathcal{F}$ by ensuring that the following $E^s$-parameterization $\mathcal{E}_0 = ((\mathcal{E}_0)_f \oplus (\mathcal{E}_0)_c)$. First, we set $V^*_{(\mathcal{E}_0)_f}(x) = 0$ for all $x \in Fix(H^{\#k}) \setminus Fix(H)$,  otherwise $V_{(\mathcal{E}_0)_f}(x^k) = V_{\mathcal{E}_f}(x)$, where $\mathcal{E}_f$ is free part of the $E^s$-parameterization $\mathcal{E}_f$ corresponding to $\mathcal{F}$; with the $C_{\ell k}$-equivariant structure on $\mathcal{E}_0$ coming from the $C_\ell$-equivariant structure on $\mathcal{E}$ via the quotient map $C_{\ell k}/C_k \to C_\ell$. Second, we require that the vector spaces associated to all elements of $S_{\bar{A}_k} \setminus \rho(S_{k\bar{A}})$ by the canonical part $(\mathcal{E}_0)_c$ of $\mathcal{E}_0$ parameterization are zero (Definition \ref{def:canonical-parameterization}), while the vector spaces associated to the elements of $\rho(S_{k \bar{A}})$ are the vector spaces associated to  the corresponding elements of $S_{\bar{A}}$ by the canonical part $\mathcal{E}_c$ of $\mathcal{E}$. With this choice, $V_{\mathcal{F}_0}(x^k,y^k) = V_{\mathcal{F}}(x,y)$, and the formulae in \eqref{eq:induced-system-of-perturbation-data} make sense, with the proviso that the bundle maps $\lambda_{V_\mathcal{F}(x^k,y^k)}$ are defined only over $U(x^k, y^k)^Z \times M$, and because $\lambda_{V_\mathcal{F}(x',y')}$ is not defined for non-$C_k$ invariant $x'$ or $y'$. 

We now have the following 
\begin{lemma}
\label{lemma:change-of-perturbation-datum-and-restratification}
    There exist open subsets with compact closure $U^c(x^k,y^k) \subset U(x^k,y^k)$ satisfying \eqref{eq:preservation-by-domain-composition-equation} with $U(x',y')$ replaced with $U^c(x', y')$ everywhere, as well as the property that 
    \begin{equation}
        \label{eq:disjointness-of-convenient-opens-from-new-strata} \overline{U^c(x^k, y^k)} \cap \overline{im \iota_{x^zy^k}} = \emptyset \text{ for any } z \in Fix(H_s^{\#k}) \setminus Fix(H_s)
    \end{equation}
    such that if any perturbation $\mathcal{P}' = (\bar{\mathfrak{A}}_k, \mathcal{F}', \{\lambda_{V_{\mathcal{F}'}(x,y)})$ for $\CC(H_s^{\#k}, J_s^{\#k})$ satisfies 
    $(\mathcal{F}')^{C_k} = \mathcal{F}$ (see Definition \ref{def:restriction-and-fixed-points-of-parameterization}) and 
    \[ \lambda_{V_{\mathcal{F}'(x^k, y^k)}}(a, v) = \lambda_{V_\mathcal{F}(x,y)}(m(a), v)\]
    for $a \in V_{\mathcal{F}'}(x,y)^{C_k}$ and $a \in U^c(x,y)$ then there is a virtual topological flow category from a restriction of $\CC'(H_s^{\#k},J_s^{\#k})_{\mathfrak{P}'}^{C_k}$ to a restriction of $(\CC'(H_s,J_s)_{\mathfrak{P}'})_r$, where $r$ is the restratification associated to the extension of integralization data $k \mathfrak{A} \to \mathfrak{A}_k.$  
\end{lemma}
\begin{proof}
    Write $T(x,y)$ for the thickenings associated to $\mathfrak{P}$, $T(x,y)_r$ for the thickenings associated to $\CC'(H_s,J_s)_{\mathfrak{P}'})_r$, and $T'(x,y)$ for the thickenings associated to $P'$. We will first produce maps 
    \[ q_{xy}: \CC(x,y) \to U(x,y), \iota_{xyz} \circ( q_{xy} \times q_{yz}) = q_{xz}\circ \bar{f}_{xyz};\]
    each $q_{xy}$ will extend to some open neighborhood of the zero locus in $T(x,y)_r$, and we will define $U^c(x,y)$ to be some system of open subsets of $U(x,y)$ with compact closure satisfying \eqref{eq:preservation-by-domain-composition-equation} such that $U^c(x,y) \supset q_{xy}(\CC(x,y))$ (such a system always exists by a straightforward inductive argument). Given $u \in \CC(x,y)$, we have a corresponding element $(u, a, v) \in T(x,y)$. By condition \ref{eq:modified-action-is-never-constant} there unique collection of additional marked $s$-coordinates $\{s_A\}_{A \in \mathfrak{A}_k \setminus k \mathfrak{A}}$ that we can add to $\paramspace{M}(x,y)^Z_a$ such that $\mathcal{A}(s_A) = A$. The corresponding marked curve defines the element $q_{xy}(u) \in U(x,y)$. The same definition clearly extends to a map $q_{xy}: T(x,y) \to U(x,y)$; we will explain how to extend it to a neighborhood $T(x,y)_r^o$ of $T(x,y)$ in $T(x,y)_r = T(x,y) \times (-1,1)^{\mathfrak{A}_k \setminus k \mathfrak{A}}. $ Namely, we shrink the second factor to $(-\epsilon, \epsilon)^{\mathfrak{A}_k \setminus k \mathfrak{A}}$, where $\epsilon$ is smaller than a third of the minimum distance between any two elements of $\mathfrak{A}_k \cup \mathcal{A}(Ob(\CC(H^{\#k}_s, J^{\#k}_s))$; and we define $q(u,a,v, (a_j))$ to be element corresponding to marking $\paramspace{M}(x,y)^Z_a$ with additional marked $s$-coordinates $\{s_A\}_{A \in \mathfrak{A}_k \setminus k \mathfrak{A}}$ such that $A - \mathcal{A}(u_{s_A)} = a_j$. We can always find such a set of points when $\epsilon$ is chosen as above, and again by condition \eqref{eq:modified-action-is-never-constant} this choice is unique (see Figure \ref{fig:description-of-marking-fixing-condition}); so we have extended $q$ as desired. By the fact that the projection $T(x,y) \to \paramspace(x,y)$ was a topological submersion, we see that $q$ is a topological submersion as well; in particular $q(T(x,y)_r^o)$ is open.

    Let us now pick $U^c(x,y)$ as specified previously with the additional condition that $U^c(x,y) \subset q(T(x,y)_r^o)$ as well as the condition \eqref{eq:disjointness-of-convenient-opens-from-new-strata} (which is straightforward since $q(T(x,y)_r^o) \cap \overline{im \iota_{x^kz y^k}}$ for the newly introduced objects $z$) and let replace the sets $T(x,y)_r^o$ by their subsets $q^{-1}(U^c(x,y))$; this defines the desired restriction of $(\CC'(H_s,J_s)_{\mathfrak{P}'})_r$. The desired restriction of $\CC'(H_s^{\#k},J_s^{\#k})_{\mathfrak{P}'}^{C_k}$ is given by taking the thickenings to be the subsets given by the preimages of $U^c(x,y)$ under the projections of the thickenings to $\paramspace_k(x,y)$.  

    Suppose now that $\mathfrak{P}'$ exists as in the lemma. The desired homeomorphism on thickenings is 
    \begin{equation}
        \label{eq:identification-of-thickening-with-restratification}
        T(x,y)_r^o \ni (u, a,v, (a_j)) \mapsto (u_{q(a)}, q(a), v) = (u', b, w) \in  T_k(x^k, y^k).
    \end{equation} 
    where $u_{q(a)}$ is the composition of $u$ with identification $\paramspace_k{M}(x,y)^Z_{q(a)} \to \paramspace{M}(x,y)^Z_a$ obtained by dropping the marked $s$ coordinates added by $q$. This is certainly a continuous map; to see that it is a homeomorphism, we note that since $w \in V_{\mathcal{F}}(x,y)$ and $b \in U_c(x,y)$, we know the perturbation of the Floer equation that $u$ solves as well as its domain; but since the perturbation functions are pulled back under projection under $\pi$, in fact $u'$ it comes from a unique $u$, with the location of the extra markings exactly specified by $a_j$. This discussion defines an inverse map. These homeomorphsims extend to an isomorphism of virtual topological flow categories by requiring that the bundle maps covering the above homeomorphisms simply act by identity on the fiber.
\end{proof}

Finally, we have 
\begin{lemma}
\label{lemma:change-of-perturbation-datum-and-restratification-exists}
    A perturbation $\mathfrak{P}^k$ as in Lemma \ref{lemma:change-of-perturbation-datum-and-restratification} exists.  
\end{lemma}
\begin{proof}
    We rerun the argument proving Proposition \ref{prop:compatible-perturbation-data-exist} twice, with minor modifications similar to the proof of Proposition \ref{prop:perturbation-data-continuation-map-exist}:
    \begin{itemize}
        \item We rerun the double induction of Proposition \ref{prop:compatible-perturbation-data-exist}, first while fixing $\mathcal{E}_\mathcal{P} = \mathcal{E}|_\mathcal{P}$ and $V_{\mathcal{P}; \mathcal{Q}}(x,y) = V_{\mathcal{F}}|_{\mathcal{P}}(x,y)$ where $\mathcal{F}$. As we run this induction, we always set perturbation functions so that they satisfy the condition
        \[\lambda_{V_{\mathcal{P}; \mathcal{Q}}(x^k,y^k)}(a, v) = \lambda_{V_\mathcal{F}(x,y)}(m(a), v) \text{ for } a \in U^c(x,y);\]
        with $U^c(x,y)$ as in Lemma \ref{lemma:change-of-perturbation-datum-and-restratification} above; such choices exist by \eqref{eq:disjointness-of-convenient-opens-from-new-strata}. 
        \item The second time we rerun the induction to stabilize the previously chosen perturbation datum by a new perturbation datum $\mathfrak{P}_k^{not-inv}$ to a regular perturbation datum, such that the underlying $F$-parametriztaion $\mathcal{F}^{not-inv}$ of $\mathfrak{P}_k^{not-inv}$ satisfies $(\mathcal{F}^{not-inv})^{C_p} = 0$; this is possible because the previously chosen perturbation datum already makes all the $C_k$-invariant Floer trajectories regular. 
    \end{itemize}
\end{proof}

\begin{proof}[Proof of Lemma \ref{lemma:extend-perturbation-data-across-restratifications}]
One follows precisely the same argument as the one above, but replacing $\paramspace_k(x,y)$ with $\paramspace'(x,y)$, the parameter spaces associated to $\CC(H_s, J_s)$ (or to $\CC(H, J)$) by the integralization data $\mathfrak{A}'$, replacing $\mathfrak{A}_k$ with $\mathfrak{A}'$ everywhere, removing all mention of $C_k$ fixed points and all quotients $C_{\ell k} \to C_\ell$, and removing all superscripts from notation $x^k$; where the condition  \eqref{eq:disjointness-of-convenient-opens-from-new-strata} in Lemma \ref{lemma:change-of-perturbation-datum-and-restratification} is vacuous, and the second step in the proof of Lemma \ref{lemma:change-of-perturbation-datum-and-restratification-exists} does not need to be implemented. 
\end{proof}

\paragraph{Compatible smoothings.}
\begin{lemma}
\label{lemma:smoothings-cyclotomic-compatibility}
    Let $\mathcal{E}$ be an $E^s$-parameterization such that there is a smooth structure on $\CC'(H_s, J_s)_{\mathcal{E}}$. There is an $E^s$-parameterization $\mathcal{E}_k$ such that there is a smooth structure on $\CC'(H_s^{\#k}, J_s^{\#k})_{\mathcal{E} \oplus \mathcal{E}'}$ such that the corresponding smooth structure on $\CC'(H_s^{\#k}, J_s^{\#k})_{\mathcal{E} \oplus \mathcal{E}'}^{C_k}$ agrees (possibly after restriction) with that induced on 
    $\CC'(H_s, J_s)_{\mathcal{E} \oplus \mathcal{E}'}$ from that chosen on $\CC'(H_s, J_s)_{\mathcal{E}}$ (also possibly both under restriction). 
\end{lemma}
\begin{proof}
    This follows from a very slight modification of the proof of Proposition \ref{prop:smoothing-a-flow-category}. First one notes that the subspaces
    \[ (T_{\CC'(H_s^{\#k}, J_s^{\#k})})_{\mathcal{E} \oplus \mathcal{E}_{\mathcal{P}, \mathcal{Q}}}(x,y)^{C_k} \subset (T_{\CC'(H_s^{\#k}, J_s^{\#k})})_{\mathcal{E} \oplus \mathcal{E}_{\mathcal{P}, \mathcal{Q}}}(x,y)\]
    are nice submanifolds. One chooses a compatible system of elegant metrics on the thickenings (such that the maps on thickenings associated to composition are isometries); then, using Lemma \ref{??}, one sees that one has for every $x,y$, a system of tubular neighborhoods of the inclusions above which one can use to deform the standard microbundle lift to one which restricts to the smooth microbundle lift coming from the chosen elegant metrics $TT_{\CC'(H_s^{\#k}, J_s^{\#k})})_{\mathcal{E} \oplus \mathcal{E}_{\mathcal{P}, \mathcal{Q}}}(x,y)^{C_k}$, mapping the latter to $T_{\CC'(H_s^{\#k}, J_s^{\#k})})_{\mathcal{E} \oplus \mathcal{E}_{\mathcal{P}, \mathcal{Q}}}(x,y)^{C_k}$. One then runs the induction of Proposition \ref{prop:smoothing-a-flow-category} using these nice submanifolds with their elegant metrics in all propositions invoked. 
\end{proof}

\subsection{Composition of continuation maps}
\label{sec:continuation-homotopy-geometry}

In this section, we explain how to generalize the construction of Section \ref{sec:continuation-maps-geometry} to produce virtually smooth flow categories witnessing the existence of homotopies between compositions of continuation maps. Specifically, suppose we are given 
\begin{itemize}
    \item four choices of admissible Floer data $((H_1, J_1), (H_2, J_2), (H_3, J_3), (H_4, J_4))$; 
    \item four choices of admissible continuation data $(H_{ab}, J_{ab})$ from $(H_b, J_b)$ to $(H_a, J_a)$ , where $(a, b) \in \{(1, 2), (1, 3), (2, 4), (3, 4)\}$, and 
    \item A \emph{(admissible) continuation homotopy datum} $(H_{h}, J_{h})$ defined as follows. Note that there is a diffeomorphism $\phi: \overline{Conf}_2 \simeq \R_+$: the boundary point corresponds to a curve with two components, each with a marked point, and as we move away from the boundary, the two components merge into one and the two marked points move closer together. Let us write  $\overline{Conf}_{2, h, 0}$ and $\overline{Conf}_{2, h, 1}$ for two copies of $\phi^{-1}([0,1/2])$. We can defined structure of a $\langle 1 \rangle$-manifold on
    \[\overline{Conf}_{2, h} = \overline{Conf}_{2, h, 0} \cup_{\phi^{-1}(1/2)} \overline{Conf}_{2, h, 1} \]
    by requiring that the boundary of $\overline{Conf}_{2, h}$ is the union of the boundaries of $\overline{Conf}_{2, h, 0}$ and of $\overline{Conf}_{2, h, 1}$, and more generally that there is a diffeomorphism 
    \[ \tilde{\tau}: \overline{Conf}_{2, h} \to [0,1], \tilde{\tau}(x)|_{\overline{Conf}_{2, h, 0}} = \phi(x),\tilde{\tau}(x)|_{\overline{Conf}_{2, h, 1}} = 1-\phi(x). \]
    Over  $\overline{Conf}_{2, h}$ lie the universal curves $\overline{Conf}_{2, h}^Z, \overline{Conf}_{2, h}^{\bar{Z}},$ and so forth, which are the union of the corresponding universal curves lying over $\overline{Conf}_{2, h, 0}$ and $\overline{Conf}_{2, h, 1}$. 

    The chosen continuation data $(H_{ab}, J_{ab})$ in fact define smooth maps 
    \[H_h: \partial \overline{Conf}_{2, h}^{\bar{\R}} \times M \to \R, \]
    \[J_h: \partial \overline{Conf}_{2, h}^{\bar{\R}}\to \mathcal{J}(M)\]
    where $\mathcal{J}(M)$ is the space of almost complex structures on $M$. These maps are defined by requiring that on $\partial (\overline{Conf}_{2, h}^{\bar{\R}})_{\tilde{\tau}^{-1}(0)}$, $H_h$ is $H_{12}$ on the first component and $H_{24}$ on the second component, under the diffeomorphisms sending the corresponding component of the universal curve to $Z$ with the marked coordinate sent to $s=0$. Similarly, on $\partial(\overline{Conf}_{2, h}^{\bar{\R}})_{\tilde{\tau}^{-1}(1)}$, we require that $H_h$ agrees with $H_{13}$ on the first component and $H_{34}$ on the second component under the corresponding diffemorphisms with $Z$. Finally we define $J^h$ in the same way, replacing $H^h$ with $H^h$ and $H_{ab}$ with $H^{ab}$ in the discussion above.

    The data of the continuation homotopy datum $(H_h, J_h)$ \emph{is the data of a smooth extension of the maps defined above to the domains} $ \overline{Conf}_{2, h}^{\bar{Z}} \times M$. 
\end{itemize}
We can symbolically describe this data by the diagram 
\[ \begin{tikzcd} H_1 \ar[r, "H_{12}"] \ar[d, "H_{13}"] & H_2 \ar[d, "H_{24}"] \ar[dl, Leftrightarrow] \\ 
H_3 \ar[r, "H_{43}"]& H_4.\end{tikzcd}\]
where the double arrow corresponds to $H_h$. 

From this data, one can produce a pre-flow category $\CC(H_h, J_h)$ with objects the disjoint union of the objects of $\CC(H_a, J_a)$ for $a=1,2,3,4$, and such that the full subcategories on the objects of $\CC(H_a, J_a)$ and $\CC(H_b, J_b)$ for $(a,b) \in \{(1,2), (1,3), (2,4), (3,4)\}$ agree with $\CC(H_{ab}, J_{ab})$, by setting, for $x_- \in Fix(H_1)$ and $x_+ \in Fix(H_4)$, 
\[\CC(H_h, J_h)(x_-, x_+) = \overline{\{(u, a), u: (\overline{Conf}_{2, h}^{Z})_a \to M, a\in \overline{Conf}_{2, h} : \partial_s u + J_h(a, s) \partial_t u = \Grad H_h^(a), \lim_{s \to \pm \infty} u(s, t) = x_\pm(s),\}}\]
where for $a \in \tilde{\tau}^{-1}(\{0, 1\})$, the condition is imposed on each component of the universal curve, and two boundary conditions are imposed at the two opposite ends of the two components of the map from the universal curve in the natural way; and the overline denotes the Gromov-Floer compactification of this parameterized moduli space.

In this section, we explain how to adapt the construction of Section \ref{sec:continuation-maps-geometry} to produce a flow category lift of $\CC(H_h, J_h)$, and a virtual smoothing of this flow category lift, in a way that is compatible with given virtual smoothings of $\CC(H_{ab}, H_{ab})$ and also with $k$-iteration operations which relate $(H_a, J_a)$ to $(H^{\#k}_a, J^{\#k}_a)$
We will follow the general strategy of \cite{rezchikov2022integral} to produce this virtual smoothing. We will explain the geometric setup for this virtual smoothing in some detail; however, because the proofs of the required properties are almost identical to those given in Section \ref{sec:continuation-maps-geometry}, we will only write out a summary of the argument, with few comments indicating small modifications that may need to be made from proofs given earlier. 

\paragraph{Integralization data for continuation homotopies. }

Note that continuation homotopy data $(H_h, J_h)$ define a map 
\[ \mathcal{A}_h: \overline{Conf}_{2, h}^{\bar{\R}} \times LM \to \R.\]

using the modified actions from $\bar{\mathfrak{A}}_{12}$ and $\bar{\mathfrak{A}}_{24}$ over $\tilde{\tau}^{-1}(0)$ and the modified actions from $\bar{\mathfrak{A}}_{13}$ and $\bar{\mathfrak{A}}_{34}$ over $\tilde{\tau}^{-1}(1)$.

\begin{definition}
\label{def:integralization-datum-continuation-homotopy}
    An \emph{integralization datum} for $\CC(H_h, J_h)$ is a tuple $(\mathcal{A}_h, \mathfrak{A}_{h}, \{A^1_c, A^2_c\})$ such that 
    \begin{enumerate}[(A)]
        \item There exists a function $c: \overline{Conf}_{2, h}^{\bar{\R}} \to \R$  which takes constant values near the nodes and which has negative nonzero derivative at all $s$ coordinates for which $\partial_s H_h$ is nonzero, such that
        \item The quantity $\mathcal{A}_h$ is a function 
            \[ \mathcal{A}_h: \overline{Conf}_{2, h}^{\bar{\R}} \times LM \to \R \]
            given by 
            \[ \mathcal{A}_h:  = D \mathcal{A}_{H_h} + c(s)\]
            for some constant $D > 0$; in addition, 
    \item For every pair of objects $x>y$ of $\CC(H_h, J_h)$ with $x \in \CC(H_1, J_1)$ and $y \in \CC(H_4, J_4)$, 
    for every $(u, a) \in \CC(H_h, J_h)(x, y)$, and every $s \in (\overline{Conf}_{2, h}^{\R})_a$, we have that
    \begin{equation}
    \label{eq:modified-action-is-never-constant}
        \frac{d}{ds} \mathcal{A}_h(u_s) < 0 \text{ for all } s. 
    \end{equation} 
    \item The function $\mathcal{A}_h$ defines  four associated functions 
    \[ \mathcal{A}_{h, ab}: \R \times LM \to \R \]
    for $(a,b) \in \{(1, 2), (1,3),(2,4), (1,4)\}$, by restricting $\mathcal{A}_h$ in the components of the universal curves over $(\overline{Conf}_{2,h}^{Z})_0$ and $(\overline{Conf}_{2,h}^{Z})_0$, so $\mathcal{A}_{h, ab}$ is related to $\mathcal{A}_{H_{ab}}$ by an affine-linear reparameterization of its codoomain. We require that for each $(a,b)$ as above, the tuple $(\mathcal{A}_{h, ab}, \mathfrak{A}_{h, ab}, \{A^{j_{ab}}_c\})$ is an integralization datum for $\CC(H_{ab}, J_{ab})$ where 
    \begin{equation}
        \label{eq:disjoint-actions-convenience}
        \begin{gathered}
        j_{12} = j_{13} = 1, j_{24} = j_{34} = 2, \\
    \mathfrak{A}_{12} = \mathfrak{A}'_1 \cup \mathfrak{A}'_2, \mathfrak{A}_{13} = \mathfrak{A}'_{1} \cup \mathfrak{A}'_{2} \cup \mathfrak{A}'_{3},\\
     \mathfrak{A}_{24} = \mathfrak{A}'_2\cup \mathfrak{A}'_3 \cup \mathfrak{A}'_4, \mathfrak{A}_{34} = \mathfrak{A}'_3 \cup \mathfrak{A}'_4. 
     \end{gathered}
     \end{equation}
    \end{enumerate}
\end{definition}

\begin{lemma}
\label{lemma:integralization-data-exist-continuation-homotopy}
    Suppose that we are given integralization data $\bar{\mathfrak{A}}_{ab}= (\mathcal{A}_{ab}, \mathfrak{A}_{ab}, \{\rho^x\})$ for $\CC(H_{ab}, J_{ab})$ with $(a, b) \in \{(1,2), (1,3), (2,4), (3,4)\}$, such that the induced pairs of integralization data on $\CC(H_i, J_i)$ agree for each of $i=1, 2,3, 4$. Write
    \[ \mathcal{A}_{ab} = D \mathcal{A}_{H_ab} + c_{ab}, \lim_{s \to - \infty} c_{ab} = C_a, \lim_{s \to + \infty} = C_b \]
    and
    \[\mathfrak{A}_{ab} = (\mathfrak{A}_a + C_a) \sqcup (\mathfrak{A}_b + C_b)\] 
    for the corresponding decompositions of marked actions. Suppose further that the following conditions are satisfied: 
    \begin{equation}
    C_a - C_b > \left(\max_{x \in \CC(H_c, J_c)} D \mathcal{A}_{H_c}(x) +C_c \right) -   \left( \min_{x \in \CC(H_a, J_a)} D \mathcal{A}_{H_a}(x) + C_a\right)
    \end{equation}
    \begin{equation}
     C_a - C_c > \left(\max_{x \in \CC(H_b, J_b)} D \mathcal{A}_{H_b}(x) +C_b \right) -   \left( \min_{x \in \CC(H_a, J_a)} D \mathcal{A}_{H_a}(x) + C_a\right)
    \end{equation}
    \begin{equation}
        C_a - C_d > \sup_{\tau \in (0,1)} \int_{(\overline{Conf}_{2,h}^Z)_{\tilde{\tau}^{-1}(\tau)}} \max_{(x,t) \in M \times S^1}|\partial_s H_h(s, x, t)| ds.
    \end{equation}
    Then there exists an integralization datum $\bar{\mathfrak{A}}_h$ such that the associated integralization data for each $\CC(H_{ab}, J_{ab})$ are $\bar{\mathfrak{A}}_{ab}$. 
\end{lemma}
\begin{proof}
   The first two conditions above mean that the extension of $\mathfrak{A}_{ab}$ by the new marked actions $\mathfrak{A}_{c} + C_b$ to the marked actions is still an integralization datum from $\CC(H_{ab}, J_{ab})$, and similarly with $b$ swapped with $c$ everywhere. The only challenge is thus to define the function $c$ (Definition \ref{def:integralization-datum-continuation-homotopy} (A)); this can be defined by the gluing construction for points projecting to points near the boundary of $\overline{Conf}_{2,h}$, and the last condition in the lemma allows us to extend $c$ to the remaining locus.  
\end{proof}

\begin{lemma}
\label{lemma:integralization-data-exist-integralization-homotopy-cyclotomic}
    Suppose we are given  $\bar{\mathfrak{A}}_h = (D \mathcal{A}_{H_h} + c, \mathfrak{A}_h, \
    \{A^1_c, A^2_c\})$ for $\CC(H_h, J_h)$. Write $C_a, C_b, C_c, C_d$ and so forth as in the statement of Lemma \ref{lemma:integralization-data-exist-continuation-homotopy}, and suppose that 
    \begin{itemize}
        \item for each $(\alpha, \beta)$ in $((a,b), (a,c), (b,c), (b,d))$, the data $(C_\alpha, C_\beta, A^{j_{ab}}_c)$ satisfy the condition required in Lemma \ref{lemma:cyclotomic-compatibility-for-integralization-data-continuation-map} with respect to the integralization datum $\bar{\mathfrak{A}}_{ab}$ associated by $\mathfrak{A}_h$ to $\CC(H_{ab}, J_{ab})$, and 
        \item The condition
        \[ C_a - C_d > DC_{H_a} + DC_{H_d} + max_{\tau \in (0,1)} \int_{(\overline{Conf}_{2,h}^Z)_{\tilde{\tau}^{-1}(\tau)}} \max_{(x,t) \in M \times S^1}|\partial_s H(s, x, t)| ds,\]
        is satisfied
        where $C_{H_a}$ and $C_{H_d}$ are the constants associated by Lemma \ref{lemma:apriori-bound-on-action-values} to $H_a$ and $H_d$, respectively. 
    \end{itemize}
    Then there exists a set $\mathfrak{A}^k_{h, new} \subset \R \setminus k \mathfrak{A}_h$ and such that 
    \[ \bar{\mathfrak{A}}_h^k = (D \mathcal{A}_{H_h^{\#k}} + k c, k \mathfrak{A}_h \cup \mathfrak{A}^k_{h, new}, \{kA^1_c, kA^2_c\}) \]
    is an integralization datum for $\CC(H_h^{\#k}, J_h^{\#k})$. 
\end{lemma}
\begin{proof} The argument is essentially identical to that of Lemma \ref{lemma:cyclotomic-compatibility-for-integralization-data-continuation-map}. 
\end{proof}

\paragraph{Stratifying the pre-flow category $\CC(H_h, J_h)$ using integralization data}
We let $\bar{\mathfrak{A}}$ be an integralization datum for $\CC(H_h, J_h)$; we will keep the notation as in Definition \ref{def:integralization-datum-continuation-homotopy}. We set $\underline{\mathfrak{A}} = \mathfrak{A} \cup \{ A^1_c, A^2_c\}$. As before we have $\underline{\mathfrak{A}}(x) = \{A \in \underline{\mathfrak{A}} : A < \mathcal{A}(x)\}$, where $\mathcal{A}$ is the augmented action underlying $\bar{\mathfrak{A}}$. The integral action for $\CC(H_h, J_h)$ is
\[ \bar{A}(x) = \# \underline{\mathfrak{A}}(x)-1\]
as usual; the definition of the  integral action is the same as before but with the new notation.

\paragraph{Parameter spaces and perturbation data for continuation homotopy.}
Let us assume that $x \in Fix(H_1)$ and $y \in Fix(H_4)$. 
We now define analogs of the spaces $\overline{Conf}^c_{\mathfrak{A}(x,y)}$ defined previously which are appropriate to this new setting. Let $(\overline{Conf}_{\mathfrak{A}(x,y)})^{\R^2}$ be the fiber product of $(\overline{Conf}_{\mathfrak{A}(x,y)})^\R$ over $(\overline{Conf}_{\mathfrak{A}(x,y)})$ with itself. Similarly, to before, we
we define a open subset $\overline{Conf}^o_{\underline{\mathfrak{A}}(x,y)} \subset (\overline{Conf}_{\mathfrak{A}(x,y)})^{\R^2}$, and then an alternative compactification $\overline{Conf}^c_{\underline{\mathfrak{A}}(x,y)}$ of this set. The set $(\overline{Conf}_{\underline{\mathfrak{A}}(x,y)})^{\R^2}$ is an element $a$ of $(\overline{Conf}_{\mathfrak{A}(x,y)})$ together with two additional marked  $s$ coordinates, which we will call $s_{A_c^1}$ and $s_{A_c^2}$, on the universal curve over $a$. The subset $\overline{Conf}^o_{\underline{\mathfrak{A}}(x,y)}$ consists of those elements such that $(a, s_{A^1_c}, s_{A^2_c})$ satisfy the following properties:

\begin{enumerate}[(A)]
    \item The element $a$ lies in a stratum where there is a component of the universal curve containing marked points labeled by actions from both $\mathfrak{A}_1$ and from $\mathfrak{A}_2 \cup \mathfrak{A}_3$ (this component must be unique if it exists)and $s_{A^1_c}$ is on that component;
    \item Similarly there is a (unique) component of the universal curve containing marked points from $\mathfrak{A}_2 \cup \mathfrak{A}_3$ and from $\mathfrak{A}_4$, and $s_{A^2_c}$ also is on that component;
    \item Finally, we require that $s_{A^1_c}$ lies to the left of $s_{A^2_c}$. 
\end{enumerate}

We now define a stratification of $\overline{Conf}^o_{\underline{\mathfrak{A}}(x,y)}$ making it into a $\langle \underline{\mathfrak{A}}(x,y)'\rangle$-manifold. Recall as before that we identify $\underline{\mathfrak{A}}'(x,y)$ with the regions in between the elements of $\underline{\mathfrak{A}}(x,y)$; given an element $\tilde{a} \in \overline{Conf}^o_{\underline{\mathfrak{A}}(x,y)}$  we put it in the stratum $S(\tilde{a}) =\underline{\mathfrak{A}}'(x,y) \setminus S^c(\tilde{a})$, where $S^c(\tilde{a})$ contains the elements corresponding to the regions where there are break-points in between marked coordinates labeled by elements of $\mathfrak{A}(x,y)$. (See Figure \ref{fig:parameter-spaces-for-continuation-map}.)

We now modify the space $\overline{Conf}^o_{\underline{\mathfrak{A}}(x,y)}$ by ``doubling'' it in a way analogous to the fashion that $\overline{Conf}_{2, h}$ is produced from $\overline{Conf}_2$. Namely, by forgetting all the marked points except for the markings $s_{A_c^1}$ and $s_{A_c^2}$ and collapsing all the remaining unsable components, we get a map 
\[\pi_h: \overline{Conf}^o_{\underline{\mathfrak{A}}(x,y)} \to \overline{Conf}_2.  \]
As before we will define two halves
\[  \overline{Conf}^{h,0}_{\underline{\mathfrak{A}}(x,y)} = \pi_h^{-1}(\tilde{\tau}^{-1}([0, 1/2]))\]
\[ \overline{Conf}^{h,1}_{\underline{\mathfrak{A}}(x,y)} = \pi_h^{-1}(\tilde{\tau}^{-1}([0, 1/2])),\]  \[\overline{Conf}^{h,\partial}_{\underline{\mathfrak{A}}(x,y)} = \pi_h^{-1}(\tilde{\tau}^{-1}(1/2))\]
and define, as a topological $\langle k \rangle$-manifold,
\[ \overline{Conf}^h_{\underline{\mathfrak{A}}(x,y)} = \overline{Conf}^{h,0}_{\underline{\mathfrak{A}}(x,y)} \cup_{\overline{Conf}^{h,\partial}_{\underline{\mathfrak{A}}(x,y)}} \overline{Conf}^{h,1}_{\underline{\mathfrak{A}}(x,y)}\]
We note now that $\overline{Conf}^{h,\partial}_{\underline{\mathfrak{A}}(x,y)}$ is in fact a codimension $1$ $\langle k \rangle$-submanifold, i.e. there are local charts to $\R^n \times \R_+^{[k] \setminus S}$ covering an open neighborhood of $\overline{Conf}^{h,\partial}_{\underline{\mathfrak{A}}(x,y)}$ such that the latter takes the form 
\[\R^{n-1} \times \R^{[k] \setminus S}_+ \subset \R^{n} \times \R^{[k] \setminus S}_+ . \]

In order to put a smooth structure on $\overline{Conf}^{h}_{\underline{\mathfrak{A}}(x,y)}$ we must choose collar neighborhoods 
\begin{equation}
\label{eq:collar-neighborhood-of-boundary-in-continuation-homotopy-parameter-space}
\phi_i: (1/2-\epsilon, 1/2] \times \overline{Conf}^{h,i}_{\underline{\mathfrak{A}}(x,y)},\; i=0,1,
\end{equation}
sending $\{1/2\} \times \overline{Conf}^{h,\partial}_{\underline{\mathfrak{A}}(x,y)}$ to $\overline{Conf}^{h,\partial}_{\underline{\mathfrak{A}}(x,y)}$. We will choose $\phi_1=\phi_2$ such that $\pi_1\phi_i^{-1} = \tilde{\tau}\pi_h$ on the image of $\phi_i$, where $\pi_1$ is projection to the first factor. To define $\phi_i$, we simply think of a point $a$ in $\overline{Conf}^{h,\partial}_{\underline{\mathfrak{A}}(x,y)},$ as taking the curve $\pi_h(a) \in (\overline{Conf}_2^Z)_{1/2}$, which has two marked $s$-coordinates, which we identify with $s_{A_1}$ and $s_{A_2}$, and adding on some additional marked coordinates; fixing the locations of these new marked coordinates and moving the originally chosen marked coordinates $s_{A_c^1}$, $s_{A_c}^2$ coordinates apart defines $\phi_i$. The two collars $\phi_i$ glue together to a smooth map 
\begin{equation}
    \label{eq:full-collar-neighborhood-of-boundary-in-continuation-homotopy-parameter-space}
\phi: (1/2-\epsilon, 1/2+\epsilon) \times \overline{Conf}^{h,\partial}_{\underline{\mathfrak{A}}(x,y)} \to \overline{Conf}^{h}_{\underline{\mathfrak{A}}(x,y)}
\end{equation}
which is a diffeomorphism onto its image by definition.

Now, for $x>y$ $\CC(H_t, J_t)$ with $x \in Fix(H_1)$ and $y \in Fix(H_4)$, we define 
\[ \paramspace(x,y) = \overline{Conf}^h_{\underline{\mathfrak{A}}(x,y)}. \] 

Otherwise we define the spaces $\paramspace(x,y)$ as before. We define the universal curves over over $\overline{M}(x,y)$ to be the pullback of the corresponding universal curves from the first factor. One has maps  $\iota_{xzy}$ covered by maps $\iota'_{xzy}$ identifying  universal curves as in \eqref{eq:identification-of-universal-curve-domains}, with the clarification that if $x \in Fix(H_1)$, $z \in Fix(H_2)$ and $y \in Fix(H_4)$ then one includes $\paramspace(x,z) \times \paramspace(z,y)$ into $\paramspace(x,y)$ via the inclusion into $\overline{Conf}^{h, 0}(x,y)$; and otherwise if $x \in Fix(H_1)$, $z \in Fix(H_3)$ and $y \in Fix(H_4)$ one includes $\paramspace(x,z) \times \paramspace(z,y)$ into $\paramspace(x,y)$ via $\overline{Conf}^{h, 1}(x,y)$. Systems of perturbation data and their compatibility are defined as earlier but with the new notation.

We have maps $\pi_h: \paramspace(x,y)^Z \to \overline{Conf}_{2, h}^{\bar{\R}}$ which come from forgetting all the markings except for $s_{A_c^1}$ and $s_{A_c^2}$, collapsing the components that are now unstable to points, and subsequently quotienting by the $S^1$ action.  Given a point on $\paramspace(x,y)^Z$ we can write $\tilde{s}$ for its image under $\pi_h$.

We define
\begin{equation}
\label{eq:thickenings-for-floer-contiuation-homotopy}
    T(x_-, x_+) = \left\{ \begin{array}{c|r} 
    (u,a,v) : v \in V^1(x_-, x_+), &  \partial_s u(s,t) + J(\tilde{s})_t \partial_t u(s,t) -\Grad H^h(\tilde{s}, u(s,t))= \lambda_V(((s,t), a), v),\\
 a \in \paramspace(x,y), u: \paramspace(x,y)^Z_a \to M &   \lim_{s \to -\infty} u_1(s, \cdot) = x_-, \lim_{s \to \infty} u_r(s, \cdot) = x_+\\ 
   u = u_1 \cup \ldots \cup u_r  & \lim_{s \to + \infty} u_j(s, \cdot) = \lim_{s \to - \infty} u_{j+1}(s, \cdot), j=1, \ldots, r-1
 \end{array}\right\}.
\end{equation}

We define the thickenings, obstruction bundles, obstruction sections, as in the section on the continuation maps but now with the new notation.

\paragraph{Properties.}

The proofs of the next two propositions are essentially identical to the proofs of the corresponding propositions in Section \ref{sec:continuation-maps-geometry}, and are thus omitted. 
\begin{proposition}
\label{prop:perturbation-data-continuation-homotopy-exist}
    Suppose that we have compatible regular systems of perturbation data $\mathfrak{P}_i$ over $\CC(H_i, J_i)$, and compatible regular systems of perturbation data $\mathfrak{P}_{ab}$ over $\CC(H_{ab}, J_{ab})$ for $(a,b) \in \{(1,2), (2,4), (1,3), (3,4)\}$ such that the restriction of $\mathfrak{P}_{ab}$ to $\CC(H_a, J_a)$ is a stabilization of $\mathfrak{P}_a$ by an $E^s$-parameterization equipped with trivial perturbation data, and similarly for its restriction to $\CC(H_b, J_b)$; and such that there is an integralization datum $\mathfrak{A}$ for $\CC(H_h, J_h)$ with underlying integralization data for $\CC(H_{ab}, J_{ab})$ agreeing with those underlying $\mathfrak{P}_{ab}$. Then there exists a regular perturbation datum $\mathfrak{P}_h$ with underlying integralization datum $\mathfrak{A}_h$ such that its restriction to $\CC(H_{ab}, J_{ab})$ is a stabilization of $\mathfrak{P}_{ab}$ by an $E^s$-parameterization equipped with trivial perturbation data. $\blacksquare$
\end{proposition}

\begin{proposition}
\label{prop:continuation-homotopy-compatible-smoothings}
 The category  $\CC'(H_h, J_h)$ defined using $\mathfrak{P}$ produced by Proposition \ref{prop:perturbation-data-continuation-homotopy-exist},  is topologically smooth. 

 Moreover, suppose we have $E^s$-parameterizations $\mathcal{E}_a$ over $\CC(H_a, J_a)$ and $\mathcal{E}_{ab}$ over $\CC'(H_{ab}, J_{ab})$, such that $(\CC'(H_{ab}, J_{ab})_{\mathcal{E}_{ab}})|_{\CC(H_{a}, J_a}$ is a semi-free stabilization (by some $E$-parameterization) of $\CC'(H_{a}, J_a)_{\mathcal{E}_a}$ (and similarly with $(H_a, J_a, \mathcal{E}_a)$ replaced with $(H_b, J_b, \mathcal{E}_b)$); and moreover such that there are smooth structures on $\CC'(H_{ab}, J_{ab})_{\mathcal{E}_{ab}}$ and $\CC'(H_{a}, J_a)_{\mathcal{E}_a}$ with the previously mentioned isomorphisms being smooth. Then there is an $E^s$-parameterization $\mathcal{E}$ such that $\CC'(H_h, J_h)_{\mathcal{E}}$ has a smooth structure, and its restriction to $\CC(H_{ab}, J_{ab})$ is a stabilization by some $E^s$-parameterization of $(\CC'(H_{ab}, J_{ab})_{\mathcal{E}_{ab}})$ with the corresponding smooth structure.
\end{proposition}

Similarly, using arguments which can be copied verbatim from the previous sections, we can choose perturbation data in a manner that is compatible with the $C_k$-symmetries on continuation homotopy data.
\begin{proposition}
\label{prop:cyclotomic-compatibility-continuation-homotopy}
    The analog of Proposition \ref{prop:cyclotomic-compatibility-all-perturbation-data-choices} and Lemma \ref{lemma:smoothings-cyclotomic-compatibility} hold with $(H_s, J_s)$ replaced with $(H_h, J_h)$, with the condition involving Lemma \ref{lemma:cyclotomic-compatibility-for-integralization-data-continuation-map} replaced with the corresponding condition involving Lemma \ref{lemma:integralization-data-exist-integralization-homotopy-cyclotomic}. $\blacksquare$
\end{proposition}

\section{Index theory}
\label{sec:index-theory}
In order to produce the framings of flow categories defined using Floer theory that are needed for the constructions of Section \ref{sec:flow-categories}, and to verify that the same method can be used to produce the \emph{equivariant framings} of Section \ref{sec:equivariant-flow-categories}, we will make extensive use of family index theory. As such, in this section, we will review the aspects of equivariant family index theory that we will need to subsequently produce framings in Section \ref{sec:producing-framings}. 
\subsection{Review of $K$ theory}
\label{sec:k-theory-review}
There is an extraordinary cohomology theory called $KO$ theory which defines a contravariant functor 
\[ X \mapsto KO^*(X)\]
from the homotopy category $hTop$ to the category of graded algebras. On spaces $X$ with the homotopy type of a finite $CW$ complex, the degree-zero part of this functor has an elementary description in terms of \emph{virtual vector bundles}. Namely, let $X$ be a paracompact Hausdorff space, and let $Vect(X)$ denote the set of isomorphism classes of finite-dimensional real vector bundles over $X$. Then, if $X$ has the homotopy type of a finite $CW$ complex, we have 
\begin{equation}
    \label{eq:def-virtual-vector-bundles}
    KO^0(X) = \{ (V_1, V_2) \in Vect(X) \times Vect(X)\}/((V_1, V_2) \sim (V_1 \oplus W, V_2 \oplus W) \text{ for } W \in Vect(X)) 
\end{equation}
where we take the quotient by the minimal equivalence relation generated by the relation above. One writes the element of $KO^0(X)$ represented by $(V_1, V_2)$ as $[V_1] - [V_2]$. 

Recall now the \emph{Atiyah-Janich theorem}:
\begin{theorem}[Atiyah\cite{Atiyah1994-ar}, J\"anich (thesis)]
\label{thm:atiyah-janich}
Let $X$ be a finite $CW$ complex. Then 
    \begin{equation}
        \label{eq:atiyah-janich}
        KO^0(X) = [X, \Fred \mathcal{H}]
    \end{equation} 
    where on the right, $\Fred \mathcal{H}$ is the space of Fredholm operators acting on a separable Hilbert space $\mathcal{H}$, equipped with the operator norm topology. 
\end{theorem}

The ideas behind the proof of this result are used repeatedly when studying moduli spaces arising in Floer homology. As such, we will recall the method of proof. We first recall the fundamental theorem of Kupier:
\begin{theorem}[Kupier \cite{KUIPER196519}]
Let $\mathcal{H}_\R$ and $\mathcal{H}_\C$ be real and complex infinite-dimensional separable Hilbert spaces, respectively. Let $GL(\mathcal{H}_\R)$ and $GL(\mathcal{H}_\C)$ be the spaces of  real- and complex-linear bounded automorphisms with bounded inverses, respectively. Then, writing 
\begin{equation}
\label{eq:orthogonal-automorphisms-definition}
    O(\mathcal{H}_\R) \subset GL(\mathcal{H}_\R), U(\mathcal{H}_\C) \subset GL(\mathcal{H}_\C)
\end{equation}
    for inclusions of real and complex-linear isometries, we have that every space in \eqref{eq:orthogonal-automorphisms-definition} above is contractible. 
\end{theorem}

\begin{proof}[Proof sketch of Thm \ref{thm:atiyah-janich}, see \cite{Atiyah1994-ar}]
We write $\underline{\R^n}$ for the trivial $\R^n$ bundle over $X$. Similarly, write $\underline{\mathcal{H}}$ for the trivial Hilbert bundle $\mathcal{H} \times X$ over $X$.

Every virtual vector bundle $[V] - [W]$ is equivalent to virtual vector bundle $[V'] - [\underline{\R^n}]$ for some $n$. Embed $V'$ into the Hilbert bundle $\mathcal{H} \times X$ as a subbundle; using Kupier's theorem, one can then produce a bundle isomorphism  
\[\bar{\mathcal{F}}: \underline{\mathcal{H}}/V \to \underline{\mathcal{H}'},\] where $\mathcal{H}' \subset \mathcal{H}$ is a complement of a fixed embedding $W \subset \mathcal{H}$. One defines the family of Fredholm operators associated to $[V] - [W]$ to be 
\[ \mathcal{F}(x) = \mathcal{H} \to \mathcal{H}/V_x \xrightarrow{\bar{\mathcal{F}}} \mathcal{H}' \to \mathcal{H}.\]
One then shows that the homotopy class of this family is independent of the choices made via Kupier's theorem, and that this class is independent of the representative $(V, W)$ of the original class via an explicit homotopy of Fredholm opreators which adds an $x$-independent trivial vector space to the kernel and the cokernel of each operator $\mathcal{F}_x$.

Let $\mathcal{F}: X \to Fred \mathcal{H}$ be a map, and write $\mathcal{F}_x: \mathcal{H} \to \mathcal{H}$ for the operator $\mathcal{F}(x)$. 
Then, for every $x \in X$, there exists a finite dimensional subspace $W_x \subset \mathcal{H}$ such that $W_x$ is transverse to $Im(\mathcal{F}_x)$, i.e. that $W_x + Im(\mathcal{F}_x) = \mathcal{H}$. But then $W_x$ is also transverse to $Im(\mathcal{F}_y)$ for all $y$ in some open neighborhood $U_x$ of $x$. Since $X$ is compact, there exists a finite collection of points $x_1, \ldots, x_r \in X$ such that $W := W_{x_1} + \ldots + W_{x_r}$ is transverse to $Im(\mathcal{F}_y)$ for all $y \in X$. As such, there is a continuous family of Fredholm operators
\[ \mathcal{F}^W: X \to Hom(\mathcal{H} \oplus W, \mathcal{H})\]
\[ \mathcal{F}^W_x(v, w) = \mathcal{F}_x(v) + w\]
which are all surjective. Therefore, the dimensions of $\ker \mathcal{F}^W_x$ are independent of $X$ and are given by the index $Ind \mathcal{F}^W_x$ for any $x$. As such, there is a vector bundle $\ker \mathcal{F}^W$ over $X$ with the fiber $\ker \mathcal{F}^W_x$ given by the kernel of the corresponding operator. We set the element of $KO^0(X)$ associated to $\mathcal{F}$ to be $[\ker \mathcal{F}^W_x] - [\underline{W}]$.

Given any two choices $W_1, W_2$ as above with $W_1 \subset W_2$, one sees that 
\[ \ker \mathcal{F}^{W_2} = \ker \mathcal{F}^{W_1} \oplus \underline{W_2/W_1}. \]
As such, the assignment above gives a well defined element of $KO^0(X)$ for any continuous map $\mathcal{F}$. To show that this assignment does not change upon a homotopy of maps $\mathcal{F}^t$, one uses the same argument to build a vector bundle $\ker \mathcal{F}^{W', t}$ over $X \times [0,1]$ for some sufficiently large finite-dimensional $W' \subset \mathcal{H}$ associated to such a homotopy of maps, and then uses the fact that this implies that $\ker \mathcal{F}^{W', t}|_{X \times \{0\}} \simeq \ker \mathcal{F}^{W', t}|_{X \times 1}$.
\end{proof}

There is another space, $\Z \times BO$, such that for spaces with the homotopy type of a finite $CW$ complex we have 
\begin{equation}
\label{eq:ko-theory-classifying-space}
KO^0(X) \simeq [X, \Z \times BO].
\end{equation}
We briefly review the construction of this space in Appendix \ref{sec:classifying-spaces-for-k-theory}. More generally, one defines $KO^0(X) = [X, \Z \times BO]$ for \emph{any space} with the homotopy type of a CW complex; this set turns out \emph{not} to agree with the set of virtual vector bundles on $X$ for infinite CW complexes $X$ \cite{anderson1968k}. Nonetheless, Theorem \ref{thm:atiyah-janich} still establishes a homotopy equivalence $\Z \times BO \simeq \Fred \mathcal{H}$ intertwining the isomorphisms  \eqref{eq:atiyah-janich} and \eqref{eq:ko-theory-classifying-space}, as we show in Lemma \ref{lemma:atiyah-janich-in-general}. 

\subsection{Index theory for paths of self-adjoint operators}
\label{sec:index-theory-for-paths}
The Atiyah-Janich theorem is only the first of a richer series of idenfications which identify classifying spaces of the cohomology theory $KO^*$, which were computed by Bott \cite{bott1956application}, with spaces of operators on $\mathcal{H}$.  In particular, we have
\begin{theorem}[Atiyah-Singer \cite{atiyah1969index}]
\label{thm:atiyah-singer-self-adjoint}
For $X$ compact Hausdorff, we have 
\[ KO^1(X) = [X, U/O] = [X, \mathcal{A}(\mathcal{H})]\]
where 
\begin{equation}
    U/O = \colim_n U(n)/O(n)
\end{equation} and $\mathcal{A}(\mathcal{H})$ is the connected component of the space $\mathcal{A}^{all}(\mathcal{H})$ of real bounded self-adjoint operators on a separable Hilbert space $\mathcal{H}$ containing operators with infinitely-many positive and infinitely-many negative eigenvalues. The other two components are $\mathcal{A}^+(\mathcal{H})$ and $\mathcal{A}^-(\mathcal{H})$, consisting of operators with finitely-many negative and finitely-many positive eigenvalues, respectively; each of $\mathcal{A}^\pm(\mathcal{H})$ is contractible.
\end{theorem}

As such, since $\Omega U/O \simeq \Z \times BO$ by Bott periodicity, one expects that there is a natural map 
\begin{equation}
    \Omega \mathcal{A}(\mathcal{H}) \to \Fred \mathcal{H}.
\end{equation}
Such a map indeed exists, and figures heavily in Floer-theoretic constructions. It is helpful to recall there the following theorems:
\begin{theorem}[Corollary of \cite{palais1966homotopy}]
\label{thm:banach-manifold-cw-complex}
Any Banach manifold has the homotopy type of a CW complex.
\end{theorem}


In fact, an equivariant version of the above is known:
\begin{theorem}
Any smooth $G$-Banach manifold has the homotopy type of a CW complex.
\end{theorem}
\begin{proof}
\label{thm:g-banach-manifold-cw-complex}
    This follows from \cite{antonyan1985equivariant}, 
    \cite{antonyan2007homotopy}, and  \cite{antonyan2009equivariant}.
\end{proof}

The following is standard:
\begin{theorem}
    The spaces $Fred(\mathcal{H}), \mathcal{A}(\mathcal{H})$ are smooth Banach manifolds for a real or complex Hilbert space $\mathcal{H}$. If $\mathcal{H}$ is a $G$-representation for a compact Lie group $G$, then all of these spaces are smooth $G$-Banach manifolds. 
\end{theorem}
\begin{proof}
    The space of $B(\mathcal{H})$ of bounded operators on $H$ is a Banach space under the operator norm; $Fred(\mathcal{H})$ is an open subset of the latter \cite{lax}.  and $\mathcal{A}(\mathcal{H})$ is manifestly a Banach submanifold of $Fred(\mathcal{H})$. The smoothness of the $G$-action in the operator norm topology is obvious.
\end{proof}

Fix two operators $A_\pm \in \mathcal{A}(\mathcal{H})$, and constants $C, D > 0$. Let the model for $\Omega \mathcal{A}(\mathcal{H})$ be 
\begin{equation}
\Omega_{A_-, A_+} \mathcal{A}(\mathcal{H}) = \{ A_t: \R \to \mathcal{A}(\mathcal{H}): |A_t - A_\pm| \leq C e^{\mp D t} \text{ as } t \to \pm \infty, A_t \text{ is smooth }. \}. 
\end{equation}
There are natural separable Hilbert spaces 
\begin{equation}
\mathcal{H}' = L^{1,2}(\R_s, \mathcal{H}), \; \mathcal{H}'' = L^2(\R_s, \mathcal{H})
\end{equation} of $\mathcal{H}$-valued functions on $\mathcal{R}$. There is a natural continuous map
\begin{equation}
\label{eq:family-self-adjoint-index-map}
\begin{gathered}
\Omega_{A_-, A_+} \mathcal{A}(\mathcal{H}) \to Fred(\mathcal{H}', \mathcal{H}'') \\
A_t \mapsto (\partial_s + A_s) \in Fred(L^{1,2}(\R_s, \mathcal{H}), L^2(\R_s, \mathcal{H})).
\end{gathered}
\end{equation}
The operators in the image of this map are Fredholm by Robbin-Salamon \cite{robbin1995spectral} (see also \cite{rabier2004robbin}).

\begin{remark}
Unfortunately, the author does not know any proof that the above map is a homotopy equivalence, nor that it in some sense agrees with the Bott periodicity map $\Omega U/O \simeq Z \times BO$ under some naturally defined equivalences 
\[ \mathcal{A}(\mathcal{H}) \simeq U/O, \Fred \mathcal{H} \simeq Fred(\mathcal{H'}, \mathcal{H}'') \simeq \Z \times BO \]
associated to Theorems \ref{thm:atiyah-janich} and \ref{thm:atiyah-singer-self-adjoint}. First of all $\mathcal{H'} \neq \mathcal{H}$, so the operators $A_s$ are not \emph{bounded} operators, and \cite{atiyah1969index} does not actually apply; secondly, even if we have $\mathcal{H}' = \mathcal{H}$, the comparison with the Bott periodicity map proven in \cite{atiyah1969index} is somewhat different. It seems to the author that in the literature (e.g. \cite{cohen2007floer}) some comparison of the above kind is implicitly used. 
\end{remark}

Below, while we will not provide a comparison with a more classically defined Bott-periodicity map, we will explain how to extend \cite{atiyah1969index} to the Floer theoretic setting (where the operators are unbounded) and subsequently how the map \eqref{eq:family-self-adjoint-index-map} can be used to produce compatible (and \emph{equivariant}!) framings of Floer theoretic moduli spaces. 

\subsection{Equivariant Operators and Classifying Spaces}
Now let $\mathcal{H}$ be a separable Hilbert space equipped with an orthogonal $G$-action. Suppose that $\mathcal{H}$ is the \emph{completion of a complete $G$-universe}: it contains an infinite number of pairwise orthogonal copies of the regular $G$-representation whose sum is dense in $\mathcal{H}$. 

Using $\mathcal{H}$, we can perform a number of constructions analogous to those of Section \ref{sec:index-theory}. Namely, there is a functor $KO^0_G$ on $hG-Top$ which sends spaces $X$  with the homotopy type of a finite $G$-CW complex to the set of isomorphism classes of \emph{virtual $G$-vector bundles on $X$}, which are defined like in \eqref{eq:def-virtual-vector-bundles} with $Vect$ replaces with $Vect_G$, the set of isomorphism classes of finite-dimensional $G$-equivariant vector bundles. Similarly, the space of Fredholm operators $\Fred \mathcal{H}$ is now a $G$-space, with the $G$-action given by $gA = gAg^{-1}$. The $G$-fixed points of $\Fred H$ are then the $G$-equivariant Fredholm operators. By replacing the definition of $W$ with $W := Span(G(W_{x_1} \oplus \ldots \oplus W_{x_r}))$ in the proof of Theorem \ref{thm:atiyah-janich}, one gets a functorial isomorphism
\[ [X, \Fred \mathcal{H}]_G \to KO_G^0(X)\]
on all finite $G$-CW complexes. 
Similarly, there are classifying spaces $BO \in hG-Top$ for $KO_G^0$ (see Appendix \ref{sec:classifying-spaces-for-k-theory}).

\subsection{Polarization classes}
\label{sec:polarization-classes}
As we have seen earlier, the Hessians of the symplectic action functional formally define a family of formally-self-adjoint operators \eqref{eq:hessian-operator-form} parameterized by points on $LM$, which are unbounded operators on the Hilbert spaces
\begin{equation}
\label{eq:l2-tangent-spaces}
    L^2(x^*TM) = (\mathcal{T}^{0,2}LM)_x
\end{equation}
which are densely defined on the subspaces
\begin{equation}
W^{1,2}(x^*TM) = (\mathcal{T}^{1,2}LM)_x \subset (\mathcal{T}^{0,2}LM)_x.
\end{equation}
These spaces piece together into Hilbert bundles $\mathcal{T}^{L^2}LM$ over $LM$, or at least over the homotopy equivalent subspace $L^{1,2}M \subset LM$ of $W^{1,2}$ loops in $M$, and the operators give a continuous bundle map 
\begin{equation}
\label{eq:hessians-bundle-maps}
    Hess(\mathcal{A}_H): \mathcal{T}^{1,2}LM \to \mathcal{T}^{0,2}LM.
\end{equation}

This is closely tied to the discussion of Section \ref{sec:index-theory}, where we recalled that the classifying space for $KO^1$ is the nontrivial component of the space of bounded self-adjoint operators on a separable Hilbert space. 

By Kupier's theorem on the contractibility of the space of invertible bounded linear endomorphisms $GL(\mathcal{H})$ of a separable Hilbert space $\mathcal{H}$, it is tempting to produce a class in $KO^1(LM)$ from the family of Hessians \eqref{eq:hessians-bundle-maps} by trivializing the bundle $\mathcal{T}^{0,2}LM$, and applying \ref{thm:atiyah-singer-self-adjoint}. However, this cannot be done directly as $LM$ is noncompact and the operators \eqref{eq:hessian-operator-form} are unbounded on $\mathcal{T}^{0,2}LM$, the idea is essentially correct, and can be carried forward with a few minor technical modifications. 

The first challenge is the non-compactness of $LM$; this, however, can be fixed using Palais' theorem, as described in Appendix \ref{sec:classifying-spaces-for-k-theory}.

The second challenge is to write \eqref{eq:hessian-operator-form} as a family of bounded operators. First, write 
\begin{equation}
    \mathcal{H} = L^2(S^1, \C^n), \mathcal{H}_1 = W^{1,2}(S^1, \C^n).
\end{equation}
Then the inclusion 
\begin{equation}
    \iota: \mathcal{H}_1 \subset \mathcal{H}
\end{equation}
is a bounded map of Hilbert spaces, and in fact a compact operator. There is also a bounded operator without kernel
\begin{equation}
i \partial_t + 1: \mathcal{H}_1 \to \mathcal{H}
\end{equation}
As such, $(i \partial_t + 1)$ has an adjoint operator $(i \partial_t + 1)^T$, and the composition $\mathcal{K} = \iota (i \partial_t + 1)$ acts via 
\begin{equation}
    \mathcal{K}: \mathcal{H} \to \mathcal{H}_1, \mathcal{K}(e_j e^{2 \pi i k t}) = \frac{-2 \pi k + 1}{1 + (2 \pi k)^2} e_j e^{2 \pi i k t}
\end{equation}
where $e_j$ is the $j$-th unit coordinate vector of $\C^n$, and $k \in \Z$. The operator $\mathcal{K}$ is a compact self-adjoint operator without kernel, has image equal to $\mathcal{H}_1$, and has \emph{mild spectrum} in the sense of \cite[Definition~33.1.1]{kronheimer2007monopoles}. The closure of its positive and negative eigenvectors in $\mathcal{H}$ and $\mathcal{H}_1$ will be denoted by $\mathcal{H}^\pm$ and $\mathcal{H}_1^\pm$, respectively.

There are spaces of operators
\begin{equation}
\begin{gathered}
    B_\R(\mathcal{H}, \mathcal{H}_1) = \{ A \in B(\mathcal{H}): A(\mathcal{H}_1) \subset \mathcal{H}_1, A^*(\mathcal{H}_1) \subset \mathcal{H}_1\}, \|A\|_{B(\mathcal{H}, \mathcal{H}_1)} = \max(\|A\|_{\mathcal{H}}, \|A\|_{\mathcal{H}_1}, \|A^T\|_{\mathcal{H}_1});\\
    B_\C(\mathcal{H}, \mathcal{H}_1)= \{A \in B_\R(\mathcal{H}, \mathcal{H}_1) : A \text{ is } \C\text{-linear}\}\\
    \mathcal{O}(\mathcal{H}; \mathcal{H}_1) = \{O \in B_\R(\mathcal{H}, \mathcal{H}_1): O^TO=1\}\\
       \mathcal{U}(\mathcal{H}; \mathcal{H}_1) = \{U \in B_\C(\mathcal{H}, \mathcal{H}_1): U^*U=1\}
        \end{gathered}
\end{equation}
with the latter two spaces naturally having the strucutre of Banach-Lie Groups (this follows from standard spectral theory of unitary operators \cite{lax}). 
There are further spaces of bounded Fredholm operators
\begin{equation}
\begin{gathered}
    \mathcal{S}^{\C}(\mathcal{H}; \mathcal{H}_1) = \{ A \in B_\C(\mathcal{H}_1, \mathcal{H}) | A\text{ Fredholm of index }0, \langle A \psi, \phi \rangle_\C = \langle \psi, A \phi \rangle \} \\
     \mathcal{S}(\mathcal{H}; \mathcal{H}_1) = \{ A \in B_\R(\mathcal{H}_1, \mathcal{H}) | A \text{ Fredholm of index }0,  \langle A \psi, \phi \rangle = \langle \psi, A \phi \rangle \} 
\end{gathered}
\end{equation}
where $B_\C$, $B_\R$, denote $\C$- and $\R$-linear bounded maps of Banach spaces, respectively. These spaces, which are all Banach manifolds,  contain subspaces
\begin{equation}
        \mathcal{S}^\C_{*}(\mathcal{H}; \mathcal{H}_1) \subset \mathcal{S}^\C(\mathcal{H}; \mathcal{H}_1),\;  \mathcal{S}_*(\mathcal{H}; \mathcal{H}_1) \subset  \mathcal{S}(\mathcal{H}; \mathcal{H}_1)
\end{equation}
consisting of those operators $A$ such that there exists $U \in \mathcal{U}(\mathcal{H}; \mathcal{H}_1)$ or $O \in \mathcal{O}(\mathcal{H}; \mathcal{H}_1)$ respectively, sending $\mathcal{H}^\pm$ and $\mathcal{H}_1^\pm$ to the closures of the spans of the non-negative (for $+$) and negative (for $-$) eigenvectors of $A$ (which are all guaranteed to lie in $\mathcal{H}_1$) in $\mathcal{H}$ and in $\mathcal{H}_1$, respectively. 

For every $x \in LM$, (after possibly restricting the bundle $x^*TM$ to a homotopy-equivalent subspace $LM' \subset LM$ consisting of sufficiently regular loops, which we will denote by $LM$ for the purposes of the rest of this proof) can be differentiably trivialized as a finite-dimensional Hermitian bundle, and a choice of trivialization defines a trivialization of $(x')^*TM$ for all $x' \in LM$ near $x$. As such, the bundle $\mathcal{T}^{0,2}LM \to LM$ has structure group $U(\mathcal{H})$, and the inclusion
\begin{equation}
    \bar{\iota}: \mathcal{T}^{1,2}LM \to \mathcal{T}^{0,2}LM
\end{equation}
defines a restriction of structure group of $\mathcal{T}^{0,2}LM$ from $\mathcal{U}(\mathcal{H})$ to $\mathcal{U}(\mathcal{H} ; \mathcal{H}_1)$. The corresponding $\mathcal{U}(\mathcal{H} ; \mathcal{H}_1)$-principal bundle has associated Hilbert bundles $\mathcal{T}^{0,2}LM$ as well as the bundle $\mathcal{T}^{1,2}LM$; the fiberwise-norms induced on these bundles from their identification with associated bundles of this principal bundle are equivalent to their usual fiber-wise inner products induced from $TM$.The operator $Hess(\mathcal{A}_H)$ \eqref{eq:hessian-operator-form} can be thought of as a $C^1$ section of the associated bundle $\mathcal{S}(\mathcal{T}^{0,2}LM, \mathcal{T}^{1,2}LM)$ corresponding to the action of $\mathcal{U}(\mathcal{H} : \mathcal{H}_1)$ on $\mathcal{S}(\mathcal{H}; \mathcal{H}_1)$. 

In fact, it is canonically  homotopic as $C^1$ sections to an operator which lies in $\mathcal{S}^\C_*(\mathcal{T}^{0,2}LM, \mathcal{T}^{1,2}LM) \subset \mathcal{S}(\mathcal{T}^{0,2}LM, \mathcal{T}^{1,2}LM)$, the bundle associated to the $\mathcal{U}(\mathcal{H} : \mathcal{H}_1)$-action on $\mathcal{S}^{\C}_*(\mathcal{H}; \mathcal{H}_1)$. To see this,  we first note that there is a canonical homotopy of the operators \eqref{eq:hessian-operator-form} through self-adjoint Fredholm operators to the operator $J_t\Grad_t$, which in turn under unitary trivializations of $x^*TM$ become  operators of the form  $i\partial_t + A^x(t)$ where $A^x(t)$ is a family of symmetric matrices. Deforming the connection on $M$ to a unitary connection induces a further deformation of operators to one such that the matrices $A^x(t)$ are all Hermitian. As such,the eigenvectors of this operator can be written as $e^{2 \pi i k t}v_j(t)$ for some finite collection of functions $v_j$, $j= 1, \ldots 2n$, corresponding to the eigenvectors of the operator $u_x$ associated to the value at $t=1$ of the solution to the equation $i\partial_tU(t) + A^x(t)U(t)=0$ with $U(0) = 0$ for $U \in U(N)$  (see \cite[Proposition 33.2.2]{kronheimer2007monopoles}). In particular, we see that these operators all lie in $\mathcal{S}^\C_*(\mathcal{H}_1; \mathcal{H})$, and so $J_t\Grad_t$ defines a section of $\mathcal{S}^\C_*(\mathcal{T}^{0,2}LM, \mathcal{T}^{1,2}LM)$. 

Chapter 33 of \cite{kronheimer2007monopoles} proves:
\begin{proposition}[Proposition 33.1.6, \cite{kronheimer2007monopoles}]
    The group $\mathcal{U}(\mathcal{H}; \mathcal{H}_1)$ is contractible, and the space $S^\C_*(\mathcal{H}; \mathcal{H}_1)$ admits a map to $U$ which is a homotopy equivalence. 
\end{proposition}
Applying this proposition to the operators $J_t\Grad_t$ thus defines an element of $[LX, U]$, which defines an associated element of $[LX, U/O] = KO^1(LX)$.

\begin{proposition}
\label{prop:analytic-polarization-class}
    The group $\mathcal{O}(\mathcal{H}; \mathcal{H}_1)$ is contractible, and the space $S_*(\mathcal{H}; \mathcal{H}_1)$ admits a map to $U/O$ which is a homotopy equivalence. 
\end{proposition}
\begin{proof}
The proof of this proposition follows directly analogously to the proof of Proposition 33.1.6 of \cite{kronheimer2007monopoles} except for one key detail, where one uses the geometry of the Lagrangian Grassmannian. We will start off assuming that $\mathcal{H}$ has no complex structure, and complexify it shortly. The contractibility of $\mathcal{O}(\mathcal{H}; \mathcal{H}_1)$ is proven in the same way as the contractibility of $\mathcal{U}(\mathcal{H}; \mathcal{H}_1)$, via the Eilenberg swindle. Let $\beta: \R \to [-1, 1]$ be a monotone function with $\beta(\lambda) = 1$ for $\lambda > \Lambda_+$ and $\beta(\lambda) = -1$ for $\lambda < \Lambda_-$. The map $A \mapsto \beta(A)$ defines a continuous map from $S(\mathcal{H}; \mathcal{H}_1)$ to 
\[ S^f(\mathcal{H}; \mathcal{H}_1) = \{ A \in B(\mathcal{H}; \mathcal{H}_1) \cap S(\mathcal{H} ;\mathcal{H}_1)| A^2-1 \text{ has finite rank}, \|A\|=1.\} \]
Its image is 
\[S^f_*(\mathcal{H}; \mathcal{H}_1) \subset S^f(\mathcal{H}; \mathcal{H}_1),\] the subspace of operators for which there are unitary isomorphisms $\ker (A\mp 1) \to \mathcal{H}^\pm$ which restrict to bounded isomorphisms $\ker (A\mp 1) \cap \mathcal{H}^\pm_1 \to \mathcal{H}_1^\pm$. The fibers of this map can be identified with the spaces of pairs of operators $(A^+, A^-)$ with 
\[ A^\pm: \ker (A \mp 1) \cap \mathcal{H}^\pm_1 \to \ker(A \mp 1)\]
a symmetric Fredholm operator of index zero with spectrum bounded from one side. Thus this map is a homotopy equivalence (by arguing analogously to \cite[Lemma 3.7]{atiyah1969index} using an inductive argument and the Dold-Thom quasifibration recognition lemma as in \cite{kronheimer2007monopoles}).  

Let \[U^f(\mathcal{H} \tensor \C; \mathcal{H}_1 \tensor \C) \subset U(\mathcal{H} \tensor \C; \mathcal{H}_1 \tensor \C)\] 
be the set of elements of such that $u-1$ is of finite rank.
Let 
\[ U^f_2(\mathcal{H} \tensor \C; \mathcal{H}_1 \tensor \C) \subset U(\mathcal{H} \tensor \C; \mathcal{H}_1 \tensor \C)\]
consisting of elements $u$ such that $u^2-1$ is of finite rank. Let 
\begin{equation}
    \label{eq:lag-grassmannian}
    Lag_f(\mathcal{H} \tensor \C; \mathcal{H}_1 \tensor \C) = U^f(\mathcal{H} \tensor \C; \mathcal{H}_1 \tensor \C)/U^f(\mathcal{H} \tensor \C; \mathcal{H}_1 \tensor \C)
\end{equation}
where the quotient is taken by the subgroup  of elements which preserve (but might not fix) the subspace 
\[i \mathcal{H} \subset \mathcal{H} \tensor \C.\] 
In other words, \eqref{eq:lag-grassmannian} parameterizes set of images of the Lagrangian subspace $\mathcal{H} \subset \mathcal{H} \tensor \C$ under elements of $U^f_2(\mathcal{H} \tensor \C; \mathcal{H}_1 \tensor \C)$, topologized with the quotient topology. Elements of $Lag_f(\mathcal{H} \tensor \C; \mathcal{H}_1 \tensor \C))$ are subspaces $\Lambda$ of $\mathcal{H} \tensor \C$ on which the symplectic form $(v, w) \mapsto Re \langle v, i w \rangle$ restricts to zero, and which have the property that $\Lambda \cap i \mathcal{H}$ is finite codimension in both $\Lambda$ and in $i\mathcal{H}$. The perpendicular to $\Lambda \cap i \mathcal{H}$ in $\Lambda$, which we will denote by $(\Lambda \cap i \mathcal{H})^\perp$, is a finite dimensional subspace of $\mathcal{H} \tensor \C$ on which the symplectic form restricts to zero, and which intersects with $i\mathcal{H}$ in $0$. As such, projection to $\mathcal{H}$ is full rank, and, writing $\pi((\Lambda \cap i \mathcal{H})^\perp)$ for its image, we have that $(\Lambda \cap i \mathcal{H})^\perp$ is the image of $\pi((\Lambda \cap i \mathcal{H})^\perp)$ under a unique symmetric unitary automorphism $U$ of $\C \pi((\Lambda \cap i \mathcal{H})^\perp)$ with eigenvalues $e^{i\pi/2 \lambda_j}$ satisfying $1 > \lambda_j > -1$ for all $j$. 

\begin{remark}
    The last claim may not be completely obvious but follows from elementary symplectic linear algebra (the product of symmetric matrices is symmetric, an element of $U(n)$ preserving $\R^n$ lies in $O(n)$, and a Lagrangian subspace that is a graph over $\R^n$ is the graph of a unique real symmetric matrix).   The eigenvalues of $U$ are the `Lagrangian angles' of the Lagrangian $(\Lambda \cap i \mathcal{H})^\perp \subset \C \pi((\Lambda \cap i \mathcal{H})^\perp)$.
\end{remark}

This implies that the continuous map
\begin{equation}
\label{eq:final-analytic-index-comparison-map}
\begin{gathered}S^f_*(\mathcal{H}; \mathcal{H}_1) \to Lag_f(\mathcal{H} \tensor \C; \mathcal{H}_1 \tensor \C)\\
A \mapsto [\exp((i \pi/2) A)]
\end{gathered}
\end{equation}
is a surjection.
The fibers of this map correspond to orbits of the decomposition $\mathcal{H} = \mathcal{H}^+ \oplus \mathcal{H}^-$ under $O(\mathcal{H}; \mathcal{H}_1)$, which are contractible. Thus, arguing as in \cite{atiyah1969index}, this map is a homotopy equivalence. (Note that the statement that the has the homotopy type of a CW complex can be proven from Lemma \ref{lemma:colimit-agrees-with-hocolim} and standard results about spectral theory of symmetric operators by filtering the domain by the rank of the matrix.)

Let $e_i$ be an orthonormal basis of eigenvectors of $K$, arranged in decreasing order of the absolute value of the eigenvalue. Let $E_n$ be the span of the first $n$ eigenvectors, and let $Lag_n \subset \Lambda \cap i \mathcal{H}$ be the subset for which $(\Lambda \cap i \mathcal{H})^\perp$ is contained in $\C E_n$. Then $Lag_n$ is  manifestly   the Lagrangian Grassmannian of subspaces in $\C E_n$, and an approximation argument \cite{palais-approximation}
shows that $\cup_n Lag_n$ is homotopy equivalent to $\Lambda \cap i \mathcal{H}$. We have thus proven the proposition. 
\end{proof}

\begin{proposition}[Equivariant analog of Proposition \ref{prop:analytic-polarization-class}]
\label{prop:equivariant-analytic-polarization-class}
    Suppose that $\mathcal{H}_1$ contains a dense complete $G$-universe. Then the map in $hTop$ constructed in the previous proposition admits a canonical lift to $G-Top$ which is also a $G$-homotopy-equivalence. 
\end{proposition}
\begin{proof}[(Sketch.)]
    The equivariant contractibility of $\mathcal{O}(\mathcal{H}, \mathcal{H}_1)$ follows from \cite{james1978equivariant}. All maps defined in the previous section are equivariant so it suffices to show that they are $G$-homotopy equivalences. The arguments of \cite{atiyah1969index} showing the maps defined above with contractible fibers are weak homotopy equivalences can be adapted to the equivariant case using the equivariant quasifibration results of \cite{waner1980equivariant}.
\end{proof}

With this background in mind, we now prove Theorem \ref{thm:condition-for-vanishing-of-equivariant-polarization-class}.

\begin{proof}[Proof of Theorem \ref{thm:condition-for-vanishing-of-equivariant-polarization-class}]
    Note that $S_*(\mathcal{H}, \mathcal{H}_1)$ has a canonical base point, namely the operator $A_0$ which acts by $\pm 1$ on $\mathcal{H}^\pm$. We will show that the map $LM \to S_*(\mathcal{H}, \mathcal{H}_1)$, $u \mapsto A_u$ produced by Proposition \ref{prop:analytic-polarization-class} from the operators $J_t \Grad'_t$ (where $\Grad'_t$ is $G$-equivariant Hermitian connection on $M$) can be equivariantly contracted to the base point. If $Spec A_u \cap \{0\} = 0$ for all $u$, the proof would be obvious: over any compact subset of $LM$, a homotopy would be given by $\beta(C(1+t) A_u)$ for some sufficiently large $C$, which would land in the equivariantly contractible fiber of \eqref{eq:final-analytic-index-comparison-map}. Thus $u \mapsto A_u$ would give the trivial map on equivariant homotopy groups and thus by Theorem \ref{thm:g-banach-manifold-cw-complex} would be equivariantly homotopic to the constant map. The hypothesis allows us to canonically `split' the zero eigenspaces of $A_u$ in this homotopy. First, note that the $t$-dependent complex structure $J_t$ can be deformed to $t$-independent one since $\pi_1(BU(n)) = 0$; second, note that we can deform $\Grad'_t$ equivariantly to a connection $\Grad''_t$ which preserves the splitting $TM = \xi \oplus i \xi$.  The matrices $U(t)$ arising in the characterization of the spectrum of $J\Grad''_t$ (see reference to \cite[Proposition 33.2.2]{kronheimer2007monopoles} above) will then actually be orthogonal matrices, and the spectrum of $J\Grad''_t$ will be symmetric under flipping the sign of the eigenvalue with an even number of eigenvectors in any symmetric interval around zero. Moreover, given a sphere $S^V$ mapping euivariantly to $LM$, choosing $\Grad''_t$ generically among equivariant connections preserving the splitting above, the same characterization shows that if we plot the spectrum of $J\Grad''_t$ as a subset $\mathcal{S}$ of $S^V \times \R$ it will be a union of local sections of the latter bundle over $S^V$. Thus there will be an $L$ such that $S^V \times ((-\infty, - L) \cup (L, \infty))$ never intersects any component of $\mathcal{S}$ which intersects the zero section. Thus by a suitable deformation $\beta_t$ one can deform $J\Grad''_t$ over $S^V$ to a family of operators  for which the spectrum is either $1, 0$ or $-1$, with the multiplicity of the zero eigenspace being even and and constant. By choosing sections of the $1$ and $-1$ eigenspaces (which are infinite-dimensional bundles over $S^V$, and thus have sections!) and deforming those eigenvectors to have eigenvalue zero, one can deform the operator to one for which the bundle of zero eigenspaces is trivial over $S^V$. A global orthormal basis of sections $e_1, \ldots, e_r, e'_1, e'_r$ of the zero eigenspace and then deforming $e_1, \ldots, e_r$ to have positive eigenvalues and $e'_1, \ldots, e'_r$ to have negative eigenvalues then gives a deformation of the family of operators to the basepoint. 
\end{proof}

\begin{remark}
There is a more direct construction of this class $\rho \in KO^1(LM)$ which has been claimed \cite{cohen2007floer} to be isomorphic to the class defined by \eqref{eq:family-self-adjoint-index-map} under the isomorphism of Theorem \ref{thm:atiyah-singer-self-adjoint}. Part of the purpose of this section is is to explain how to make sense of this claim, given that the Hessians of the Floer action functional are \emph{unbounded} self-adjoint operators, while Theorem \ref{thm:atiyah-singer-self-adjoint} is about \emph{bounded} self-adjoint operators.  (We do not give a complete argument for the claim in \cite{cohen2007floer}.) In this paper, we will use Propositions \ref{prop:analytic-polarization-class} and \ref{prop:equivariant-analytic-polarization-class} to define the polarization class, in order to establish the existence of equivariant framings, and in order to keep the discussion as unambiguous as possible. Note that \cite{cohen1995floer} also gives an argument for the non-equivariant analog of Theorem \ref{thm:condition-for-vanishing-of-equivariant-polarization-class}; an equivariant version of the argument of \cite{cohen1995floer} was not obvious to the author.
\end{remark}

\subsection{Floer homology and loop space structures on spaces of operators}
\label{sec:trivialize-one-moduli-space}
We now explain how the maps \eqref{eq:family-self-adjoint-index-map} of Section \ref{sec:index-theory-for-paths} arise in the context of Floer homology. Indeed, given a pair $x_\pm \in Fix(H)$ such that $Hess(\mathcal{A}_H)_{x_\pm}$ has no zero eigenvalues, and a smooth path of loops
\begin{equation}
    u: \R \times S^1 \to M, \lim_{s \to \pm \infty} u(s,t) = x_\pm(t)
\end{equation}
(where we require that in a chart near $x_\pm$, the above convergence is exponential in the $C^0$ norm), the linearization of the floer equation 
is 
\begin{equation}
    \label{eq:linearized-floer-equation}
    \xi \mapsto D_u\xi  = 
   \Grad_s \xi +  J_t \Grad_t \xi + (\Grad_{\xi} J_t) \partial_t u - \Grad_\xi \Grad H_t.
\end{equation}
Thinking $u(s,t) = x_s(t)$ as a path in $LM$, we have bundles of Hilbert spaces $x_s^*\mathcal{T}^{1,2}LM$ and $x_s^*\mathcal{T}^{0,2}LM$ over $\R$. Choosing a smooth trivialization $\rho$ of the $\mathcal{O}(\mathcal{H}, \mathcal{H}_1)$-principal bundle underlying $\mathcal{T}^{0,2}LM$ gives trivializations of Banach bundles
\begin{equation}
\label{eq:hilbert-bundle-trivializations}
\begin{gathered}
 \rho_s: x_s^*\mathcal{T}^{1,2}LM \to \mathcal{H}_1 \times \R, \\ \rho_s: x_s^*\mathcal{T}^{0,2}LM \to \mathcal{H} \times \R,  
\end{gathered}
\end{equation}
(where the latter is a trivialization of Hilbert bundles)
denoted by the same symbol since the latter is the completion of the former. In this trivialization, the linearized Floer operator becomes 
\begin{equation}
\label{eq:operator-in-trivialization}
    \partial_s + \rho_s Hess(\mathcal{A}_H)\rho_s^{-1} + h_s: L^{1,2}(\R_s, \mathcal{H}) \cap L^{2}(\R_s, \mathcal{H}_1) \to  L^{0,2}(\R_s, \mathcal{H})
\end{equation}
where $h_s$ is continuous family of bounded operators on $\mathcal{H}_1$. Explicitly, $h_s$ is defined by 
\begin{equation}
    \label{eq:define-error-term}
    \partial_s + h_s = \rho_s \circ \Grad_s \circ \rho_s^{-1};
\end{equation} 
so these operators are bounded and skew-adjoint on $\mathcal{H}$, and since $\partial_s u \to 0$ exponentially as $s \to \pm \infty$, we have that $\lim_{s \to \pm \infty} h_s = 0$ in operator norm. Thus, in particular the operators 
\[ A_{x_\pm} := \lim_{s \to \pm \infty} \rho_s Hess(\mathcal{A}_H)\rho_s^{-1} + h_s = \lim_{s \to \pm \infty} \rho_s Hess(\mathcal{A}_H)\rho_s^{-1}\]
are $1-ASAFOE$ operators in the sense of \cite[Section 14.1]{kronheimer2007monopoles}. Write 
\begin{equation}
    \label{eq:hessians-after-trivialization}
    A_s =  \rho_s Hess(\mathcal{A}_H)\rho_s^{-1}.
\end{equation}

The operator \eqref{eq:operator-in-trivialization} Fredholm according to Robbin-Salamon \cite{robbin1995spectral}  (see \cite{rabier2004robbin} for a further relaxation of the conditions required for Fredholmness).

Suppose now that the class $[\rho]$ defined by Proposition \ref{prop:analytic-polarization-class} is zero. 
\begin{definition} \label{def:trivialization-of-polarization-map} A trivialization of the polarization map is a homotopy 
\begin{equation}
    \label{eq:homotopy-of-operators}
    \kappa(A_{\gamma}, t): LM \to S_*(\mathcal{H}; \mathcal{H}_1) , t \in [0,1]
\end{equation}
from the family $A_\gamma$, $\gamma \in LM$ to a constant family. When the original family $A_\gamma$ is $C_{p^k}$-equivariant this homotopy must also be $C_{p^k}$-equivariant, and we will refer to $\kappa$ then as a $C_{p^k}$-equivariant trivializtion of the polarization map (the group may be dropped from the notation when it is implicit).
\end{definition}
One can use this homotopy to write a deformation of Fredholm operators from \eqref{eq:operator-in-trivialization} to an operator that does not depend on $u$, but only on the endpoints $A_{x_\pm}$. Let us explicitly write out such a deformation. 
Write $A^1 = \kappa(A, 1)$ for the relevant $\gamma$-independent operator. Write $A_s = A'_{s'(s)}$ where $s'(s) \in (0,1)$ is the orientation-preserving diffeomorphism from $\R$ to $(0,1)$ given by 
\begin{equation}
\label{eq:make-r-compact}
    s'(s) = \frac{1}{2} + \frac{1}{\pi}\tan^{-1}(s).
\end{equation} Find a smooth family of functions with parameter $t$
\begin{equation}
\label{eq:homotopy-parameterization}
    \tilde{r}_t: [0,1] \to [0,1] \times [0,1], \tilde{r}_t = (\tilde{s}_t, \tilde{t}_t)
\end{equation}
such that 
\begin{itemize}
    \item $\tilde{r}_0(s') = (s', 0)$
    \item $\tilde{r}_t(0) = (0, 0), \tilde{r}_t(1) = (1, 0)$ \text{ for all } $t \in [0,1]$; and 
    \item the path $\tilde{r}_1(s')$ lies entirely on $\partial ([0,1] \times [0,1]) \setminus (0,1) \times \{0\}$, and satisfies
    \item $\tilde{r}_1([0, 1/5]) \subset \{0\} \times [0,1]$, 
    \item $\tilde{r}_1^{-1}(0,1) = [1/5, 2/5],$, 
    \item $\tilde{r}_1^{-1}((0,1) \times \{1\}) = (2/5, 3/5)$, and
    \item $\tilde{r}_1(s') = (1,0) + \tilde{r}_1(1-s')$ for $s'\in [3/5, 1]$.
\end{itemize} 

Then there is a homotopy of operators
\[ D_{u, t}: L^{1,2}(\R_s,\mathcal{H}) \cap L^{2}(\R_s, \mathcal{H}_1) \to L^{0,2}(\R_s, \mathcal{H}). \]
where
\begin{equation}
    \label{eq:initial-homotopy}
    D_{u, t} = \partial_s + \kappa(A'_{\tilde{s}_t(s')}, \tilde{t}_t(s')) + (1-t)h_s.
\end{equation}

The space of choices of functions $\tilde{r}_t$ satisfying the above conditions is contractible, and the homotopy of families operators $h$ is well defined up to contractible choice so the above homotopy $D_{u,t}$ is unique up to contractible choice. Moreover, with little bit of effort, one can make the operator $D_{u,1}$ be independent of $u$ for any family of $u$ in a compact CW subcomplex of the space of paths in $LM$ from $x_-$ to $x_+$.  Indeed, the operator $D_{u,1}$ defined above only depends on $A_{x_\pm}$, and we will write 
\begin{equation}
    \label{eq:canonical-family}
    D_{A_{x_-}, A_{x_+}} = \partial_s + A^{x_-, x_+}_s = D_u^1.
\end{equation}
Note that there is a constant $D$ only dependent on $\kappa$, such that $A^{x_-, x_+}_s = A_{x_\pm}$ for $\pm s > D$; and there is a constant $C < D$ such that $A^{x_-, x_+}(s) = A_0$ for $|s| < C$. 
Because of this, we will be able to use the homotopies $A^t_{x}$ to find homotopies of the operators $D_{u, t}$ over all $u$ arising in a flow category, in a way compatible with gluing, and thus to produce framings of a virtually smooth flow category. 

For later reference, we record a compatiblity condition between trivializations of the polarization map:

\begin{definition}Fix a prime number $p$, and write $\phi$ for the map sending a loop to its $p$-fold cover.

    Write $\mathcal{H}^k$ for $\mathcal{H}$ thought of a $C_{p^k}$-representation, and similarly for $\mathcal{H}_1^k$. Then there is an isomorphism of $C_{p^k}\simeq C_{p^{k+1}}/C_p$-representations 
    \begin{equation}
    \label{eq:cyclotomic-iso-on-hilbert-spaces}
        \Gamma_k: (\mathcal{H}^{k+1})^{C_p} \simeq \mathcal{H}^k 
    \end{equation}
    induced by the pullback of functions along the $p$-fold cover map $S^1 \to S^1$. 
    Suppose that we have a $C_{p^k}$-equivariant map 
    \[ A^k_\gamma: LM \to S_*(\mathcal{H}^k; \mathcal{H}^k_1) \]
    and a similar map $A^{k+1}_\gamma$ with $k$ replaced by $k+1$ everywhere. 

    We say that these are \emph{cyclotomically compatible} if $(A^{k+1}_{\phi(\gamma)})^{C_p} = \Gamma_k A_\gamma \Gamma_k^{-1}$. Similarly, given a $C_{p^k}$-equivariant trivializations $\rho_k$ of the $\mathcal{O}(\mathcal{H}^k, \mathcal{H}^k_1)$ principal bundle underlying $\mathcal{T}^{0,2}LM$, and similarly $\rho_{k+1}$ which is like $\rho_k$ with $k$ replaced by $k+1$ everywhere, we say that these are cyclotomically compatible if 
    \[ (\rho_{k+1})_{\phi(\gamma)}^{C_k} = \Gamma_k (\rho_k)_\gamma.\]
    When one has cyclotomically compatible trivializations of these principal bundles, the polarization maps associated to Floer data for $(H^{\#p^k}, J^{\#p^k})$ and $(H^{\#p^{k+1}}, J^{\#p^{k+1}})$ will be cyclotomiclly compatible as well.

   A pair of $(C_{p^k}, C_{p^{k+1}})$-equivariant trivializations of polarization maps which are homotopies through cyclotomically compatible pairs of such maps is said to be a cyclotomically compatible pair of trivializations of polarization maps. 
\end{definition}

The following lemma is straightforward to prove using an inductive argument and the results of the previous sectoin:
\begin{lemma}
\label{lemma:choose-trivializations-of-polarization-map-cyclotomically}
    Suppose that the polarization class $[\rho]\in KO^1_{C_{p^{k+1}}}(LM)$ vanishes. Then, given a trivialization of the polarization map for $(H^{\#k}, J^{\#k})$, there exists a cyclotomically compatible trivialization of the polarization map for $(H^{\#k}, J^{\#k+1}$. 
\end{lemma}

\subsection{Compatible Framings}
\label{sec:producing-framings}
We now explain how to produce a compatible system of framings for the smooth $C_{p^k}$-flow categories produced in Section \ref{sec:global-charts}. We first explain how to produce a framing for the virtual smoothing $\CC'(H, J)$ defined using perturbation data $\mathcal{P} = (\bar{\mathfrak{A}}, \mathcal{F}^1, \{\lambda_{V_{\mathcal{F}^1}}(x,y)\})$ as in Proposition \ref{prop:smooth-flow-category-floer-trajectories}. These framings will depend on a choice of $C_{p^k}$-equivariant trivialization of the polarization map (Definition \ref{def:trivialization-of-polarization-map}).

First, choose once and for all a system of trivializations 
\[ T\overline{Conf}_{k+1} \simeq \underline{T\mathbb{R}^k_+}\]
which satisfy the condition that the restriction of this trivialization to $\overline{Conf}_{k_1+1} + \overline{Conf}_{k_2+1}$ for $k_1+1 + k_2+1 = k+1$ agrees with the product of the trivializations composed with the inclusion 
\[ \underline{T\mathbb{R}^{k_1}_+} \boxplus \underline{T\mathbb{R}^{k_2}_+} \simeq  \underline{T\mathbb{R}^{k_1}_+ \oplus 0 \oplus T\mathbb{R}^{k_2}_+} \subset \underline{T\mathbb{R}^{k}_+}.\]

This defines a system of trivializations of the parameter spaces 
\begin{equation}
\label{eq:parameter-space-trivializations}
    T\paramspace(x,y) \simeq  \underline{T \R^{\mathfrak{A}'(x,y)_+}}.
\end{equation} 
Thus there are partial trivializations
\begin{equation}
    \label{eq:thickening-partial-framing}
    TT(x,y) = \underline{T\R^{\mathfrak{A}'(x,y)_+}}  \oplus \ker D_u^V 
\end{equation} 
where $\ker D_u^V$ is the bundle of kernels of the operators \eqref{eq:stabilized-linearized-operator} (or the kernels of the corresponding direct sum of operators when the domain of $u$ has multiple components). 

Now, recall that the obstruction bundles of $\CC'(H, J)$ are framed with the corresponding parameterization given by
\[ \mathcal{F}^1_T = \mathcal{F}^1_{can} \oplus \mathcal{F}^1\]
with $\mathcal{F}^1_{can}$ the canonical $F'$-parameterization $\mathcal{F}^1_{can}$ corresponding to the $E'$-parameterization $\mathcal{E}^1_{can}$ \eqref{eq:canonical-obstruction-f-parameterization}.

\paragraph{Step 0: Collars.}
Our first step is to choose a \emph{compatible system of collars} on $\CC'(H, J)$. 

There is a functor $T$ from $DerMan'_{sm, no-stab}$ to $Man'$ sending a derived $\langle k \rangle$-manifold to its thickening. This is a monoidal functor; given a virtually smooth flow category $\CC'$, write $\CC'_T$ for the corresponding smooth flow category. The collaring operator $Coll$ of Definition \ref{def:collars} defines a monoidal functor $Man' \to Man'$; thus, given a smooth flow category $\CC'_T$ there is an associated smooth flow category $Coll(\CC'_T)$ obtained by applying $Coll$ to all morphism objects. The latter is $G$-equivariant if the former is.
\begin{definition}
    A \emph{compatible system of collars} on a virtually smooth flow category $\CC'$ is an (equivariant) isomorphism $Coll(\CC'_T) \to \CC'_T$ such that the maps on morphism objects preserve strata, i.e. the induced map $S_{Coll(\CC'_T)(x,y)} \to S_{\CC'_T(x,y)}$ is the identity for all $x > y$ in $\CC$. 
\end{definition}

\begin{lemma}
\label{lemma:collars}
    A virtually smooth flow category $\CC'$ admits a compatible system of collars. If $H \subset G$ is a normal subgroup then a ($G/H)$-equivariant) compatible system of collars on $(\CC')^H$ extends to an ($G$-equivariant) compatible system of collars on $\CC'$.  
\end{lemma}
\begin{proof}
    One follows an equivariant double induction; for each individual moduli space this is the content of \cite[Proposition 4.1]{albin2011resolution}. 
\end{proof}

Now, every sequence $x_0 > x_1 > \ldots x_r$ in $\CC$ has an associated stratum $S_{x_0 > \ldots > x_r}$ in $\CC'(x_0, x_r)$, and by the assumption that $H$ is $C_k$-action nondegenerate, the stratum $\CC'_T(x_0, x_r)(S_{x_0 > \ldots > x_r})$ consists of perturbed Floer trajectories from $x_0$ to $x_r$ which break on a set of Hamiltonian orbits containing a set of the form $\{g_i x_i\}_{i=1}^{r-1}$ for some $g_i \in G_{x_0x_r}$. We will write $\CC'_T(x_0, x_r)[x_0, \ldots, x_r] \subset \CC'_T(x_0, x_r)(S_{x_0 > \ldots > x_r})$ for the corresponding connected component. 

We will write $T(x,y)$ for the thickening of $\CC'(x,y)$, and we will write $T(x,y)^o= T(x,y) \setminus \partial T(x,y)$. We will write 
\[ T(x,y)^\square \text{ for the image of } T(x,y) \text{ under } T(x,y) \subset Coll(T(x,y)) \to T(x,y).\]

\paragraph{Step 1: Deforming operators.}
We now build compatible deformations of Fredholm operators over the spaces $T(x,y)^o \times [0,1]_t$. 

Let $\pi: T(x,y) \to \paramspace(x,y)$ be the projection; write \[T(x,y)^\R = \pi^*\paramspace(x,y)^{\mathbb{R}}.\] 
There is a map 
\[u_s: T(x,y)^\R \to LM \] 
which sends $s$-coordinates on the domain of a perturbed fer trajectory to the corresponding loops. There is also a section $A^u_s$ of the bundle $u_s^*\mathcal{S}_*(\mathcal{T}^{0,2}LM, \mathcal{T}^{1,2}LM)$. 

We follow the strategy of Section \ref{sec:trivialize-one-moduli-space}. 

Let us imagine that the polarization class $\rho \in KO^1_G(LM)$ vanishes for $G = C_k$ for some $k$ or $G = S_1$. 
We fix, once and for all, a $G$-trivialization of the $\mathcal{O}(\mathcal{H}, \mathcal{H}_1)$-principal bundle associated to $\mathcal{T}^{0,2}LM$, which induces trivializations 
\begin{equation}
    \label{eq:hilbert-bundle-trivializations-2}
    \begin{gathered}
        \rho_s: u_s^*\mathcal{T}^{1,2}LM \to \mathcal{H}_1 \times T^C(x,y),\\  \rho_s: u_s^*\mathcal{T}^{0,2}LM \to \mathcal{H} \times T^C(x,y).
    \end{gathered}
\end{equation}

We choose the $G$-homotopy $\kappa$ \eqref{eq:homotopy-of-operators} with the condition that that the base point base point $A^0 \in S_*(\mathcal{H}, \mathcal{H}_1)$ has no zero eigenvalues.

    We wish to find a free stabilization of the perturbation data $\mathfrak{P}$ by trivial perturbation data such that stabilized perturbation data $\mathfrak{P}' = (\bar{\mathfrak{A}}, \mathcal{F}', \lambda_{V^1(x,y)})$ (so $\mathcal{F}' = \{V^1(x,y)\}$ is a free stabilization of $\mathcal{F}$ by $\mathfrak{F}_f$) define the $\tau=0$ case of a  $1$-parameter family of \emph{compatible thickening-dependent regular perturbation data} $\mathfrak{P}'_\tau = (\bar{\mathfrak{A}}, \mathcal{F}', \{\tilde{\lambda}^\tau_{V^1(x,y)}\})$ for $\tau \in [0,1]$. Here, the maps $\tilde{\lambda}^\tau_{V^1(x,y)}$ are differentiable bundle maps from $\underline{V^1(x,y)}$ to $\mathcal{H}_1$ over $T(x,y)^\R$ which are zero near the ends of every fiber of the universal curve over $T(x,y)$, and restrict to one another under the isomorphisms $T(x,y)^\R(T(x,y)[x,z,y]) \simeq T(x,z)^\R \times T(z,y) \cup T(x,z) \times T(z,y)^\R$; the continuity condition in $\tau$ is that as $\tau$ varies these define differentiable bundle maps over $T(x,y)^\R \times [0,1]_\tau$. The perturbation datum $\mathfrak{P}'$ defines the thickening-dependent perturbation datum $\mathfrak{P}'_0$ via 
    \[ \tilde{\lambda}^0_{V^1(x,y)}(u,a,v; s) = \rho_s \circ \lambda_{V^1(x,y)}(a,v)|_{u_s}.\]

    To specify the conditions that $\mathfrak{P}'_\tau$ must satisfy we specify the following operators. 
    First, recall that for $x, y \in \CC$ we have the canonical operators $D_{A_{x_-}, A_{x_+}}$ defined in \eqref{eq:canonical-family}. Now the gluing construction allows us to define, for every $x_- = x_0, x_1, \ldots, x_{r-1}, x_+ \in \CC$, a `canonical family of operators' $D_{A_{x_0}, A_{x_1}, \ldots, A_{x_r}}=: D_{x_0, \ldots, x_r}$ over $[0,1]^{r-1}$. We outline this below; first, however we specify a related gluing construction. Define the families of operators $A^{x, b}_s$ by setting $A^{x, b}_s = A^{x, y}_{s}$ for  $s \leq -C$ (where $y$ is arbitrary), and $A^{x_-, b}_s = A^0$ for $s > -C$. Define the operator 
    \begin{equation}
        \label{eq:asymptotic-fredholm-operator-not-stabilized}
        D_x: L^{1,2}(\R_s, \mathcal{H}) \cap L^2(\R_s, \mathcal{H}_1)  \to L^{0,2}(\R_s, \mathcal{H}), D_x\xi = \partial_s \xi + A^{x, b}_s \xi.
    \end{equation} 
    The operator $D_y$ can be \emph{canonically glued} to the operator $D_{x,y}$ in the following sense. Set, for $\tau_2>D$,  
\begin{equation}
    A^{x_-, x_+}_{s, \tau_2} = \begin{cases}
        A^{x_-, x_+}_s & \text{ if } s < \tau_2 \\
    A^{x_+, b}_{s - (D+\tau_2)} & \text{ if } s \geq \tau_2 
    \end{cases}
\end{equation}
and for $0 < \tau_2 \leq D$, we set 
\begin{equation}
    A^{x_-, x_+}_{s, \tau_2} = \begin{cases}
        A^{x_-, x_+}_s & \text{ if } s < \tau_2 \\
    A^{x_+, b}_{s - (2\tau_2)} & \text{ if } s \geq \tau_2.
    \end{cases}
\end{equation}

Thus for $\tau_2 > D$, $\tau_2-D$ measures the distance between the regions where the family of operators is staying at $A_0$ is occurring; and for $0 > \tau_2 > D$ the two homotopies $\kappa(A_{x_+}, t)$ going in opposite directions are cancelling themselves out. We can then set, for $\tau_2 \in [0,\infty)$, 
\begin{equation}
    \begin{gathered}
        D_{x,y, \tau_2}:L^{1,2}(\R_s,\mathcal{H}) \cap L^{2}(\R_s, \mathcal{H}_1)  \to L^{0,2}(\R_s, \mathcal{H}).\\
     D^{\mathfrak{P'}}_x(\xi, v) =  \partial_s\xi  +  A^{x_-, x_+}_{s, \tau_2} \xi. 
    \end{gathered}
\end{equation}

As $\tau_2 \to \infty$, this operator becomes the `gluing' of $D_{x,y}$ and $D_y$. 

Choosing once and for all an orientation-reversing diffeomorphism $[0,\infty)_{\tau_2} \to (0, 1]_{\tau_2'}$ (we choose to work with the map $\tau_2 = 1/\tau'_2-1$) this defines a family of operators $D_{x,y, \tau_2'}$ over $(0,1]$. We set $D_{x,y, 0} = D_{x,y} \oplus D_y$. 

By performing a similar construction with the $A^{x_i, x_{i+1}}_s$, one defines operators 
\begin{equation}
    \begin{gathered}
        D_{x_0, \ldots, x_{r}; \tau'_2}:L^{1,2}(\R_s,\mathcal{H}) \cap L^{2}(\R_s, \mathcal{H}_1)  \to L^{0,2}(\R_s, \mathcal{H}), \tau'_2 \in (0, 1]^{[r-2]}\\
     D^{\mathfrak{P'}}_x(\xi, v) =  \partial_s\xi  +  A^{x_0, \ldots, x_r}_{s; \tau_2}. 
    \end{gathered}
\end{equation}
Here $\tau_2$, after applying the diffeomorphism $(0,1]^{[r-1]} \to [0, \infty)^{[r-1]}$ given by applying the previously chosen diffeomorphism coordinate-wise, measures the sizes of the regions where homotopies from $A_{x_i}$ to $A^0$ to $A_{x_{i+1}}$ are occurring. Thinking of $[0,1]^{[r-1]} \subset \R^{[r-1]}_+$ as a stratified spaces, note that 
\[ [0,1]^{[r-1]}([r-1] \setminus\{a_1, \ldots, a_\ell\}) = [0,1]^{\{0, \ldots, a_1-1\}} \times [0,1]^{\{a_1+1, \ldots, a_{2}-1\}} \times \ldots \times [0,1]^{\{a_{\ell-1}+1, \ldots, a_\ell\}}.\]
We set, for $\tau'_2$ in $[0,1]^{[r-1]}([r-1] \setminus\{a_1, \ldots, a_\ell\})$, 
\[D_{x_0, \ldots, x_r; \tau'_2} = D_{x_0, \ldots,x_{a_1}} \oplus D_{x_{a_1}, \ldots, x_{a_2}} \oplus \ldots \oplus D_{x_{a_{\ell}}, \ldots, x_{r}}.\]

We now define operators
\[ D_{\tilde{u}, 1}: L^{1,2}(T(x,y)^\R_{\tilde{u}}, \mathfrak{H}) \cap L^2(\paramspace(x,y)^\R_{a}, \mathfrak{H}_1) \to L^{0,2}(T(x,y)^\R_{\tilde{u}}, \mathfrak{H}) \]
where we write $\tilde{u} = (u,a,v)$ for an element in $T(x,y)$, as follows. 
 We recall two facts. First, recall that the fibers of the universal curve $T(x,y)^\R_{\tilde{u}}$ are disjoint ordered copies of $\R$, with an $s$-coordinate well defined up to translation; if $u$ lies in the interior of a stratum corresponding to the set $S_{\CC'(x,y)} \setminus S'$, then the fiber consists of $S'$+1 copies of $\R$. In a moment, we we will choose trivializations 
 \[ \phi^\R_{S}: (T(x,y)(S)^o)^\R \to \underline{\R}^{\sqcup |S_{\CC(x,y)} \setminus S|+1}. \]
 Second, recall that the compatible system of collars defines maps 
 \[ \Gamma: T(x_0, x_r)[x_0, \ldots, x_r] \times [0,1]^{[r-1]} \to T(x_0, x_r)\]
 for every $x_0 > x_1 > \ldots > x_r$ in $\CC$ (where we can think of $[0,1]^{[r-1]}$ as $[0,1]^{\{x_1, \ldots, x_{r-1}\}}$). 
 These maps are diffeomorphisms onto subsets $U(x_0, \ldots, x_r) \subset T(x_0, x_r)$; giving these the induced stratification, we can write $\partial U(x_0, \ldots, x_r)$ for the image of $T(x_0, x_r)[x_0, \ldots, x_r] \times ([0,1]^{[r-1]} \cap \R^{[r-1]}_+)$  the map defined by the collar. We write $\epsilon_{x_0, \ldots, x_r}: U(x_0, \ldots, x_r) \to [0,1]^{[r-1]}$ for the map given by the composition of $\Gamma^{-1}$ with the projection to the second factor.
 We choose the trivializations $\phi^\R_S$ such that if we set
 \[ A^{\tilde{u}}_{s; 1} = A^{x_0, \ldots, x_r; \epsilon_{x_0, \ldots,x_r}(\tilde{u})}_{\phi_{S(x_0, \ldots, x_r}(\tilde{u}; s)} \text{ for } \tilde{u} \in U(x_0, \ldots, x_r)\]
 then there is a continuous extension of $A^{\tilde{u}}: T(x,y)^\R \to S_*(\mathcal{H}, \mathcal{H}_1)$ to $T(x,y)^{\bar{\R}}$. 
  Using these maps, we define $D_{\tilde{u}, 1}$ as 
  \[ D_{\tilde{u}, 1} = \partial_s + A^{\tilde{u}}_{s; 1} \text{ for } u \in T(x,y)^o, \]
  and more generally, writing 
  \[ T(x_0, x_r)[x_0, \ldots, x_r]^o \simeq T(x_0, x_1)^o \times \ldots \times T(x_{r-1}, x_r)^o, \tilde{u} \leftrightarrow (\tilde{u_1},\ldots, \tilde{u_r}), \]
  we have 
  \[ D_{\tilde{u},1} = D_{\tilde{u}_1, 1} \oplus \ldots \oplus D_{\tilde{u}_r, 1}. \]
  Given $\mathfrak{P}'_\tau$ we then introduce
  \[ D^{\mathfrak{P}'}_{\tilde{u}, 1}:(L^{1,2}(\R_s,\mathcal{H}) \cap L^{2}(\R_s, \mathcal{H}_1))^{\oplus r}\oplus V^1(x_0, x_r) \to L^{0,2}(\R_s, \mathcal{H}) \text{ for } \tilde{u} \in T(x_0, x_r)[x_0, \ldots, x_r]^o\]
  \[D^{\mathfrak{P}'}_{\tilde{u}, 1}((\xi)_{j=1}^r, v) = \oplus_j D_{\tilde{u}_j, 1}\xi_j + \tilde{\lambda}^0_{V^1(x,y)}(u,a,v; s) \]
  \begin{lemma}
  \label{lemma:bundle-of-kernels-1}
      Suppose that the operators $D^{\mathfrak{P}'}_{\tilde{u}, 1}$ are all surjective. Then, topologizing the space of kernels $\ker D^{\mathfrak{P}'}_{\tilde{u}, 1} =: \ker D^{\mathfrak{P}'}_{\tilde{u}; x,y, 1}$ in the Hausdorff topology on sections of the bundle $\underline{\mathcal{H} \oplus V^1(x_0, x_r)}$ over $T(x,y)^\R$ we have that the bundle of kernels of $D^{\mathfrak{P}'}_{\tilde{u}, 1}$ forms a topological vector bundle over $T(x_0, x_r)$ for every $x_0, x_r$. 
  \end{lemma}
  \begin{proof}
      This follows from the compatibility conditions for the $A^{\tilde{u}_s}$ and standard exponential decay estimates for solutions to operators of the form  $D^{\mathfrak{P}'}_{\tilde{u}, 1}$, which are simply of the form $\partial_s + A$ at the ends of $\R$.  
  \end{proof}

  We wish to interpolate between these bundles and the corresponding bundles associated to stabilizing $\mathfrak{P}$ to $\mathfrak{P}'$ as perturbation data for $\CC'$. Note that, writing $A^{\tilde{u}}_{s, 0} = \rho_sA^{\tilde{u}}_s\rho_s^{-1}$ we can identify the latter bundle with the bundle of fibers of the operators 
  \[ D_{\tilde{u}, 0}^{\mathfrak{P}'}: (L^{1,2}(\R_s,\mathcal{H}) \cap L^{2}(\R_s, \mathcal{H}_1))^{\oplus r}\oplus V^1(x_0, x_r) \to L^{0,2}(\R_s, \mathcal{H}) \text{ for } \tilde{u} \in T(x_0, x_r)[x_0, \ldots, x_r]^o \]
  \[D^{\mathfrak{P}'}_{\tilde{u}, 0}((\xi)_{j=1}^r, v) = \oplus_j (\partial_s\xi_j + A^{\tilde{u}}_{s;0} \xi_j + h^u_s \xi) + \tilde{\lambda}^0_{V^1(x,y)}(u,a,v; s).\]
  We will find families of operators $A^{\tilde{u}}_{s,\tau}: T(x,y)^{\bar{\R}} \to S_*(\mathcal{H}; \mathcal{H}_1)$ which restrict over $T(x,y)[x,z,y]^{\bar{\R}}$ continuous in $\tau \in [0,1]$ to  the corresponding operators 
  $A^{\tilde{u}_1}_{s, \tau} \sqcup A^{\tilde{u}_2}_{s, \tau}$ and are equivariant with respect to the $C_k$ action, such that, defining 
  \[D^{\mathfrak{P}'}_{\tilde{u}, \tau}((\xi)_{j=1}^r, v) = \oplus_j (\partial_s\xi_j + A^{\tilde{u}}_{s;\tau} \xi_j + (1 - \tau) h^u_s \xi) + \tilde{\lambda}^0_{V^1(x,y)}(u,a,v; s).\]
  (where this has the same domain and codomain as $D^{\mathfrak{P}'}_{\tilde{u}, 1}$) the spaces of kernels of these operators, topologized as in Lemma \ref{lemma:bundle-of-kernels-1}, define topological vector bundles over $T(x,y) \times [0,1]_\tau$. 

  \begin{lemma}
      \label{lemma:bundle-of-kernels-2}
      Homotopies $A^{\tilde{u}}_{s, \tau}$ as above exist such that if $\mathfrak{P}'_\tau$ is such that the above operators $D^{\mathfrak{P}'}_{\tilde{u}, \tau}$ are all surjective then the spaces of kernels of these operators, topologized as above, define (equivariant) topological vector bundles over $T(x,y) \times [0,1]_\tau$.
  \end{lemma}
   \begin{proof}(Sketch.)
       We will define this homotopy $A^{\tilde{u}}_{s, \tau}$ in two steps, first over $\tau \in [0, 1/2]$ by deforming the operators $A^{\tilde{u}}_{s, 0}$ in certain regions where they are exponentially close to $A_{x}$ to actually \emph{agree} with $A_x$ in those regions; and then, over $\tau \in [1/2, 1]$, by inserting analogs of the homotopies defined by functions like $\tilde{r}_t$ in \eqref{eq:homotopy-parameterization} in the remaining regions to complete the deformation to $A^{\tilde{u}}_{s, 1}$. To define this deformation, we will choose functions 
       \[ \epsilon_{[z], \pm}: T(x,y) \to \R_{\geq 0}, [z] \in \CC/G, x > gz > y \text{ for some }g\]
       which define certain regions in the universal curves $T(x,y)^\R$ as follows. Note that the operators $A^{\tilde{u}}_{s, 1}$ are defined to be \emph{constant} in certain regions, with the regions determined by the trivializations $\phi^\R_S$ and the values of the collar coordinates $\epsilon_{x_0, \ldots, x_r}$. These regions are canonically identified with a collection of intervals of certain lengths; the number of intervals is between $0$ and $r-1$ on $U(x_0, \ldots, x_r)$, and is exactly $r$ on $\epsilon_{x_0, \ldots, x_r}^{-1}([0,1]^{[r-1]} \setminus [K, 1]^{[r-1]})$ for some universal constant $K= K(D)$, with one of these $r$ intervals shrinking to zero as one approaches $\epsilon_{x_0, \ldots, x_r}^{-1}([K, 1] \cap ((K, \ldots, K)+\R^{[k-1]}_+))$. The functions $\epsilon_{[x], \pm}$ will be thought of as defining \emph{subintervals} of these intervals via their offsets from the ends of these intervals, with $\epsilon_{[x], +}: T(x,y) \to \R_{\geq 0}$ also defining a half-infinite-interval corresponding to the region $(\infty, -D-\epsilon_{[x], +}(\tilde{u})]$ on the correspoding copy of $\R$, and similarly $\epsilon_{[y], -}: T(x,y) \to \R_{\geq 0}$ defining $[D + \epsilon_{[y], -}, \infty)$ on the corresponding copy of $\R$. Thus these functions are required the following three conditions:
       \begin{enumerate}[(a)]
       \item The intervals have size at least zero, and cannot have negative size (thus, for example, the functions $\epsilon_{[z], \pm}$ are zero on $U(x_0, \ldots, x_r)$ if $[z] \notin \{[x_0], \ldots, [x_r]\}$);
       \item The functions are equivariant in the sense that $\epsilon_{[z], \pm} \circ g_{xy} = \epsilon_{[z], \pm}$, and 
       \item The functions restrict to one another: on $T(x,y) \times T(y,z) \subset T(x,z)$, $\epsilon_{[w]_\pm}$ restricts either to $\pi_1^* \epsilon_{[w]_\pm}$  (if $x > gw > y$) $\pi_2^*\epsilon_{[w]_\pm}$ (if $y > gw > z$), or $0$ (if neither holds). 
       \end{enumerate}
       The remaining condition that these functions must on the intervals defined by these functions over $U(x_0, \ldots, x_r)$, one must have $u_s$ is exponentially close to the Hamiltonian orbit corresponding to that interval. Specifically, for every Hamiltonian orbit $x$, there is a  constant $C_x$ such that curves asymptotic to $x$ have $d(u_s, x) \leq C'_x e^{\pm C_x s}$ for some $C'_x>0$ (with the sign determined by whether one is asymptoting at a positive or a negative end). There are also constants $\epsilon_x$ such that the exponential map of radius $\epsilon_x$ with respect to the metric determined by $J_t$ is injective at $x(t)$, for every $t \in S^1$. Using this exponential map one can produce a trivialization of $\mathcal{T}^{1,2}(LX)$ around $x$, and in this trivialization one has that if $d(x_s, x) < \epsilon_0$ then one can write $A_{x_s}$ as $A^x + h(x_s)$ with $h(x_s)$ a compact operator on $\mathcal{T}^{1,2}(LX)_{x_s}$, with the property that $A_s + \epsilon h(x_s) \in S_*(\mathcal{T}^{1,2}(LM)_{x_s}, \mathcal{T}^{0,2}(LM)_{x_s})$ for any $\epsilon$. Thus, writing these operators in our fixed $\rho$ gives operators of the form $A^{\tilde{u}}_{s, 0} = A^{x} + h'(u_s)$ for operators $h'(u_s)$ where the operators $A^{x} + \epsilon h'(u_s)$ are in $S_*(\mathcal{H}, \mathcal{H}_1)$ for any $\epsilon$, whenever $s$ is an element of the universal curve over $\tilde{u} \in U(x_0, \ldots, x, \ldots, x_r)$ lying in the interval corresponding to $r$. Thus, by multiplying these $h'(u_s)$ by smooth cutoff functions which  $1$ at open neighborhoods of the ends of the interval and are zero neighborhoods of the subintervals that are offset by one more from either side (with no condition except continuity once the subintervals thus defined would have size less than zero) one can perform the homotopy for $\tau \in [0,1/2]$. Thus replacing the previously defined intervals by these subintervals (corresponding to simply increasing the functions $\epsilon_{[x], \pm}$ by a certain amount) one has that $A^{\tilde{u}}_{s, 1/2}$ are now constant in these subintervals, and agree with the corresponding values of $A^{\tilde{u}}_{s, 1}$. In the complements of the subintervals one then glues in canonical homotopies using the same reparameterization trick as in Section \ref{sec:trivialize-one-moduli-space} (see \eqref{eq:homotopy-parameterization}), going from one such homotopy to the corresponding homotopy for the region where some subintervals have been dropped as the sizes of the subintervals shrink from some fixed size (e.g. length $1$) to length zero. 

       The gluing theorems continue to hold throughout these homotopies because on the subintervals one simply rescales the operators $A^{\tilde{u}}_{s, 0}$, which are exponentially small perturbations of some $A^x$, to constant operators $A^x$, preserving exponential decay bounds in all possible norms (controlled by some norm-specific quantity as well as the constant $C_x$ and the length of the given interval) in these regions, which are the estimates needed to establish the relevant linear gluing theorem. 
   \end{proof}



After choosing $A^{\tilde{u}}_{s, \tau}$ via Lemma \ref{lemma:bundle-of-kernels-2}, we will choose $\mathfrak{P}'_\tau$ such that the conditions of Lemma \ref{lemma:bundle-of-kernels-2} are satisfied. We will also require that 
\begin{equation}
\label{eq:constancy-of-perturbation-function-on-interior}
\parbox{0.8\textwidth}{ the functions $\tilde{\lambda}^1_{V^1(x,y)}(u,a,v; s)$ are independent of $\tilde{u} = (u,a,v)$ for $\tilde{u} \in T(x,y)^\square$.}
\end{equation} 
With this latter condition, the operators $D_{\tilde{u}, 1}^{\mathfrak{P}'}$ are actually independent of $\tilde{u}$ for $\tilde{u}$ in the image of $T(x,y)$ under the collaring $Coll(T(x,y)) \to T(x,y)$. We will write this operator as $D_{x,y}$, and we will write the functions $\tilde{\lambda}^1_{V^1(x,y)}(u,a,v; s)$ for $\tilde{u}$ as above as $\tilde{\lambda}^1_{V^1(x,y)}(s)$. 

We will require that a little bit more data is chosen. Recall that the $F$-parameterization $\mathcal{F}'$ is associated to some $\mathcal{E}$-parameterization $\mathcal{E}^1_s = \{V^1(x)\}$. We will choose
\begin{equation}
    \label{eq:embedding-perturbation-data}
    \tilde{\lambda}_{V^1(x)}: V(x) \to C^\infty(\R, \mathcal{H}_1)
\end{equation}
and also 
\begin{equation}
    \tilde{\lambda}_{V^1(x),V^1(y), \tau_2}: V(x) \times [0, \infty)_{\tau'_2} \to C^\infty(\R, \mathcal{H}_1)
\end{equation}
such that 
\[ \tilde{\lambda}_{V(x)}(v)(s) = 0 \text{ for } |s| >> 0, \text{and}\]
\[ \tilde{\lambda}_{V^1(x),V^1(y), 0} = \tilde{\lambda}_{V^1(x)}(v) \]
and for all sufficiently large $\tau'_2 >>0$,   
\[ \tilde{\lambda}_{V^1(x),V^1(y),\tau_2}(\bar{\alpha}_{xy}(v_{xy}, v_y))(s) = \tilde{\lambda}^1_{V^1(x,y)}(v_{xy})(s) \text{ for } s < \tau_2, \tilde{\lambda}_{V^1(x),V^1(y),\tau'_2}(s) = \tilde{\lambda}^1_{V(x)}(v_y)(s-(D+\tau'_2)) \text{ for } s \geq \tau+2.\]
When the underlying flow category is $G$-equivariant, we require that 
\begin{equation}
\label{eq:embedding-perturbation-data-equivariance}
    \lambda_{V^1(gx)}(gv)= g\lambda_{V^1(x)}(v), \tilde{\lambda}_{V^1(gx),V^1(gy), \tau_2}(gv)= g  \tilde{\lambda}_{V^1(x),V^1(y), \tau_2}(v).
\end{equation}

We set for $\tau_2 \in [0, \infty)$, 
\begin{equation}
\label{eq:asymptotic-fredholm-operator-stabilized}
    \begin{gathered}
        D^{\mathfrak{P'}}_{x,y, \tau_2}, D^{\mathfrak{P}'}_x:L^{1,2}(\R_s,\mathcal{H}) \cap L^{2}(\R_s, \mathcal{H}_1) \oplus V(x) \to L^{0,2}(\R_s, \mathcal{H}).\\
     D^{\mathfrak{P'}}_x(\xi, v) =  \partial_s\xi  +  A^{x_-, x_+}_{s, \tau_2} \xi +\tilde{\lambda}_{V^1(x),V^1(y),\tau_2}(v); D^{\mathfrak{P'}}_x(\xi, v) =  \partial_s\xi  +  A^{x_-, b}_{s} \xi +\tilde{\lambda}_{V^1(x)}(v)
    \end{gathered}
\end{equation}
Note that $D^{\mathfrak{P'}}_{x,y, 1} = D^{\mathfrak{P'}}_{x}$ and that standard gluing theory implies that  $D^{\mathfrak{P'}}_{x,y, \tau_2}$ is surjective for $\tau_2 >> 0$ if $D^{\mathfrak{P'}}_{x,y}$ and $D^{\mathfrak{P'}}_{y}$ are.

\begin{lemma}
\label{lemma:perturbation-data-and-numerical-condition-for-framing}
    We can choose data $(\mathfrak{P}'_\tau, \{\tilde{\lambda}_{V^1(x)}\}, \{\tilde{\lambda}_{V^1(x), V^1(y), \tau_2}\})$ such that conditions of Lemma \ref{lemma:bundle-of-kernels-2} are satisfied, such that $\mathfrak{P}'_0$ arises from $\mathfrak{P}$ via a stabilization by an $E$-parameterization with trivial perturbation data as described above, with the perturbation functions underlying $\mathfrak{P}'_\tau$ satisfying the independence conditions \eqref{eq:constancy-of-perturbation-function-on-interior} making $\tilde{\lambda}^1_{V(x,y)}(v)$ well defined, and with the operators $D^{\mathfrak{P'}}_{x}$ and $D^{\mathfrak{P}'}_{xy, \tau_2}$ all surjective (for $\tau_2 \in [0, \infty)$). 
    
    Moreover, this data can be chosen such that the kernels of the operators $D^{\mathfrak{P'}}_{x}$ operators satisfy the following numerical condition. Write 
    \begin{equation}
        \label{eq:define-v0-parameterization-from-kernels-of-asymptotic-operators}
        V^0(x) = \ker D^{\mathfrak{P}'}_x, V^0(x,y) = \ker D^{\mathfrak{P}'}_{x,y}. 
    \end{equation}
    Then by the equivariance conditions \eqref{eq:embedding-perturbation-data-equivariance} imply that there are isometries
    \[ g^0_x: V^0(x) \to V^0(gx)\]
    satisfying \eqref{eq:equivariance-conditions-for-poset-rep-1} 
    Given vector spaces $V^*(y') \subset V^0(y')$ as $y'$ ranges over $y' \in Gy$ which are  fixed by the operators $g^0_{g'y}$ over $g \in G$, $g' \in G$, and given 
     $x > y$, write 
     \[\oplus_{y'\in G_x y} V^*(y') \in RO(G_x) \] for the corresponding $G_x$-representation isomorphic to $Ind_{G_{xy}}^{G_x} Res^{G_y}_{G_{xy}} V^*(y)$.
    We require that the inductively elements
    \begin{equation}
    \label{eq:numerical-condition-for-existence-of-e-param-2}
    V^*(x) = V^0(x) - \oplus_{\{ y \in Fix(H): y < x\}/G_x} V^*(y) \in RO(G_x)
    \end{equation}
    which is a-priori virtual, is represented  by an actual orthogonal $G_x$-representation.
\end{lemma}
\begin{proof}
    The existence of $\mathfrak{P}'_\tau$ follows by the same equivariant double induction as the proof of Proposition \ref{prop:compatible-perturbation-data-exist}. One first uses this argument to define an initial variant of $\mathfrak{P}'_1$; then running this argument again one defines $\mathfrak{P}'_\tau$ (here $\mathfrak{P}'_0$ can be chosen to arise from $\mathfrak{P}$ by a stabilization with trivial perturbation data because $\mathfrak{P}$ is already regular, and so useing $\mathfrak{P}_0$ defined from $\mathfrak{P}$ will already define $\mathfrak{P}_\tau$ which satisfies the regularity condition of Lemma \ref{lemma:bundle-of-kernels-2} for all sufficiently small $\tau$). Then one runs another equivariant double induction to define the  $\{\tilde{\lambda}_{V(x)}\}, \{\tilde{\lambda}_{V^1(x),V^1(y), \tau_2}\}$, stabilizing the previous $\mathfrak{P}'_\tau$ by a sufficiently large $E$-parameterization with trivial perturbation data at each step such first the surjectivity conditions for $\ker D^{\mathfrak{P}'}_{x,y}$ and $\ker D^{\mathfrak{P}'}_{x}$, and then afterwareds the numerical condition, are satisfied, to define the final $\mathfrak{P}'_\tau$.
\end{proof}

\begin{remark}
Note that we do not require any compatiblity conditions between $\tilde{\lambda}_{xy, \tau_2}$ and $\tilde{\lambda}_{yz, \tau_2}$; this makes these choices straightforward. We essentially only need the $\tilde{\lambda}_{xy, \tau_2}$ to be able to make sense of the numerical conditions in the lemma above, and to have gluing isomorphisms  $\ker D^{\mathfrak{P}'}_x \simeq \ker D^{\mathfrak{P}'}_{x,y} \oplus \ker D^{\mathfrak{P}'}_y$ \emph{only well defined in the homotopy category, and satisfying the relevant associativity conditions in the homotopy category}, which standard from linear gluing theory. A crucial point of the entire construction is that the complications arising from the non-associativity of gluing are systematically avoided at every step. 
\end{remark}

Lemma \ref{lemma:numerical-condition-for-existence-of-e-param} then shows that the vector spaces $V^0(x)$ and operators $g^0_x$ underly an $E^s$-parameterization $\mathcal{E}^0$.
Moreover, the equivariance of the homotopies of operators above implies that there are isomorphisms of $G_{xy}$-representations, 
    \[ V^0(x) = V^0(x,y) \oplus V^0(y). \]
In particular, picking a distinguished pair $(x,y)$ in every $G$-orbit $\{ x,y \in Fix(H)^2 | x > y\}/G$, picking an isomorphism  $\alpha_{xy}:V^0(x,y) \to (Im \bar{\alpha}_{xy})^\perp$, and fixing 
\[ \alpha_{gx gy}: V^0(x,y) \to (Im \bar{\alpha}_{xy})^\perp, \alpha_{gx gy} = g \alpha_{xy} g^{-1}\]
we see that the vector spaces $V^0(x,y)$ underly an $F^s$-parameterization $\mathcal{F}^0$ associated to the $E$-parameterization $\mathcal{E}^0$. 

The $F^s$ parameterization underlying the framing of the thickenings of the morphism objects of $\CC'(H_s, J_s)_{\mathfrak{F}_f} = \CC'(H_s, J_s)_{\mathfrak{P}'}$ that we will produce below will be precisely $\mathcal{F}^0$.

\paragraph{Step 2: Flat connections.}
We now note that since $\mathfrak{P}'$ is a stabilization of $\mathfrak{P}$ by trivial perturbation data, there are canonical projection maps 
\[ \pi_{xy}:(\CC'(H, J)_{\mathfrak{P}'})_T(x,y) \to (\CC'(H, J)_{\mathfrak{P}})_T(x,y) \]
making the former into a trivial vector bundle with fiber the vector space assigned by $\mathcal{F}_f$ to the pair $(x,y)$. Thus, taking the product of the collaring maps for $\CC'(x,y)$ with the corresponding vector space arising from $\mathcal{F}_f$ we get a compatible system of collars for $\CC'(H, J)_{\mathfrak{P}'})$, and pulling back the data $\mathfrak{P}'_\tau$ under $\pi_{xy}$ defines data satisfying the same compatibility coonditions for $\CC'(H, J)_{\mathfrak{P}'})$. In particular, we write $\ker D^{\mathfrak{P}'}_{\tilde{u}; x,y, \tau}$ for the bundles defined via this pulled-back data (which is just the same as the pullback of the corresponding bundles); and let us write $T(x,y)$ now for the thickenings of $\CC'(H, J)_{\mathfrak{P}'}$. We have isomorphisms of bundles
\begin{equation}
    \label{eq:associative-deformation-of-bundles}
    \ker D^{\mathfrak{P}'}_{\tilde{u}; x,y,\tau} \oplus \ker D^{\mathfrak{P}'}_{\tilde{u};y,z,\tau} \to \ker D^{\mathfrak{P}'}_{x,z,\tau}|_{Im T(x,y) \times T(y,z) \subset T(x,z)}
\end{equation} 
which satisfy obvious associativity and $G$-equivariance conditions, and which are continuous in $\tau \in [0,1]$. Moreover, the fact that $D^{\mathfrak{P}'}_{\tilde{u}; x,y,1}$ is independent of $\tilde{u}$ for $\tilde{u} \in T(x,y)^\square$, together with the fact that there is an isomorphism of $G_{xy}$-representations 
\[ \ker D^{\mathfrak{P}'}_x = V(x) \simeq \ker D^{\mathfrak{P}'}_{x,y} \oplus \ker D^{\mathfrak{P}'}_y \]
(due to the surjectivity of the operators $D^{\mathfrak{P}'}_{x,y; \tau_2}$ and the equivariance properties of the latter), means that we can choose $G_{x,y}$-equivariant isomorphisms 
\[ a_{xy}: \ker D^{\mathfrak{P}'}_{\tilde{u}; x,y,1}|_{T(x,y)^\square} = \underline{\ker D^{\mathfrak{P}'}_{xy}} \simeq \ker \simeq \underline{V^0(x,y)} \]
which satisfy 
\[ a_{gx, gy}(gv) = g a_{x,y}(v). \]
by simply choosing corresponding isomorphisms 
\begin{equation}
    \label{eq:framings-over-complement-of-collar}
    \ker D^{\mathfrak{P}'}_{xy} \to V^0(x,y). 
\end{equation}
Using the first run of the equivariant double induction, we choose these such that in the diagram 
\begin{equation}
    \label{eq:correct-connected-component}
\begin{tikzcd}
    V^0(x) \ar[equal, d]& \ar[l]V^0(x,y)\ar[d]  & \oplus & V^0(y)\ar[equal, d] \\
    \ker D^{\mathfrak{P}'}_{x} & \ar[l] \ker D^{\mathfrak{P}'}_{xy} & \oplus &  \ker D^{\mathfrak{P}'}_{y}
\end{tikzcd}
\end{equation}
(where the bottom map is induced by gluing of operators, and is well-defined up to homotopy)
the two maps from the top right to the bottom left agree up to homotopy of $G_{x,y}$-equivariant isometries. Now, we wish to extend the framings  $a_{xy}$ to framings

\begin{equation}
\label{eq:framing-of-vertical-tangent-bundle}
 a_{xy}: \ker D^{\mathfrak{P}'}_{\tilde{u}; x,y,1} \simeq \underline{V^0(x,y)}
\end{equation} 
defined over all of $T(x,y)$ satisfying the same equivariance conditions, such that there are commutative diagrams

\begin{equation}
\label{eq:compatible-framings-of-vertical-tangent-bundle}
\begin{tikzcd}
    \ker D^{\mathfrak{P}'}_{\tilde{u}; x,y,1} \boxplus \ker D^{\mathfrak{P}'}_{\tilde{u}; y,z,1} \ar[r]\ar[d] & \ker D^V(x,z)\ar[d] \\
    \underline{V^0(x,y)} \boxplus \underline{V^0(y,z)} \ar[r] & \underline{V^0(x,z)}
\end{tikzcd}
\end{equation}
over the inclusions of boundary strata, where the maps on the bottom come from the chosen $F$-parameterization $\mathcal{F}^0$.
Once this is done, by choosing trivializations 
\begin{equation}
    \label{eq:homotopy-of-trivializations-of-vertical-bundle}
    \xi_{xy}: \ker D^{\mathfrak{P}'}_{\tilde{u}; x,y,\tau} \simeq \pi_1^*\ker D^{\mathfrak{P}'}_{\tilde{u}; x,y,1}
\end{equation}
where $\pi_1: T(x,y) \times [0,1] \to T(x,y) \times \{1\}$ is the projection (satisfying the obvious compatiblity conditions relative to the isomorphisms \eqref{eq:associative-deformation-of-bundles} and obvious equivariance conditions) one gets a an embedded framing for $T(x,y)$ with with $F_{ss}^2$-parameterization $(\mathcal{F}^0, \mathcal{F}^1_T)$. Compatible choices of trivializations \eqref{eq:homotopy-of-trivializations-of-vertical-bundle} are straightforward to construct via an equivariant double induction so we focus on the extensions of $a_{x,y}$ to $T(x,y)$.  

To construct $a_{x,y}$, one first notes that if, for $z > t$ in $\CC$,  the condition \eqref{eq:compatible-framings-of-vertical-tangent-bundle} has been achieved for all $(x,y) \in G(z,t)$, then one has a well-defined trivialization $a_{xy}:  \ker D^{\mathfrak{P}'}_{\tilde{u}; x,y,1}|_{\partial T(x,y)} \simeq \underline{V^0(x,y)}$
by combining the previously chosen $a_{xy}$ with the maps associated to $\mathcal{F}^0$. The problem is to make sure that we can consistently glue this choice to the choice of $a_{xy}$ made over $T(x,y)^\square$. The construction of the previous section was done entirely with the purpose of making sure this is possible: the crucial property is that for $x_0 > x_1 > \ldots > x_r$ in $\CC$, over $U(x_0, \ldots, x_r) \subset T(x_0, x_r)$, we have that 
\[ \ker D^{\mathfrak{P}'}_{\tilde{u}; x,y,1}|_{U(x_0, \ldots, x_r)} = \epsilon_{x_0, \ldots, x_r}^* \ker D_{x_0, \ldots, x_r; \tau'_2} \]
where $\ker D_{x_0, \ldots, x_r; \tau'_2}$ denotes the bundle of kernels of the corresponding operatosr over $[0,1]^{[r-1]}$. In the process of the equivariant double induction, we will choose extensions of $a_{xy}$ such that $a_{x_-x_r}|_{U(x_0, \ldots, x_r)} =  \epsilon_{x_0, \ldots, x_r}^* \tilde{a}_{x_0, \ldots, x_r}$ for 
\[\tilde{a}_{x_0, \ldots, x_r}: \ker D_{x_0, \ldots, x_r; \tau'_2}  \to \underline{V^0(x_0, x_r)}\]
where the $\tilde{a}_{x_0, \ldots, x_r}$ satisfy the corresponding obvious compatiblity and equivariance conditions. The choices of $\tilde{a}_{x_0, \ldots, x_r}$ are easy to construct at each step of the inner loop of the equivariant double induction by induction up the face poset of $T(x_0, x_r)$, with $\tilde{a}_{x_0, \ldots, x_r}$ existing at each step because of the homotopies \eqref{eq:correct-connected-component}, which force the $\tilde{a}_{x_0, \ldots, x_r}|_{[0,1]^{[r-1]} \cap \partial \R^{[r-1]}_+}$ (which is forced by the compatiblity conditions and the previous choices) and $(\tilde{a}_{x_0,\ldots, x_r})_{(1, \ldots, 1)}$ (which comes from the isomorphism chosen to be the the vertical arrow in \eqref{eq:correct-connected-component}) homotopic to one another under any equivariant trivialization of $\ker D_{x_0, \ldots, x_r; \tau'_2}$.


We have established: 
\begin{proposition}
\label{prop:framings-for-floer-moduli-spaces-exist}
    Let $\CC'(H, J)$ be the virtual smoothing constructed as in Section \ref{sec:global-charts}. Then if the $S^1$-equivariant polarization class vanishes then a free stabilization of $\CC'(H, J)$ admits an equivariant embedded framing. 
\end{proposition}

\begin{remark}
    The reader may note that the whole construction of this paper never asks for any results about linear or nonlinear associativity of gluing. The statement that `the obstruction to moduli spaces being framed is controlled by the polarization classes restricted to the the interiors of the moduli spaces of Floer trajectories' is simply Proposition 12 of \cite{cohen2007floer}; however, that paper does not actually explain how to perform an inductive sequence of stabilizations of the category in order to frame the relevant flow categories (see diagram (28) therein, where rigorously translating this statement together with the triviality of the polarization class into a statement about the existence of compatible trivializations of tangent bundles involves a careful analysis of the (non)-associativity of linear Fredholm gluing, which is not entirely trivial to perform). The construction above makes this clear, and moreover makes it clear that this procedure can be done equivariantly. Below we will prove construct framings of additional flow categories satisfying various kinds of compatibility conditions using this strategy. In later works \cite{large2021spectral, cote-kartal-1}, the notion of an `isomorphism of virtual bundles' is used in the discussion of framings, but no definition is given and no definition is known to this author that makes the construction of framings completely elementary. Similarly, in this paper, we keep track of the distinction between 
    \[ F_V \Sigma^V X \text{ and } X\]
    which are canonically weakly equivalent but in fact denote distinct orthogonal spectra. 
\end{remark}

\subsection{Generalization of framing constructions: continuation maps and cyclotomic compatibility.}

\paragraph{Behavior of framings of $\CC'(H, J)$ under stabilization}
\label{lemma:framings-under-stabilization}
\begin{lemma}
    Suppose that $\CC'(H, J)$ is defined by Proposition \ref{prop:smooth-flow-category-floer-trajectories}, and $\CC'(H, J)_{\mathcal{E}^f}$ is given a framing via the construction of Proposition \ref{prop:framings-for-floer-moduli-spaces-exist}. Then given any other $E'$-parameterization $\mathcal{E}$, the framing of $\CC'(H, J)_{\mathcal{E}^f \oplus \mathcal{E}}$ induced from the framing of $\CC'(H, J)_{\mathcal{E}^f}$ via stabilization by $\mathcal{E}$ can be produced as in Proposition \ref{prop:framings-for-floer-moduli-spaces-exist} as well. 
\end{lemma}
\begin{proof}
   Recall that we can realize $\CC'(H, J)_{\mathcal{E}^f \oplus \mathcal{E}}$ by replacing the perturbation data $\mathfrak{P}$ defining $\CC'(H, J)_{\mathcal{E}^f}$ with the perturbation data $\mathfrak{P}_\mathcal{E}$ given by stabilizing $\mathfrak{P}$ by the $F$-parameterization $\mathcal{F}$ associated to $\mathcal{E}$, equipped with trivial perturbation data (Definition \ref{def:trivial-stabilization-of-perturbation-data}). Writing $(\{V^0(x,y)\}, \{V^1_{can}(x,y) \oplus V^1(x,y)\}) = \mathcal{F}_0$ and $(\{V^0(x)\}, \{V^1_{can}(x) \oplus V^1(x)\})=\mathcal{E}_0$ for for the $F^s_2$ and $E^s_2$-parameterizations associated to  $\CC'(H, J)_{\mathcal{E}^f}$, and writing 
   \begin{equation}
       \label{eq:f-param-decomposition}
       V^i_{\mathcal{E}}(x,y) = V^i(x,y) \oplus V_{\mathcal{F}}(x,y), i=0,1
   \end{equation}
   \begin{equation}
       \label{eq:e-param-decomposition}
       V^i_{\mathcal{E}}(x) = V^i(x) \oplus V_{\mathcal{F}}(x), i=0,1
   \end{equation} for the vector spaces associated to the $F$ and $E$- parameterizations $(\mathcal{F}', \mathcal{E}')$ coming from stabilizing each of $(\mathcal{F}_0, \mathcal{E}_0)$ by $\mathcal{E}$, we can define collars for $\CC'(H, J)_{\mathcal{E}}$ by taking products of previously chosen collars with the corresponding vector space associated to $\mathcal{F}$, setting
   \[ \tilde{\lambda}^\tau_{V^1_{\mathcal{E}}(x,y)}(\tilde{u})(v,w) = \tilde{\lambda}^\tau_{V^1(x,y)}(\pi_{xy}(\tilde{u})(v) \]
   \[\tilde{\lambda}_{V^1_{\mathcal{E}}(x)}(v, w) = \lambda_{V^1(x)}(v), \tilde{\lambda}_{V^1_{\mathcal{E}}(x), V^1_{\mathcal{E}}(y)}(v, w) = \lambda_{V^1(x), V^1(y)}(v)\]
   (here we use the decompositions \eqref{eq:f-param-decomposition} and \eqref{eq:e-param-decomposition}).
   With this choice the vector spaces defined in \eqref{eq:define-v0-parameterization-from-kernels-of-asymptotic-operators} are exactly  $V^0_{\mathcal{E}}(x)$ and $V^1_{\mathcal{E}}(x,y)$, as desired; and we can pick the $E$ parameterization subsequently produced by Lemma \ref{lemma:numerical-condition-for-existence-of-e-param} to be the first factor of $\mathcal{E}'$. It remains to choose the isomorphisms \eqref{eq:framing-of-vertical-tangent-bundle}; here, in the inductive procedure for constructing these, we simply choose $a_{xy}$ and $\xi_{xy}$ to be products of the previous choices with the corresponding identity map on the vector space arising from $\mathcal{F}$. 

\end{proof}

\paragraph{Behavior of framings of $\CC'(H, J)$ under extension of integralization data}

\begin{lemma}
\label{lemma:framings-under-restratification}
    Suppose we are in the situation of Lemma \ref{lemma:extend-perturbation-data-across-restratifications}: we have an integralization datum $\bar{\mathfrak{A}}$ for $\CC'(H, J)$ and an integralization datum $\bar{\mathfrak{A}}_{new}$ which is an extension of $\bar{\mathfrak{A}}$, and there are perturbation data $\mathfrak{P}$ and $\mathfrak{P}_{new}$ with underlying integralization data $\bar{\mathfrak{A}}$ and $\bar{\mathfrak{A}}_{new}$ such that (a restriction of the) the virtually smooth flow category $\CC'(H, J)_{\mathfrak{P}_{new}}$ is is isomorphic to (a restriction of) the restratification $\CC'(H, J)_{\mathfrak{P}, r}$ of $\CC'(H, J)_{\mathfrak{P}}$ associated to the extension of integralization data, with the isomorphism given as in Lemma \ref{lemma:extend-perturbation-data-across-restratifications}. Suppose we are also given an $E$-parameterization $\mathcal{E}^f$ such that $(\CC'(H, J)_{\mathfrak{P}})_{\mathcal{E}^f}$ has a framing via the construction of Proposition \ref{prop:framings-for-floer-moduli-spaces-exist}. Then there is a preferred associated framing for $(\CC'(H, J)_{\mathfrak{P}_{new}})_{\mathcal{E}^f}$. 
\end{lemma}
\begin{proof}
    In the proof of Lemma \ref{lemma:extend-perturbation-data-across-restratifications}, we see that the data $\mathfrak{P}$ defines a perturbation datum $\mathfrak{P}'_0$ over certain sets $U(x,y)^Z \subset \paramspace'(x,y)^Z$, where $\paramspace'(x,y)$ are the parameter spaces for $\bar{\mathfrak{A}}_{new}$ (see \eqref{eq:induced-system-of-perturbation-data}). The final perturbation datum produced by Lemma \ref{lemma:extend-perturbation-data-across-restratifications} agrees with this one on $U(x,y)^Z$, and the argument of that lemma then shows that the preimages of $U(x,y)^Z$ in $\CC'(H,J)_{\mathfrak{P}_{new}}$ agree with a restriction of the corresponding restratification $\CC'(H, J)_r$ of $\CC'(H, J)$; moreover, we have a canonical identification of bundles
    \begin{equation}
        \label{eq:identification-of-vertical-bundles-under-restrat}
        ker D^{V_{\mathfrak{P}}}_u \to \ker D^{V_{\mathfrak{P}'_{new}}}_{w(u)}
    \end{equation} covering the identification $w$ \eqref{eq:identification-of-thickening-with-restratification} of thickenings. 

    Now, given a compatible system of collars for $\CC'$, there is an obvious associated compatible system of collars for any  restratification of $\CC'$; using this latter compatible system of collars for $\CC'(H, J)_{\mathfrak{P}_{new}}$ and setting
    \[ \tilde{\lambda}^\tau_{V^1_{\mathfrak{P}'_{new}}(x,y)}(w(\tilde{u}), v) = 
    \tilde{\lambda}^\tau_{V^1_{\mathfrak{P}}(x,y)}(\pi(\tilde{u}),v)\]
    (where $\pi$ is the projection on thickenings from the restratification to the original category)
    Here we use that the vector spaces $V^1_{\mathfrak{P}'_{new}}(x,y)$ and $V^1_{\mathfrak{P}}(x,y)$ are simply equal. This extends the previous canonical identification \eqref{eq:identification-of-vertical-bundles-under-restrat} to one for all $\tau \in [0,1]$.

    Similarly, we define the functions $\epsilon_{[x], \pm}$ of Lemma \ref{lemma:bundle-of-kernels-2} for $\CC'(H, J)_{\mathfrak{P}'_{new}}$, and subsequently the homotopies $A^{\tilde{u}}_{s, \tau}$, by pulling them back from the projections from the restratification of $\CC'(H, J)_{\mathfrak{P}}$ to $\CC'(H, J)_{\mathfrak{P}}$ and subsequently utilizing the identifications $w$ (which are covered by maps of universal curves);  keeping $\tilde{\lambda}_{V^1(x)}$ and $\tilde{\lambda}_{V^1(x), V^1(y)}$ unchanged; these choices satisfy all required conditions for $\tilde{\lambda}^\tau_{V^1_{\mathfrak{P}'_{new}}(x,y)}$ and $\tilde{\lambda}_{V^1(x)}$, and \eqref{eq:define-v0-parameterization-from-kernels-of-asymptotic-operators} defines the same $V^0(x,y)$ and $V^1(x,y)$ as for $\CC(H, J)$. The final $E^s$-parameterization underlying the framing of the thickenings that we will give to $\CC(H, J)_{r, \mathcal{E}_1}$ is the same as the previously chosen parameterization. The remaining the remaining choice needed to be in Step 2 of Proposition \ref{prop:framings-for-floer-moduli-spaces-exist} are those of the trivializations \eqref{eq:framing-of-vertical-tangent-bundle} and \eqref{eq:homotopy-of-trivializations-of-vertical-bundle}, which we make using the canonical identifications \eqref{eq:identification-of-vertical-bundles-under-restrat} (or rather, their corresponding analogues for $\tau \in [0,1]$).
\end{proof}

\begin{remark}
    Note that we do \emph{not} state in Lemma \ref{lemma:framings-under-restratification} that the framing on the flow category $\CC'(H, J)_{\mathfrak{P}_{new}}$ arises from the framing on $\CC'_{\mathfrak{P}}$ by transport of structure from the framing on $\CC'(H, J)_{\mathfrak{P}, r}$  (with the framing on the latter coming from Lemma \ref{lemma:restratifications-exist}). Indeed, these framings will not agree; their restrictions to the vertical tangent bundles $\ker D_u^{V_{\mathfrak{P}}}$ will agree, but the components coming from the framing of the parameter spaces will disagree. However, we will use the embedding coming from Lemma \ref{lemma:restratifications-exist} to make sense of the Floer homotopy type of $\CC'(H, J)_{\mathfrak{P}_{new}}$; since the framings of the obstruction bundles of $\CC'(H, J)_{\mathfrak{P}_{new}}$ and $\CC'(H, J)_{\mathfrak{P}, r}$ \emph{do} agree under transport of structure and since the floer homotopy type only depends on the embedding, we will be able to use the conclusion of Proposition \ref{prop:simple-restrat-effect-on-homotopy-type} in this setting (see e.g.  Proposition \ref{prop:continuation-homotopy-geometry}). Thus, this discrepancy does not affect our ability to prove the invariance statements that we require in the geometric setting of Section \ref{sec:floer-homotopy}. A somewhat tedious geometric argument can be used to show that these two framings are in fact connected by an (essentially canonical) homotopy of framings, but we will not use this in the paper; see also Remark \ref{rk:how-do-framings-of-thickenings-play-a-role}.  
\end{remark}

\paragraph{Framings for continuation map moduli spaces.}
We now explain how to extend the method of the previous section produce framings of the categories $\CC'(H_s, J_s)$ associated to a convex continuation datum. The basic setup is that of Proposition \ref{prop:continuation-map-compatible-smoothings}: there exist $E^s$-parameterizations $\mathcal{E}_\pm$ of $\CC(H_\pm, J_\pm)$ such that $\CC'(H_s, J_s)|_{\CC(H_\pm, J_\pm)} \simeq \CC'(H_s, J_s)_{\mathcal{E_\pm}}$ as a virtually smooth flow category. Now Proposition \ref{prop:framings-for-floer-moduli-spaces-exist}, via Lemma \ref{lemma:perturbation-data-and-numerical-condition-for-framing}, defines $E$-parameterizations $\mathcal{E}_\pm^f$ of $\CC(H_\pm, J_\pm)$ such that $\CC'(H_\pm, J_\pm)_{\mathcal{E}_\pm^f}$ are given framings. 

\begin{lemma}
\label{lemma:framings-continuation-map}
    With notation as in the above paragraph, there is an $E^s$-parameterization $\mathcal{E}^f$ of $\CC(H_s, J_s)$ such that $\CC'(H_s, J_s)_{\mathcal{E}^f}$ has a framing, and moreover we can write $\mathcal{E}^f = \mathcal{E}^f_0 \oplus \mathcal{E}^f_+ \oplus \mathcal{E}^f_-$ where $\mathcal{E}^f_\pm$ are some $E^s$-parameterizations of $\CC(H_s, J_s)$ which agree with the previously mentioned $\mathcal{E}^f_\pm$ when restricted to $\CC(H_\pm, J_\pm)$, and under this isomorphism we have that the induced framing on $\CC'(H_s, J_s)_{\mathcal{E}^f}|_{\CC(H_\pm, J_\pm)}$ is agrees with the framing induced via stabilization by $\mathcal{E}^f_0\oplus \mathcal{E}^f_\mp \oplus \mathcal{E}_\pm$ of $\CC'(H_\pm, J_\pm)_{\mathcal{E}_\pm^f}$. 
\end{lemma}
\begin{proof}
    We must explain how we are framing $\paramspace(x,y)$ in this setting. After making this choice, the rest of the construction
    follows the construction of Proposition \ref{prop:framings-for-floer-moduli-spaces-exist}, using Lemma \ref{lemma:framings-under-stabilization} in the process of the double inductions of Lemma \ref{lemma:perturbation-data-and-numerical-condition-for-framing} and Step 2 of Proposition \ref{prop:framings-for-floer-moduli-spaces-exist} in the proof of Lemma \ref{lemma:perturbation-data-and-numerical-condition-for-framing} to define all choices made whenever we are considering a pair $(x,y)$ where both $x$ and $y$ are in $\CC(H_+, J_+)$ or in $\CC(H_-, J_-)$.

    Let $x \in Fix(H_-)$ and $y \in Fix(H_+)$. There is a submersion
    \[\pi_f: \paramspace(x,y) = \overline{Conf}^o_{\underline{\mathfrak{A}}(x,y)} \subset \overline{Conf}_{\mathfrak{A}(x,y)}^{\R} \to \overline{Conf}_{\mathfrak{A}(x,y)}\]
    given by composing the inclusion with the map from the universal curve to the parameter space. We have already chosen trivializations of $T\overline{Conf}_{|\mathfrak{A}(x,y)}|$ earlier, and the vertical tangent bundle of $\pi_f$ is canonically trivialized (one can identify each fiber with a copy of $\R$ or $\R \sqcup \R$ up to $s$-translation by moving the point on the corresponding component of the universal curve; the vertical tangent bundle is trivialized by saying that the tangent vector corresponding to moving $s_{A_c}$ to the right with speed $1$ is sent to $1 \in \R$.). We then get a trivialization 
    \[ T\overline{Conf}^o_{\underline{\mathfrak{A}}(x,y)} = \pi_f^* T\overline{Conf}_{\mathfrak{A}(x,y)} \oplus T^v_{\pi_f} \simeq \underline{T\R^{\mathfrak{A}'(x,y)}} \oplus \underline{\R} \simeq \underline{T\R^{\underline{\mathfrak{A}'}(x,y)}}\]
    where the last isomorphism identifies the $\underline{\R}$ factor with the $\underline{T\R^{\{A_c\}}_+}$ factor. This convention makes the framings satisfy the conditions of Definition \ref{def:framing-of-flow-category}, which uses the ordering on the set $S_{\bar{A}}$ in the compatibility condition \eqref{eq:framing-diagram-1}.
\end{proof}

\paragraph{Cyclotomic compatibility of framings.}
Now, we explain how to extend framings compatibly upon replacing the (potentially $C_\ell$-invariant) data $(H_s, J_s)$ with ($C_{k \ell}$-invariant data) $(H_s^{\#k}, J_s^{\#k})$, as in Proposition \ref{prop:cyclotomic-compatibility-all-perturbation-data-choices}. As in that Proposition we have that $\CC'(H_s^{\#k}, J_s^{\#k})$ is a virtually smooth flow category, and $\CC'(H_s^{\#k}, J_s^{\#k})^{C_k} \simeq \CC'(H_s, J_s)_{r, \mathcal{E}_k}$ as virtually smooth $C_{k\ell}/C_k$-equivariant flow categories, where $\CC'(H_s, J_s)_{r, \mathcal{E}_k}$ is the stabilization by an $E^s$-parameterization $\mathcal{E}_k$ of a restratification $\CC'(H_s, J_s)_r$ of $\CC'(H_s, J_s)$ associated to an extension of integralization data.

\begin{lemma}
\label{lemma:framings-cyclotomic-compatibility}
    Suppose that there is an $E^s$-parameterization $\mathcal{E}^f$ of $\CC(H_s, J_s)$ such that $\CC'(H_s, J_s)_{\mathcal{E}^f}$ is given a framing via the construction of Lemma \ref{lemma:framings-continuation-map}. Then there exists an $E^s$-parameterization $\mathcal{E}^f_k = \mathcal{E}^f_{k, 0} \oplus \mathcal{E}^f$ of $\CC(H_s^{\#k}, J_s^{\#k})$, where $\mathcal{E}^f$ here refers to some $E^s$-parameterization which agrees with the previously specified $\mathcal{E}^f$ when restricted to $\CC(H_s, J_s)$, and an equivariant framing of $\CC'(H_s^{\#k}, J_s^{\#k})_{\mathcal{E}^f_k}$ such that the induced equivariant framing of $\CC'(H_s^{\#k}, J_s^{\#k})_{\mathcal{E}^f_k}^{C_k}$ agrees, under its isomorphism with $\CC'(H_s, J_s)_{r, \mathcal{E}_k \oplus  \mathcal{E}^f_{k, 0} \oplus \mathcal{E}^f}$, with the framing on the latter produced by applying Lemma \ref{lemma:framings-under-restratification} to produce a framing of $\CC'(H_s, J_s)_{r, \mathcal{E}^f}$ from the given framing of $\CC'(H_s, J_s)_{\mathcal{E}^f}$, and subsequently taking the framing produced via stabilization by $\mathcal{E}_k \oplus  \mathcal{E}^f_{k, 0}$. 
\end{lemma}
\begin{proof}
    This follows identically to that of Lemma \ref{lemma:framings-under-restratification}, just as the proof of Lemma \ref{lemma:extend-perturbation-data-across-restratifications} is a simpler variant of the proof of Proposition \ref{prop:cyclotomic-compatibility-all-perturbation-data-choices}. One gets  collars on $\CC'(H^\#k_s, J^\#k_s)^{C_k}$ from those for the corresponding restratification of $\CC'(H_s, J_s)$, one then extends the collars to all of $\CC'(H^\#k_s, J^\#k_s)$ via Lemma \ref{lemma:collars}. Subsequently one first defines that the functions $\epsilon_{[x], \pm}$ in the proof of Lemma \ref{lemma:bundle-of-kernels-2}, and subsequently the homotopies $A^{\tilde{u}}_{s, \tau}$ of Lemma \ref{lemma:bundle-of-kernels-2} can be extended to the thickenings of $\CC'(H^\#k_s, J^\#k_s)$, over the thickenings of $\CC'(H^\#k_s, J^\#k_s)^{C_k}$, and then extends each of these to all of $\CC'(H^\#k_s, J^\#k_s)$. One then produces $\mathfrak{P}'_\tau$ for $\CC'(H^\#k_s, J^\#k_s)$ by pulling back those for $\CC'(H_s, J_s)$ to the restratification and then transporting them to $\CC'(H^\#k_s, J^\#k_s)^{C_k}$ to get $\mathfrak{P}'_{\tau, pre}$; subsequently, an equivariant double induction as in Lemma \ref{lemma:change-of-perturbation-datum-and-restratification-exists} allows us to extend these choices to the compatible system of thickening-dependent perturbation data $\mathfrak{P}'_{\tau}$ defined over all of $\CC'(H^\#k_s, J^\#k_s)$, taking care to arrange that (A) the independence condition \eqref{eq:independence-condition} is still satisfied, and (B) that the new value of $\tilde{\lambda}_{V^1(x,y)}(s)$ is obtained from the previous value by setting it to be zero on the vector space that $V^1(x,y)$ is stabilized by when going from $\mathfrak{P}'_{\tau, pre}$ to $\mathfrak{P}'_{\tau}$. One gets new $\tilde{\lambda}_{V^1(x^k)}$ and $\tilde{\lambda}_{V^1(x^k), V^1(y^k)}$ by stabilizing previous ones sufficiently by trivial data, and defines the newly required $\tilde{\lambda}_{V^1(x), V^1(y)}$ whenever one of $x$ or $y$ is in $Fix(H_s^{\#k}) \setminus Fix(H_s)$ such that the conditions of Lemma \ref{lemma:perturbation-data-and-numerical-condition-for-framing} are still satisfied.  One transports the trivializations $\xi_{xy}$ to $\CC'(H^\#k_s, J^\#k_s)^{C_k}$ and then extends them to $\CC'(H^\#k_s, J^\#k_s)$; to choose the $a_{xy}$ \eqref{eq:framing-of-vertical-tangent-bundle} one notes that one only needs to make new choices of isomorphisms $a^0_{xy}$ \eqref{eq:framings-over-complement-of-collar} for one of $x,y$ in $Fix(H_s^{\#k}) \setminus Fix(H_s)$, which one chooses using the same criterion; then subsequently choosing $\tilde{a}_{x^k_0, \ldots, x^k_r}$ to be the sums of the previously chosen $\tilde{a}_{x_0, \ldots, x_r}$ with the identity map on the corresponding vector space associated with the stabilization, and otherwise noting that as before, the diagram \eqref{eq:correct-connected-component} means that $\tilde{a}_{x_0, \ldots, x_r}$ will still exist for the remaining choices (i.e. whenever one of the $x_j \in Fix(H^{\#k}_s) \setminus Fix(H_s)$. )
    
    Here it is crucial that data needed to produce the the part of framing coming from the framings of $\paramspace(x,y)$ never interacts with the data needed to produce the framings of the vertical bundles of $T(x,y) \to \paramspace(x,y)$.
\end{proof}

\paragraph{Framings for homotopies of continuation map moduli spaces.}

Finally, in this section we explain how to construct framings for the smoothings of the flow categories associated to thomotopies of continuation maps $(H_h, J_h)$. We will also explain the analog of the cyclotomic compatibility statement (Lemma \ref{lemma:framings-cyclotomic-compatibility}), since that is what we will use to prove invariance of the cyclotomic structure on symplectic homology. 

\begin{lemma}
\label{lemma:framings-continuation-map-homotopy}
    Suppose we are in the situation of Proposition \ref{prop:cyclotomic-compatibility-all-perturbation-data-choices}. Namely, there is a perturbation datum $\mathfrak{P}_h$ for the homotopy datum $(H_h, J_h)$ with underlying integralization datum $\bar{\mathfrak{A}}$  Moreover, write $(H_a, J_a)_{a \in \{1, 2, 3, 4\}}$ and $(H_{ab}, J_{ab})_{(a,b) \in \{(1,2), (1,3), (2,4), (3,4)}\})$ for the underlying Floer data and continuation data associated to $(H_h, J_h)$. Write $\CC'(H_{a}, J_{a})$ and $\CC'(H_{ab}, J_{ab})$ for the virtual smoothings associated using the restruction of $\mathfrak{P}_h$ to $\CC(H_a, J_a)$ and $\CC(H_{ab}, J_{ab})$. Suppose that we have
    \begin{itemize}
        \item $E$-parameterizations $\mathcal{E}^f_a$ for $\CC(H_a, J_a)$ such that $\CC(H_a, J_a)_{\mathcal{E}^f_a}$ is given a framing via Proposition \ref{prop:framings-for-floer-moduli-spaces-exist}, with fixed extensions to $\CC(H_h, J_h)$, together with 
        \item $E^s$-parameterizations $\mathcal{E}^f_{ab}$ for $\CC(H_{ab}, J_{ab})$ constructed as in Lemma \ref{lemma:framings-continuation-map} (with associated decompositions $\mathcal{E}^f_{ab} = \mathcal{E}^f_{0, ab} \oplus \mathcal{E}^f_a \oplus \mathcal{E}^f_{b}$ such that $\CC'(H_{ab}, J_{ab})_{\mathcal{E}^f_{ab}}$ is given a framing as in Lemma \ref{lemma:framings-continuation-map}. Say we also fix extensions of $\mathcal{E}^f_{ab}$ to $\CC(H_h, J_h)$. 
    \end{itemize}
    Let us stabilize $\CC'(H_h, J_h)$ to $\CC'(H_h, J_h)_{\R}$ by a canonical parameterization (Definition \ref{def:canonical-parameterization})  $f(A)=\delta_{A,\#{A \in \mathfrak{A} : A < A^2_c}}$; in other words, one stabilizes all moduli spaces containing curves with $s_{A^2_C}$ by a rank one vector space, and otherwise no stabilization is performed. Then there is an $E^s$-parameterization $\mathcal{E}^f_{h, 0}$ of $\CC'(H_h, J_h)$ such that, writing 
    \begin{equation}
        \label{eq:decomposition-e-parameterization-continuation-map-homotopy}
        \mathcal{E}^f_{h}  = \mathcal{E}^f_{h, 0} \oplus \bigoplus_{(a, b) \in \{(1,2), (1,3), (2,4), (3,4)} \mathcal{E}^f_{ab, 0} \oplus \bigoplus_{i \in \{1, 2, 3, 4\}} \mathcal{E}^f_i 
    \end{equation}
    we have that $\CC'(H_h, J_h)_{\R \oplus \mathcal{E}^f_h}$ has a framing such that the framing associated to $\CC'(H_h, J_h)_{\R \oplus \mathcal{E}^f_h}|_{\CC'(H_{ab}, J_{ab})}$ are induced from the given framings of $\CC'(H_{ab}, J_{ab})_{\mathcal{E}^f_{ab}}$ via stabilization for $(a, b) \in \{(1,2), (1,3), (2, 4) \}$, but for $(a, b) = (3, 4)$, the framing $\CC'(H_h, J_h)_{\R \oplus \mathcal{E}^f_h}|_{\CC'(H_{ab}, J_{ab})}$ differs from the framing induced from $\CC'(H_{ab}, J_{ab})_{\mathcal{E}^f_{ab}}$ via a \emph{sign}: while  the obstruction bundles of $\CC'(H_h, J_h)_{\R \oplus \mathcal{E}^f_h}|_{\CC'(H_{34}, J_{34})}$ are framed as induced from the stabilization, the framing of the tangent spaces to the thickenings of $\CC'(H_h, J_h)_{\R \oplus \mathcal{E}^f_h}|_{\CC'(H_{ab}, J_{ab})}$ differs by a flip in the $\R$ factor (associated to the canonical stabilization by $\R$) from the one associated to the stabilization. 
\end{lemma}

\begin{figure}[h!]
    \centering
    \includegraphics[width=\textwidth]{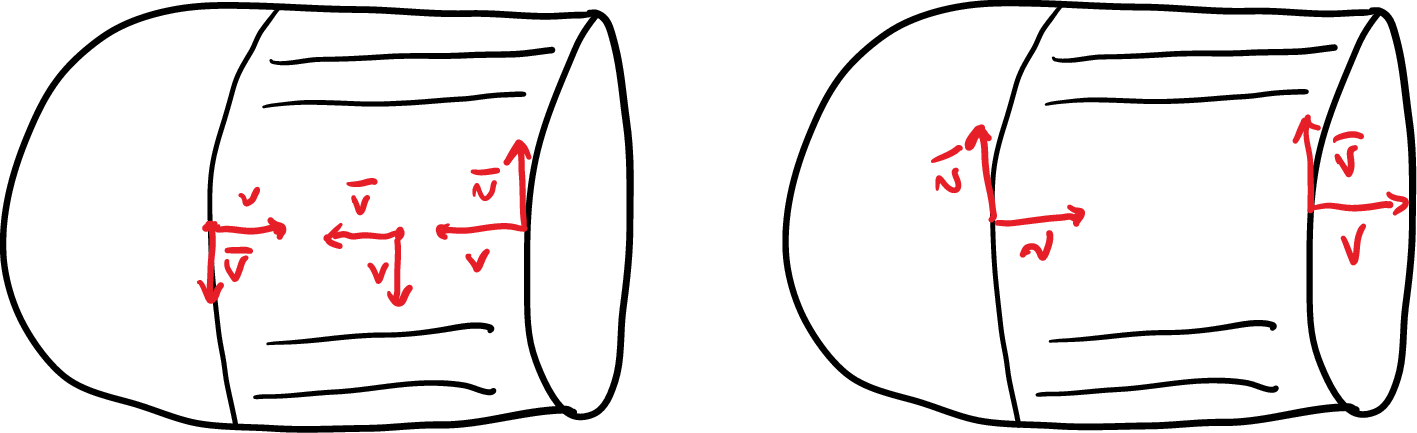}
    \caption{Given two framed manifolds with boundary, we take the product of each manifold with $\R$ and use the extra $\R$ coordinate together with a normal coordinate vector field to insert a rotation in the framing such that we can glue these two manifolds along their common boundary stratum to a new framed manifold. This is the construction used in Lemma \ref{prop:cyclotomic-compatibility-all-perturbation-data-choices}. }
    \label{fig:rotate-framing}
\end{figure}

\begin{proof} 
 Define $\paramspace'(x,y)$ to be $\R \times \paramspace(x,y)$ whenever $y \in Fix(H_4)$ and $x \notin Fix(H_4)$, and otherwise define the parameter spaces to be as before. These parameter spaces manifestly have embedding maps $\iota_{xzy}$ like in Section \ref{sec:parameter-spaces-for-global-charts}; Pulling back universal curves and perturbation data to these new parameter spaces under the projection and replacing previous uses of $\paramspace(x,y)$ with $\paramspace'(x,y)$ exactly replaces $\CC'(H_h, J_h)$ with the canonically stabilized $\CC'(H_h, J_h)_\R$ as in the statement of the lemma. We must explain how we are framing of $\paramspace'(x,y)$; the rest proceeds as in Lemma \ref{lemma:framings-continuation-map}.

First, we choose framings of $\overline{Conf}^o_{\bar{A}(x,y)} \subset \overline{Conf}_{\mathfrak{A}(x,y)}^{\R_2}$ (as defined in Section \ref{sec:continuation-homotopy-geometry}) in the same fashion that we framed $\overline{Conf}^o_{\bar{A}(x,y)}$ in the proof of Lemma \ref{lemma:framings-continuation-map}.

We can restrict the framing $\overline{Conf}^o_{\bar{A}(x,y)}$ to $\pi_h(\tau^{-1}([0, 1/2])$, which frames $\overline{Conf}^{h, i}_{\bar{A}(x,y)}$ for $i=0,1$. However, these two framings do not glue to a framing of $\overline{Conf}^{h}_{\bar{A}(x,y)} = \paramspace(x,y)$ and we need to twist these framings along the collars $\phi_i$ (\eqref{eq:collar-neighborhood-of-boundary-in-continuation-homotopy-parameter-space}) so that that after twisting these regions glue along the image of $\phi$ \eqref{eq:full-collar-neighborhood-of-boundary-in-continuation-homotopy-parameter-space}. This is where the stabilization $\paramspace(x,y) \rightsquigarrow \paramspace'(x,y)$ comes in handy.  Below let $\xi$ denote the product of the previously chosen framing of $\overline{Conf}^{h, i}_{\bar{A}(x,y)}$ with a possible identity map $\R \to \R$ whenever $\paramspace(x,y) \neq \paramspace'(x,y)$.

Let $t$ be the coordinate on $(1/2-\epsilon, 1/2]$, and let $v = \phi_*(\partial_t)/|\phi_*\partial_t|$, where we define the length of $\phi_*(\partial_t)$ using the previously chosen framings $\xi$. We then define the framing $\xi'$ on $\phi_0((1/2-\epsilon, 1/2] \times \overline{Conf}^{h,\partial}_{\bar{A}(x,y)}$ to agree with $\xi$ on the perpendicular bundle to $v, v'$, and to be 
\[ \xi'(v) = (cos \theta)\xi(v) + (\sin \theta) \xi(v'), \xi'(v) = -(\sin \theta)\xi(v) + (\sin \theta) \xi(v')\]
where $\theta$ is a function of $\tilde{\tau} \circ \pi_r$ that interpolates from $0$ to $\pi$ as one goes from $\tilde{\tau} \circ \pi_r = 1/2-\epsilon$ to $\tilde{\tau} \circ \pi_r = 1/2$ and stays at $\pi$ afterwards.  On the complement of the image of $\phi_1$, one leaves the framing uncahnged, and on $\overline{Conf}^{h, 1}_{\bar{A}(x,y)}$ one sets the framing to be $\xi' = \bar{\xi}$ where $\bar{\xi}$ agrees with $\xi$ on the complement of the span of span of $v'$ and otherwise acts by $\bar{\xi}(v') = -v'$. These define the framing of $\paramspace(x,y)$. See Figure \ref{fig:rotate-framing}. These framings extend to a compatible system of framings of the parameter spaces $\paramspace'(x,y)$ one one replaces the previously defined framings $\xi$ with corresponding framings $\bar{\xi}$ which reverse the framing of $v'$ for all $x,y$ such that $x \in Fix(H_3)$ and $y \in Fix(H_4)$. 
\end{proof}

\begin{lemma}
\label{lemma:framings-continuation-map-homotopy-cyclotomic-compatibility}
    Suppose we are in the situation of Proposition \ref{prop:cyclotomic-compatibility-continuation-homotopy}: we have perturbation data $\mathfrak{P}_k$ for $\CC(H^{\#k}_h, J^{\#k}_h)$ and perturbation data $\mathfrak{P}$ for $\CC(H_h, J_h)$ such there is a $C_{k\ell/C_k}\simeq C_\ell$-equivariant isomorphism virtually smooth flow categories between a restriction of $\CC(H^{\#k}_h, J^{\#k}_h)^{C_k}$ and a restriction of $\CC'(H_h, J_h)_{r, \mathcal{E}_k}$ where $\mathcal{E}_k$ is some $E^s$-parameterization and $\CC'(H_h, J_h)_r$ is a restratification of $\CC'(H_h, J_h)$ associated to the extension of perturbation data as in Lemma \ref{lemma:restratifications-exist}. 

    Suppose further that $\CC'(H_h, J_h)_{\mathcal{E}^f_h}$ has been given a framing as in Lemma \ref{lemma:framings-continuation-map-homotopy}. Write $\mathcal{E}^f_{h, 0, ab}$ for the direct sum of all the terms in \eqref{eq:decomposition-e-parameterization-continuation-map-homotopy} not contributing to $\mathfrak{E}^f_{ab}$. Then, fixing an extension of $\mathcal{E}^f_h$ to $\CC'(H^{\#k}_h, J^{\#k}_h)$, there is an $E^s$-parameterization $\mathcal{E}^f_{h, k, 0}$ of $\CC'(H^{\#k}, J^{\#k})$ such that, writing $\mathcal{E}^f_{h, k} = \mathcal{E}^f_{h, k, 0} \oplus \mathcal{E}^f_h$, there is a framing of $\CC'(H^{\#k}_h, J^{\#k}_h)_{\mathcal{E}^f_{h, k}}$ such that the induced framing of 
    $\CC'(H^{\#k}_h, J^{\#k}_h)_{\mathcal{E}^f_{h, k}}^{C_k}$ agrees under the previous isomorphism with the framing of $\CC'(H_h, J_h)_{r, \mathcal{E}_k \oplus \mathcal{E}^f_{h, k, 0} \oplus \mathcal{E}^f_h}$ induced by applying Lemma \ref{lemma:framings-under-restratification} to the framing of $\CC'(H_h, J_h)_{\mathcal{E}^f_h}$ to get a framing of $\CC'(H_h, J_h)_{r,\mathcal{E}^f_h}$ and subsequently stabilizing. Moreover, when restricting to the objects of $\CC(H_{ab}, J_{ab})$, the framing of $\CC'(H^{\#k}_h, J^{\#k}_h)_{\mathcal{E}^f_{h, k}}$ is simply a the stabilization of the framing of $\CC'(H^{\#k}_{ab}, J^{\#k}_{ab})_{\mathcal{E}^f_{h, k}}$ produced by applying Lemma \ref{lemma:framings-cyclotomic-compatibility} with $\mathcal{E}^f_{k, 0}$ in that lemma given by $\mathcal{E}^f_{h,0, ab} \oplus \mathcal{E}^f_{h, k, 0}$. 
\end{lemma}
\begin{proof}
    This follows completely analogously to Lemma \ref{lemma:framings-cyclotomic-compatibility}: we first choose a sufficiently large initial $\mathcal{E}^f_{h, k, 0}$ such that all the $E^s$-parameterizations of $\mathcal{V}^1_{h, k, ab}$ underlying the parameterizations of the framings of the obstruction bundles of $\CC(H^{\#k}, J^{\#k}_{ab})_{r, \mathcal{E}^f_{h, k}\oplus \mathcal{E}^f_{h, 0, ab} \oplus \mathcal{E}^f_{h, k, 0}}$ are sub-$E^s$-parameterizations of this; we then do the induction of Lemma \ref{eq:numerical-condition-for-existence-of-e-param-2} arranging that the vector space assigned to an object always contains the vector space assigned to that object by $\mathcal{E}^f_{h, k, 0}$, and otherwise requiring that the data $\lambda_{V(x,y)}$ and $\lambda_{V^1(x,y)}$ satisfy the conditions of the lemma whenever the lemma statement in fact fixes these data.
\end{proof}

\section{Background on Homotopy Theory}
\label{sec:homotopy-theory-background}

\subsection{Conventions}

For the entire paper we fix a Grothendieck Universe $\overline{\mathcal{U}}$. All model categories will consider will have all objects elements of $\overline{\mathcal{U}}$.

All groups denoted $G$ will be compact Lie groups; in most cases they are finite groups. In Section \ref{sec:nikolaus-scholze-comparison} we will consider the group $C_{p^\infty} = (\varprojlim C_{p^n}) \subset S^1$ in a general sense, but we will never explicitly consider the action of this group on any spaces. 

We will write $Top$ for the category of compactly generated weak-Hausdorff (CGWH) spaces, $G-Top$ for the category of (CGWH) spaces with continuous $G$-action, and $Top_*, G-Top_*$ for their based variants. All of these categories are enriched over themselves, are symmetric monoidal, and admit all small limits and colimits.

\subsection{Model Categories}
In this section, we briefly review what we need to know about model categories.  A reference for most of what we need about model categories is \cite{hovey2007model}; note that we can always replace conditions on existence of all colimits with existence of all $\overline{\mathcal{U}}$-small colimits. 

A model category $\CC$ comes with three classes of distinguished morphisms: the \emph{cofibrations}, \emph{fibrations}, and \emph{weak equivalences}. These are each closed under composition. Model categories always admit initial objects $0$ and final objects $1$; objects admitting a cofibration from $0$ are called \emph{cofibrant}, and objects admitting a fibration to $1$ are called \emph{fibrant}. A \emph{trivial (co)fibration} is a morphism which is a (co)fibration and also a weak equivalence.

There is a functor from a model category $\CC$ to its homotopy category $Ho(\CC)$; this functor exhibits $Ho(\CC)$ as a localization of $\CC$ at the weak equivalences. In particular, there is a functor $\CC \to Ho(\CC)$, and given a weak equivalence $f: X \to Y$ in $\CC$, its image in $Ho(\CC)$ has a canonical inverse. 

Model categories $\CC$ admit \emph{weak factorization systems}: given any morphism $f: X \to Y$, there is a factorization $X \to Z_1 \to Y$ where $X \to Z_1$ is a cofibration and $Z_1 \to Y$ is a trivial fibration, and also a factorization $X \to Z_2 \to Y$ where $X \to Z_2$ is a trivial cofibration and $Z_2 \to Y$ is a fibration. Moreover, these weak factorization systems are \emph{functorial} in the sense that if we think of $X \xrightarrow{f} Y$ as a diagram, then a map of diagrams extends to a map of factorization diagrams $X \to Y \to Z$. 

In particular, applying the weak factorization to the canonical maps $0 \to X$ gives diagrams $0 \to cX \to X$, functorial in $X$, where the first map is a cofibration and the second is a trivial fibration. Thus there is a cofibrant replacement functor $c: \CC \to \CC$ sending objects to cofibrant objects, and a natural transformation $c \to id$ with all defining maps being trivial fibrations. Dually there is fibrant replacement functor $f: \CC \to \CC$ and a natural transformation $id \to f$ with all defining maps trivial cofibrations. 

We will repatedly use the following 
\begin{lemma}
    Cofibrant replacement preserves fibrancy of objects, and vice versa. 
\end{lemma}
\begin{proof}
    The object $X$ is fibrant if $X \to 1$ is a fibration; but then $cX \to X$ is a trivial fibration, so $cX \to 1$ is a fibration because fibrations are closed under composition. The proof for the other statement is similar. 
\end{proof}

There is a notion of \emph{left homotopy} and \emph{right homotopy} on morphisms from $X$ to $Y$ in $\CC$. It is a fundamental theorem (e.g. \cite[Theorem 1.2.10]{hovey2007model}) that if $X$ is cofibrant and $Y$ is fibrant then these notions coincide, homotopic maps are sent to the same morphism in the homotopy category, and the map 
$\CC(X, Y)/\sim \to Ho(X, Y)$ is an isomorphism, where the equivalence relation on the domain is homotopy of maps. 

\subsection{Orthogonal Spectra}
\label{sec:orthogonal-spectra-detailed-review}
In this section, we give a more detailed review of the properties of orthogonal spectra that we need. Throughout, $G$ will be a finite group. Helpful references include \cite{schwede2022orthogonal, schwede2019lectures} and the foundational \cite{mandell2002equivariant}, as well as \cite{hill2016nonexistence}. We will use the model structure on orthogonal spectra defined in \cite{mandell2002equivariant}, although we will largely avoid any discussion of universes, as explained below. 

\paragraph{Universes.} A $G$-universe is a direct sum of orthogonal $G$-representations which, if it contains a copy of an irreducible representation $V$ of $G$, contains infinitely many copies of $V$. A complete $G$-universe is one that contains the regular representation. 

Throughout this paper, we will fix a single complete $S^1$-universe $\mathcal{U} \supset \R^{\infty}$ (where the inclusion is equivariant and $G$ acts trivially on $\R^\infty$). Whenever we consider $C_{p^k}$-equivariant stable homotopy theory, we will use the complete $C_{p^k}$-universe obtained from $\mathcal{U}$ by restricting along the inclusion map $C_{p^k} \subset S^1$. 

For every closed subgroup $G \subset S^1$ we fix a copy of an embedding of the regular representation $\rho_G \subset \mathcal{U}$.

\paragraph{Orthogonal spectra.} Given pair of finite-dimensional dimensional sub-$G$-representations $V, W \subset \mathcal{U}$, we have an associated $G$-space $O_G(V, W)$, which is the Thom space of the bundle $\xi(V, W) \to L(V, W)$, where $L(V,W)$ is the space of linear isometric inclusions from $V$ to $W$ and the fiber of $\xi(V, W)$ of over $\phi$ is $(Im \phi)^\perp$. These assemble into into a topological category $O_G$; one definition of an orthogonal $G$-spectrum \cite[Definition 4.1]{mandell2002equivariant} is an  functor $X: O_G \to G-Top_*$ enriched over $G-Top_*$, the structure maps
\[ O_G(V, W) \wedge X(V) \to X(W) \]
are $G$-equivariant. A map of $G$-spectra (with respect to the universe $\mathcal{U}$) is natural transformation of enriched functors.

Let $O$ be the category of $O_G$ consisting only of the representations $\R^n = \R^n \oplus 0^\infty \subset \R^\infty \subset \mathcal{U}$. The functor $X$ can functorially and canonically recovered from the functor $X|_O: O \to G-Top_*$ up to a canonical isomorphism, via the map
\[ L(\R^n, W)_* \wedge_{O^{\R^n}} X(\R^n) \simeq X(W);\]
which arises from a factorization of the structure map $O_G(\R^n, W) \wedge X(\R^n) \to X(W)$ using the fact that $O_G(\R^n, W) = L(\R^n, W)$ \cite[Remark 2.7]{schwede2019lectures}.

We will \emph{always} use the latter definition of an orthogonal $G$-spectrum as an enriched functor $O \to G-Top$. The category of orthogonal $G$-spectra in this sense will be denoted $G-OSp$. This category is a symmetric monoidal category.

However, when computing various quantities like fixed points we will define the associated functor $O_G \to G-Top_*$ (where we \emph{always use the one $G$-universe $U$}) to do the computation, and then if the final object is an $H$-spectrum, we will restrict the resulting object back to $O$ to produce a functor $O \to H-Top$. This convention will allow us to make clean statements about the interaction of various functors and change of groups. 

In particular, given a closed inclusion of Lie groups $H \subset G$, an $H$-spectrum $X$, and a $G$-spectrum $Y$, it makes sense to talk about an $H$-equivariant map $f: X \to Y$ by implicitly viewing $Y$ as an $H$-spectrum via the forgetful functor $G-Top \to H-Top$. The associated map on spectra, called `group restriction', is denoted 
\[ R^G_H: G-OSp \to H-OSp.\]
Thus, in the situation above, when we wish to be more precise we will write $f$ as a map 
\[ f: X \to R^G_H Y.\]
However, due to a number of helpful properties of the group restriction functor will allow us to drop many instances of $R^G_H$ from our notation, which will help both with typesetting and conceptual understanding.

\paragraph{Equivariant Homotopy Groups.}
For $j \geq 0$ have by definition that 
    \[ \pi^G_j(X) := \colim_n [\SS(\rho_G^{\oplus n} \oplus \mathbb{R}^j), X(\rho_G^{\oplus n})]_G\]
    where $\rho_G$ is the regular representation of $G$ and $[\cdot, \cdot]_G$ denotes equivariant homotopy classes of maps. Here $\SS$ is the \emph{sphere spectrum}, reviewed below in the paragraph on free spectra.

For $j<0$ we have that 
    \[ \pi^G_j(X) := \colim_n [S(\rho_G^{\oplus n}), X(\rho_G^{\oplus n} \oplus \mathbb{R}^j)]_G,\]

\paragraph{Equivalences and model structures.}
A \emph{level equivalence} of orthogonal $G$-spectra is a map $X \to Y$ such that the associated maps $X(\R^n) \to Y(\R^n)$ are weak $G$-homotopy equivalences. 

A \emph{stable equivalence} of orthogonal $G$-spectra is a map of orthogonal $G$-spectra such that the associaed maps on $H$-stable homotopy groups $\pi^H_k$ are isomorphisms for all closed subgroups $H \subset G$ and all $k \in \Z$.

The stable equivalences of orthogonal $G$-spectra form the the weak equivalences of several model structures on $G-OSp$. We will use the one defined in \cite[Theorem 2.4]{mandell2002equivariant}. This model structure has the property that the objects $F_VX$ for all $G$-CW complexes $X$ are always cofibrant, where $F_V$ the free spectrum on the $G$-representation $V \subset \mathcal{U}$, a notion reviewed below. More generally, the cofibrations (called $q$-cofibrations in \cite{mandell2002equivariant}) are generated by $F_Vi$, where $i: (G/H \times S^{n-1})_* \to (G/H \times D^n)_*$ are the cell cofibrations. 

Given a functor $F$ on $G-OSp$ with values in a category with weak equivalences, we will say that it is \emph{homotopical on $X$} if we would get a diagram of weak equivalences relating $F(X)$ and $F(cf X)$ from the factorization diagram. Similar language will be used for bifunctors; and if the target is not a category with an obvious model-category structure, we simply mean that the weak equivalences of the codomain of $F$ are exactly the isomorphisms. A \emph{derived} functor for $F$ is a functor on $Ho(F)$ that is naturally isomorphic to $Ho(Fcf)$.

\paragraph{Basics.}
Given $X \in G-OSp$ and $K \in G-Top_*$ one can form $X \wedge K$ by smashing $X$ with $K$ levelwise, i.e. 
\[ (X\wedge K)(V) := X(V) \wedge K, \]
and also form the mapping spectrum \cite[Example 5.2]{schwede2019lectures}
\[ X^K(V) := G-Top_*(K, X(V)).\]
These form an adjoint pair:
\[ G-OSp(X, Y^K) \simeq G-OSp(A \wedge X, Y).\]
Similarly, given $Z \in G-OSp$, there is a mapping space
\[ map(X,Z) \in G-Top_*\]
given by the set of natural transformations topologized as a subspace of $\prod_{n} G-Top_*(X(\R^n), Y^(\R^n))$. 

We write $\Omega^V X = X^{S^V}$. 

For the model structure of  \cite[Theorem 2.4]{mandell2002equivariant}, given a cofibrant $X \in G-OSp$ and a fibrant $Y \in G-OSp$, a pair of maps $f_0, f_1: X \to Y$ are homotopic in the sense of model categories if and only if there is a map 
\[ h: f: X \wedge [0,1]_* \to Y\]
such that its restrictions to $X \wedge \{i\}_*$ agree with $f_i$ for $i=0,1$.

\paragraph{Free Spectra.}
For $V \in O_G$, one defines \emph{free spectrum} 
\[ F_V(W) = O_G(V, W), \text{ i.e. } F_V = O_G(V, \cdot). \]
(Note that the kinds of free spectra allowed \emph{depend} on the universe being considered; but for each $G \subset S^1$ we have fixed the universe that we are using, so this notation is unambiguous.)

There are canonical isomorphisms of orthogonal spectra \cite[Eq. 5.8]{schwede2019lectures}
\[ F_V \wedge F_W  \simeq F_{V \oplus W}. \]

More generally, for $K \in G-Top_*$, one sets
\[ F_V K := F_V \wedge K, \]
and we have isomorphisms of orthogonal spectra
\[ F_V K_1 \wedge F_V K_2 \simeq F_{V \oplus W} K_1 \wedge K_2.\]

The functor 
\[ G-Top_* \to G-OSp, K \mapsto F_V K\]
is left adjoint to the evaluation functor at $V$; thus there is a natural map 
\[ F_V A \to \Omega^V(\Sigma^\infty A)\]
adjoint to the adjunction unit (in $G-Top_*$) $A \to \Omega^V( A \wedge S^V) = \Omega^V(\Sigma^\infty A)(V)$. This map turns out to be a stable equivalence \cite[Proposition 5.12]{schwede2019lectures}.

The sphere spectrum is $\SS := F_0$. This makes sense as an object in $G-OSp$ for any $G$, and for any closed Lie subgroup $H \subset G$ we have $R^G_H\SS = \SS$. We write 
\[ \Sigma^\infty X = F_0 X \in G-OSp\]
for a $G$-space $X$.

There are also canonical maps of orthogonal $G$-spectra 
\begin{equation}
    \label{eq:desuspension-maps}
    F_{W \oplus V}\Sigma^V X \to F_WX. 
\end{equation}
which are weak equivalences when $X$ and $Y$ are $G$-CW complexes.

\paragraph{Mapping Spectra.}

Given $V \in O_G$ and $X \in G-OSp$, there is a \emph{shift spectrum} 
\[ sh^V X (W) = X(V \oplus W)\]
with structure maps induced from the inclusions $O_G(W) \subset O_G(V \oplus W)$ which are natural in $W$.

Given $X, Y \in G-OSp$, there is an internal mapping spectrum $F(X,Y)$  \cite[5.11]{schwede2019lectures} with 
\[ F(X,Y)(V) = map(X,sh^V Y). \]

The internal function spectrum makes $G-OSp$ into a closed symmetric monoidal category with unit $\SS$ \cite[Theorem 3.1]{mandell2002equivariant}. As usual, the composition maps 
\[ F(X,Y) \wedge X \to Y \]
are the counit of the adjunction 
\[ F(Y \wedge X, Z) \simeq F(Y, F(X, Z))\]
(valid for each $X$).


The mapping spectrum $F(X,Y)$ is homotopical whenever $X$ is cofibrant and $Y$ is fibrant. This follows from \cite[Proposition 5.17]{schwede2019lectures} and the characterization of the fibrant objects as $G-\Omega$ spectra \cite[Theorem 2.4]{mandell2002equivariant}, namely those orthogonal $G$-spectra $X$ such that for all closed subgroups $H \subset G$ and all pairs of $H$-representations in the universe $\mathcal{U}$, the map $X(V) \to \Omega^W X(V oplus W)$ adjoint to the structure map $X(V) \wedge S^W \to X(V \oplus W)$ is a weak $H$-equivalence. 

\paragraph{Geometric Fixed Points.}

Given a closed \emph{normal} subgroup $H \subset G$, there is a monoidal geometric fixed point functor 
\[ \Phi^H: G-OSp \to H-OSp.\]
There are several differing definitions which all give isomorphic functors on homotopy categories. We always use variant of the definition \cite{mandell2002equivariant} which gives a monoidal functor; 
this definition agrees with the definition in \cite{hill2016nonexistence} under a nontrivial isomorphism, and differs from the definition in \cite{schwede2019lectures}. 

We have the helpful propositions:
\begin{proposition}[Proposition 4.5 \cite{mandell2002equivariant} ]
For $X \in G-Top_*$, there is a natural isomorphism in $H-OSp$
\begin{equation}
    \label{eq:geometric-fixed-points-and-free-spectra}
    \Phi^H(F_V X) \simeq F_{V^H}X^H. 
\end{equation} 
which preserves ($q$-)cofibrations and acyclic ($q$-)cofibrations.
\end{proposition}

The argument of \cite[Proposition 4.5]{mandell2002equivariant} shows also that 
\begin{lemma}
    The functor $\Phi^H$ commutes with the formation of mapping telescopes (see below.)
\end{lemma}

\begin{proposition}[Proposition 4.5\cite{mandell2002equivariant}]
For $X, Y \in G-OSp$, there is a natural map in $G/H-OSp$
\[ \Phi^H X \wedge \Phi^H Y \to \Phi^H(X \wedge Y)\]
which is an isomorphism if one of $X$ and $Y$ are cofibrant. 

\end{proposition}

The functor $\Phi^H X$ is homotopical whenever $X$ is cofibrant  \cite[Prop. 4.12, Corr. 4.13]{mandell2002equivariant}.

\paragraph{Interactions with change of group. }

We now compile several results about change of group. 

\begin{lemma}[Theorem IV.1.2 \cite{mandell2002equivariant}]
\label{lemma:group-retriction}

Let $H \subset G$ be a closed subgroup of a compact Lie group $G$, and let  $R^G_H$ be the restriction functor from orthogonal $G$-spectra to orthogonal $H$-spetra. Then $R^G_H$ preserves stable equivalences, $q$-fibrations, and $q$-cofibrations. In particular, they preserve fibrant and cofibrant objects.
\end{lemma}

It follows immediately from the definitions above there is an equality of bifunctors
\begin{equation} 
    \label{eq:restriction-and-function-spectra}
F(R^G_H X, R^G_H Y) = R^G_H F(X, Y).
\end{equation}

It follows immediately from \cite[Definition 4.1]{mandell2002equivariant} that we have an equality of functors
\begin{equation}
\label{eq:restriction-and-geometric-fixed points}
    \Phi^H R^{G}_{G'} = R^G_{G'}\Phi^H.
\end{equation}

\paragraph{Mapping Telescopes.}

We will collate some results about mapping telescopes below. First, a mapping telescope of a family of spaces is 
\[ X_0 \xrightarrow{f_0} X_1 \xrightarrow{f_1} \cdots \xrightarrow{f_{n-1}} X_n\]
is the pushout
\[ Tel(X_0 \to \cdots \to X_n) := (X_0(V) \times [0,1]) \cup_{f'_1} (X_1(V) \times [0,1]) \cup_{f'_2} \cdots \cup_{f'_{n-1}} (X_n \times \{0\}) \]
where $f'_k: X_k \times \{1\} \to X_{k+1} \times \{0\}$ is is just the action of $f_k$ on the first factor. The mapping telescope of an infinite sequence 
\begin{equation}
    \label{eq:sequence}
    X_0 \to X_1 \to \cdots
\end{equation} 
is the colimit of the mapping telescope of the finite subsequences, i.e. 
\[ Tel(X_0 \to X_1 \to \cdots) = \colim_n (X_0 \to X_1 \to \cdots \to X_n)\]
under the natural closed embeddings of these spaces. When all the $X_i$ are $G$-spaces and $f_i$ are $G$-maps then mapping telescopes inherit a $G$-action, and if all spaces involved are $G$-CW complexes and $G$-CW maps then the mapping telescope is a $G$-CW complex. 

There are collapse maps 
\[ coll_{0, 1}: Tel(X_0 \to \cdots \to X_n) \to Tel(X_1 \to \cdots \to X_n)\]
induced by the collapse map $X_0 \times [0,1] \to X_1$, and they are (equivariant) homotopy equivalences, with explicit homotopy-inverse the closed inclusion 
\begin{equation}
\label{eq:inclusion-of-mapping-telescopes}
    i_{1,0}: Tel(X_1 \to \cdots \to X_n) \subset Tel(X_0 \to \cdots \to X_n)
\end{equation} 
Taking colimits in $n$ shows that the analogous collapse and inclusion maps (denoted with in the same way) on infinite mapping telescopes are $(G-)$-homotopy equivalences as well. One defines $coll_{a, b} = coll_{a, a+1}coll_{a+1, a+2} \cdots coll_{b-1, b}$ and similarly for $i_{b,a}$. 

Given a diagram like \eqref{eq:sequence} of orthogonal $G$-spectra, the mapping telescope is defined by applying the above constructions level-wise. One thus immediately sees that if all $X_i$ are cofibrant $G$-spectra then so is the mapping telescope. Taking mapping telescopes manifestly commutes with group restriction, and the collapse maps and their canonical homotopy inverses when applied level-wise define corresponding maps of mapping telescopes of orthogonal $G$-spectra which are also level equivalences. 

The following lemma is standard:
\begin{lemma}
\label{lemma:colimit-agrees-with-hocolim}
    Consider a diagram of $X_i \in G-Top_*$
    \[ X_0 \subset X_1 \subset X_2 \subset \ldots \subset X = \colim_i X_i.\]
    Then if all inclusions above are closed inclusions then the canonical map from the mapping telescope $Tel(X_0 \to \ldots) \to X$ is a weak equivalence.
\end{lemma}

We have the following useful standard result:
\begin{lemma}
    Given a diagram of orthogonal $G$-spectra equipped with levelwise nondegenerate basepoints
    \[ Y_0 \to Y_1 \to \cdots\]
    (with either finitely or infinitely many terms) where all maps are level-wise closed embeddings, we have a functorial isomorphism
    \begin{equation}
    \label{eq:colimit-and-stable-homotopy}
        \colim \pi_G^j(Y_i) \simeq \pi_G^j(\colim Y_i)
    \end{equation}
    of equivariant homotopy groups. 
\end{lemma}
\begin{proof}
    This is proven exactly like its analog for $G=1$, which is explained in \cite[Proposition 1.24]{schwede2022orthogonal}. For $j \geq 0$, by \cite[Prop 2.4.2]{hovey2007model}, any map from $S(\rho_G^{\oplus n} \oplus \mathbb{R}^j)$ to $\colim_i Y_i(\rho_G^{\oplus n})]_G$ factors through some $Y_i(\rho_G^{\oplus n})]_G$, as does any homotopy between a pair of such maps. Thus the canonical map 
    \[\colim_i [S(\rho_G^{\oplus n} \oplus \mathbb{R}^j), Y_i(\rho_G^{\oplus n})]_G \to [S(\rho_G^{\oplus n} \oplus \mathbb{R}^j), \colim_i 
 Y_i(\rho_G^{\oplus n})]_G\]
    is an isomorphism. 
    Passing to colimits in $n$ we conclude that the canonical map \eqref{eq:colimit-and-stable-homotopy} is an isomorphism. 

    For $j < 0$, 
    and a directly analogous argument goes through.
    
\end{proof}

This result easily implies the following helpful lemmata on mapping telescopes of equivariant orthogonal spectra:
\begin{lemma}
    Given a (strictly) commutative diagram 
    \begin{equation}
    \label{eq:commutative-diagram-of-spectra-a}
        \begin{tikzcd}
            X_0 \ar[r] \ar[d] & X_1 \ar[r] \ar[d] & X_2 \ar[d] \ar[r] & \cdots \\
            X'_0 \ar[r] & X'_1 \ar[r] & X'_2 \ar[r] & \cdots 
        \end{tikzcd}
    \end{equation}
    of orthogonal $G$-spectra with all the vertical arrows stable equivalences, the induced map on mapping telescopes 
    \[ Tel(X_0 \to X_1 \to \cdots) \to Tel(X'_0) \to Tel(X'_0 \to X'_1 \to \cdots) \]
    is a stable equivalence.
\end{lemma}
\begin{proof}
    A stable equivalence of orthogonal $G$-spectra is one for which the induced maps on $\pi^H_k$ are isomorphisms for all subgroups $H \subset G$ and all $k \in \Z$. If the diagram \eqref{eq:commutative-diagram-of-spectra-a} has finitely many entries and terminates with a morphism $X_n \to X'_n$ for some $n$, then there is a commutative diagram 
    \begin{equation}
        \begin{tikzcd}
            Tel(X_0 \to X_1 \to \cdots)\ar[d] & X_n \ar[d] \ar[l] \\
            Tel(X'_0 \to X'_1 \to \cdots)& X'_n \ar[l]
        \end{tikzcd}
    \end{equation}
    where the horizontal arrows are the inclusion maps (and are thus level equivalences), the left vertical arrow is the induced map on mapping telescopes, and the right vertical arrow is the map in the diagram \eqref{eq:commutative-diagram-of-spectra-a}. In the case where \eqref{eq:commutative-diagram-of-spectra-a} has infinitely many entries, the statement for finite-size diagrams implies that the map 
    \[ \colim_n \pi^H_j(Tel(X_0 \to \cdots \to X_n)) \to \colim_n \pi^H_j(Tel(X'_0 \to \cdots \to X'_n))\]
    is an isomorphism for all $H, j$; the previous lemma then implies that the map between mapping telescopes is a stable equivalence. 
\end{proof}

\begin{lemma}
\label{lemma:weak-map-of-mapping-telescopes}
    Given a diagram 
     \begin{equation}
     \label{eq:diagram-homotopy-commutative}
        \begin{tikzcd}
            X_0 \ar[r, "a_0"] \ar[d, "f_0"] & X_1 \ar[r, "a_1"] \ar[d, "f_1"] & X_2 \ar[d, "f_2"] \ar[r, "a_2"] & \cdots \\
            X'_0 \ar[r, "b_0"] & X'_1 \ar[r, "b_1"] & X'_2 \ar[r, "b_2"] & \cdots 
        \end{tikzcd}
    \end{equation}
    of orthogonal $G$-spectra which is \emph{not} commutative, but such that one has chosen homotopies $h_i: X_i\wedge [0,1]_* \to X'_{i+1}$ from $f_{i+1}a_i$ to $b_if_i$, there is a map 
    \begin{equation}
        \label{eq:induced-map-on-telescopes-from-homotopy-commutative-diagram}
        Tel(X_0 \to X_1 \to \cdots) \to Tel(X'_0) \to Tel(X'_0 \to X'_1 \to \cdots) 
    \end{equation}
    defined by as the induced map on colimits of finite mapping telescopes by maps which, when evaluated on the representation $V$, map $X_i(V) \times\{0\}$ to $X'_i(V) \times \{0\}$ via $f_i(V)$, and more generally map $X_i(V) \times [0,1]$ to 
    $X'_i(V) \times [0,1] \cup_{b_i(V)} X'_{i+1} \times \{0\}$ via by sending $(x, t) \mapsto (h_i(x, 2t-1), 0)$ for $t \geq 1/2$, and otherwise $(x, t) \mapsto (f_i(x), t/2)$ for $t \leq 1/2$. 

    If all the vertical arrows of \eqref{eq:diagram-homotopy-commutative} are stable equivalences then the map \eqref{eq:induced-map-on-telescopes-from-homotopy-commutative-diagram} is a stable equivalence.
    
\end{lemma}
\begin{proof}
    The proof of this lemma is \emph{identical} to the proof of the previous lemma. 
\end{proof}

\section{Floer Homotopy and its Properties}
\label{sec:floer-homotopy}
In this section, we apply the tools that we have developed in this paper to construct the genuine equivariant Floer homotopy type,  as well as the cyclotomic structure on symplectic cohomology. First, in Section \ref{sec:geometric-results-review}, we will combine the results we have established so far and use these to define the (genuine equivariant) Floer homotopy and cohomotopy types $CF_\bullet(H, J; \SS)$ and $CF^\bullet(H, J; \SS)$ for admissible Floer data $(H, J)$. A-priori these depend on many choices, which we will refer to as the \emph{homotopy definition data}; in this same section, we describe a stabilization of homotopy definition data, and explain how the geometric results we have established allow us to establish maps and relations between maps between these Floer homotopy types. Then, in Section \ref{sec:comparison-with-floer-homology}, we describe some general results about flow categories (which in fact are used in Section \ref{sec:geometric-results-review}), and use these and other arguments regarding the framings of our virtually smooth flow categories $\CC'(H, J)$ to show that on homology, the spectra and maps of spectra defined in Section  \ref{sec:comparison-with-floer-homology} define, when the homotopy definition data are chosen appropriately (and in particular when the Floer data are sufficiently regular), spectra and maps of spectra which correspond to the Floer chain complex and the Floer continuation maps on Floer homology. Subsequently, in Section \ref{sec:invariance-of-equivariant-floer-homotopy}, we use these aforementioned results to prove invariance of the full genuine equivariant Floer homotopy types $CF_\bullet(H, J; \SS)$ from all choices made. 

We then proceed with the discussion of the cyclotomic structure on symplectic cohomology. In order to construct symplectic cohomology, we imitate the usual construction of symplectic cohomology as a colimit of Hamiltonian Floer homology groups. 

Let us recall the construction of symplectic cohomology. Suppose that $(H, J)$ is a \emph{strongly nondegenerate} convex linear Floer datum on $(M, \omega, \lambda)$.  Thus, the Floer data $(H^{\# k}, J^{\# k})$ are also convex linear nondegenerate Floer data, with the slope of $H^{\# k}$ being $k$ times the slope of $H$. Fix a prime $p$. If we choose a convex continuation datum $(H_s, J_s)$ from $H_1$ to $H^{\# p}$, then $(H^{\#p^k}_s, J^{\# p^k}_s)$ is a convex continuation datum from $(H^{\#p^k}, J^{\# p^k})$ to $(H^{\# p^{k+1}}, J^{\# p^{k+1}})$. If all the Floer data and continuation data are also regular, then we have a diagram of Floer cohomology groups 
\[ CF^*(H, J) \to CF^*(H^{\# p}, J^{\# p}) \to \cdots. \]
The symplectic cohomology $SH^*(M)$ is the colimit of the corresponding cohomology groups 
\begin{equation}
    \label{eq:sh-colimit-on-homology}
    HF^*(H, J) \to HF^*(H^{\# p}, J^{\#p}) \to \cdots. 
\end{equation}
It turns out to be independent of all the data chosen above up to canonical isomorphism \cite{abouzaid-sh}. 

To construct spectral symplectic cohomology $SH^\bullet(M; \SS)$ one imitates this construction. The geometric constructions of the paper (reviewed in Section \ref{sec:geometric-results-review} below) allow us to define Floer cohomotopy types
\[CF^\bullet(H^{\# k}, J^{\# k}; \SS)_{\mathfrak{H}^k} := F(|\CC(H^{\# k}, J^{\#k})|_{\mathfrak{H}^k}, \SS). \]
where $|\CC(H^{\#k}, J^{\#k})|_{\mathfrak{H}^k}$ is the orthogonal spectrum associated to a virtual smoothing of the proper flow category $\CC(H^{\#k}, J^{\#k})$, together with a framing, an embedding of the framing, and so forth, as needed in Sections \ref{sec:flow-categories} - \ref{sec:index-theory}; all the choices made in the construction are denoted by $\mathfrak{H}^k$, and are referred to as the \emph{homotopy definition data}.  Here, the function spectrum $F( \cdot, \SS)$ is taken in the derived sense. The constructions of this paper allow us to define a diagram of orthogonal spectra
\begin{equation}
    \label{eq:spectral-lift-of-floer-homology-colimit}
CF^\bullet(H, J; \SS)_{\mathfrak{H}^1} \to CF^\bullet(H^{\# p}, J^{\# p}; \SS)_{\mathfrak{H}^p} \to CF^\bullet(H^{\# p^2}, J^{\# p^2}; \SS)_{\mathfrak{H}^{p^2}} \to \cdots 
\end{equation}
These maps lift the maps \eqref{eq:sh-colimit-on-homology} in the sense that upon taking (derived) smash product with $H\Z$ and subsequently taking stable homotopy groups one recovers the maps \eqref{eq:sh-colimit-on-homology}. Thus, the spectral symplectic cohomology reduces to the usual symplectic cohomology upon smashing with $H\Z$ and taking homotopy groups as well. 

Constructing the cyclotomic structure is somewhat more challenging due to the potentially large numbers of compatibility conditions involved. As stated in the introduction, we use the definition of a \emph{genuine $p$-cyclotomic spectrum}, originally defined by Hesselholt-Madsen \cite{hesselholt2004rham}; we use closely related variant of this definition introduced by Nikolaus-Scholze \cite{nikolaus-scholze} which still uses genuine equivariant spectra, but uses $\infty$-categories to keep track of some of the coherence data required. We wish to show that $SH^\bullet(M; \SS)$ is the underlying spectrum of a genuine $p$-cyclotomic spectrum for every $p$. (We do not work out the compatibility conditions between different primes $p$ in this paper, although it is essentially a tedious exercise given the tools developed in this paper.) The notion of a genuine $p$-cyclotomic spectrum, reviewed in Section \ref{sec:nikolaus-scholze-comparison}, is purely $\infty$-categorical, while our constructions are concrete and geoometric. Now, the symmetries of the geometric situation translate furnish part of the data needed for the $\infty$-categorical notion. Indeed, a a genuine $p$-cyclotomic spectrum $X$ is, for each $\ell$, a $C_{p^\ell}$-spectrum $X_\ell$, such that $X_{\ell-1}$ is weakly equivalent to $F^{C_{p^\ell}}_{C_{p^{\ell-1}}}X_\ell$, together with cyclotomic structure equivalences 
\begin{equation}
    \label{eq:cyclotomic-structure-equivalence-formal}
    X_{\ell-1} \to \Phi^{C_p} X_\ell
\end{equation}

such that all this data is compatible up to homotopy as $k$ varies. The geometry of the virtual smoothings of the flow categories of Floer moduli spaces that we have established in the earlier portions of this paper gives us each piece of this structure. To define $X_\ell$, we consider the diagram \eqref{eq:spectral-lift-of-floer-homology-colimit} over all $k \geq \ell$, and note that the relevant flow categories are naturally $C_{p^\ell}$-equivariant, and all the smoothings, framings, etc. can be chosen to be so as well. As such, the construction of this paper give us a lift of that part of the diagram \eqref{eq:spectral-lift-of-floer-homology-colimit} thought of as a diagram in the homotopy category of spectra to a diagram of $C_{p^\ell}$-spectra; taking the homotopy colimit in this latter category defines $X_\ell$. The equivalences \eqref{eq:cyclotomic-structure-equivalence-formal} are then defined by taking the map on homotopy colimits induced by the weak equivalences 
\[ CF^\bullet(H^{\# p^{k-1}}, J^{\# p^{k-1}}; \SS)_{\mathfrak{H}^{p^{k-1}}} \simeq \Phi^{C_p} CF^\bullet(H^{\# p^k}, J^{\# p^k}; \SS)_{\mathfrak{H}^{p^k}}\]
which arise essentially for free from Theorem \ref{thm:geometric-fixed-points-of-flow-category} once one has done all constructions equivariantly, as detailed earlier in this paper. 

Thus the geometric constructions give the input to the data of a genuine $p$-cyclotomic spectrum. The challenge is to verify that the compatibility conditions betwen these various maps required by the definition of a genuine $p$-cyclotomic structure hold. All geometric choices to be made are made in the straightforward Section  \ref{sec:cyclotomic-compatibility-geometry} (which combines many geometric results in the paper). This data is then carefully repackaged in Sections \ref{sec:cyclotomic-structure-cofibrancy} and Section \ref{sec:nikolaus-scholze-comparison}. The first basic issue is that one is working with many different categories of genuine equivariant spectra (specifically, orthogonal $C_{p^k}$ spectra for varying $k$) \emph{for which the notions of cofibrant and fibrant spectra are all distinct}. Using summary of the geometric results of Section \ref{sec:geometric-results-review}, we explain, in Section \ref{sec:cyclotomic-structure-cofibrancy}, how to take a sequence of cofibrant-fibrant replacements of the resulting objects and maps in the varying model categories of orthogonal spectra, as well as explicit models of homotopy colimits of the resulting diagrams, in order to define certain relatively concrete representatives of the objects $X_\ell$ in the categories of $C_{p^\ell}$-orthogonal spectra. Subsequently, using the explicit nature of the homotopy colimits and the sequences of fibrant-cofibrant replacements taken, we show in Section \ref{sec:nikolaus-scholze-comparison} that this process also gives the compatibility data needed for these maps to define a lift of $SH^\bullet(M; \SS)$ to the category of genuine $p$-cyclotomic spectra. 

We now proceed with the construction. 

\subsection{Basic Results on Floer Homotopy}
\label{sec:geometric-results-review}

In this section we use the earlier geometric geometric constructions of this paper in order to define Floer homotopy types and define relationships between them. We will work with two kinds of homotopy types: \emph{homological} and \emph{cohomological}. One will be the Spanier-Whitehead dual of the other, and the relationship between them be analogous to the relationship between Floer \emph{homology} and Floer \emph{cohomology}. The dualization from homological Floer homotopy types to cohomological Floer homotopy types will be introduced to reverse the directions of various maps defined by virtually smooth flow categories, which is important when defining symplectic cohomology since in that setting the continuation maps used in the definition of symplectic cohomology only go one way, in the direction of increasing slope of the Hamiltonian.

Throughout, we are working with an admissible symplectic manifold. 

\begin{remark}
Below, whenever we consider free spectra $F_V$ for $V$ a $C_{p^k}$-representation, we find an $S^1$-representation $V'$ which gives $V$ upon restriction, and fix an embedding of $V'$ into $\mathcal{U}$. Often we will have pairs of representations $V_1 \subset V_2$ where $V_1$ is a $C_{p^k}$-representation, $V_2$ is a $C_{p^\ell}$ representation for $\ell \geq k$, and the embedding is $C_{p^\ell}$-equivariant. In that case, if we have already picked a lift of $V_1$ to an $S^1$-represetation $V'_1$ and an embedding $V'_1 \subset \mathcal{U}$, we proceed by choosing a lift of $V_2$ an $S^1$-representation $V'_2$ which contains $V'_1$ as a sub-representation, and extend of $V'_1$ to an embedding of $V'_2$. In the cases we consider it will be straightforward to see that such embeddings exist and that their choices of embeddings do not affect the results up to isomorphism.
\end{remark}

\begin{proposition}
\label{prop:homotopy-type-geometry}
    Let $(H, J)$ be an admissible Floer datum on $M$. Then, given a choice of integralization data $\mathfrak{A}$, there is a proper flow category $\CC(H, J)$. There exists a compatible regular system of perturbation data $\mathfrak{P}$ with underlying integralization data $\mathfrak{A}$ and a smooth structure on the virtual smoothing $\CC'(H, J)$ of $\CC(H, J)$ associated to $\mathfrak{P}$. If the polarization class $\rho \in KO^1(LM)$ of Section \ref{sec:producing-framings} vanishes, then, by choosing a trivialization of the polarization map $\tilde{\rho}$, the virtual smoothing can be  framed, and admits extension $\bar{\mathbf{E}}$ of an embedding $\mathbf{E}$. If $(H, J)$ are $C_k$ equivariant, and the polarization class $\rho \in KO^1_{C_k}(LM)$ of Section \ref{sec:producing-framings} vanishes, then all of the above data can be chosen to be equivariant. The associated stable homotopy type is the orthogonal spectrum 
    \begin{equation}
        |\CC'(H, J)| = F_{\tilde{V}^1}\{\CC'(H, J)\}
    \end{equation}
    where $\tilde{V}^1$ is the \emph{shift space} associated to $\CC'$ (Remark \ref{rk:shift-spaces-for-categories} and Definition \ref{def:shift-spaces}) and $\{\CC'(H, J)\}$ is a $C_k$-space with the equivariant homotopy type of a $C_k$-CW complex. This orthogonal spectrum depends on all of the data described above; we call the data $\mathfrak{H} = (\mathfrak{A}, \mathfrak{P}, \tilde{\rho}, \mathbf{E}, \bar{\mathbf{E}})$ the \emph{homotopy definition data}, and we will refer to $\tilde{V}^1 = \tilde{V}^1_\mathfrak{H}$ and $\{\CC'(H, J)\}_{\mathfrak{H}}$ when we wish to keep track of the dependence on the homotopy definition data.
\end{proposition}

We will define the \emph{Floer homotopy type} associated to $\mathfrak{H}$ to be 
\[ CF_\bullet(H, J)_{\mathfrak{H}} := |\CC'(H, J)|_{\mathfrak{H}} := F_{\tilde{V}^1_{\mathfrak{H}}}\{\CC'(H, J)\}_{\mathfrak{H}}.\] The \emph{Floer cohomotopy type} associated to $\mathfrak{H}$ is the (equivariant) Spanier-Whitehead dual of the above object in the (equivariant) stable homotopy category, and is represented in the equivariant stable homotopy category by the function spectrum 
\[CF^\bullet(H, J)_{\mathfrak{H}} := F(|\CC'(H, J)_{\mathfrak{H}}|, f\SS),\]
where $f$ denotes fibrant replacement.

\begin{proof}[Proof of Proposition \ref{prop:homotopy-type-geometry}]
    This follows from a combination of the results we have established in this paper. Namely, we invoke Propositions \ref{prop:compatible-perturbation-data-exist},
\ref{prop:smooth-flow-category-floer-trajectories}, and \ref{prop:framings-for-floer-moduli-spaces-exist}
to produce the virtually smooth framed flow category $\CC'(H, J)$. To produce the embedding and its extension we invoke Propositions \ref{prop:embeddings-exist} and Lemma \ref{lemma:extensions-of-embeddings-exist}. The definition of $\{\CC'(H, J)\}$ is given in Definitions \ref{def:floer-homotopy-type} and \ref{def:equivariant-stable-homotopy-type}, and the statement that it has the homotopy type of a $G$-CW complex is Proposition \ref{prop:homotopy-type-is-G-CW-complex}. 
\end{proof}

In order to relate these homotopy types, we will use the continuation map construction, for which virtual smoothings are constructed in Section \ref{sec:continuation-maps-geometry}. We first need to introduce the notion of a stabilization of a homotopy definition datum:

\begin{definition}
\label{def:stabilization-of-homotopy-definition-data}
Let $\mathfrak{H}'$ and $\mathfrak{H}$ be homotopy definition data for $\CC(H, J)$. We say that $\mathfrak{H}'$ is a stabilization of $\mathfrak{H}$ if the following conditions hold:

\begin{enumerate}
    \item The integralization datum  $\bar{\mathfrak{A}}'$ underlying $\mathfrak{H}'$ is an extension of the integralization datum $\bar{\mathfrak{A}}$ underlying $\mathfrak{H}$ in the sense of Lemma \ref{lemma:extend-perturbation-data-across-restratifications}, and
    \item Writing $\mathfrak{P}$ and $\mathfrak{P}'$ for the perturbation data underlying $\mathfrak{H}$ and $\mathfrak{H}'$, respectively, we have that there is an isomorphism $f$ of virtually smooth equivariant flow categories  from $\CC'(H, J)_{\mathfrak{P}'}$ to (a restriction of) $\CC'(H, J)_{\mathfrak{P}, r, \mathcal{E}}$, which is the stabilization by an $E^s$-parameterization $\mathcal{E}$ of the restratification $\CC'(H, J)_{\mathfrak{P}, r}$ of $\CC'(H, J)_{\mathfrak{P}}$ associated to the fact that $\bar{\mathfrak{A}}'$ is an extension of $\bar{\mathfrak{A}}$
    \item There is an isomorphism $f^E: \mathcal{E}' \to \mathcal{E}^s$, of $E^s$-parameterizations, where $\mathcal{E}'$ is the $E^s_2$-parameterization underlying  $\CC'(H, J)_{\mathfrak{P}'}$, and $\mathcal{E}^s$ is the $E^s_2$-parameterization underlying the framing of $\CC'(H, J)_{\mathfrak{P}, r, \mathcal{E}}$ induced from that of $\CC(H, J)'_{\mathfrak{P}}$ by Lemmas \ref{lemma:framings-under-restratification} and \ref{lemma:framings-under-stabilization}. Moreover,
    under the identifications $f$ and $f^E$, the framing of $\CC'(H, J)_{\mathfrak{P}'}$ agrees with the framing of $\CC'(H, J)_{\mathfrak{P}, r, \mathcal{E}}$ induced by these same Lemmas. 
    \item Furthermore,the embeddings underlying $\mathfrak{H}'$ agree (under $f$ and the map $f^F$ on $F^s_2$-parameterizations induced by $f^E$) with the embeddings induced from the embeddings underlying $\mathfrak{H}$ via  Lemma \ref{lemma:restratifications-exist} and \ref{lemma:stabilizing-embeddings-by-parameterizations}, and similarly the extension of the embedding underlying $\mathfrak{H}'$ is given (under $f^F$) by applying Lemmas \ref{lemma:restratifications-exist} and \ref{lemma:stabilizing-embeddings-by-parameterizations}.

    The same definition applies to make sense of stabilizations of homotopy definition data for $\CC(H_s, J_s)$ or $\CC(H_h, J_h)$, where $(H_s, J_s)$ are convex continuation data and $(H_h, J_h)$ are convex continuation homotopy data. 
\end{enumerate}
\end{definition}

\begin{proposition}
\label{prop:effect-of-stabilization-on-floer-homotopy-type}
Say $\mathfrak{H}'$ is a stabilization of of the homotopy definition datum $\mathfrak{H}$. Using the notation of Definition \ref{def:stabilization-of-homotopy-definition-data}, write $V^1_{\mathcal{E}}$ and $\mathcal{V}^1_{\mathcal{E}^s}$ for the shift spaces underlying $\mathcal{E}$ and $\mathcal{E}^s$, and write 
\[ \Delta \bar{A}_{max} = |\mathfrak{A}'| \setminus |\mathfrak{A}|\]
where $(\mathfrak{A}, \mathfrak{A}')$ are the marked actions underlying the integralization data $(\bar{\mathfrak{A}},\bar{\mathfrak{A}}')$, respectively. Then there is an equivariant homeomorphism 
\begin{equation}
    \label{eq:effect-of-stabilization-on-floer-homotopy-space}
    \{\CC(H, J)\}_{\mathfrak{H}'} = \Sigma^{V^1_{\mathcal{E}} \oplus \R^{\Delta \bar{A}_{max}}+M} \{\CC(H, J)\}_{\mathfrak{H}}
\end{equation}
and an isomorphism of orthogonal spectra
\begin{equation}
    \label{eq:effect-of-stabilization-on-floer-homotopy-type}
    |\CC(H, J)|_{\mathfrak{H}'} \simeq F_{V^1_{\mathcal{E}}  \oplus \R^{\Delta \bar{A}_{max}}} \Sigma^{V^1_{\mathcal{E}} \oplus \R^{\Delta \bar{A}_{max}}+M} |\CC(H, J)|_{\mathfrak{H}},
\end{equation}
where $M$ is the restratification embedding shift of Lemma \ref{lemma:restratifications-exist}.
Moreover we can replace $(H, J)$ with $(H_s, J_s)$ or with $(H_h, J_h)$ in any of the statements above. 
\end{proposition}
\begin{proof}
The homeomorphism \eqref{eq:effect-of-stabilization-on-floer-homotopy-space} follows immediately from Property (2) of Definition \ref{def:stabilization-of-homotopy-definition-data} together with Proposition \ref{prop:simple-restrat-effect-on-homotopy-type} and Lemma \ref{lemma:effect-of-stabilizing-by-parameterization-on-floer-homotopy-type}. 

The isomorphism \eqref{eq:effect-of-stabilization-on-floer-homotopy-type} follows from \eqref{eq:effect-of-stabilization-on-floer-homotopy-space} together with the statement that, writing $V^1_{\mathcal{E}_0}$ for the shift space associated to $\CC'(H, J)_{\mathfrak{P}}$, we have an equivariant isomorphism 
\begin{equation}
    \label{eq:computation-of-change-of-shift-space}
    V^1_{\mathcal{E}_0} \oplus V^1_{\mathcal{E}} \oplus \R^{\Delta \bar{A}_{max}} \simeq V^1_{\mathcal{E}^s}. 
\end{equation}
But this follows from Property (3) of Definition \ref{def:stabilization-of-homotopy-definition-data} together with the fact that in the context 
of Lemma \ref{lemma:framings-under-restratification} , the shift space associated to the part of the $E^s$-parameterization $V^1_{can}$ of that lemma changes exactly by $\R^{\Delta \bar{A}_{max}+M} \oplus V^1_{\mathcal{E}_1}$ for some semi-free parameterization $\mathcal{E}_1$ (see Remark \ref{rk:can-incorporate-shift-in-canonical-param}, in the context of which $M =  \Delta \mathfrak{A}_{lower}$) while in the context of Lemma \ref{lemma:framings-under-stabilization}, we simply add $\mathcal{V}^1_{\mathcal{E}_2}$ to the shift space for some semi-free parameterization $\mathcal{E}_2$; and in the end we have $\mathcal{E} = \mathcal{E}_1 \oplus \mathcal{E}_2$. Combining these statements gives exactly \eqref{eq:computation-of-change-of-shift-space}. The proofs with $(H, J)$ replaced by $(H_s, J_s)$ or $(H_h, J_h)$ are identical. 
\end{proof}

Below we summarize the geometric data associated to continuation maps:
\begin{proposition}
\label{prop:continuation-map-geometry}
    Let $(H_\pm, J_\pm)$ be a pair of admissible floer data for $M$, and let $(H_s, J_s)$ be admissible continuation data from $(H_-, J_-)$ to $(H_+, J_+)$. If the polarization class $\rho \in KO^1(LM)$ vanishes, then, given homotopy definition data $\mathfrak{H}_\pm$ for $|\CC'(H_\pm, J_\pm)|$, there exists a homotopy definition datum $\mathfrak{H}$ defining a virtual smoothing $\CC'(H_s, J_s)$ of $\CC(H, J)$, such that the restrictions of $\mathfrak{H}$ to $\CC(H_\pm, J_\pm)$ are stabilizations of $\mathfrak{H}_\pm$, and such that the integral actions are given by shifts of the original integral actions for $\CC(H_\pm, J_\pm)$ as in \eqref{eq:shifted-integral-action-around-contiunuation-map}. Moreover, if either of $\mathfrak{H}_\pm$ are $C_k$-equivariant,  and the polarization class $\rho \in KO^1_{C_k}(LM)$ vanishes then the  restriction of $\mathfrak{H}$ to the corresponding flow subcategory $\CC(H_\pm, J_\pm)$ can be chosen to carry a $C_k$ equivariant structure,  with the stabilization from $\mathfrak{H}_\pm$ to $\mathfrak{H}|_{\CC(H_\pm, J_\pm)}$ an equivariant stabilization.In particular, the shift spaces $(\tilde{V}^1_+, \tilde{V}^1_-, \tilde{V}^1)$ associated to $(\mathfrak{H}_-, \mathfrak{H}_+, \mathfrak{H})$ satisfy $\tilde{V}^1 = \tilde{V}^1_\pm \oplus \tilde{V}^1_{\pm, s}$, where this decomposition is a decomposition into $C_k$-representations if $(H_\pm, J_\pm)$ is $C_k$-equivariant. 

    We have an inclusion 
    \begin{equation}
        \label{eq:inclusion-of-lower-category-space}
        \Sigma^{V^1_{+, s}}\{\CC(H_+, J_+)\}_{\mathfrak{H}_+} \subset \{\CC(H_s, J_s)\}_{\mathfrak{H}}
    \end{equation}  
    which is an equivariant homeomorphism onto its image, and the quotient by the inclusion is equivariantly homeomorphic to a shift of the space underlying $|\CC(H_-, J_+)|$: 
    \begin{equation}
        \label{eq:identification-of-upper-category-space}
        \{\CC(H_s, J_s)\}_{\mathfrak{H}}/\left(\Sigma^{V^1_{+, s}}\{\CC(H_+, J_+)\}_{\mathfrak{H}_+}\right) \simeq \Sigma^{\R \oplus V^1_{-, s}} \{\CC(H_-, J_-)\}_{\mathfrak{H}_-} \simeq \{\CC(H_-, J_-)\}_{\mathfrak{H}}. 
    \end{equation}
\end{proposition}

\begin{proof}
    Write $\mathfrak{P}_\pm$ for the perturbation data underlying $\mathfrak{H}_\pm$. We choose integralization data via Lemma \ref{lemma:choose-integralization-data-for-continuation-map}. We then initially choose perturbation data via Lemma \ref{prop:perturbation-data-continuation-map-exist}, and then subsequently a trivial stabilization of these perturbation data which gives the associated topological smoothing a smooth structure via Lemma \ref{prop:continuation-map-compatible-smoothings} (see Definition \ref{def:trivial-stabilization-of-perturbation-data}). This virtual smoothing is framed via Lemmas \ref{lemma:choose-trivializations-of-polarization-map-cyclotomically} and \ref{lemma:framings-continuation-map}. Choosing  $\CC_1 = \CC(H_-, J_-)$ and $\CC_2 = \CC(H_+, J_+)$ in Proposition \ref{prop:relative-variant-of-embedding-construction} and using Lemma \ref{lemma:stabilizing-embeddings-by-parameterizations} to define $\mathbf{E}_i$ and $\bar{\mathbf{E}}_i$ in this Proposition defines the embedding and the extension of the embedding. 

    We now apply Lemma \ref{lemma:filtrations-of-floer-homotopy-types} below; the inclusion \eqref{eq:inclusion-of-lower-category-space} is simply \eqref{eq:subspace-inclusion-flow-subcategory} applied to this setting. Similarly, the homeomorphism  \eqref{eq:identification-of-upper-category-space} is simply \eqref{eq:cone-on-subspace-subcategory} combined with Proposition \ref{prop:effect-of-stabilization-on-floer-homotopy-type} (here $M=1$, as per Remark \ref{rk:can-incorporate-shift-in-canonical-param}; see also \eqref{eq:computation-of-change-of-shift-space}).

\end{proof}

We can use this geometric construction to produce a map of orthogonal spectra:

\begin{proposition}
\label{prop:map-associated-to-continuation-floer-homotopy-type}
    The data $\mathfrak{H}$ as in Proposition \ref{prop:continuation-map-geometry} determines a map $f$, unique up to (equivariant) homotopy, of the form 
    \begin{equation} 
    \label{eq:map-between-floer-homology-spaces}
    \Sigma^{\R \oplus V^1_-, s} \{ \CC(H_-, J_-)\}_{\mathfrak{H}_-} \xrightarrow{f} \Sigma^{\R \oplus V^1_+, s} \{\CC(H_+, J_+)\}_{\mathfrak{H}_+},
    \end{equation}
    which in turns defines the map $\bar{f} = F_{\R \oplus V^1_\mathfrak{H}} f$  on orthogonal spectra sitting in the diagram below:
    \begin{equation}
        \label{eq:map-between-floer-homology-spectra}
        \begin{gathered}
        |\CC(H_-, J_-)|_{\mathfrak{H}_-} \leftarrow F_{\R \oplus V^1_{-, s}} \Sigma^{\R \oplus V^1_{-, s}} |\CC(H-, J_-)|_{\mathfrak{H}_-} \simeq \\
        F_{\R \oplus V^1} \Sigma^{\R \oplus V^1_-, s} \{ \CC(H_-, J_-)\}_{\mathfrak{H}_-} \xrightarrow {\bar{f}}
        F_{\R \oplus V^1} \Sigma^{\R \oplus V^1_+, s} \{\CC(H_+, J_+)\}_{\mathfrak{H}_+}     
         \\\simeq F_{\R \oplus V^1_{+, s}}\Sigma^{\R \oplus V^1_{+, s}}|\CC(H_+, J_+)|_{\mathfrak{H}_+} \to  |\CC(H_+, J_+)|_{\mathfrak{H}_+}.
    \end{gathered}
    \end{equation}
    Here the top arrow pointing left and the bottom arrow pointing right denote the canonical weak equivalences \eqref{eq:desuspension-maps}. In particular, the above diagram defines a map in $[f]$ in the homotopy category from $|\CC(H_-, J_-)|_{\mathfrak{H}_-} \to |\CC(H_+, J_+)|_{\mathfrak{H}_+}$.
\end{proposition}

\begin{proof}
    The existence of the map $f$ follows from Lemma \ref{lemma:maps-of-floer-homotopy-types} stated and proven below. To produce $f$, we apply that Lemma with $\{\CC\} = \{\CC(H_s, J_s)\}_{\mathfrak{H}}$, so $\{\CC_2\} = \{\CC(H_+, J_+\}_{\mathfrak{H}}$, and thus $\Sigma^{V_2}\{\CC_2\} = \Sigma^{V^1_+, s} \{\CC(H_+, J_+)\}_{\mathfrak{H}_+}$. Meanwhile, one gets that $\{\CC_1\} = \{\CC(H_-, J_-)\}_{\mathfrak{H}} \simeq \Sigma^{\R \oplus V^1_{-, s}} \{\CC(H_-, J_-)\}_{\mathfrak{H}_-}$ by \eqref{eq:identification-of-upper-category-space}. Thus the map $f: \{\CC_1\} \to \Sigma^{\R \oplus V_2}\{\CC_2\}$ produced by that lemma takes the the form of \eqref{eq:map-between-floer-homology-spaces}.

    To conclude it suffices to explain the isomorphisms on the top right and the bottom left of \eqref{eq:map-between-floer-homology-spectra}. The top right follows from the fact that $|\CC(H_-, J_-)|_{\mathfrak{H}_-} = F_{V^1_-}\{\CC(H_-, J_-)\}_{\mathfrak{H}_-}$ and the fact that $V^1_{-, s} \oplus V^1_- = V^1$; and the one on the bottom left follows by a similar argument.
\end{proof}

Finally, we can use the constructions regarding continuation homotopy data of Section \ref{sec:continuation-homotopy-geometry} to produce commutative diagrams associated to maps of Floer homotopy types. Below we summarize the geometric results proven in the paper. 

\begin{proposition}
\label{prop:continuation-homotopy-geometry}
    Suppose that the polarization class $\rho \in KO^1_{C_k}(LM)$ vanishes, and we are given $C_k$-equivariant homotopy continuation data $(H_h, J_h)$ as well as 
    \begin{itemize}
        \item Homotopy definition data $\mathfrak{H}_a$, $a \in \{1, 2,3,4\}$ for the associated floer data $(H_a, J_a)$, and 
        \item homotopy definition data $\mathfrak{H}_{ab}$, $(a,b) \in \{(1,2), (1,3), (2, 4), (3,4)\}$ for the associated continuation data $(H_{ab}, J_{ab})$, such that for every such $(a, b)$ the homotopy definition data $\mathfrak{H}_{ab, a}$ obtained from $\mathfrak{H}_{ab}$ by restriction to $\CC(H_a, J_a)$ are stabilizations of $\mathfrak{H}_a$, and similarly the homotopy definition data $\mathfrak{H}_{ab, b}$ obtained from $\mathfrak{H}_{ab}$ by restriction to $\CC(H_b, J_b)$ are stabilizations of $\mathfrak{H}_b$, and 
        \item The integralization data underlying the $\mathfrak{H}_{ab}$ satisfy the conditions of Lemma \ref{lemma:integralization-data-exist-continuation-homotopy}.
    \end{itemize}
    Moreover, all of the above are assumed to be $C_k$-equivariant. Then there exists a homotopy definition datum $\mathfrak{H}_h$ for $\CC(H_h, J_h)$ which such that the homotopy definition data $\mathfrak{H}_{h, ab}$ for $\CC(H_{ab}, J_{ab})$ obtained from $\mathfrak{H}_h$ by restriction are all stabilizations of $\mathfrak{H}_{ab}$, for $(a, b) \in \{(1,2), (1,3), (2, 4), (3,4)\}$. Write $(V^1_h, \{V^1_{ab}\}, \{V^1_a\})$ for the shift spaces assoicated to $(\mathfrak{H}_h, \{\mathfrak{H}_{ab}\}, \{\mathfrak{H}_a\})$. Then there are $G$-vector spaces $(\{V^1_{ab, s}\}, \{V^1_{a, \tilde{s}}\}, \{V^1_{a, s_{ab}}\})$ and decompositions 
    \[ V^1_h = V^1_{ab} \oplus V^1_{ab, s} = V^1_a \oplus V^1_{a, \tilde{s}}, V^1_{ab} = V^1_a \oplus V^1_{a, s_{ab}}\]
    (the last decomposition are the decompositions of shift spaces described in Proposition \ref{prop:continuation-map-geometry}). This forces additional isomorphisms 
    \[ V^1_{a, \tilde{s}} = V^1_{a, s_{ab}} \oplus V^1_{ab, s} = V^1_{a, s_{ab'}} \oplus V^1_{ab', s}\]
    \[V^1_{b, \tilde{s}} = V^1_{b, s_{ab}} \oplus V^1_{ab, s} = V^1_{a', s_{a'b}} \oplus V^1_{a'b, s}.\]
    Writing $\CC_i$ for $\CC(H_i, J_i)$, and writing $\CC_{abc}$ for the subcategory on the objects of $\CC(H_i, J_i)$ as $i=a,b,c$, we have an inclusion 
    \begin{equation}
        \label{eq:continuation-homotopy-exact-sequence-1}
    \Sigma^{V^1_{4, \tilde{s}}} \{\CC_4\}_{\mathfrak{H}_4} \to \{\CC_{234}\}_{\mathfrak{H}}, \{\CC_{234}\}_{\mathfrak{H}}/\Sigma^{V^1_{4, \tilde{s}}} \{\CC_4\}_{\mathfrak{H}_4} \simeq \Sigma^{\R \oplus V^1_{2, \tilde{s}}} \{\CC_2\}_{\mathfrak{H}_2} \vee \Sigma^{\R \oplus V^1_{3, \tilde{s}}} \{\CC_3\}_{\mathfrak{H}_3}.
    \end{equation} 
    Moreover there is an inclusion $\{\CC_{234}\}_{\mathfrak{H}} \to \{\CC(H_h, J_h)\}_{\mathfrak{H}}$ and composing the previous inclusion with this one gives an inclusion
    \begin{equation}
        \label{eq:continuation-homotopy-exact-sequence-2}
    \Sigma^{V^1_{4, \tilde{s}}} \{\CC_4\}_{\mathfrak{H}_4} \to \{\CC(H_h, J_h)\}_{\mathfrak{H}}, \{\CC(H_h, J_h)\}_{\mathfrak{H}}/\Sigma^{V^1_{4, \tilde{s}}} \{\CC_4\}_{\mathfrak{H}_4} \simeq \{ \CC_{123}\}_{\mathfrak{H}};
    \end{equation}
    In particular by the previous two claims there is an inclusion 
    \begin{equation}
        \label{eq:continuation-homotopy-exact-sequence-3}
     \Sigma^{\R \oplus V^1_{2, \tilde{s}}} \{\CC_2\}_{\mathfrak{H}_2} \vee \Sigma^{\R \oplus V^1_{3, \tilde{s}}} \{\CC_3\}_{\mathfrak{H}_3} \xrightarrow{i_{23}} \{ \CC_{123}\}_{\mathfrak{H}}, \{ \CC_{123}\}_{\mathfrak{H}}/Im i_{23} \simeq \Sigma^{\R \oplus \R \oplus V^1_{1, \tilde{s}}}\{\CC_1\}_{\mathfrak{H}_1}. 
     \end{equation}
    
\end{proposition}

\begin{proof}
    We choose the integralization data via
    Lemma \ref{lemma:integralization-data-exist-continuation-homotopy}, the perturbation data via
    Proposition \ref{prop:perturbation-data-continuation-homotopy-exist}, the smoothings via
    Proposition \ref{prop:continuation-homotopy-compatible-smoothings}, and the framings via 
    Lemma \ref{lemma:choose-trivializations-of-polarization-map-cyclotomically} and Lemma \ref{lemma:framings-continuation-map-homotopy}.

    Choosing and $\CC_1 = \CC(H_1, J_2)$ and $\CC_2 =\CC(H_{24}, J_{24}) \cup \CC(H_{34}, J_{34})$ in Proposition \ref{prop:relative-variant-of-embedding-construction} and using Lemma \ref{lemma:stabilizing-embeddings-by-parameterizations} to define $\mathbf{E}_i$ and $\bar{\mathbf{E}}_i$ in this Proposition defines the embedding and the extension of the embedding.

   All homeomorphisms are produced (as above in the proof of Proposition \ref{prop:continuation-homotopy-geometry}) via Lemma \ref{lemma:filtrations-of-floer-homotopy-types}, Proposition \ref{prop:effect-of-stabilization-on-floer-homotopy-type}, and Lemma \ref{lemma:effect-of-stabilizing-by-parameterization-on-floer-homotopy-type}, with the restratification embedding shifts in Lemma \ref{lemma:effect-of-stabilizing-by-parameterization-on-floer-homotopy-type} giving rise to the $\R$ factors in these homeomorphisms above. The computation of the restratification embedding shifts comes coming from the fact that $V^1_{can}$ (see \eqref{eq:obstruction-bundle-for-my-stabilization-continuation-map}; the definition is the same in this case) does not contain coordinates corresponding to the marked continuation map points $s_{A^1_c}$ and $s_{A^2_c}$.
\end{proof}

\begin{lemma}
\label{lemma:model-for-summing-maps}
    Given two maps of pointed Hausdorff spaces $g_1, g_2: X \to Y$, when $\overline{g_1^{-1}(*)^c} \cap \overline{g_2^{-1}(*)^c} = \emptyset$, write $g_1 \underline{\vee} g_2 X \to Y_*$ for the map defined by $g_1 \underline{\vee} g_2(x) = g_i(x)$ if $x \in g_i^{-1}(*)^c$ for $i=1,2$, and otherwise $g_1 \underline{\vee} g_2(x) = *$; by the non-intersection condition this is a well-defined continuous map. When $X$ and $Y$ are pointed CW complexes, this map represents the class of $g_1 + g_2$ in the stable homotopy category.  
\end{lemma}
\begin{proof}
    This occurs because of a canonical homotopy. Indeed, given any two maps $g_1, g_2: X \to Y$, there one can define the class of the map $g_1 + g_2$ concatenating the maps $\Sigma g_1$ and $\Sigma g_2$ along the suspension coordinate, i.e. by composing $\Sigma g_1 \wedge \Sigma g_2$ with the natural map $\Sigma X \to \Sigma X \wedge \Sigma X$; and then formally desuspending this map. But this latter map is homotopic to $\Sigma g_1 \underline{\vee} \Sigma g_2$ whenever the latter makes sense.
\end{proof}

\begin{lemma}
\label{lemma:analyzing-continuation-map-choices-in-commutative-square}
  We can choose the maps associated by Proposition \ref{prop:map-associated-to-continuation-floer-homotopy-type} to the data \eqref{eq:continuation-homotopy-exact-sequence-1}  to be $\Sigma^{V^1_{24, s}}f_{24} \vee (-\Sigma^{V^1_{34, s}} f_{34})$, where the $(-\Sigma^{V^1_{34, s}} f_{34})$ denotes the composition of reflection in $\R$ shift coordinate of $\Sigma^{R \oplus V^1_{3, \tilde{s}}}\{\CC_3\}_{\mathfrak{H}_3}$ with $(\Sigma^{V^1_{34, s}} f_{34})$. Moreover, $\{\CC_{123}\}_{\mathfrak{H}}$ is homotopy equivalent to the cone on $\Sigma^{V^{1}_{12,s}} f_{12} \underline{\vee} \Sigma^{1}_{13, s} f_{13}$, and the map $\{\CC_{123}\}_{\mathfrak{H}} \to \Sigma^{V^1_{4, \tilde{s}} \{\CC_4\}_{\mathfrak{H}_4}}$ associated by Proposition \ref{prop:map-associated-to-continuation-floer-homotopy-type} to \eqref{eq:continuation-homotopy-exact-sequence-2} can be taken to agree with $\Sigma^{V^1_{24, s}}f_{24} \vee \Sigma^{V^1_{34, s}} f_{34}$ on the images of $\Sigma^{\R \oplus V^1_{2, \tilde{s}}} \{\CC_2\}_{\mathfrak{H}_2} \vee \Sigma^{\R \oplus V^1_{3, \tilde{s}}} \{\CC_3\}_{\mathfrak{H}_3}$ in $\{\CC_{123}\}_{\mathfrak{H}}$.

\end{lemma}
\begin{proof}
    The first claim follows from the following reasoning. Given a pushout diagram 
    \[ \begin{tikzcd}
        X \ar[r, "a"] \ar[d, "b"] & Y_1 \ar[d] \\
        Y_2 \ar[r] & Z
    \end{tikzcd}\]
    with the map to $X \to Z$ called $C$ and both $a$ and $b$ being cofibrations, if we have chosen homotopy equivalences $h_1: Y_1/X \to Cone(a)$ and $h_2: Y_2/X \to Cone(b)$ then the wedge of these defines an equivalence $h_3: Y_1/X \wedge Y_2/X \to Cone(c)$. Moreover, using this definition, the map give by the compositio of $h_3$ with the Puppe map is manifestly the wedge products associated to the compositions of the maps associated to $h_i$ and the corresponding Puppe maps for $i=1,2$.

    The next claim goes as follow. If we consider $\{\CC_{123}\}_{\mathfrak{H}}/ \Sigma^{\R \oplus V^1_{2, \tilde{s}}} \{\CC_2\}_{\mathfrak{H_2}}$ then this is homeomorphic to $\Sigma^{\R \oplus V^1_{13, s}} \{\CC_{13}\}_{\mathfrak{H}_{13}}$; by the fact that every triple of maps in a long exact puppe sequence is homotopy equivalent to $X \xrightarrow{f} Y \to C(f)$, we have that $\Sigma^{\R \oplus V^1_{13, s}} \{\CC_{13}\}_{\mathfrak{H}_{13}}$ is homotopy equivalent to $\Sigma^{V^1_{13, s}} f_{13}$. This statement together with the corresponding statemtn with $2$ and $3$ exchanged implies the claim.

    The last claim goes as follows: if we extend the earlier diagram as 
    \[ \begin{tikzcd}
        X \ar[r, "a"] \ar[d, "b"] & Y_1 \ar[d]& \\
        Y_2 \ar[r]& Z \ar[r] & W
    \end{tikzcd}\]
    where the map $Z \to W$ is a cofibration, and we call the map $X \to W$ as $d$, then there is an inclusion $Cone(c) \to Cone(d)$, and we can choose the equivalence $W/X \to Cone(d)$ to restrict to the previously chosen equivalence $Z/X = Y_1/X \wedge Y_2/X \to Cone(c)$. 
\end{proof}

\begin{proposition}
\label{prop:map-associated-to-continuation-homotopy-floer-homotopy-type}
  Using the notation of Proposition \ref{prop:continuation-homotopy-geometry}, let
     $f_{ab}: \Sigma^{\R \oplus V^1_{a, s_{ab}}} \{\CC_a\}_{\mathfrak{H}_a} \to \Sigma^{\R \oplus V^1_{b, s_{ab}}} \{\CC_b\}_{\mathfrak{H}_b}$ be the maps associated by Proposition \ref{prop:map-associated-to-continuation-floer-homotopy-type} to the homotopy definition data $\mathfrak{H}_{ab}$. 
    Then we have a homotopy  
    \begin{equation}
        \label{eq:homotopy-of-continuation-maps-on-spaces}
        \Sigma^{V^1_{24}, s}f_{24} \Sigma^{V^1_{12}, s} f_{12} \underline{\vee} (-\Sigma^{V^1_{34}, s}f_{34}) \Sigma^{V^1_{13}, s} f_{13}. 
    \end{equation}

    With notation as in Proposition \ref{prop:map-associated-to-continuation-floer-homotopy-type}, we have equality in the stable homotopy category
    \begin{equation}
        \label{eq:stable-homotopy-category-homotopy-commutative}
        [f_{24}][f_{12}] = -[f_{34}][f_{13}].
    \end{equation}
\end{proposition}
\begin{proof}
    First, we will explain why the domains and codomains of the maps in \eqref{eq:homotopy-of-continuation-maps-on-spaces} match appropriately so that this equation make sense. Indeed, the map $\Sigma^{V^1_{ab}, s}f_{ab}$ defines a map 
    \begin{equation}
    \label{eq:reinterpreted-suspended-continuation-maps}
    \Sigma^{\R \oplus V^1_{a, \hat{s}}} \{C_a\}_{\mathfrak{H}_a} \simeq \Sigma^{V^1_{ab, s}} \Sigma^{\R \oplus V^1_{a, s_{ab}}}\{C_a\}_{\mathfrak{H}_a} 
    \xrightarrow{\Sigma^{V^1_{ab}, s} f_{ab}}
    \Sigma^{V^1_{ab, s}} \Sigma^{\R \oplus V^1_{b, s_{ab}}}\{C_b\}_{\mathfrak{H}_b} \simeq \Sigma^{\R \oplus V^1_{b, \hat{s}}} \{C_b\}_{\mathfrak{H}_b}, 
    \end{equation}
    using the identities $V^1_{a, s_{ab}} \oplus V^1_{ab, s} = V^1_{a, \hat{s}}$; we use these isomorphisms to make sense of the compositions in \eqref{eq:homotopy-of-continuation-maps-on-spaces}.

    Then, we will explain how \eqref{eq:homotopy-of-continuation-maps-on-spaces} implies \eqref{eq:stable-homotopy-category-homotopy-commutative}. Indeed, the identity \eqref{eq:stable-homotopy-category-homotopy-commutative} holds if the identity 
    \begin{equation}
        \label{eq:suspended-identity}
        \Delta [f_{24}]\Delta[f_{12}] = \Delta[f_{34}]\Delta[f_{13}] 
    \end{equation} 
    holds, where $\Delta$ denotes the functor homotopy category $X \to F_{\R \oplus V^1_h}\Sigma^{\R \oplus V^1_h}X$ on the stable homotopy category. This latter identity holds in turn of we an find representing maps for each of the $\Delta f_{ab}$ in the category of orthogonal spectra and a homotopy between the composed maps. Now, we can apply the functor $F_{\R \oplus V^1_{h}}$ to the maps in \eqref{eq:reinterpreted-suspended-continuation-maps} and apply the isomorphism 
    \[F_{\R \oplus V^1_{h}}\Sigma^{\R \oplus V^1_{a, \hat{s}}} \{C_a\}_{\mathfrak{H}_a} \simeq F_{\R \oplus V^1_{a, \hat{s}}}\Sigma^{\R \oplus V^1_{a, \hat{s}}} |\CC(H_a, J_a)|_{\mathfrak{H}_a}\]
    to interpret these maps potentially being representing maps for the $\Delta [f_{ab}]$; given this claim, the homotopy \eqref{eq:homotopy-of-continuation-maps-on-spaces} establishes \eqref{eq:suspended-identity} and thus \eqref{eq:stable-homotopy-category-homotopy-commutative}. To see this latter claim about representing maps, we simply note that applying $F_{V^1_{ab, s}}\Sigma^{V^1_{ab, s}}$ to the map $\bar{f}_{ab}$ of \eqref{eq:map-between-floer-homology-spectra} and applying isomorphisms in \eqref{eq:reinterpreted-suspended-continuation-maps} to the maps inside the desuspension functors, and the isomorphisms $V^1_{ab, s} \oplus V^{1}_{ab} = V^1_h$ gives exactly $F_{\R \oplus V^1_{h}}$ applied to the maps in \eqref{eq:reinterpreted-suspended-continuation-maps}. This proves the claim about representing maps.

    It remains to justify \eqref{eq:homotopy-of-continuation-maps-on-spaces}. But this follows from Lemma \ref{lemma:analyzing-continuation-map-choices-in-commutative-square}. Indeed, given maps $X \xrightarrow{f} Y \xrightarrow{g} Z$, an extension of $g$ to a map $Cone(f) \to Z$ specifies a homotopy of $g \circ f \sim *$. This lemma states that we are precisely in this position with $f = \Sigma^{V^1_{12, s}}f_{12} \underline{\vee} \Sigma^{V^1_{13, s}} f_{13}$ and $g = \Sigma^{V^1_{24,s}} f_{24} \wedge \Sigma^{V^1_{34, s}} f_{34}$. So we have a homotopy  $\Sigma^{V^1_{24}, s}f_{24} \Sigma^{V^1_{12}, s} f_{12} \underline{\vee} (-\Sigma^{V^1_{34}, s}f_{34}) \Sigma^{V^1_{13}, s} f_{13} \sim *$, as desired.
\end{proof}

\subsection{Comparison with Floer Homology}
\label{sec:comparison-with-floer-homology}

\begin{lemma}
\label{lemma:filtrations-of-floer-homotopy-types}
    Let $\CC$ be a virtually smooth semi-freely framed proper flow category together with an embedding $\mathbf{E}$ and its extension $\bar{\mathbf{E}}$. If $\CC_2 \subset \CC$ is a downwards-closed subposet of $\CC$ that is sent to itself by the $G$-action, then, letting $\CC_2$ also denote the full virtually smooth flow subcategory on the objects of $\CC_2$ equipped with the restrictions of $\mathbf{E}$ and $\bar{\mathbf{E}}$, we have a q-cofibration
    \begin{equation}
    \label{eq:subspace-inclusion-flow-subcategory}
    \Sigma^{V_2}\{\CC_2\} \subset \{\CC\}
    \end{equation}
    where $V_2$ is the perpendicular to the image of the shift space associated to $\CC_2$ in the shift space associated to $\CC$. Similarly, we have a level $q$-cofibration \cite[Theorem~2.4]{mandell2002equivariant}
    \begin{equation}
        \label{eq:subspectrum-inclusion-flow-subcategory}
        F_{V_2} \Sigma^{V_2} |\CC_2| \subset |\CC|.
    \end{equation} 
    Moreover, we have that there is an (equivariant) homeomorphism
    \begin{equation}
        \label{eq:cone-on-subspace-subcategory}
        \{\CC\}/\Sigma^{V_2}\{\CC_2\} \simeq \{\CC_1\}
    \end{equation} 
    where $\CC_1$ is the full flow subcategory of $\CC$ on the objects $Ob(\CC_2) \setminus Ob(\CC_2)$, and is also equipped with the restricted embedded framings, embeddings, and extensions of embeddings. 
    
    Thus, in particular, given 
    a grading  $\alpha: \CC \to \Z$ on the poset $\CC$, i.e. an integer-value function on the objects of $\CC$ such that $\alpha(x) \geq \alpha(y)$ if $x > y$ and $\alpha(gx) = \alpha(x)$ for all $x \in \CC$ and all $g \in G$, then there is a filtration by sub-spectra, i.e. there is are diagrams of $G-CW$-subspectra
    \begin{equation} 
    \label{eq:filtration-by-subcomplexes-floer-homotopy-space}
    \{\CC\}_{\leq k-1} \subset  \{\CC\}_{\leq k} \subset \{\CC\}_{\leq k+1} \subset \ldots, 
    \end{equation}
    and a diagram of orthogonal $G$-spectra
    \begin{equation}
    \label{eq:filtration-floer-homotopy-type}\cdots \to |\CC|_{\leq k-1} \to |\CC|_{\leq k} \to |\CC|_{\leq k+1} \to \cdots 
    \end{equation}
    where each map is a level $q$-cofibration, and 
    \[ \{\CC\} = \cup_k \{\CC\}_{\leq k}, |\CC| = \colim_k |\CC|_{\leq k}.\]
\end{lemma}

\begin{proof}
    The inclusion \eqref{eq:subspace-inclusion-flow-subcategory} follows from the construction, since, writing $V^1_{\CC}$ for the shift space associated to $\CC$, we have that $V^1_{\CC_2}(x)^\perp \oplus V_2 = V^1_{\CC}(x)^\perp$; thus the attaching maps \eqref{eq:fundamental-attaching-map} for the the cells of $\{\CC\}$ coming from objects of $\CC_2$ are therefore just suspensions of the corresponding attaching maps in $\{\CC_2\}$. The inclusion \eqref{eq:subspectrum-inclusion-flow-subcategory} then follows immediately from the fact that $V_2 \oplus V^1_{\CC_2} = V^1_\CC$. The claim \eqref{eq:cone-on-subspace-subcategory} follows from the fact that when defining $\{\CC_1\}$ one collapses all parts of the boundaries of the cells associated to objects of $\CC_1$ that would have mapped to the cells associated to objects of $\CC_2$ to basepoint instead (see Condition F of Definition \ref{def:extension} for extensions of embeddings). The remaining claims follow immediately from the definitions. 
\end{proof}

\begin{lemma}
\label{lemma:maps-of-floer-homotopy-types}
    In the setting of Lemma \ref{lemma:filtrations-of-floer-homotopy-types}, there is a map defined up to contractible choice
    \[ f: \{\CC_1\} \to \Sigma^{\R \oplus V_2} \{\CC_2\}\]
    as well as 
    \[ \bar{f} = F_Vf: |\CC_1| \to F_{V_2}\Sigma^{\R \oplus V_2} |\CC_2|. \]
\end{lemma}
\begin{proof}
This follows from \eqref{eq:subspectrum-inclusion-flow-subcategory}, \eqref{eq:cone-on-subspace-subcategory} by composing the last map in the Puppe sequence 
\[ \Sigma^{V_2} \{\CC_2\} \to \{\CC\} \to Cone(\Sigma^{V_2} \{\CC_2\} \to \{\CC\}) \to \Sigma^{\R V_2} \{\CC_2\} \]
with an equivariant homotopy inverse of the equivalence \eqref{eq:cone-on-subspace-subcategory}.
\end{proof}

\begin{lemma}
\label{lemma:filtered-maps-of-floer-homotopy-types}
    Let $\CC_2$ be a downwards-closed flow-subcategory of a virtually smooth flow category $\CC$ with an $E^s_2$ framing and an extension of an embedding; equip $\CC_2$ with the restrictions of all of this data. Let $\CC_1$ be the full subcategory of $\CC$ on $Ob(\CC) \setminus Ob(\CC_2)$, and equip it with the restrictions of all the data associated to $\CC$. Suppose that we are given gradings $\alpha_1$ and $\alpha_2$ of the posets $\CC_1$ and $\CC_2$, and $\CC$ \emph{respects} these gradings in the sense that if $x_i \in \CC_i$ for $i=1,2$ and $x_i > x_2$ in $\CC$ then $\alpha_1(x_1) \geq \alpha_2(x_2)$. Then the map 
    \[f: \{\CC_1\} \to \Sigma^{\R \oplus V_2} \{\CC_2\}
    \]
    defined by Lemma \ref{lemma:maps-of-floer-homotopy-types}
    where $V_2$ as in that lemma, respects the filtrations on the domain and codomain induced by suspending the filtrations of Equation \ref{eq:filtration-by-subcomplexes-floer-homotopy-space}, Lemma \ref{lemma:filtrations-of-floer-homotopy-types}. Moreover, the space of choices of such $f$ which respect these filtrations is contractible.

\end{lemma}
\begin{proof}
Define $\alpha: \CC \to \Z$ by $\alpha(x) = \alpha_1(x)$ if $x \in \CC_2$, and otherwise $\alpha(x) = \alpha_2(x)$ if $x \in \CC_2$. 
    Then we have a filtration of $\{\CC\}$  \eqref{eq:filtration-by-subcomplexes-floer-homotopy-space} associated to $\alpha$, and the restriction of this filtration to $\Sigma^{V_2} \{\CC_2\}$ is the suspension of the filtration by $\alpha_2$ of $\{\CC_2\}$. This defines a filtration by $\Z$ of $Z = Cone(\Sigma^{V_2} \{\CC_2\} \to \{\CC\})$ and the projection \eqref{eq:cone-on-subspace-subcategory} sends this filtration onto the filtration of $\{\CC_1\}$ associated to $\alpha_1$. Thus the inverse map associated to \eqref{eq:cone-on-subspace-subcategory} can be resepects  these filtrations as well, and it is an inverse to an explicit homeomorphism so there are no choices to make. One then sees that the Puppe map $Z \to \Sigma^{\R \oplus V_2} \{\CC_2\}$ respects filtrations as well (with the latter filtered by  the suspension of the filtration associated to $\alpha_2$). The space of choices of a Puppe map is contractible, thus proving the lemma, 
\end{proof}

Recall that the homology groups of a cofibrant spectrum $X$ as $H_*(X, \Z) = \pi_*(X \wedge H\Z)$, where $H\Z$ is the Eilenberg-MacLane spectrum.   Moreover, if $X$ is a colimit over $k$ of q-cofibrant spectra $X_{\leq k}$ along q-cofibrations, then there is a spectral sequence 

\begin{equation}
    \label{eq:spectral-sequence-assoc-to-filtration}
    E^1_{s,t} = H_{s+t}(Cone(X_{\leq s-1} \to X_{\leq s}), \Z) \Rightarrow H_*(X, \Z),
\end{equation} 
which converges if the colimit is finite. 

An elementary computation shows that, using the notation of Lemma \ref{lemma:filtrations-of-floer-homotopy-types} if $X = |\CC|$, $X_{\leq k} = |\CC|_{\leq k}$ as in \eqref{eq:filtration-floer-homotopy-type},  and $\alpha$ is \emph{stict} (i.e. if $x > y$ then $\alpha(x) > \alpha(y)$), then there is a level equivalence
\begin{equation}
\label{eq:homotopy-equivalence-of-cones} 
Cone(X_{\leq s-1} \to X_{\leq s}) \simeq \bigvee_{x \in \CC : \alpha(x) = s} \SS^\mu(x)
\end{equation}
where $\mu(x)$ is the grading on $\CC$. Indeed, we have that $X_{\leq s} = F_{V^1_{max}} \{\CC\}_{\leq s}$, and cone on the inclusion
\[ \{\CC\}_{\leq s-1} \to \{\CC\}_{\leq s}\]
is 
\[ \bigvee_{x \in \CC : \alpha(x) = s} \SS^{\R^{\bar{A}(x)+1} \oplus V^0(x) \oplus V^1(x)^\perp}; \]
by the grading condition \eqref{eq:index-condition-for-embedding-of-framing} we have that 
\[\mu(x) = \dim V^0(x) - \dim V^1(x) + (\bar{A}(x)+1) = \bar{A}(x)+1 + \dim V^0(x) + \dim V^1(x)^\perp - \dim V^1_{max}\]
which proves the equivalence \eqref{eq:homotopy-equivalence-of-cones}. 

Now, given an orientation of $V^0(x)$, $V^1(x)$, and $V^1_{max}$, we get an isomorphism of $\tilde{H}_*(\SS^{\R^{\bar{A}(x)+1} \oplus V^0(x) \oplus V^1(x)^\perp})$ (where $\tilde{H}$ denotes reduced homology) with a copy of $\Z$ in the relevant degree.

Thus, after choosing \emph{orientation data} $\mathfrak{o}$, which consists of an orientation of $V^1_{max}$ as well as orientations of $V^0(x)$ and $V^1(x)$ for each $x$,  the $E^1$ pages of the associated spectral sequence takes the form 
\begin{equation}
\label{eq:spectral-sequence-of-filtered-flow-category}
E^1_{s,t} = \oplus_{x \in \CC: \alpha(x) = s, \mu(x) = s+t} \Z.
\end{equation}

If we are in the situation of Lemma \ref{lemma:filtered-maps-of-floer-homotopy-types}, the map $f$ induces a map $\bar{f} = F_{V^1_{max}} f: |\CC_1| \to F_{V_2} \Sigma^{\R \oplus V_2} |\CC_2|$, which in turn induces a map of spectral sequences 
\begin{equation}
    \label{eq:map-of-spectral-sequence-of-filtered-flow-category}
    \bar{f}^1_{s,t}: (E^1_{\CC_1})_{s,t} \to (E^1_{\CC_2})_{s,t+1}.
\end{equation} 
(We will find that this shift by $1$ in the grading cancels out when applying this computation to flow categories associated to continuation maps.) 

A choice of orientation data in turn gives orientations of $V^i(x,y)$ for $x > y$ (by requiring that the map $V^i(x,y) \oplus V^i(y) \to V^i(x)$ is orientation preserving). In particular, if $\mu(x) = \mu(y)+1$, then \eqref{eq:index-condition-for-embedding-of-framing} implies that the map $e_{xy}$ associated to the embedding $\mathbf{E}$ of $\CC$ is a map between equidimensional manifolds. The definition of the attaching map $e_x$ then shows:
\begin{lemma}
\label{lemma:identify-counts-with-degrees}
    Suppose that for all $s$ and all $x \in \CC_1$ and $y\in \CC_2$ with $\alpha_1(x) = \alpha_2(y) = s$, the map $e_{xy}$ associated to the embedding of $\CC$ has isolated zeros in its interior. Then 
    the map $\bar{f}^1_{s,t}$ sends the fundamental class of $\SS^{\mu(x)}$ to a sum of fundamental classes of $\SS^{\mu(y)+1}$ for $y \in \CC_2$ such that $\mu(y) = \mu(x)-1$ and $\alpha_1(x) = \alpha_1(y)$, with the coefficient of the fundamental class of $\SS^{\mu(y)+1}$ being exactly the signed number of zeros of $e_{x,y}$ on the interior, with the signs defined via the standard orientation of $\R^{\bar{A}(x,y)}$ and the orientations of $V^i(x,y)$ assigned above. $\blacksquare$
\end{lemma}

Now, let us suppose that we are in a geometric setting, i.e. $\CC = \CC(H, J)$ for some convex nondegenerate Floer data $(H, J)$, or $\CC = \CC(H_s, J_s)$ for convex nondegenerate continuation data $(H_s, J_s)$. Then,  the orientation data $\mathfrak{o}$ actually trivializes the \emph{orientation lines} associated to the objects of $\CC$ (\cite{abouzaid-sh, seidel2008fukaya, floer1993coherent}); in other words, it defines orientation data in the sense of Section \ref{sec:grading-conventions}), which lets us associate signs $s_u, s_{u'} \in \pm$ to regular Floer trajectories $u \in \CC(H, J)(x,y)$ or regular continuation map trajectories $u' \in \CC(H_s, J_s)(x', y')$ via usual coherent orientations procedure \cite{floer1993coherent} outline in Section \ref{sec:grading-conventions}. Indeed, the orientation lines $o(x)$ of that section are canonically isomorphic to the determinant lines of the operators $D_x$ \eqref{eq:asymptotic-fredholm-operator-not-stabilized}; since $D_x^{\mathfrak{P}'}$ \eqref{eq:asymptotic-fredholm-operator-stabilized} is a stabilization of $D_x$ by $V^1(x)$, one has that 
\[ o(x) \simeq \det D_x \simeq \det D_x^\mathfrak{P} \tensor \det V^1(x) = \det V^0(x) \tensor \det V^1(x). \]
We are now able to explain how to compare usual computations in Floer homology with the maps on homology of Floer homotopy defined above: 

\begin{theorem}
\label{thm:comparison-with-floer-homology}
    Suppose that the admissible Floer data $(H, J)$ are regular. In the setting of Lemma \ref{lemma:filtrations-of-floer-homotopy-types}, let us choose $\alpha$ to be the grading $\mu$ by the Conley-Zehnder index $\tilde{CZ}$ (see Section \ref{sec:grading-conventions}). Then we can choose homotopy definition data $\mathfrak{H}$ and orientation data $\mathfrak{o}$ for $\CC(H, J)$ such that the differential on the first page of the first pace of the spectral sequence \eqref{eq:spectral-sequence-of-filtered-flow-category} agrees with the Floer differential. 

    Similarly, suppose that we have a continuation map category $\CC(H_s, J_s)$ such that for some strict filtrations $\alpha_i$ of $\CC(H_i, J_i)$ (with $(H_-, J_-) = (H_1, J_1)$ and $(H_+, J_+) = (H_2, J_2)$) with $\CC(H_s, J_s)$ respecting this strict filtration, we have that for every $x \in \CC(H_1, J_1)$ and $y \in \CC(H_2, J_2)$ with $\alpha_1(x) = \alpha_2(y)$ and $\mu^{CZ}(x) = \mu^{CZ}(y)$,  every unbroken continuation map trajectory $u' \in \CC(H_s, J_s)(x,y)$ is regular. Then the coefficient of the fundamental class of $\SS^{\mu(y)+1}$ in  $\bar{f}^1([\SS^{\mu(x)}])$ is exactly $\sum_{u'}  s_{u'}$ where the sum runs over these trajectories (up to a global sign).
\end{theorem}
\begin{proof}
    We will prove the statement about the Floer differential; the proof of the second statement is similar. Write $\CC'(x,y) = (T(x,y), V(x,y), \sigma(x,y))$. There is a map $\phi: \CC(H, J)(x,y) \to T(x,y)$ which is a homeomorphism onto the zero-set of $\sigma(x,y)$, and thus a map $\CC(x,y) \to \paramspace(x,y)$, $u \mapsto a(\phi(u))$. In particular, regular index $1$ Floer trajectories $u$ are mapped by $\phi$ to isolated zeros of $\sigma(x,y)$, and to the top open stratum of $\paramspace(x,y)$. We will show the claim in the theorem under the condition that the perturbation map $\lambda_{V^1(x,y)}$ \eqref{eq:perturbation-data-2} is identically zero in an open neighborhood of $a(\phi(u))$, for every (regular) index $1$ Floer trajectory $u$; this can be achieved in the proofs of Propositions \ref{prop:compatible-perturbation-data-exist} (and \ref{prop:perturbation-data-continuation-map-exist} in the continuation map case). The choice of orientation data means that there is an integer assigned to the zero of $e_{x,y}$ corresponding to $u$; this is the same integer as one assigned to the zero $\phi(u)$ of $\sigma(x,y)$, when we give $V(x,y)$ and $TT(x,y)$ the orientations coming from the framings of these spaces and the orientation data. We will first argue that this integer $s'_u$ is $\pm 1$; we will then argue that in fact $s'_u$ agrees with $s_u$, the sign assigned to $u$ by Floer homology. 
    
    Let us argue the first point. By Equation \ref{eq:thickening-partial-framing} and the condition on the triviality of $\lambda_{V^1(x,y)}$ near $u$ we have an isomorphism 
    \begin{equation}
    \label{eq:summarizing-decompositions-of-tangent-bundles}
        \begin{gathered}TT(x,y)_{\phi(u)} = T(x,y)_{a(\phi(u))} \oplus V^1(x,y)_u \oplus \ker D_u = 
        \R^{\mathfrak{A}'(x,y)} \oplus V^1(x,y) \oplus \R
        \end{gathered}
    \end{equation} 
    where $D_u$ is the linearization of the \emph{unperturbed} Floer equation. Now one can produce a convenient chart around $\phi(u) = (\tilde{u}, a(\phi(u)), 0)$ with domain $V^1(x,y) \oplus \R^{\bar{A}(x,y)+1}$ by: translating the map $\tilde{u}$ in the $s$ direction, moving the $\bar{A}(x,y)+1$ marked points on its domain (up to global simultaneous translation), or deforming $0 \in V^1(x,y)$ to a nonzero element. In this chart, using the definition of $\sigma(x,y)$ together with the description \begin{equation}
    \label{eq:summarizing-decomposition-of-obstruction-bundle}
        (V(x,y)^1_T)_u = V^1(x,y) \oplus V^1_{can}(x,y)
    \end{equation} 
    the $V^1(x,y)$ component of $\sigma(x,y)$ is simply the composition of the projection to $V^1(x,y)$ with the identity map $V^1(x,y) \to V^1(x,y)$, and the second component of $\sigma(x,y)$ is the projection to the $\R^{\bar{A}(x,y)+1}$ component of the chart together with a map $\R^{\bar{A}(x,y)+1} \to \R^{\bar{A}(x,y)+1}$ determined by $\delta A$. This latter map is always manifestly orientation-preserving; thus, 
    we have that $s'_u = \pm 1$, with the sign given by the difference in orientations of $V^1(x,y)$ coming from the framing of the obstruction bundles, and of $V^1(x,y)$ induced from the inclusion $V^1(x,y) \oplus \R \simeq \ker D_u^{\mathfrak{P}'}$. 
    
    We will now compare this sign with $s_u$. To do that we must consider how the framings we produce orient $T(x,y)_u$ and $V(x,y)_u$. The latter are produced from the orientation data and the construction of Section \ref{sec:producing-framings}, for which $V^1_{can}(x,y)$ and the component $\R^{\bar{A}(x,y)}$ (this is all but the last two terms of \eqref{eq:summarizing-decompositions-of-tangent-bundles}) are framed in a canonical way while the terms $\ker D_u^{\mathfrak{P}} = V^1(x,y)_u \oplus \ker D_u$ and $V^1(x,y)_u$ get induced orientations (from the framing) via the procedure described in that section, starting from the orientations of of the $\ker D_x^{\mathfrak{P}}$ and $V^1(x)$ -- i.e., the orientation data. On the other hand, in Floer homology, one has that the orientation line 
    \[ o(x) \simeq \det \ker D_x^{\mathfrak{P}} \tensor \det V^1(x)\]
    receives an orientation from the orientation data, and the sign $s_{xy}$ is the comparison sign from comparing the orientation of $o(x)$ produced from orientation data from the orientation of $o(x)$ produced from the canonical orientation of the determinant line $\ker D_u$ of $u$ and the orientation of $o(y)$ via the linear gluing map $\ker D_u \tensor o(y) \to o(x)$. But we have that the framings $a_{xy}$ of Proposition \ref{prop:framings-exist} are chosen precisely such they induce the orientations on $\ker D_u^{\mathfrak{P}'}$ induced via the stabilized gluing maps $D^{\mathfrak{P}'}_{V^1(x), V^1(y)}$ of \eqref{eq:asymptotic-fredholm-operator-stabilized}; this shows that the Floer comparison sign and the sign $s'_u$ agree. The proof for continuation maps maps is similar.
\end{proof}

\begin{proposition}
\label{eq:compute-for-admissible-Floer-data}
    Let $(H, J)$ be admissible regular Floer data. Then there exists homotopy definition data $\mathfrak{H}$ such that if we define the Floer homotopy type $|\CC(H, J)|$ using $\mathfrak{H}$ we have an isomorphism 
    \[ \pi_*(|\CC(H, J)| \wedge H\Z) \simeq HF_*(H, J).\]
    In fact, recall that the Eilenberg-MacLane spectrum $H\Z$ has a single nonzero homotopy group, namely $\pi_0(H\Z) = \Z$; thus given the chain complex $CF_*(H, \Z)$ and choosing orientation data $\mathfrak{o}$ so that we get isomorphisms $CF_k(H, \Z) = \Z^{a_k}$, there is an associated $H\Z$-module spectrum which we denote $|CF_*(H, J)|$ which is built by attaching wedges of suspensions of $H\Z$ via the maps given by the coefficients of the differential, and one has a weak equivalence 
    \[ |\CC(H, J)| \wedge H\Z \simeq |CF_*(H, J)|.\]
\end{proposition}
\begin{proof}
    Choosing $\alpha$ to be the Conley-Zehnder index, this follows by induction  from \eqref{eq:homotopy-equivalence-of-cones}  and \eqref{thm:comparison-with-floer-homology}. 
\end{proof}

\subsection{Choosing data compatibly with the cyclotomic action}
\label{sec:cyclotomic-compatibility-geometry}

We will show how to imitate the colimit \eqref{eq:sh-colimit-on-homology} in the setting of Floer homotopy types in the setting where none of the Floer data or continuation data are required to be regular. First, fix a \emph{complete $S^1$-universe} $\mathfrak{U}$; by restriction of representations this universe is also a complete $C_k$-universe for all $k$. Assume that the $S^1$-equivariant polarization class $\rho \in KO^1_{S^1}(LM)$ vanishes, and inductively choose a sequence of cyclotomically compatible trivializations of the polarization map for all $(H^{\#p^k}, J^{\#p^k})$ (see Lemma \ref{lemma:choose-trivializations-of-polarization-map-cyclotomically}) once and for all.


By repeatedly applying Proposition \ref{prop:continuation-map-geometry} , we produce a sequence of homotopy definition data $\mathfrak{H}_k^{cont}$ for $\CC(H_s^{\#p^k}, J_s^{\#p^k})$ which have the following property. Write
\[\mathfrak{H}^+_k = \mathfrak{H}_k^{cont}|_{\CC(H^{\#p^k}, J^{\#p^k})}, \]
\[\mathfrak{H}^-_k = \mathfrak{H}_k^{cont}|_{\CC(H^{\#p^{k+1}}, J^{\#p^{k+1}})}, \]
\[\mathfrak{H}_k = \mathfrak{H}^-_{k-1} \text{ for } k \geq 1, \mathfrak{H}_0 = \mathfrak{H}^+_0. \]

Proposition \ref{prop:continuation-map-geometry} allows us to require that for $k \geq 1$, $\mathfrak{H}_k^+$ is a $C_k$-equivariant stabilization of $\mathfrak{H}_{k-1}^-$ (where $\mathfrak{H}_{k-1}$ was chosen such that $\mathfrak{H}_{k-1}^-$ has a $C_k$-equivariant structure). Write 
\begin{equation}
    \label{eq:shift-space-decomposition-in-colimit}
    \tilde{V}(\mathfrak{H}^{cont}_k) = \tilde{V}(\mathfrak{H}_k^-) = \tilde{V}(\mathfrak{H}_{k-1}^-) \oplus \tilde{V}(\mathfrak{H}_{k-1}^-, s)
\end{equation}
for the corresponding $C_k$-equivariant decomposition of shift spaces. 


Using \eqref{eq:inclusion-of-lower-category-space} and \eqref{eq:identification-of-upper-category-space}, the map associated to a Puppe sequence defines, for every $k \geq 0$, a $C_{p^{k}}$-equivariant map
\begin{equation}
    \label{eq:continuation-maps-on-floer-homotopy-types-in-colimit}
     g_1^k: \Sigma^{\R} \{\CC(H^{\# p^{k+1}}, J^{\#p^{k+1}})\}_{\mathfrak{H}^-_k} \to \Sigma \{\CC(H^{\# p^k}, J^{\# p^k})\}_{\mathfrak{H}^+_k} \simeq \Sigma^{\R \oplus \tilde{V}(\mathfrak{H}_{k-1}^-, s)} \{\CC(H^{\# p^{k}}, J^{\#p^{k}})\}_{\mathfrak{H}_{k-1}}.
\end{equation}
where the second equivalence only makes sense for $k \geq 1$.

Using \eqref{eq:desuspension-maps}, the decompositions \eqref{eq:shift-space-decomposition-in-colimit} define associated maps 
\begin{equation}
    \label{eq:destabilizing-map-on-floer-homotopy-types-in-colimit}
     g_2^k: F_{ \R \oplus \tilde{V}(\mathfrak{H}^-_k)}\Sigma^{\R \oplus \tilde{V}(\mathfrak{H}_{k-1}^-, s)} \{\CC(H^{\# p^{k}}, J^{\#p^{k}})\}_{\mathfrak{H}^-_{k-1}}. \to  F_{ \R \oplus \tilde{V}(\mathfrak{H}_k-1)}\Sigma^\R\{\CC(H^{\# p^{k}}, J^{\#p^{k}})\}_{\mathfrak{H}^-_{k-1}}. 
\end{equation}
Define 
\[ |\CC'(H, J)|'_{\mathfrak{H}} = F_{\R\oplus\tilde{V}^1_\mathfrak{H}}\Sigma^\R\{\CC'(H, J)\}\]

The composition of \eqref{eq:continuation-maps-on-floer-homotopy-types-in-colimit} and \eqref{eq:destabilizing-map-on-floer-homotopy-types-in-colimit} is a map 
\begin{equation}
    \label{eq:joint-map-on-floer-homotopy-types-in-colimit}
    g^k = g_2^k \circ (F_{\R \oplus \tilde{V}(\mathfrak{H}^-_k)} g_1^k): |\CC(H^{\# p^k}, J^{\# p^k})|'_{\mathfrak{H}_k} \to |\CC(H^{\# p^{k-1}}, J^{\# p^{k-1}})|'_{\mathfrak{H}_k};
\end{equation}
this map makes sense even for $k=1$ by simply avoiding \eqref{eq:destabilizing-map-on-floer-homotopy-types-in-colimit}. 

Dualizing, we get $C_{k-1}$-equivariant maps
\begin{equation}
    f_k: F(|\CC(H^{\# p^{k-1}}, J^{\# p^{k-1}})|'_{\mathfrak{H}_k}, \SS) \to F(|\CC(H^{\# p^k}, J^{\# p^k})|'_{\mathfrak{H}_k}, \SS). 
\end{equation}
for all $k \geq 1$.

Now, we have one more geometric proposition:
\begin{proposition}
\label{prop:compatible-homotopy-continuation-data}
    The homotopy continuation data $\mathfrak{H}_k$ can be chosen such that 
    \[ \{\CC(H_s^{\#p^k}, J_s^{\#p^k})\}^{C_p} \simeq \Sigma^{\tilde{V}(\mathfrak{H}^-_{k-1}, s)^{C_p}} \{\CC(H_s^{\#p^{k-1}}, J_s^{\#p^{k-1}})\}.\]
\end{proposition}

Using this proposition, we can choose the maps $g_1^k$ (which come from Puppe sequences) such that
\[ (g_1^k)^{C_p} = \Sigma^{\tilde{V}(\mathfrak{H}^-_{k-1}, s)^{C_p}}g_1^{k-1}.\]

Writing
\[|H^k| :=|\CC(H^{\# p^{k-1}}, J^{\# p^{k-1}})|'_{\mathfrak{H}_k}, \]

and
\[ \SS|H^k| := F(|\CC(H^{\# p^{k-1}}, J^{\# p^{k-1}})|'_{\mathfrak{H}_k}, \SS) \]
for short, we have a diagram of maps 
\begin{equation} 
\label{eq:unreplaced-diagram}
\SS|H^{\# p^0}| \xrightarrow{t_0} \SS|H^{\# p^1}| \xrightarrow{t_1} \SS|H^{\# p^2}| \xrightarrow{t_2} \cdots.
\end{equation}

Here the map $t_k$ is $C_{p^k}$-equivariant. 

\begin{proof}[Proof of Proposition \ref{prop:compatible-homotopy-continuation-data}]
    We simply make sure that the integralization data for $\CC(H_s, J_s)$ is chosen such that it satisfies the conditions of Lemma \ref{lemma:cyclotomic-compatibility-for-integralization-data-continuation-map}; this is always possible by choosing $c(s)$ to decrease a sufficient amount. This statement is proven by induction on $k$ with the base case at $k=0$ being Proposition \ref{prop:continuation-map-geometry}. Given the case for $k-1$, the case for $k$ follows from
     follows from Proposition \ref{prop:cyclotomic-compatibility-all-perturbation-data-choices}, Lemma \ref{lemma:smoothings-cyclotomic-compatibility}, Lemma \ref{lemma:framings-cyclotomic-compatibility},  Proposition \ref{prop:cyclotomic-variant-of-embedding-construction} (after first using Lemma \ref{lemma:restratifications-exist} and Lemma \ref{lemma:stabilizing-embeddings-by-parameterizations} to induce an embedding and the extension of the embedding for $\CC'(H_s^{\# p^k}, J_s^{\# p ^k})^{C_p}$ from that chosen for $\CC'(H_s^{\# p^{k-1}}, J_s^{\# p ^{k-1}})$), and subsequently using, Theorem \ref{thm:geometric-fixed-points-of-flow-category}, Lemma \ref{lemma:effect-of-stabilizing-by-parameterization-on-floer-homotopy-type}, and Proposition \ref{prop:simple-restrat-effect-on-homotopy-type}. One simply notes that now the perturbation data satisfy the hypotheses of Lemma \ref{lemma:cyclotomic-compatibility-for-integralization-data-continuation-map} for $\CC(H_s^{p^k},J_s^{\# p^k})$, and one proceeds onwards with the induction.
\end{proof}

\subsection{Invariance of equivariant Floer homotopy}
\label{sec:invariance-of-equivariant-floer-homotopy}
\begin{proposition}
\label{prop:invariance-of-floer-homotopy}
    Suppose that $H_\pm$ are parts of admissible Floer data $(H_\pm, J_\pm)$, with $H_\pm$ having the same slope at infinity if $M$ has contact boundary. Then, choosing an admissible continuation datum $(H_s, J_s)$ between $(H_\pm, J_\pm)$ (assuming the relevant condition on polarization classes) choosing homotopy definition data $\mathfrak{H}_s$ as in Proposition \ref{prop:continuation-map-geometry}, the map $f$ defined by Proposition \ref{prop:map-associated-to-continuation-floer-homotopy-type} on Floer homotopy types is an equivalence of (nonequivariant) spectra.
\end{proposition}
\begin{proof}
Choose a continuation homotopy datum $(H_h, J_h)$ such that $(H_{1}, J_1) = (H_3, J_3) = (H_4, J_4) =  (H_-, J_-), (H_2,J_2) = (H_+, J_+)$, and $(H_{13}, J_{13})$ and $(H_{34}, J_{34})$ are $s$-independent.  By Proposition \ref{prop:continuation-homotopy-geometry} we can find a homotopy definition datum $\mathfrak{H}_h$ for $\CC(H_h, J_h)$ such that 
\begin{itemize}
    \item The associated homotopy definition datum $\mathfrak{H}_{h, 12}$ for $\CC(H_{12}, J_{12})$ is a stabilizations of $\mathfrak{H}_{s}$, 
    \item For every $x \in Fix(H_-)$, we get three distinct objects $(x_1, x_3, x_4) \in \CC(H_h, J_h)$ with $x_i \in \CC(H_i, J_i)$ corresponding to the Hamiltonian trajectory $x$ of $H_i = H_-$. Because the actions of $x_i$ are all the same, both $\CC(H_{13}, J_{13})(x_1, x_3)$ and $\CC(H_{34}, J_{34})(x_3, x_4)$ are each a single point corresponding to the (regular) constant continuation map trajectory at $x$.  We can thus arrange for the perturbation function $\lambda_{V^1(x_i, x_j)}$ for the perturbation data $\mathfrak{P}_h$ underlying $\mathfrak{H}_h$ associated to $(x_1, x_3)$ and to $(x_3, x_4)$ to be zero; this is possible because in the proof of Proposition \ref{prop:continuation-homotopy-geometry}  (or more saliently in the proof of Proposition \ref{prop:perturbation-data-continuation-homotopy-exist}) there are no constraints arising from boundary strata to $\CC(x_3,x_4)$ or $\CC(x_1,x_3)$ in the inductive procedure for choosing perturbation data. 
    \item The homotopy definition data $\mathfrak{H}_{i}$ associated to $\CC(H_{i}, J_i)$ for $i=1,3,4$ are all stabilization of the same homotopy definition datum $\mathfrak{H}_-$.
\end{itemize}
Choose a strict grading $\alpha_-$ on $\CC(H_-, J_-)$, i.e. a map $\alpha_-: \CC(H_-, J_-) \to \Z$ such that if $x > y$ then $\alpha_-(x) > \alpha_-(y)$. By Lemma \ref{lemma:filtrations-of-floer-homotopy-types} and Equation \ref{eq:homotopy-equivalence-of-cones} this defines filtrations of
 $|\CC(H_i, J_i)_{\mathfrak{H}, x}|$.  An inductive argument shows that the posets $\CC(H_{13}, J_{13}) \simeq \CC(H_{34}, J_{34}) \simeq \CC(H_s, J_s)$ respect these filtrations in the sense of Lemma \ref{lemma:filtered-maps-of-floer-homotopy-types}.
 Write
\[ f_{ij}: F_{\R \oplus V^1_{i, s_{ij}}}\Sigma^{\R \oplus V^1_{i, s_{ij}}}|\CC(H_i, J_i)_{\mathfrak{H}_i}| \to F_{\R \oplus V^1_{j, s_{ij}}}\Sigma^{\R \oplus V^1_{j, s_{ij}}}|\CC(H_j, J_j)_{\mathfrak{H}_j}| \]
for a choice of such maps associated to $\mathfrak{H}_{ij}$ by Lemma \ref{prop:map-associated-to-continuation-floer-homotopy-type}; by Lemma \ref{lemma:filtered-maps-of-floer-homotopy-types} these can be chosen to preserve the filtrations associated to the gradings. In particular, we see that for both $f_{13}$ and $f_{34}$, the maps \eqref{eq:map-of-spectral-sequence-of-filtered-flow-category} induced in the spectral sequences \eqref{eq:spectral-sequence-assoc-to-filtration} are isomorphisms, since these satisfy the condition Theorem \ref{thm:comparison-with-floer-homology} and the only trajectories contributing to the maps on spectral sequences are the constant continuation map trajectories. Thus the maps induced by $-F_{V^1_{34, s}}\Sigma^{V^1_{34, s}}f_{34}F_{V^1_{12, s}}\Sigma^{V^1_{12, s}}f_{13}$ (see Proposition \ref{prop:map-associated-to-continuation-homotopy-floer-homotopy-type}) are isomorphisms on integral homology; by Proposition \ref{prop:map-associated-to-continuation-homotopy-floer-homotopy-type}, that implies that the map induced by $F_{V^1_{24,s}}\Sigma^{V^1_{24,s}}f_{24}F_{V^1_{12,s}}\Sigma^{V^1_{12, s}}{f_{12}}$ is also an isomorphism on integral homology. So $f_{12}$ induces an isomorphism on integral homology, and thus a weak equivalence of spectra. 

A similar argument but now with a choice of continuation homotopy datum $H_h$ and homotopy definition datum $\mathfrak{H}_h$ such that $(H_1, J_1) = (H_3, J_3) = (H_4, J_4) = (H_+, J_+)$, and $(H_2, J_2) = (H_-, J_-)$, with $\mathfrak{H}_{34}$ a stabilization of $\mathfrak{H}_s$, shows that the map $f$ is a surjection on integral homology. So $f$ is an isomorphism on integral homology; but a map of finite spectra which is an isomorphism on integral homology is a weak equivalence. 
\end{proof}

\begin{proposition}
\label{prop:equivariant-floer-homotopy-well-defined-up-to-equivalence} Let $p$ be a prime.
        Suppose that we are in the setting of Proposition \ref{prop:invariance-of-floer-homotopy}, but now all data chosen are $C_{p^k}$-equivariant (and the corresponding assumption on polarization classses holds). Then the map of $C_{p^k}$-spectra  $f$ defined by Proposition \ref{prop:map-associated-to-continuation-floer-homotopy-type} on Floer homotopy types is an equivalence of $C_{p^k}$-spectra.
\end{proposition}
\begin{proof}
    Write $(H_{0, \pm}, J_{0, \pm}, H_{s, 0}, J_{s,0})$ for the data which upon $p^k$-iteration give $(H_\pm, J_\pm, H_s, J_s)$, respectively. We have that $\CC'(H_s, J_s)^{C_{p^r}}$ for $0 < r \leq k$ is exactly $\CC'(H_{s, 0}^{\#p^{k-r}}, J_{s, 0}^{\#p^{k-r}})$, where the latter is obtained by taking the $C_{p^r}$ invariant parts of the integralization data perturbation data in the natural fashion, with framing as in Lemma \ref{lemma:framing-of-invariant-category}; the latter category is equipped with embeddings and extensions of embeddings by taking $C_{p^r}$-invariants. Similar results hold for $(H_s, J_s, H_{s,0}, J_{s,0})$ replaced by $(H_\pm, J_\pm, H_{0, \pm}, J_{0, \pm})$, respectively, and we see that $V_{\pm, s}^{C_{p^r}}$ is exactly the $C_{p^{k-r}}$ representation $V_{\pm, s}$ arising in Proposition \ref{prop:map-associated-to-continuation-floer-homotopy-type} when applied to $\CC'(H_{s, 0}^{\#p^{k-r}}, J_{s, 0}^{\#p^{k-r}})$. Thus see in this setting that $f^{\Phi C_{p^r}}$ is a map that would be assigned by Proposition \ref{prop:map-associated-to-continuation-floer-homotopy-type} to $\CC'(H_{s, 0}^{\#p^{k-r}}, J_{s, 0}^{\#p^{k-r}})$. 
    
    Thus the claim follows, from the above, together with Proposition \ref{prop:invariance-of-floer-homotopy} applied to each of the maps $f^{\Phi C_{p^r}}$ (which is sensible by the previous observation), together with the fact that a map of $C_{p^k}$-spectra that is a weak equivalence on ordinary spectra of all geometric fixed point spectra is a weak equivalence of $C_{p^k}$-spectra \cite{schwede2019lectures}.
\end{proof}
\begin{remark}
    We can easily upgrade this proposition to its $C_k$-equivariant analog by choosing equivariant perturbation data in an inductive fashion indexed by the lattice of subgroups in $C_k$, but we have not done so here.
\end{remark}

\begin{proof}[Proof of Theorem \ref{thm:equivariant-floer-spectra-exist}]
This follows immediately from Proposition \ref{prop:equivariant-floer-homotopy-well-defined-up-to-equivalence} and Proposition \ref{eq:compute-for-admissible-Floer-data}.
\end{proof}

\begin{proof}[Proof of Theorem \ref{thm:fixed-point-thm}]
This is implicit in the proof of Proposition \ref{prop:equivariant-floer-homotopy-well-defined-up-to-equivalence}. 
\end{proof}

\subsection{Constructing symplectic cohomology and the cyclotomic structure maps.}
\label{sec:cyclotomic-structure-cofibrancy}

We now wish to make sense of the homotopy colimit  of the diagram \eqref{eq:spectral-lift-of-floer-homology-colimit} where all spectra and maps are defined with the constructions of Section \ref{sec:geometric-results-review}, choosing the homotopy definition data $\mathfrak{H}_k$ such that the statements in Proposition \ref{prop:compatible-homotopy-continuation-data} hold.  This would be the spectral analog of the definition of symplectic cohomology. To to this we should take a  ``fibrant-cofibrant'' replacement of this diagram, and then subsequently form the mapping telescope on it. However, this desire is complicated by the fact that we are repeatedly changing groups throughout the diagram, so the diagram does not sit in any one model category; and moreover that none of the function spectra we have written are homotopical because $\SS$ is generally not fibrant. We will make the desired ``cofibrant-fibrant replacement'' precise below.

 Let $c_k$, $f_k$ denote the cofibrant and fibrant replacement functors in the category of orthogonal $C_{p^k}$-spectra, and for $k > \ell$, let  $R^k_\ell$ be the restriction functors from orthogonal $C_{p^k}$-spectra to orthogonal $C_{p^\ell}$-spectra. Explicitly, we will replace the \eqref{eq:unreplaced-diagram} with a digram for which the $k$-th term is $c_k f_k F(|H^k|,f_k \SS)$, and then define the symplectic cohomology to be be a mapping telescope assocaited to this diagram.

 The basic challenge for performing this replacement 
\[ c_{k-1}R^{k}_{k-1} \neq R^{k}_{k-1}c_k, f_{k-1}R^{k}_{k-1} \neq R^{k}_{k-1}f_k.\]

However, let us recall Lemma \ref{lemma:group-retriction}: let $X$ be an $H$-spectrum, $Y$ be a $G$-spectrum, and $t: X \to Y$ be an $H$-map of spectra. For precision we will denote this by 
\[ t: X \to R^{G}_H Y.\]
Then we can form the following diagrams of spectra:
\begin{equation}
\label{eq:cofibrant-fibrant-relative-group-replacement}
    \begin{tikzcd}
        c_HX \ar[r, "c_H t"] & 
        c_HR^G_H Y & 
        cHR^G_H c_G Y \ar[l] \ar[r] &
        R^G_H c_G Y, \\
        f_HX \ar[r, "f_H t"] &
        f_H R^G_H Y \ar[r] & 
        f_H R^G_H f_G Y & 
        R^G_H f_G Y \ar[l]
    \end{tikzcd}
\end{equation}
Here $(c_H, f_H, c_G, f_G)$ denote the cofibrant and fibrant replacement functors for orthogonal $H$ and $G$-spectra. Moreover, there is a diagram 
\begin{equation}
\label{eq:c-to-cf}
    \begin{tikzcd}
        c_H f_H X \ar[r] & 
        c_H f_H R^G_H Y\ar[r] & 
        c_H f_H R^G_H f_G Y & 
        c_H f_H R^G_H c_G f_GY \ar[l] \ar[r] & 
        f_H R^G_H c_G f_GY  &
        R^G_H c_G f_G Y \ar[l]  \\
        c_H X \ar[u] \ar[r] &
        c_H R^G_H Y \ar[u] \ar[r, equal]&
        c_H R^G_H Y \ar[u]&
        c_H R^G_H c_G Y \ar[u] \ar[l] \ar[r] &
        R^G_H c_G Y \ar[u] \ar[r, equal] &
        R^G_H c_G Y \ar[u].
    \end{tikzcd}
\end{equation}

\pagebreak
\newgeometry{bottom=0.1cm}
\begin{sidewaysfigure}
\relsize{-1}
  \centering\thisfloatpagestyle{floatpage}%
\begin{equation*}
        \begin{tikzcd}
    \Phi^{C_p} c_k F(|H^k|, f_k\SS) \ar[r]\ar[d] &
    \Phi^{C_p} c_k F(|H^k|, f_k f_{k+1}\SS)\ar[d]&
    \Phi^{C_p} c_k F(|H^k|, f_{k+1}\SS) \ar[r] \ar[l]\ar[d] & 
    \Phi^{C_p} c_k F(|H^{k+1}|, f_{k+1}\SS) \ar[d]& 
    \Phi^{C_p} c_k c_{k+1} F(|H^{k+1}|, f_{k+1}\SS) \ar[l] \ar[d] \ar[r] &
    \Phi^{C_p} c_{k+1} F(|H^{k+1}|, f_{k+1}\SS) \ar[r] \ar[d]&
    \cdots\\
    F(\Phi^{C_p}|H^k|, \Phi^{C_p}f_k\SS) \ar[r] \ar[d]&
    F(\Phi^{C_p}|H^k|, \Phi^{C_p}f_k f_{k+1}\SS) \ar[d]&
    F(\Phi^{C_p}|H^k|, \Phi^{C_p}f_{k+1}\SS) \ar[r] \ar[l] \ar[d] & 
    F(\Phi^{C_p}|H^{k+1}|, \Phi^{C_p}f_{k+1}\SS) \ar[d]& 
    F(\Phi^{C_p}|H^{k+1}|, \Phi^{C_p}f_{k+1}\SS) \ar[l] \ar[r] \ar[d]&
    F(\Phi^{C_p}|H^{k+1}|, \Phi^{C_p}f_{k+1}\SS) \ar[r] \ar[d]&
    \cdots\\
    F(\Phi^{C_p}|H^k|, f_{k-1}\Phi^{C_p}f_k\SS) \ar[r] &
    F(\Phi^{C_p}|H^k|, f_{k-1}f_{k}\Phi^{C_p}f_k f_{k+1}\SS) &
    F(\Phi^{C_p}|H^k|, f_k\Phi^{C_p}f_{k+1}\SS) \ar[r] \ar[l] & 
    F(\Phi^{C_p}|H^{k+1}|, f_k\Phi^{C_p}f_{k+1}\SS) & 
    F(\Phi^{C_p}|H^{k+1}|, f_k\Phi^{C_p}f_{k+1}\SS) \ar[l] \ar[r] &
    F(\Phi^{C_p}|H^{k+1}|, f_k\Phi^{C_p}f_{k+1}\SS) \ar[r] & \cdots
    \\
    c_{k-1} F(|H^{k-1}|, f_{k-1}\SS) \ar[r] \ar[u]& 
    c_{k-1}F(|H^{k-1}|, f_{k-1} f_{k}\SS) \ar[u]& 
    c_{k-1} F(|H^{k-1}|, f_k\SS) \ar[r] \ar[l] \ar[u]& 
    c_{k-1} F(|H^k|, f_k\SS) \ar[u]& 
    c_{k-1} c_{k} F(|H^{k}|, f_{k}\SS) \ar[l] \ar[r]\ar[u] &
    c_{k} F(|H^{k}|, f_{k}\SS) \ar[r] \ar[u]& \cdots 
    \end{tikzcd}
\end{equation*}
  \caption{A large commutative diagram. The diagram continues on the right as it starts on the left.}
  \label{fig:large-diagram}
\end{sidewaysfigure}
\thispagestyle{empty}
\pagebreak

\restoregeometry

Below, we will discuss `diagrams' several times; however, each time, we are considering objects and arrows in a various homotopy categories satisfying a relation that can be depicted in a diagrammatic fashion, rather than an actual functor from some simple indexing category to a given category. We will explain in each case what the relations are. 
\begin{theorem}
\label{thm:sh-equivariant-data}
The diagram \eqref{eq:unreplaced-diagram} canonically defines associated objects and morphisms, depicted
\begin{equation}
    \label{eq:diagram-in-homotopy-category}
    F(|H^0|, \SS) \xrightarrow{[t_0]} F(|H^1|, \SS)\xrightarrow{[t_1]} \cdots
\end{equation} 
where $F(|H^k|, \SS)$ is thought of an object in $Ho(C_{p^k}-OSp)$, and $[t_k] \in Ho(C_{p^k}-OSp)$ is a map $F(|H^k|, \SS) \to R^{k+1}_k F(|H^k|, \SS)$. In fact all objects in the above diagram have canonical representatives $c_kf_k F(|H^k|,f_k\SS) \in C_{p^k}-OSp$.

Similarly, there are objects in morphisms depicted as

\begin{equation}
\label{eq:commutative-diagram-in-homotopy-category}
\begin{tikzcd}
    \Phi^{C_p} F(|H^1|, \SS)\ar[r, "\Phi^{C_p}{[}t_1{]}"] \ar[d, "{[}\phi_1{]}"] &
    \Phi^{C_p}F(|H^2|, \SS)\ar[r, "\Phi^{C_p}{[}t_2{]}"] \ar[d,"{[}\phi_2{]}"] &
    \cdots \\
     F(|H^0|, \SS)\ar[r, "{[}t_0{]}"] &
    F(|H^1|, \SS)\ar[r, "{[}t_1{]}"]  &
    \cdots 
\end{tikzcd}
\end{equation}

where the horizontal maps are induced from the maps \eqref{eq:diagram-in-homotopy-category} functorially, vertical maps $[\phi_{k+1}]$ are all isomorphisms in $Ho(C_{p^k}-OSp)$, and 
\begin{equation}
    \label{eq:cyclotomic-homotopies}
    [\phi_{k+1}] \Phi^{C_p}[t_{k+1}] = [t_k][\phi_{k+1}]\text{ in }Ho(C_{p^k}-OSp).
\end{equation}

Choosing representatives $t'_k$ of $[t_k]$ and defining $SH(M, \SS)_{\geq r}$ to be the mapping telescope 
and choosing representatives of $\phi_k$ of $[\phi_k]$ and homotopies representing \eqref{eq:cyclotomic-homotopies}, one has stable equivalences 
\begin{equation}
    (\tilde{\phi}^r: \Phi^{C_p}SH(M, \SS)_{\geq r} \to SH(M, \SS)_{\geq r-1}) \in C_{p^{r-1}}-OSp. 
\end{equation} 
Moreover, there are objects and morphisms, depicted pictorially as

\begin{equation}
\label{eq:compatibility-of-cyclotomic-maps-and-change-of-groups}
        \begin{tikzcd}
            \Phi^{C_p}SH(M, \SS)_{\geq 1} \ar[d] &
            \Phi^{C_p}SH(M, \SS)_{\geq 2} \ar[l] \ar[d] & 
            \Phi^{C_p}SH(M, \SS)_{\geq 3} \ar[l] \ar[d] &
            \cdots \ar[l] \\
            SH(M, \SS)_{\geq 0} &
            SH(M, \SS)_{\geq 1}\ar[l]  &
            SH(M, \SS)_{\geq 2} \ar[l] &
            \cdots  \ar[l]
        \end{tikzcd}
\end{equation}

where the vertical arrows are the $\tilde{\phi}^k$, the horizontal arrows are the inclusion maps $i_{k+1, k}$ on mapping telescopes (see \eqref{eq:inclusion-of-mapping-telescopes}). This diagram commutes strictly in the sense that $\tilde{\phi}^ki_{k+1, k} = i_{k, k-1}\tilde{\phi}^{k+1}$ in $C_{p^{k-1}}-OSp$.

\end{theorem}

\begin{proof}
The second diagram in \eqref{eq:cofibrant-fibrant-relative-group-replacement} defines maps 
\begin{equation}
    \label{eq:sphere-fibrant-replacement}
    f_k \SS \to f_k f_{k+1} \SS \to f_{k+1} \SS.
\end{equation}
Above, and for the rest of this proof (including Figure \ref{fig:large-diagram}) we will adopt the convention that 
\begin{itemize}
    \item We will drop all symbols denoting group restriction. Thus $f_kf_{k+1}\SS$ is a $C_{p^k}$-spectrum, while $f_{k+1}\SS$ is a $C_{p^{k+1}}$-spectrum.
    \item All maps are ``as equivariant as possible'', in the sense that if we have a map $f: X' \to Y'$ where $X'$ is produced by some cofibrant and fibrant replacements from $X'$ and similarly $Y'$ from $Y$, then $f$ is equivariant with respect to the maximal subgroup of $S^1$ that it could be given the kinds of spectra that $X'$ and $Y'$ are.
\end{itemize}

All ambiguities regarding insertion of $R^G_H$ into formulae can be resolved by Equations \eqref{eq:geometric-fixed-points-and-free-spectra}, \eqref{eq:restriction-and-function-spectra}, and \eqref{eq:restriction-and-geometric-fixed points}.

The diagram \eqref{eq:sphere-fibrant-replacement} lets us replace each map in \eqref{eq:unreplaced-diagram} with a zig-zag of maps:

\begin{equation}
    \begin{tikzcd}
    F(|H^k|, f_k \SS) \ar[r] & F(|H^k|, f_k f_{k+1} \SS)  F(|H^k|, f_{k+1} \SS)  \ar[l] \ar[r] & F(|H^{k+1}|, f_{k+1} \SS).
    \end{tikzcd}
\end{equation} 
Using the top diagram of \eqref{eq:cofibrant-fibrant-relative-group-replacement} then lets us produce the bottom row of Figure \ref{fig:large-diagram}. Applying $\Phi^{C_p}$ to the bottom row and shifting the indices over produces the top row of Figure \ref{fig:large-diagram}.

Now there are canonical maps
    \[ \Phi^H cF(X, Y) \to F(\Phi^H X, \Phi^H Y)\]
    adjoint to
    \[ \Phi^H cF(X, Y) \wedge \Phi^H \to \Phi^H (cF(X, Y) \wedge Y) \to \Phi^H(F(X,Y) \wedge Y) \to \Phi^H(Y). \]
    These maps are stable equivalences whenever $X$ is $F_V Y$ for some \emph{finite} $G$-CW complex $Y$ and $Y$ is fibrant. These maps produce maps to the second row of Figure \ref{fig:large-diagram}, and the maps from  second row to the third row are produced by maps to fibrant replacements. 

    The maps from the fourth row to the third row are built from the strictly commutative diagram 
   \begin{equation}
   \label{eq:diagram-2}
        \begin{tikzcd}
            F(|H^{\# p^0}|, \SS) \ar[r] \ar[d] & F( |H^{\# p^1}|, \SS)  \ar[r] \ar[d] & F( |H^{\# p^2}|, \SS) \ar[r] \ar[d] & \cdots \\
            F(\Phi^{C_p} |H^{\# p^1}|, \SS) \ar[r]  & F(\Phi^{C_p} |H^{\# p^2}| , \SS) \ar[r] &F(\Phi^{C_p} |H^{\# p^3}| , \SS)    \ar[r] & \cdots
        \end{tikzcd}
    \end{equation}
    by performing the same operation to the bottom row of \eqref{eq:diagram-2} as is used to produce the bottom row of Figure \ref{fig:large-diagram} from the top row of \eqref{eq:diagram-2}, and then applying maps from cofibrant and to fibrant replacements.

     The fact that there is a diagram like \eqref{eq:diagram-2} follows from Proposition \ref{prop:compatible-homotopy-continuation-data} and some careful index checking. Indeed, Proposition \ref{prop:compatible-homotopy-continuation-data} implies that there is a natural transformation of diagrams from the diagram with maps $(F_{\R \oplus \tilde{V}(\mathfrak{H}^-_k)} g_1^k)^{\Phi C_p}$ to the diagram with maps 
    $(F_{\R \oplus \tilde{V}(\mathfrak{H}^-_{k-1}) \oplus \tilde{V}(\mathfrak{H}^-_{k-1}, s)^{C_p}} \Sigma^{\tilde{V}(\mathfrak{H}^-_{k-1}, s)^{C_p}}g_1^{k-1})$, and the maps \eqref{eq:desuspension-maps} give a natural transformation from the latter diagram to the diagram with maps $(F_{\R \oplus \tilde{V}(\mathfrak{H}^-_{k-1}) } g_1^{k-1})$.  Similarly, Proposition \ref{prop:compatible-homotopy-continuation-data} gives a natural transformation from the diagram with maps $(g_2^k)^{\Phi^k_2}$ to the diagram with maps 
    \[(F_{\tilde{V}(\mathfrak{H}_{k-1}, s)^{C_p}} \Sigma^{\tilde{V}(\mathfrak{H}_{k-1}, s)^{C_p}})g^{k-1}_2 \] 
    (the diagram where we insert the extra representation into the free spectrum and the suspension) and then again \eqref{eq:desuspension-maps} lets us remove the copies of $F_{\tilde{V}(\mathfrak{H}_{k-1}, s)^{C_p}} \Sigma^{\tilde{V}(\mathfrak{H}_{k-1}, s)^{C_p}}$.  This pair of natural transformations fits together to define a natural transformation from the diagram with maps $(g^k)^{\Phi C_p}$ to the diagram with maps $g^{k-1}$. One then gets the diagram by applying $F(\cdot, \SS)$.

    It is tedious but straightforward to check that the diagram Figure \ref{fig:large-diagram} commutes strictly, with all vertical arrows weak equivalences. Using \eqref{eq:c-to-cf}, after extending the right most column of Figure \ref{fig:large-diagram} with three more identical columns, one gets a bigger diagram where the bottom row maps into a row like the current bottom row of 
    Figure \ref{fig:large-diagram} but extended like the top row of \eqref{eq:c-to-cf} such that all objects are cofibrant fibrant, and similarly the top row of the new diagram admits a map to $\Phi^H$ applied to the bottom row. The whole diagram still strictly commutes, and the objects on the top corners are cofibrant and the objects on the bottom corners are fibrant. Looking at this diagram in the homotopy category of $G$-spectra produces the diagrams \eqref{eq:diagram-in-homotopy-category}, \eqref{eq:commutative-diagram-in-homotopy-category}. 

    Choosing representatives $t'_k$ as in the statement of the lemma defines representatives $\Phi^{C_p}t'_k$ of $\phi^{C_p}[t'_k]$. Choosing representatives of $[\phi_k]$ and the homotopies as in the lemma, and using the fact that 
    mapping telescopes commute with geometric fixed points, defines the maps $\tilde{\phi}^k$. These are weak equivalences by Lemma \ref{lemma:weak-map-of-mapping-telescopes}. 
    
    The fact that these fit into \eqref{eq:compatibility-of-cyclotomic-maps-and-change-of-groups} follows immediately from basic facts about the mapping telescope reviewed in Section \ref{sec:orthogonal-spectra-detailed-review}. 
   
\end{proof}

\subsection{Comparison with Nikolaus-Scholze.}
\label{sec:nikolaus-scholze-comparison}

Let $N$ denote the nerve of a $1$-category and let $N^{hc}$ denote the homotopy-coherent nerve of a simplicially-enriched category. \cite[\href{https://kerodon.net/tag/00KM}{Tag 00KM}]{kerodon}. Given a category $\CC$, there is a simplicially enriched category associated to $\CC$ with underlying category $\CC$ and the morphisms from $x$ to $y$ the simplicial set obtained by taking the constant simplicial object on $\CC(x,y)$; we will also denote this category by $\CC$. There is a map of simplicial sets \cite[\href{https://kerodon.net/tag/00KZ}{Tag 00KZ}]{kerodon}.
\begin{equation}
\label{eq:nerve-is-hc-nerve}
    N(\CC) \simeq N^{hc}(\CC).
\end{equation} 

 Given a model category  $\CC$ with objects contained in the Grothendieck universe $\overline{\mathcal{U}}$, we will denote the associated hammock localization by $L^{H}(\CC)$ \cite{dwyer1980calculating, hinich2017lectures}. Thus, writing $L^{H}(\CC)^f$ for the Begner fibrant replacament \cite{bergner2007model}, there is a simplicial set $N^{hc}(L^H(\CC)^f)$ which is an $\infty$-category, and we call the \emph{$\infty$-category associated to $\CC$.}. 

There is thus a functor of $\infty$-categories from $N(\CC)$ to the $\infty$-category associated to $\CC$. In fact, the latter is a localization of the $\infty$-category $N(\CC)$ at the weak equivalences \cite{dwyer1980calculating, hinich2013dwyer}. In particular, our notion of of the $\infty$-category $C_{p^n}Sp$ underlying the homotopy theory of orthogonal $C_{p^n}$-spectra agrees with that of \cite{nikolaus-scholze}.

\begin{definition}(Nikolaus-Scholze \cite{nikolaus-scholze})
    The $\infty$-category of genuine $C_{p^\infty}$-spectra is the $\infty$-category 
    \[ C_{p^\infty}Sp := \varprojlim C_{p^k}Sp\]
    where the inverse limit is taken along the group restriction functors $F^{C_{p^{k+1}}}_{C_{p^k}}$ (or rather the functors obtained from applying  $N^{hc}(L^H(\cdot)^f)$ to these functors). 

    The functor $\Phi^{C_p}$ induces functors 
    \[ \Phi^{C_p}: C_{p^k}Sp \to C_{p^k}/{C_p}Sp \simeq C_{p^{k-1}}Sp.\]
    for all $p$, and commutes with $F^{C_{p^{k+1}}}_{C_{p^k}}$. It thus induces a functor 
    \[ C_{p^\infty}Sp \to (C_{p^\infty}/C_p)Sp\simeq C_{p^\infty}Sp\]

    The $\infty$-category of genuine $p$-cyclotomic spectra is the equalizer of the diagram 
    \begin{equation}
    \label{eq:equalizer-diagram}\begin{tikzcd}
        C_{p^\infty}Sp \ar[r, yshift=1ex, "Id"]
  \arrow[r, yshift=-1ex,swap,"\Phi^{C_p}"] &C_{p^\infty}Sp.
    \end{tikzcd} 
    \end{equation}
\end{definition}

\begin{theorem} 
\label{thm:compare-w-nikolaus-scholze}
    The data of Theorem \ref{thm:sh-equivariant-data} define an associated genuine $p$-cyclotomic spectrum in the sense of Nikolaus-Scholze. 
\end{theorem}
\begin{proof}
    This is essentially an exercise in standard $\infty$-categorical combinatorics. We first recall two notions and some useful theorems. 

    First, recall that a category can be regarded as a simplicial category with the simplicial sets of morphisms being discrete; we will use this identification to speak about functors from categories to simplicial categories. Equivalently, a functor from a category to a simplicial category is a functor from the category to the underlying category of the simplicial category.

    The simplicial category $QCat$ has objects all small $\infty$-categories and morphisms $QCat(X,Y)$ the core of the $\infty$-category of functors $Fun(X,Y)$. The $\infty$-category of $\infty$-categories $\mathcal{QC}$ is the homotopy-coherent nerve of $QCat$.

    Suppose we have a category $\CC$ and a functor $\mathcal{F}: \CC \to QCat$. To this data, one can associate
the weighted nerve $N_{\bullet}^{\mathcal{F}(\CC)}$ \cite[\href{https://kerodon.net/tag/025X}{Tag 025X}]{kerodon}, which admits a map 
        \[ U: N_{\bullet}^{\mathcal{F}}(\CC) \to  \mathcal{N}(\CC)\]
        which is a cocartesian fibration  Moreover, objects of the $\infty$-category $N_{\bullet}^{\mathcal{F}(\CC)}$ are tuples $(C, x)$, where $C \in Ob(\CC)$ and $x$ is an object of $\mathcal{F}(C)$, and $1$-morphisms are tuples $(f, e): (C,x) \to (D, y)$ where $f: C \to D$ is a morphism in $\CC$ and $e: \mathcal{F}(f)(x) \to y$ is a 1-morphism in $\mathcal{F}(D)$. Finally, $(f,e)$ is $U$-cocartesian if and only if $e$ is an isomorphism in the $\infty$-category $\mathcal{F}(D)$. \cite[\href{https://kerodon.net/tag/046Y}{Tag 046Y}]{kerodon}.
        
    The weighted nerve is an $\infty$-category, and there is a map of $\infty$-categories 
    \[N_{\bullet}^\mathcal{F}(\CC) \xrightarrow{\theta} \int_{\mathcal{N}(\mathcal{F})} \mathcal{N}^{hc}(F). \]
    which is an equivalence of $\infty$-categories and gives a bijection on $U$-cocartesian morphisms \cite[\href{https://kerodon.net/tag/027J}{Tag 027J}]{kerodon}.

    Now, Theorem \ref{thm:sh-equivariant-data} produces the data 
    \begin{equation}
        \label{eq:data-1}
        SH(M, \SS)_0 \xleftarrow{i_{01}} SH(M, \SS)_1 \xleftarrow{i_{12}} SH(M, \SS)_2 \xleftarrow{i_{23}} \cdots 
    \end{equation}
    where $SH(M, \SS)_r \in C_{p^r}-OSp$ and \[i_{ab}: F^{C_{p^b}}_{C_{p^a}} SH(M, \SS)_a \to SH(M, \SS)_b\] 
    for $b>a$ are the corresponding inclusions of mapping telescopes, and are all stable equivalences of orthogonal spectra. 

    Write $\Z_{\geq 0}$ for the  category  corresponding to the poset of nonnegative integers ordered by magnitude There is a functor 
    \[ \mathcal{F}: \Z_{\geq 0}^{op} \to 1-Cat\]
    mapping $n$ to $C_{p^n}-OSp$ and mapping $n \to m$ to $F^{C_{p^n}}_{C_{p^m}}$. This defines an associated functor of $\infty$-categories
    \[ N(F): \Z_{\geq 0}^{op} \to QCat\]
    which sends each $1$-morphism $n \to m$ to $N(F^{C_{p^n}}_{C_{p^m}})$.
    There is another functor 
    \[ N^{hc}(L^H(\mathcal{F})^h): Z^{op} \to QCat\]
    which sends each $1$-morphism $n \to m$ to $N^{hc}(L^H(F^{C_{p^n}}_{C_{p^m}})^f)$ where we are are using the fact that $F^{C_{p^n}}_{C_{p^m}}$ sends $C_{p^n}-$stable equivalences to $C_{p^m}$-stable equivalences to define $L^H(F^{C_{p^n}}_{C_{p^m}})$. Note that 
    \[N^{hc}(L^H(\mathcal{F}^h))(n) = C_{p^n}Sp,\]
    the $\infty$-category underlying the homotopy theory of genuine $C_{p^n}$-spectra. There is a natural transformation of functors 
    \begin{equation}
        \label{eq:nerve-to-coherent-nerve}
        N(\mathcal{F}) \to  N^{hc}(L^H(\mathcal{F}^h))
    \end{equation} 
    induced by the maps 
    \[ N(\CC) \to N^{hc}(L^H(\CC)^f).\]
    This defines a functor 
    \begin{equation}
    \label{eq:map-of-nerves}
        N_\bullet^{N(\mathcal{F})} \to N_\bullet^{N^{hc}(L^H(\mathcal{F}^h)}.
    \end{equation}

    The data \eqref{eq:data-1} produces an element of $Fun(N(Z^{op}), N_{\bullet}^{N(\mathcal{F})})$ which sends $n \to m$ to $i_{nm}$. There is a corresponding element of 
    \[Fun(N(Z^{op}),N_\bullet^{N^{hc}(L^H(\mathcal{F}^h)})\]
    sending $n \to m$ to $i_{nm}$, related to the previous element by the functor above between weighted nerves.
    
    Each of these morphisms is $U$-cocartesian for 
    \[ U: N_{\bullet}^{N^{hc}(L^H(\mathcal{F}^h)}) \to N(Z^{op}) \]
    since it is an isomorphism in the $\infty$-category $C_{p^n}-Sp$. Applying $\theta$, we have produced an element of 
    \begin{equation}
        \label{eq:model-of-Cpinfty-spectra}
        \underleftarrow{\text{holim}}(\mathcal{F}) := Fun^{U-coCart}_{\mathcal{N}(Z^{op})}(N(Z^{op}), N^{hc}(N^{hc}(L^H(\mathcal{F})^h).
    \end{equation} 

    Now, $\underleftarrow{\text{holim}}(\mathcal{F})$ is a limit of the diagram $N(N^{hc}(L^H(\mathcal{F}^h)))$ \cite[\href{https://kerodon.net/tag/03AB}{Tag 03AB}]{kerodon}; we have produced an element of $C_{p^\infty}Sp$. 

    To argue that the data of the cyclotomic structure maps enhances this to an genuine $p$-cyclotomic spectrum, we argue similarly. We first detail the construction of the functor 
    \[ \Phi^{C_p}: C_{p^\infty}Sp \to C_{p^\infty}Sp.\]
    Let $\mathcal{F}': Z_{\geq 0}^{op} \to 1-Cat$ be the unique functor sending $0$ to $\Delta^0$, and otherwise $\mathcal{F}(n \to m) = N^{hc}(L^H(\mathcal{F}')^h)(n+1 \to m+1)$ for all $n, m \in \Z$. There is a natural transformation 
    \[ N^{hc}(L^H(\mathcal{F})^h) \to \mathcal{F}'\] given by applyin $N^{hc}(L^H(\Phi^{C_p})^f): C_{p^n}Sp \to C_{p^{n-1}}Sp$ for all $n \geq 1$, and applying the unique map  $C_{p}Sp \to \Delta^0$ for $n=0$. This induces a functor $N^{N^{hc}(L^H(\mathcal{F})^h)}_\bullet \to N^{\mathcal{F}'}_\bullet$  preserving $U$-cocartesian morphisms, and by postcomposition this defines a map 
    \[ \underleftarrow{\text{holim}}(\mathcal{F}) \to \underleftarrow{\text{holim}}(\mathcal{F'}) \simeq \underleftarrow{\text{holim}}(\mathcal{F})\]
    (see \cite[\href{https://kerodon.net/tag/03AC}{Tag 03AC}]{kerodon}); here the last isomorphism uses the fact that, writing $Z_{\geq 1} \subset Z_{\geq 0}$ for the full subcategory on the strictly positive integer, cocartesian sections of the cocartesian fibration $N^{\mathcal{F}'}_\bullet$ over $N(Z_{\geq 1}^{op})$ extend uniquely to $N(Z_{\geq 0}^{op})$, and this induces an isomorphism of simplicial sets.   

    We now have the strictly commutative diagram of $\infty$-categories 
    \eqref{eq:equalizer-diagram}, and we use the same method to produce an element of its limit. Explicitly, we can produce a genuine $p$-cyclotomic spectrum by picking an element $X \in C_{p^\infty}Sp$ and an isomorphism $\Phi^C_p X \to X$ in $C_{p^\infty}Sp$. Specifically, the fact that there are equivalences $\Phi^{C_p}(SH(M, \SS)_{\geq r}) \to \Phi^{C_p}(SH(M, \SS)_{\geq r-1})$ such that $\Phi^{C_p}i_{ab} = i_{a-1, b-1}$ as in Theorem \ref{thm:sh-equivariant-data} satisfying the commutative diagram \eqref{eq:compatibility-of-cyclotomic-maps-and-change-of-groups}
    defines a map in the category of sections $Fun_{/N(Z_{\geq 0}^{op})}(N(Z_{\geq 0}^{op}), N^{N(\mathcal{F})})$, which becomes an isomorphism in $C_{p^\infty}-Sp$ after applying the map \eqref{eq:map-of-nerves} \cite[\href{https://kerodon.net/tag/035R}{Tag 035R}]{kerodon}. 

\end{proof}

Using this data, one can define the $p$-typical topological cyclic homology groups $TR(SH(M, \SS), p), TC(SH(M, \SS), p)$ as in \cite{bokstedt1993cyclotomic} (see \cite[Definition II.4.4]{nikolaus-scholze}). There is a functor from genuine $p$-cyclotomic spectra to $p$-cyclotomic spectra in the sense of Nikolaus-Scholze which is an equivalence on the subcategory of bounded below objects; however, in general symplectic cohomology is bounded neither above nor below.

\begin{proof}[Proof of Theorem \ref{thm:sh-is-cyclotomic}]
It suffices to compute the homotopy groups of $H\Z \wedge SH^\bullet(M, \SS)$ as an ordinary non-equivariant spectrum. But it is constructed via the homotopy colimit \eqref{eq:diagram-in-homotopy-category};  smashing commutes with homotopy colimits, and on homology, by the functorial isomorphisms $H\Z \wedge F(X, \SS) \simeq F(X, H \Z)$ for finitely-built spectra $X$, Proposition \ref{prop:invariance-of-floer-homotopy} identifies $\pi_*H\Z \wedge F(|H^i|, \SS)$ with $HF^{-*}(\widetilde{H^{\#p^k}}; \Z)$, where $\widetilde{H^{\#p^k}}$ is a non-equivariant perturbation of $H^{p^k}$ which is now nondegenerate, and the new homotopy definition data for $F(|widetilde{H^{\#p^k}}|, \SS)$ are chosen such that that proposition holds. To identify the action of $H\Z \wedge [t_i]$ on homotopy groups with the corresponding action of continuation maps on $HF^{-*}(H^{\#p^k}; \Z)$, one constructs commutative diagrams via Proposition \ref{prop:continuation-homotopy-geometry} and Proposition \eqref{prop:map-associated-to-continuation-homotopy-floer-homotopy-type} where the new horizontal maps are induced by \emph{regular, non-equivariant} admissible continuation data, which satisfy the regularity hypotheses of Theorem \ref{thm:comparison-with-floer-homology} letting us make the comparison with the usual Floer-theoretic continuation map.  
\end{proof}

\begin{remark}
    While this document does not contain a complete proof that the cyclotomic spectrum $SH^\bullet(M, \SS)$ is independent of the choices made, verifying this should simply be a slightly tedious exercise in the technology of this paper. One first argues that $SH^\bullet(M, \SS)$ is unchanged up to equivalence if one compatibly stabilizes all homotopy definition data (a fairly simple type of stabilization, for which one allows only stabilization of perturbation data by parameterizations with trivial perturbation functions, is sufficient). One then uses Proposition \ref{prop:continuation-homotopy-geometry}, and an analog of Proposition \ref{prop:compatible-homotopy-continuation-data} for homotopy definition data to build a `map of diagrams' of the form  \eqref{eq:unreplaced-diagram}. Performing the functorial replacements of Proposition \ref{thm:sh-equivariant-data} then defines maps which give maps of diagrams \eqref{eq:diagram-in-homotopy-category} satisfying compatibility conditions with respect to the geometric fixed point functors, as in \eqref{eq:cyclotomic-homotopies}. These maps are pieced together to define a map of the mapping telescopes defining equivariant maps of the spectra $SH(M, \SS)_{\geq 0}$, and this should define exactly a map in the corresponding infinity category of genuine $p$-cyclotomic spectra. Checking that it is an equivalence is then done on homotopy groups, noting that the cyclotomic structure and comparison with symplectic homology (Theorem \ref{thm:sh-is-cyclotomic}) immediately implies that it is an equivalence on homology after taking geometric fixed points for any $C_{p^k}$, and is thus an a map of genuine $p$-cyclotomic spectra which is an equivalence of underlying $C_{p^\infty}$-spectra, i.e. an equivalence of genuine $p$-cyclotomic spectra. 
\end{remark}

\pagebreak

\section{Appendix}

\subsection{Classifying spaces for $K$ theory}
\label{sec:classifying-spaces-for-k-theory}

\paragraph{(Equivariant) Grassmannians.}
We first recall the construction of the space $BO$ and $BU$, as well as their equivariant analogs. Given a map $A: V \to W$ between real inner product spaces, let $A^T: W \to V$ denote its transpose. If instead $(V, W)$ are Hermitian inner product spaces and $A$ is complex linear, let $A^\dagger: W \to V$ denote its Hermitian transpose. For a real inner product space $V$ and a Hermitian inner product space $W$, let
\begin{equation}
EO(k, V) = \{ A\in Hom(\R^k, V)  | A^TA = id\}, \; EU(k, W) = \{ A \in Mat_\C(\C^n, W) | A^\dag A = id \}.
\end{equation}
For $k \leq \dim V$ and $k \leq \dim_\C W$, there are fibrations
\begin{equation}
    EO(k, V) \to EO(k+1, V \oplus \R) \to S^n, EU(k, W) \to EU(k+1, W \oplus \C) \to S^{2n+1}
\end{equation}
show inductively that for any $C$ there is an $N= N(C, k)$  such that for $\dim V > N$ or $\dim_\C W > N$, we have $\pi_i(EO(k, V)) = \pi_i(EU(k, W) = 0$ for $0 \leq i \leq C$. There is a free $O(k)$ action on $EO(k, n)$ by precomposition and a similar free $U(k)$ action on $EU(k, n)$. Thus 
\begin{equation}
\label{eq:define-steifel-manifolds}
    EO(k) = \colim_n EO(k, \R^n), EU(k) = \colim_n EU(k, \C^n)
\end{equation}
are contractible spaces (which can be given a CW stucture, since all above maps are embeddings of manifolds) with free actions of $O(k)$ and $U(k)$. Therefore 
\[ BO(k) = EO(k)/O(k) = Gr(k, \R^\infty), BU(k) = EU(k)/U(k) = Gr_{\C}(k, \C^\infty)\] 
are classifying spaces for $O(k)$ and $U(k)$, respectively. They manifestly parameterize $k$-planes in $\R^\infty$ and $\C^\infty$, respectively. In particular, there are  $\R^k$ and $\C^k$ bundles $\eta_k^\R \to BO(k), \eta_k^\C \to BU(k)$ called universal bundles.  

The inclusions 
\begin{equation}
\label{eq:inf-dim-vector-spaces-inclusions}
\R^\infty \simeq 0 \oplus \R^\infty \subset \R^\infty, \C^\infty \simeq 0 \oplus \C^\infty \subset \R^\infty
\end{equation}
induce closed inclusions
\begin{equation}
    BO(k) \subset BO(k+1), BU(k) \subset BU(k+1) 
\end{equation}
covered by inclusions 
\begin{equation}
    EO(k) \subset EO(k+1), EU(k) \subset EU(k+1)
\end{equation}
where we send $A \in EO(k)$ to the map sending first $k$ coordinate vectors $e_1, \ldots, e_k$ to the images of $A(e_1), \ldots A(e_k)$ under \eqref{eq:inf-dim-vector-spaces-inclusions}, and send $e_{k+1}$ to the first coordinate vector in $\R^\infty$ or $\C^\infty$ respectively. These inclusions equivariant with respect to the inclusions of groups 
\begin{equation}
    O(n) \subset O(n+1), U(n) \subset U(n+1)
\end{equation}
given on matrices by $A \mapsto A \oplus 1$. 

by taking quotients by $O(k) \subset O(k+1)$, $U(k) \subset U(k+1)$ thought of as matrices with the top left $k \times k$ entries nonzero. In particular, the colimits
\begin{equation}
    EO = \colim_k EO(k), EU = \colim_k EU(k)
\end{equation}
are contractible spaces (which can be given a CW structure, using the filtrations arising from \eqref{eq:define-steifel-manifolds}) with free actions of the topological groups 
\begin{equation}
    O = \colim_n O(n), U = \colim_n U(n)
\end{equation}
and thus $BO = EO/O = \colim_k BO(k)$ and $BU = EU/U = \colim_k BU(k)$ are classifying spaces (also CW complexes) for $O$ and $U$. Manifestly, any map from $f:K \to BO$ from a compact $CW$ complex $k$ factors through some $BO(k)$, and thus the pullback $f^*\eta^\R_k$ to $K$ will be a vector bundle; the indeterminacy of the choice of $k$ implies that the virtual vector bundle $[f^*\eta^\R_k] - [\underline{\R}^k]$ is well defined. The fact that vector bundles admit complements then implies that $\Z \times BO$ and $\Z \times BU$ are classifying spaces for $KO^0$ and $KU^0$, respectively.

An equivariant analog of the entire construction is straightforward. We will only do the orthogonal case as the unitary case is completely analogous. Fix a compact Lie group $G$, and let $\mathcal{U}$ be a complete $G$-universe. The previous discussion are $(EO(k), EU(k))$ are free contractible $G$-CW complexes with free actions of $O(k)$ and $U(k)$ commuting with the $G$-actions.

Fix an isomorphism
\begin{equation}
    \mathcal{U} \simeq V_1^{\Z} \oplus V_2^{\Z} \oplus \ldots
\end{equation}
where the $\{V_i\}_{i=1}^\infty$ are a choice of irreps of $G$ such that each irreducible representation occurs exactly once. Write 
\[ \mathcal{U}_{a_1, \ldots, } \simeq V_1^{\Z_{\geq -a_1}} \oplus V_2^{\Z_{\geq -a_2}} \oplus \ldots \subset \mathcal{U}.\]
for $(a_i)_{i=1}^\infty$ any sequence of nonnegative integers with only finitely many nonzero terms.

Let $\mathcal{Z}$ be the poset of such sequences with $(a_i)_{i=1}^\infty \leq (b_i)_{i=1}^\infty$ if $a_i \leq b_i$ for all $i$. Then the inclusions $\mathcal{U}_{(a_i)} \subset \mathcal{U}_{(b_i)}$ for $(a_i) \leq (b_i)$ define a functor 
\[ F: \mathcal{Z} \to G-Top, F((a_i)) = E(\R^{a_1} \oplus \R^{a_2} \oplus \ldots, \mathcal{U}_{(a_i)})\]
with all morphisms closed inclusions (which can be refined to inclusions of $G$-CW complexes). The colimit of this functor over this filtered diagram is a space that is homeomoprhic to $EO$. This colimit is what we mean by $EO \in G-Top$, and it is straightforward to show that it is a classifying space for $KO^0_G$.

\begin{lemma}
    \label{lemma:atiyah-janich-in-general}
    Let $\mathcal{H}$ be the completion of a complete $G$-universe. There is a $G$-homotopy equivalence 
    \[ \Z \times BO \to \Fred H\]
    which induces a natural isomorphism between the functors \eqref{eq:atiyah-janich} and \eqref{eq:ko-theory-classifying-space}. 
\end{lemma}

\begin{proof}
    The argument for $G=1$ is straightforward: the proof of Theorem \ref{thm:atiyah-janich} allows one to build a map $c(Fred H) \to \Z \times BO$ which is a weak equivalence. The result then follows by Theorem \ref{thm:banach-manifold-cw-complex}. For general $G$-the equivariant analog of the proof of Theorem \ref{thm:atiyah-janich} and Theorem \ref{thm:g-banach-manifold-cw-complex}.
\end{proof}

\printbibliography

\end{document}